# Some 4-manifold geometry from hyperbolic knots in $S^3$

Clifford Henry Taubes[†]

Department of Mathematics
Harvard University
Cambridge, MA 02138

(chtaubes@math.harvard.edu)

ABSTRACT: A 4-manifold is constructed with some curious metric properties; or maybe it is many 4-manifolds masquerading as one, which would explain why it looks curious. Anyway, knots in the 3-sphere with complete finite volume hyperbolic metrics on their complements play a role in this story.

[†]Supported in part by the National Science Foundation

# 1. Introduction

The purpose of this paper is to advertise an intriguing set of oriented 4 dimensional manifolds, all homeomorphic to $(\#_3\mathbb{CP}^2)\#(\#_{23}\overline{\mathbb{CP}}^2)$ and each admitting Riemannian metrics with anti-self dual Weyl curvature. This set of manifolds might be interesting from the point of view of 4-manifold differential topology, or it might be interesting from the point of view of 4-manifold differential geometry (or perhaps from both points of view). Here is why: All of the manifolds from this set have vanishing Seiberg-Witten invariants, so it *might* be that all are diffeomorphic to one and the same manifold (perhaps $(\#_3\mathbb{CP}^2)\#(\#_{23}\overline{\mathbb{CP}}^2)$). *But*, if there are infinitely many that are diffeomorphic to one and the same manifold, then this manifold has the following curious property: An infinite set of topologically distinct 4-manifolds appear as Gromov-Hausdorff limits of sequences of metrics on this one manifold with each metric having anti-self dual Weyl curvature, fixed volume and an a priori bound on the $L^2$ norm of its Riemann curvature tensor.

These 4-manifolds with their associated anti-self dual Weyl curvature metrics are constructed from a chosen knot in $S^3$ via a version of Fintushel-Stern knot surgery [FS]. The knot in question (denoted henceforth by K) is chosen to satisfy some topological conditions, the most important one being that its complement in $S^3$ has a complete, finite volume hyperbolic metric. The knot K then labels a non-compact part of the moduli space of anti-self dual Weyl metrics on the surgered 4-manifold. In particular, there are sequences of metrics in this part of the moduli space with the afore-mentioned Gromov-Hausdorff limit being $(S^3-K)\times S^1$ with the product metric given by the complete, finite volume hyperbolic metric on $S^3-K$ with sectional curvature -1 and the length $2\pi$ metric on $S^1$. (This product metric on $(S^3-K)\times S^1$ is locally conformally flat so it has vanishing Weyl curvature.)

Pretend for the moment that the 4-manifolds for different knot choices are not diffeomorphic to one and the same manifold. The appearance of the complete, finite volume hyperbolic metric on $S^3-K$ in the Gromov-Hausdorff limit suggests (tantalizingly?) that some property of the moduli space of anti-self dual metrics that manifests itself in this case as a property of the hyperbolic metric on $S^3-K$ will distinguish the smooth structures on these manifolds. But maybe not…

The construction of the anti-self dual Weyl curvature metrics to be described momentarily uses a differential equation approach that follows for the most part Kovalev-Singer [KS]. There is also likely an equivalent twistor space construction for these anti-self dual Weyl curvature metrics along the lines of LeBrun-Singer [LS]. In fact, the existence theorem asserted in this paper would follow directly from theorems in [KS] or [LS] were there not cokernel obstructions to deal with. These cokernel issues are analyzed using ideas of Ache-Viaclovsky [AV] and they are dealt with by a strategy that was first used by Donaldson-Friedman [DF].

By way of some history, the existence theorems for anti-self dual Weyl curvature metrics on connect sums in all of the afore-mentioned papers are descendants of a gluing construction that was pioneered many years ago by Andreas Floer [Fl]. Even so, Floer's work doesn't address the case of $(\#_3 \mathbb{CP}^2) \# (\#_m \overline{\mathbb{CP}}^2)$. Anti-self dual Weyl curvature metrics on the standard smooth version of $(\#_3 \mathbb{CP}^2) \# (\#_m \overline{\mathbb{CP}}^2)$ with m ≥ 42 were first constructed by Claude LeBrun [Le]; subsequently Yann Rollin and Michael Singer [RS] found metrics on these manifolds for the cases when m ≥ 30. As a parenthetical remark, the connect sum of the K3 manifold with any n ≥ 3 copies of $\overline{\mathbb{CP}}^2$ has metrics with anti-self dual Weyl curvature (constructed by Donaldson/Friedman [DF]); this n-fold blow up of the K3 manifold is homeomorphic to, but not diffeomorphic to $(\#_3 \mathbb{CP}^2) \# (\#_{19+n} \overline{\mathbb{CP}}^2)$. Peter Braam [Br] was likely the first to construct non-product 4-manifolds with anti-self dual Weyl curvature using hyperbolic 3-manifolds. Braam's construction and the subsequent generalizations by [Ki] (see also [AKO]) use infinite volume manifolds without cusp like ends; and in any event, they lead to non-simply connected 4-manifolds.

The remainder of this paper is organized as follows: Section 2 describes the manifolds in question and Section 3 describes families of metrics on these manifolds with very small anti-self dual Weyl curvature. The existence theorem for the anti-self dual metrics in question is in Section 4; it is stated formally as Theorem 4.3. This theorem is a direct consequence of two preliminary propositions in Section 4. The first proposition is essentially an application of theorems from Kovalev-Singer [KS]: This first proposition reduces the existence question to a question about the vanishing or not of a vector in a finite dimensional obstruction vector space. The second proposition (Proposition 4.2) says when this obstruction vector can be guaranteed to be zero. Sections 5 and 6 set the stage for the proof of Proposition 4.2 which is in Section 7. Section 8 talks about the Gromov-Hausdorff limits of the metrics that are described in Theorem 4.3. There is an appendix to this article whose first section outlines the proofs of two propositions from Section 5. The next three sections supply the background for the proofs of the two propositions. The fifth section of the appendix proves these two propositions.

SECTION 2: The 0-framed knot surgery manifolds

SECTION 3: Metrics with very small anti-self dual Weyl curvature.

SECTION 4: Deformation to metrics with zero anti-self dual Weyl curvature.

SECTION 5: The obstruction space.

SECTION 6: The obstruction vector.

SECTION 7:  Proof of Proposition 4.2.

SECTION 8:  Gromov-Hausdorff limits.

APPENDIX:  The obstruction vector space.

This introduction hereby ends with my sincere thanks to Selman Akbulut and Ron Fintushel for sharing their thoughts about the constructions that are described below.

**2. The 0-framed knot surgery manifolds**

The first part of this section describes the 0-framed knot surgery construction. The subsequent parts of this section say something about the topology of the resulting manifolds.

**a) The construction**

The manifolds that are constructed in what follows are labeled by a knot in $S^3$. The knot will be denoted by K and the resulting manifold will be denoted by $X_K$. The construction of $X_K$ has three parts.

*Part 1*:  This first part of the construction reviews the Kummer description of the K3 manifold. To start, introduce by way of notation $\mathbb{T}$ to denote the 4-torus, written as the quotient of $\mathbb{R}^4$ by the lattice of points whose coordinates are integer multiples of $2\pi$. An involution to be denoted by $\iota$ acts on $\mathbb{T}$ by sending any given 4-tuple of coordinates $(t_1, t_2, t_3, t_4) \in \mathbb{R}^4/(2\pi\mathbb{Z}^4)$ to the point $(-t_1, -t_2, -t_3, -t_4)$. It has 16 fixed points, these being the points where the Euclidean coordinate functions $\{t_k\}_{k=1,2,3,4}$ are equal to 0 mod $\pi$. A neighborhood of any given singular point in the quotient space $\mathbb{T}/\iota$ is a cone on $\mathbb{RP}^3$. Since the boundary of the unit disk bundle in $T^*S^3$ is diffeomorphic to $\mathbb{RP}^3$, a smooth manifold is obtained by removing a cone neighborhood of each singular point in $\mathbb{T}/\iota$ and replacing it with a suitable radius disk bundle in $T^*S^2$. This manifold is the K3 manifold.

*Part 2*:  Fix $t_* \in (0, \frac{1}{4}\pi]$ and let $T \subset \mathbb{T}$ denote the 2-torus where the Euclidean coordinates $t_3$ and $t_4$ are both $t_*$ (mod $2\pi$). Note that the coordinate functions $(t_1, t_2)$ restrict to give $\mathbb{R}^2/(2\pi\mathbb{Z}^2)$ coordinates on T. The torus T has a product neighborhood in $\mathbb{T}$ where $|t_3 - t_*|^2 + |t_4 - t_*|^2 < \frac{1}{16} t_*^2$. Let N denote this neighborhood. Note in particular that N is disjoint from the fixed point set of $\iota$ and N is disjoint from $\iota(N)$. It is useful to use polar coordinates $(\rho, \theta_N)$ on the constant $t_1, t_2$ disks in N. The coordinate $\rho$ takes values in [0, $\frac{1}{4}t_*$] and with $\theta_N$ being $\mathbb{R}/2\pi\mathbb{Z}$ valued. They are defined so that $t_3 = t_* + \rho\cos\theta_N$ and $t_4 = t_* + \rho\sin\theta_N$.

*Part 3*: Supposing that K denotes a knot in $S^3$, let $N_K$ denote for the moment a suitably small radius solid torus neighorhood of K in $S^3$. The plan is to replace the copies of N and $\iota(N)$ in $\mathbb{T}$ with copies of $(S^3-N_K)\times S^1$ so that the involution $\iota$ extends from the $\mathbb{T}-(N\cup\iota(N))$ part of the resulting manifold to the whole manifold. Replacement surgery of this sort whereby a product neighborhood of a 2-torus in a 4-manifold is replaced by the product of the complement of solid torus in $S^3$ with the circle is called *knot surgery* (see [FS]). Note however that the version of knot surgery that is defined below does not use the same boundary identification as the one used in [FS]. Knot surgery with the identifications used here will be called *0-framed knot surgery*. The result of the forthcoming 0-framed knot surgery on $\mathbb{T}_O$ is denoted by $\mathbb{T}_K$.

To describe the 0-framed knot surgery operation, identify N as before with the product of the radius $\frac{1}{4}t_*$ disk about the origin in $\mathbb{R}^2$ and the 2 dimensional torus T. Use radial coordinates $(\rho,\theta_N)$ for the disk in $\mathbb{R}^2$ and coordinates $(t_1, t_2)$ for the $\mathbb{R}/(2\pi\mathbb{Z})$ coordinates on T.

A neighborhood in $S^3$ of a given knot K can be chosen so that the complement of K in this neighborhood has coordinates $(s, \tau_1, \tau_2)$ of the following sort: What is denoted by s is a proper function on $S^3-K$ with values in $[-2, \infty)$ with no critical values in $[-1, \infty)$. The knot K is approached in $S^3-K$ as $s \to \infty$. In particular, each constant $s \in [-1, \infty)$ surface is the boundary of solid torus neighborhoods of K. Meanwhile, $\tau_1$ and $\tau_2$ are $\mathbb{R}/(2\pi\mathbb{Z})$ valued. They are chosen so that each $\tau_1$ = constant annulus in the $s \geq -1$ part of $S^3-N_K$ (which is parametrized by the pair $(s, \tau_2)$) is the $s \geq -1$ part of some Seifert surface for K in $S^3-K$. Since the product of $S^3-K$ with $S^1$ is used in the upcoming surgery, a coordinate for the $S^1$ factor is required. This $\mathbb{R}/(2\pi\mathbb{Z})$-valued coordinate is denoted by $\theta$.

Fix a $R > 100(1 + |\ln\frac{1}{4}t_*|)$ and remove the $\rho < e^{-R}$ part of N; then replace it with the $s < R + \ln(\frac{1}{4}t_*)$ part of $(S^3-K)\times S^1$ using the identifications:

$$(t_1, t_2) = (\tau_1, \tau_2) \ and \ \rho = e^{s-R} \ and \ \theta_N = \theta.$$

(2.1)

This identification defines the 0-framed knot surgery operation. The manifold $\mathbb{T}_K$ is obtained from $\mathbb{T}$ by doing two instances of this knot surgery operation. The first operation replaces N from $\mathbb{T}$ with $(S^3-N_K)\times S^1$ using (2.1). Here and henceforth, $N_K$ denotes the specific solid torus neighborhood of K in $S^3$ whose boundary is the locus where $s = R + \ln(\frac{1}{4}t_*)$. The second instance of 0-framed knot surgery replaces $\iota(N)$ with second copy of $(S^3-N_K)\times S^1$ using the $\iota(N)$ version of (2.1). The replacement of $\iota(N)$ in the second surgery can and should be made so that the involution $\iota$ on $\mathbb{T}-(N\cup\iota(N))$ extends to $\mathbb{T}_K$ so as to interchange the two copies of $(S^3-N_K)\times S^1$.

*Part 4*: Since $\iota$ interchanges the two copies of $(S^3-N_K)\times S^1$, its fixed points on $\mathbb{T}_K$ lie in the complement of these two subspaces. Meanwhile, the preceding constructions identify the complement of the two copies of $(S^3-N_K)\times S^1$ in $\mathbb{T}_K$ with the complement of $N\cup\iota(N)$ in $\mathbb{T}$; and the involution $\iota$ acts on this $\mathbb{T}-(N\cup\iota(N))$ part of $\mathbb{T}_K$ as if it were part of $\mathbb{T}$. It follows as a consequence that $\iota$ has precisely sixteen fixed points on $\mathbb{T}_K$ which are the points in the $\mathbb{T}-(N\cup\iota(N))$ part of $\mathbb{T}_K$ and that are the fixed points of $\iota$'s action on $\mathbb{T}$. By way of a reminder, these are the points where the functions $\{t_k\}_{k=1,2,3,4}$ are 0 mod $\pi$.

As noted in Part 1, each of the resulting sixteen singular points in $(\mathbb{T}-(N\cup\iota(N)))/\iota$ has a neighborhood modelled by a cone on $\mathbb{RP}^3$ neighborhood. Therefore, these singular points can be resolved as in Part 1 by replacing suitable $\mathbb{RP}^3$ cone neighborhoods by copies of a constant radius disk subbundle in $T^*S^2$. The manifold that results from $\mathbb{T}_K/\iota$ by replacing these sixteen $\mathbb{RP}^3$ cones by disk subbundles in $T^*S^2$ is the manifold $X_K$.

**b) The topology of $X_K$**

The following propostion summarized what can be said about the topology of $X_K$.

**Proposition 2.1**: *Fix a knot in $S^3$ to be denoted by K and let $X_K$ denote the 4-manifold that is constructed from K using the rules from the previous subsection.*
- *The manifold $X_K$ is homeomorphic to $(\#_3 \mathbb{CP}^2)\#(\#_{19}\overline{\mathbb{CP}}^2)$.*
- *The Seiberg-Witten invariants of $X_K$ are zero.*

By way of a parethetical remark, Bauer and Furuta [BF], [Ba] have introduced what they call a stable cohomotopy refinement of the Seiberg-Witten invariants; but it is likely that Bauer-Furuta invariants are also trivial for $X_K$.

*Proof of Proposition 2.1*: Let O denote an unknotted circle in $S^3$, the *unknot*. As explained to the author by Ron Fintushel and by Selman Akbulut, the manifold $X_{K=O}$ is diffeomorphic to $(\#_3 \mathbb{CP}^2)\#(\#_{19}\overline{\mathbb{CP}}^2)$. The Seiberg-Witten invariants of the latter manifold vanish because it has a metric with positive scalar curvature.

Now consider the case when K is not the unknot. The 0-framed knot surgery operation can be viewed as a 2-step procedure, the first step producing $X_O$ according to the rules in Section 2a and the second step replacing the $(S^3-N_O)\times S^1$ part of $X_O$ with $(S^3-N_K)\times S^1$ by identifying the respective $s\geq 0$ parts of $S^3-N_O$ and $S^3-N_K$ via their respective $(s,\tau_1,\tau_2)$ coordinate charts. Since $S^3-N_O$ and $S^3-N_K$ have isomorphic homology and cohomology groups, the respective homology groups of $X_K$ and $X_O$ will be isomorphic via an isomorphism that intertwines their cup product operations. Thus, if it is the case that $X_K$ is simply connected, then it follows from Freedman's theorem [Fr] that $X_K$ is homeomorphic to $(\#_3 \mathbb{CP}^2) \# (\#_{19} \overline{\mathbb{CP}}^2)$.

The fact that $X_K$ is simply connected follows via the Seifert-van Kampen theorem from the two observations that follow directly. To state the first observation, let $Z_K$ denote the part of $X_K$ that comes by way of $\mathbb{T}_K/\iota$ from the $(S^3 - N_K) \times S^1$ part of $\mathbb{T}_K$. The first observation is that $X_K - Z_K$ is simply connected. This is so because $X_K - Z_K$ is $X_O - Z_O$; the latter is the complement in the K3 manifold of a torus fiber when K3 is viewed as an elliptic fibration over $\mathbb{CP}^1$; and the complement of such a fiber in the K3 manifold is known to be simply connected (see, e.g. [FS].) To state the second observation, let $\vartheta: T^2 \to S^3 - N_K$ denote the inclusion of the boundary torus. The second observation is that the group $\pi_1(S^3 - N_K)$ is the least normal subgroup in $\pi_1(S^3 - N_K)$ containing $\vartheta_*(\pi_1(T^2))$. This is to say that any element in $\pi_1(S^3 - N_K)$ can be written as a product of elements that have the form $k \, i \, k^{-1}$ with $i$ being in $\vartheta_*(\pi_1(T^2))$ and with k being in $\pi_1(S^3 - N_K)$. This second observation follows from the Seifert-van Kampen theorem also because $S^3 = (S^3 - K) \cup N_K$ and because $\pi_1(S^3) = 1$.

The proof that the Seiberg-Witten invariants of $X_K$ are all zero invokes the gluing theorem for these invariants in [T] with the fact that the Seiberg-Witten invariants of $X_O$ are zero. A lemma about the second homology of $Z_K$ and $X_K$ is needed for this. The statement of this lemma introduces by way of notation $X_+$ to denote $X_K - Z_K$. This notation does not indicate the knot K because $X_K - Z_K$ has a canonical identifications with $X_O - Z_O$ for each knot K. These canonical identifications are implicit in what follows.

**Lemma 2.2**: *The second homology of $Z_K$, $X_+$ and $X_K$ have the following properties:*
- $H_2(Z_K; \mathbb{Z}) \approx \mathbb{Z}$. *Moreover, the generator comes from $H_2(\partial Z_K; \mathbb{Z})$ via the homomorphism that is induced by the inclusion of $\partial Z_K$ into $Z_K$.*
- *The inclusion of $Z_K$ into $X_K$ induces a monomorphism from $H_2(Z_K; \mathbb{Z})$ to $H_2(X_K; \mathbb{Z})$.*
- *The inclusion of $X_+$ into $X_K$ induces a homomorphism from $H_2(X_+; \mathbb{Z})$ to $H_2(X_K; \mathbb{Z})$ whose kernel is independent of* K.

This lemma is proved momentarily.

The three steps that follow use this lemma to prove that the Seiberg-Witten invariants of $X_K$ are zero. There is a fourth step that sketches a possible proof of an assertion to the effect that $X_K$ has trivial the Bauer-Furuta invariant.

<u>Step 1</u>: Let $\mathbb{Z} H_2(X_K; \mathbb{Z})$ denote the free $\mathbb{Z}$-module generated by the classes in $H_2(X; \mathbb{Z})$. This is to say that an element in $\mathbb{Z} H_2(X_K; \mathbb{Z})$ is a finite sum of the form $\sum_{\mu \in H_2} a_\mu \mu$ with the coefficients $\{a_\mu\}_{\mu \in H_2}$ being integers and with only finitely many of them being non-zero. It proves useful in what follows to define a multiplicative structure on $\mathbb{Z}H_2(X_K; \mathbb{Z})$ by the rule whereby

$$(\sum_{\mu \in H_2} a_\mu \mu) \cdot (\sum_{\mu' \in H_2} b_{\mu'} \mu') = \sum_{\mu'' \in H_2} a_\mu b_{\mu'} \delta_{\mu''-\mu-\mu'} \mu''$$

(2.2)

with $\delta_{(\cdot)}$ being 1 on the class 0 and zero otherwise.

The Seiberg-Witten invariants for $X_K$ can viewed as defining an element in $\mathbb{Z}H_2(X_K; \mathbb{Z})$. This element has the form $\sum_{\mu \in H_2} a_\mu \mu$ with the coefficient $a_\mu$ defined as follows: The Seiberg-Witten invariants of $X_K$ as originally defined associate an integer to each $\text{Spin}_\mathbb{C}$ structure on $X_K$ (see [W] or [M].) Since there is no 2-torsion in the second homology of $X_K$, any given $\text{Spin}_\mathbb{C}$ structure on $X_K$ is determined by the first Chern class of a certain associated, complex line bundle. Thus, Poincaré duality can be invoked to label $\text{Spin}_\mathbb{C}$ structures by elements in $H_2(X_K; \mathbb{Z})$. Supposing that $\mu \in H_2(X_K; \mathbb{Z})$, the coefficient $a_\mu$ is the Seiberg-Witten invariant of the $\text{Spin}_\mathbb{C}$ structure labeled by $\mu$.

<u>Step 2</u>: It follows from the first two bullets of Lemma 2.2 that there is a class in $H^2(X_K; \mathbb{Z})$ with non-zero restriction to $H^2(\partial X_+; \mathbb{Z})$. Given such a class, the constructions in [T] can be used to define a Seiberg-Witten invariant for $X_+$ that takes values in the free $\mathbb{Z}$ module that is generated by the classes in $H_2(X_+; \mathbb{Z})$. This $\mathbb{Z}$-module is denoted in what follows by $\mathbb{Z}H_2(X_+; \mathbb{Z})$ and the Seiberg-Witten element is denoted by $SW_+$.

By the same token, there is also a Seiberg-Witten element in $\mathbb{Z}H_2(Z_K; \mathbb{Z})$. Because $H_2(Z_K; \mathbb{Z}) \approx \mathbb{Z}$, a given elements in the $\mathbb{Z}$-module $\mathbb{Z}H_2(Z_K, \mathbb{Z})$ can be depicted as a finite Laurent polynomial in a single variable (this denoted by $t_K$) by using $t_K^n$ for $n \in \mathbb{Z}$ to denote the n'th power of a chosen generator. As explained in [MT], the Seiberg-Witten element of $Z_K$ is $\Delta_K(t_K^2)$ with $\Delta_K(\cdot)$ denoting the symmetrized Alexander polynomial of the knot K.

<u>Step 3</u>: Let $i_K: Z_K \to X_K$ and $j_K: X_+ \to X_K$ denote the respective inclusion maps. The corresponding homorphisms from $H_2(Z_K; \mathbb{Z})$ and $H_2(X_+; \mathbb{Z})$ to $H_2(X_K; \mathbb{Z})$ induce in turn respective homomorphisms from $\mathbb{Z}H_2(Z_K; \mathbb{Z})$ and $\mathbb{Z}H_2(X_+; \mathbb{Z})$ into $\mathbb{Z}H_2(X_K; \mathbb{Z})$. The latter are denoted by $i_{K*}$ and $j_{K*}$. The main theorem in [T] asserts that Seiberg-Witten element in $\mathbb{Z}H_2(X_K; \mathbb{Z})$ can be written as

$$SW_K = j_{K*}(SW_+) \cdot i_{K*}(SW_{Z_K}) \ .$$

(2.3)

with the multiplication rule given by (2.3).

To exploit this identity, consider first the case when K is the unknot O. In this case, the left hand side of (2.3) is zero and thus so is the right hand side. To see the implications, let [0] denote for the moment the trivial class in $H_2(Z_O; \mathbb{Z})$. Since the Alexander polynomial of the unknot is 1, it follows that $SW_{Z_O} = [0]$. Since Lemma 2.2 says that $i_{O*}$ is injective, it follows from the multiplication rules in (2.2) that the right hand side of (2.3) is zero if and only if $j_{O*}(SW_+) = 0$. Now let $A_+$ denote the kernel of the homomorphism induced by $j_O$ from $H_2(X_+; \mathbb{Z})$ to $H_2(X_O; \mathbb{Z})$. The vanishing of $j_{O*}(SW_+)$ requires that $SW_+$ be in the multiplicative ideal generated by elements of the form $[A] - [0]$ with [0] denoting the trivial class and with [A] denoting a class $A_+$. (Multiplication is again defined using (2.2).)

Now let K be any given knot. Since the kernel of the homomorphism from $H_2(X_+; \mathbb{Z})$ to $H_2(X_K; \mathbb{Z})$ is $A_+$, it follows from what was said in the preceding paragraph that $j_{K*}(SW_+) = 0$. Thus, the right hand side of (2.3) is zero for K; and so $SW_K = 0$ also.

<u>Step 4</u>: As remarked above, it is likely that the stable homotopy refinement of the Seiberg-Witten invariants that is defined by Bauer and Furuta [BF] is also trivial. What follows is a sketch of a possible proof of this assertion: Since $X_O$ is the connect sum of a manifold with $\mathbb{CP}^2$ and since the signature of $\mathbb{CP}^2$ is 1 mod(4), Proposition 4.5 in [Ba] can be invoked to see that the Bauer-Furuta invariant is trivial for $X_O$. Meanwhile, there is likely a gluing theorem for the Bauer-Furuta invariant along an embedded three dimensional torus if the assumptions for the gluing theorem in [T] are met (which is the case by virtue of Lemma 2.2.) This new gluing theorem should assert that the Bauer-Furuta invariant of $X_K$ is the product in some sense of Bauer-Furuta invariants for $X_+$ and $Z_K$, thus the analog of (2.3). An argument much like that given above in Steps 1-3 using this hypothetical gluing theorem would presumably prove that the Bauer-Furuta invariants for $X_K$ are also zero since they are zero for $X_O$. Granted the analysis and a priori estimates in [T], the proof of this hypothetical stable homotopy gluing theorem along three dimensional tori should not be much different from the proof of Theorem 1.1 in [Ba] which asserts a gluing theorem along an embedded three dimensional sphere. The key point in both gluing theorems is that the moduli space of solutions of the relevant version of the Seiberg-Witten equations on the sphere or torus has but a single point. (In the case of the sphere, the equations are the standard three dimensional analog of the Seiberg-Witten equations with the metric being the constant curvature round metric. In the case of the torus and under the assumptions in [T], the equations are a certain deformation of the standard three dimensional Seiberg-Witten equations that admit a unique solution when the metric on the torus is flat.)

***Proof of Lemma 2.2***: The proof has six steps. Step 1 computes $H_2(Z_K; \mathbb{Z})$ and the remaining steps prove that the $i_K: Z_K \to X_K$ and $j_K: X_+ \to X_K$ induce homomorphism on the second homology with the asserted properties.

<u>Step 1</u>: The manifold $Z_K$ is $(S^3 - N_K) \times S^1$. The assertion that $H_2((S^3 - N_K) \times S^1; \mathbb{Z})$ is $\mathbb{Z}$ follows from the fact that $H_1(S^3 - N_K; \mathbb{Z}) \approx \mathbb{Z}$ which follows in turn from the Mayer-Vietoris sequence for the decomposition of $S^3$ as the union of $S^3 - K$ and $N_K$.

Keeping in mind that $N_K$ is a solid torus neighborhood of K, this same Mayer-Vietoris sequence implies that a generator for $H_2((S^3 - N_K) \times S^1; \mathbb{Z})$ can be taken to be the class of the product $\gamma \times S^1 \subset (S^3 - N_K) \times S^1$ with $\gamma \in S^3 - N_K$ being a loop that has linking number 1 with the knot K. In particular, a suitable loop $\gamma$ is given in the part of $S^3 - K$ parametrized by the coordinates $(s, \tau_1, \tau_2)$ by taking $s$ = constant and $\tau_2$ = constant. This loop has linking number 1 with K because it has intersection number 1 (or -1) with Seifert surfaces for the knot K. Keep in mind in this regard that the annuli in $S^3 - K$ where $\tau_1$ is constant are in the $s \geq 0$ part of some Seifert surfaces for K.

<u>Step 2</u>: It follows from what is said in Step 1 that the generator of $H_2(Z_K; \mathbb{Z})$ is in the image of the homomorphism from $H_2(\partial Z_K; \mathbb{Z})$ to $H_2(Z_K; \mathbb{Z})$. Since $\partial Z_K$ is geometrically the boundary of $X_+$, the homomorphism $i'_K: H_2(\partial X_+; \mathbb{Z}) \to H_2(Z_K; \mathbb{Z})$ induced from the inclusion map is an monomorphism. The kernel of this monomorphism $i'_K$ is 2-dimensional because $H_2(\partial X_+; \mathbb{Z}) \approx \mathbb{Z}^3$. As explained next, this kernel is independent of K. To see why this is so, note first that the boundary of $X_+$ is the boundary of N and so it can be parametrized by the coordinates $(t_1, t_2, \theta_N)$ for the boundary of N. The kernel of $i'_K$ in for any knot K can be generated by the fundamental classes of the $\theta_N = 0$ torus and the $t_{1\text{out}} = t_*$ torus. Indeed, the $\theta_N = 0$ torus appears when written using the $(\tau_1, \tau_2, \theta)$ coordinates of (2.2) as the torus where $\theta = 0$. This torus in the boundary of $(S^3 - N_K) \times S^1$ is the boundary of $(S^3 - N_K) \times \{0\}$. Meanwhile, the torus where $t_{1\text{out}} = t_*$ appears when written using the $(\tau_1, \tau_2, \theta)$ coordinates as the torus where $\tau_1 = t_*$. This is the boundary of $\Sigma \times S^1$ with $\Sigma$ being the $S^3 - N_K$ part of some Seifert surface for the knot K.

<u>Step 3</u>: Suppose for the moment that the homomorphism $i_{K*}$ from $H_2(Z_K; \mathbb{Z})$ to $H_2(X_K; \mathbb{Z})$ is an isomorphism (this is the assertion of the second bullet of Lemma 2.2). The third bullet of the lemma is shown in this step to follow from this assumption and from what is said in Step 2 about the kernel of $i'_K$ being independent of K.

The argument for the third bullet invokes the Mayer-Vietoris sequence for the decomposition of $X_K$ as the union of the closures of $X_+$ and $Z_K$:

$$\cdots \to 0 \to H_2(\partial X_+; \mathbb{Z}) \to H_2(Z_K; \mathbb{Z}) \oplus H_2(X_+; \mathbb{Z}) \to H_2(X_K; \mathbb{Z}) \to \cdots .$$
(2.4)

The 0 on the far left of this sequence is $H_3(X_K; \mathbb{Z})$ and the arrows from $H_2(\partial X_+; \mathbb{Z})$ and from $H_2(Z_K; \mathbb{Z})$ and $H_2(X_+; \mathbb{Z})$ are induced by the relevant inclusion maps. In particular, the respective arrows from $H_2(Z_+; \mathbb{Z})$ and $H_2(X_+; \mathbb{Z})$ are $i_{K*}$ and $j_{K*}$. (Technically, the Mayer-Vietoris sequence holds for a decomposition of $X_K$ as the union of two sets whose interiors cover $X_K$, but since $X_+$ and $Z_K$ are smooth manifolds with boundary in $X_K$, this technicality has no bearing on the subsequent discussion.)

Let $j_+$ denote the homomorphism from $H_2(\partial X_+; \mathbb{Z})$ to $H_2(X_+; \mathbb{Z})$ that is induced by the inclusion map. Supposing that the homomorphism $i_K$ is injective (as asserted by Lemma 2.2's second bullet), then the composition of first $\iota_K$ and then $i_K$ isn't zer because $\iota_K$ is an monomorphism (according to Step 1.) Let $A_+$ denote the kernel of $\iota_K$ which is independent of K (according to Step 2). As a consequence of this and the fact that (2.4) is exact, the homomorphism $j_+$ must be injective and the kernel of $j_+$ must be $j_+(A_+)$.

Step 4: The remaining steps proves the assertion of the lemma's second bullet to the effect that the homomorphism $i_{K*}: H_2(Z_K; \mathbb{Z}) \to H_2(X_K; \mathbb{Z})$ is injective. The proof is along the following lines: As noted in Step 1, the identification of $Z_K$ with $(S^3 - N_K) \times S^1$ gives a generated for $H_2(Z_K; \mathbb{Z})$ of the form $[\gamma \times S^1]$ with $\gamma$ being some s=constant and $\tau_2$ = constant loop in the $s \geq 0$ part of $S^3 - N_K$. This homology class is non-zero in $H_2(X_K; \mathbb{Z})$ if it has positive intersection number with a surface in $X_K$. The plan for what follows is to exhibit such a surface where the intersection number is +1 (or -1). By way of a parenthetical remark, the surface in question in the case when K is the unknot will be a 2-sphere with trivial normal bundle. The latter 2-sphere can be used as a starting point to decompose $X_O$ as the connect sum of three $\mathbb{CP}^2$s and nineteen $\overline{\mathbb{CP}}^2$s.

Step 5: The loop on $\partial X_+$ where $\theta_N = 0$ and $t_1 = t_*$ is the boundary of an embedded surface in $Z_K$ because this loop appears in the $(\tau_1, \tau_2, \theta)$ coordinates as the circle where $\theta = 0$ and $\tau_1 = t_*$; and this is the boundary of a surface of the form $\Sigma \times \{0\}$ in $(S^3 - N_K) \times S^1$ with $\Sigma$ being the $S^3 - N_K$ part of a Seifert surface for the knot K that contains the $s \geq 0$ and $\tau_1 = t_*$ annulus. Denote this surface in $Z_K$ by $\Sigma_+$. A corresponding surface is needed in $\iota(Z_K)$. This surface in $\iota(Z_K)$ is not $\iota(\Sigma_+)$ but rather a surface of the form $\iota(\Sigma' \times \{\pi\})$ with $\Sigma'$ being the $S^3 - N_K$ part of a Seifert surface for K that contains the $s \geq 0$ and $\tau_1 = -t_*$ annulus. Use $\Sigma_-$ to denote $\iota(\Sigma' \times \{0\}.)$

Since the complement in $\mathbb{T}_K$ of $Z_K \cup \iota(Z_K)$ is the complement in $\mathbb{T}$ of $N \cup \iota(N)$, the boundary loops of both $\Sigma_+ \subset Z_K$ and $\Sigma_-$ in $\iota(Z_K)$ can be written using the $(t_1, t_2, t_3, t_4)$ coordinates for $\mathbb{T}$:

- *The boundary of $\Sigma_+$ is the loop in $\mathbb{T}$ where $t_1 = t_*$, $t_3 = (1 + \frac{1}{4})t_*$ and $t_4 = t_*$.*
- *The boundary of $\Sigma_-$ is the loop in $\mathbb{T}$ where $t_1 = t_*$, $t_3 = (-1 + \frac{1}{4})t_*$ and $t_4 = -t_*$.*

(2.5)

As explained next, these two loops are the boundary of an annulus in the $t_1 = t_*$ locus of $\mathbb{T}-(N \cup \iota(N))$. To depict such a annulus, let T´ for the moment denote the torus with coordinates $(t_3, t_4)$. Let D denote the disk in T´ of radius $\frac{1}{4} t_*$ centered on $(t_*, t_*)$ and let $\iota(D)$ denote the disk in T´ of radius $\frac{1}{4} t_*$ centered on $(-t_*, -t_*)$. Choose a path in T´$-(D \cup \iota(D))$ that starts at the point $(t_3 = (1 + \frac{1}{4})t_*, t_4 = t_*)$ on $\partial D$ and ends at the point in $\partial \iota(D))$ where $(t_3 = (-1 + \frac{1}{4})t_*, t_4 = t_*)$. Parametrize this path by the interval $[0, 1]$ and denote it by $\eta$. The image of the map from $[0, 1] \times S^1$ into $\mathbb{T}$ that sends a pair $(\sigma, \tau)$ with $\sigma \in [0, 1]$ and $\tau \in \mathbb{R}/2\pi\mathbb{Z}$ to the point in $\mathbb{T}$ where $t_1 = t_{1*}$, $t_2 = \tau$ and $(t_3, t_4) = \eta(\sigma)$ is an embedded annulus in the $t_1 = t_*$ locus of $\mathbb{T}-(N \cup \iota(N))$ whose boundary curves are depicted in (2.5).

Denote the annulus from the preceding paragraph by $\Sigma_0$ and let A denote the surface $\Sigma$ given by the union $\Sigma_- \cup \Sigma_0 \cup \Sigma_+$. This is a closed, embedded surface. Note by the way that in the case K = O, both $\Sigma_-$ and $\Sigma_+$ can be assumed to be disks; and if they are disks, then A is an embedded 2-sphere.

<u>Step</u> <u>6</u>: The surface A is disjoint from $\iota(A)$ and thus it is mapped homeomorphically to its image in $\mathbb{T}_K/\iota$. Because this image in $\mathbb{T}_K/\iota$ has distance no less than $\frac{1}{2} t_*$ from the singular points, it also defines an embedded surface in $X_K$. The surface A also has a single transversal intersection with the torus $\gamma \times S^1$ in $Z_K$. It follows as a consequence that the homology classes in $X_K$ of both $\gamma \times S^1$ (from $Z_K$) and A generate a $\mathbb{Z}^2$ submodule in $H_2(X_K; \mathbb{Z})$.

## 3. Metrics with very small anti-self dual Weyl curvature

Suppose in what follows that K is a hyperbolic knot in $S^3$, which is to say that $S^3-K$ has a complete, finite volume hyperbolic metric. This section describes Riemannian metrics on $X_K$ with very small self-dual Weyl curvature.

### a) The curvature decomposition in dimension 4.

Let X denote for the moment a four dimensional manifold with a given Riemannian metric. The metric's Hodge star maps 2-forms to 2-forms as an involution with square 1. The ±1 eigenspaces of the Hodge star are both three dimensional; the +1 eigenspace is the bundle of self-dual 2-forms and the -1 eigenspace is the bundle of anti-self dual 2-forms. These are denoted respectively by $\Lambda^+$ and $\Lambda^-$. The Riemannian metric in any dimension defines a symmetric endomorphism of the bundle of 2-forms. In the case at hand where the dimension is 4, this endomorphism can be written in $2 \times 2$ block form with respect to the decomposition $\wedge^2 T^*X = \Lambda^+ \oplus \Lambda^-$ as follows:

$$\begin{pmatrix} \mathcal{W}_+ - \tfrac{1}{12} R \mathbb{I} & \mathcal{B} \\ \mathcal{B}^T & \mathcal{W}_- - \tfrac{1}{12} R \mathbb{I} \end{pmatrix}$$

(3.1)

with $\mathbb{I}$ denoting the identity $3 \times 3$ endomorphism, R denoting the scalar curvature, $\mathcal{B}$ encoding the traceless Ricci tensor and $\mathcal{W}_+$ and $\mathcal{W}_-$ denoting the respective self-dual and anti-self dual parts of the Weyl curvature. A Riemannian metric is said to have anti-self dual Weyl curvature when $\mathcal{W}_+ = 0$.

A Riemannian metric on a manifold is said to be locally conformally flat when it has a cover by coordinate chart such that the metric in each chart is proportional at every point to the Euclidean metric. The vanishing of the Weyl curvature tensor is a necessary and sufficient condition for this. In dimension 4, this means that both $\mathcal{W}_+$ and $\mathcal{W}_-$ are everywhere zero. By way of an example, suppose that Y is a three dimensional manifold and that $\mathfrak{g}$ on Y is a metric with constant sectional curvature equal to -1, thus a hyperbolic metric. Let $\theta$ denote the Euclidean $\mathbb{R}/2\pi\mathbb{Z}$ coordinate for $S^1$. Then the metric $\mathfrak{g} + d\theta^2$ on $Y \times S^1$ is locally conformally flat. In this case, contraction with the vector field $\tfrac{\partial}{\partial \theta}$ identifies both $\Lambda^+$ and $\Lambda^-$ with $T^*Y$; and this identification intertwines the endomorphism $\mathcal{B}$ in (3.1) with 4 times the identity endomorphism. (This is $\mathcal{B} = +4\,\mathbb{I}$, not $-4\,\mathbb{I}$.)

If $\mathfrak{g}$ is a locally conformally flat metric on a given manifold, then so is $\phi \mathfrak{g}$ with $\phi$ being any strictly positive function. This freedom to multiply a given locally conformally flat metric by a function is the reason for introducing the notion of a *locally conformally flat structure*, which is an equivalence class of locally conformally flat metrics with the equivalence relation identifying metrics when one is the product of a positive function times the other. Any locally conformally flat structure is defined by the following data:

- *A locally finite open cover, $\mathfrak{U}$, of the manifold in question.*
- *A locally conformally flat metric on each set from the cover.*
- *Supposing that $U \in \mathfrak{U}$ and $V \in \mathfrak{U}$, let $\mathfrak{g}_U$ and $\mathfrak{g}_V$ denote their associated metrics. If $U \cap V \neq \emptyset$, then $\mathfrak{g}_U = \phi_{UV} \mathfrak{g}_V$ on $U \cap V$ with $\phi_{UV}$ being a positive function on $U \cap V$.*

(3.2)

A partition of unity can be used to construct an honest locally conformally flat metric from the collection $\{\mathfrak{g}_U\}_{U \in \mathfrak{U}}$. This metric on any given set $U \in \mathfrak{U}$ will be the product of a function times the metric $\mathfrak{g}_U$.

The term *conformal structure* is used in what follows to denote a data set (indicated by brackets, as in [$\mathfrak{g}$]) that obeys the conditions in (3.2) without the conformally flat requirement for the metrics on the open sets of the given cover. Thus, [g] signifies an open cover of the manifold and an assigned metric to each open set in the cover subject to the constraint in the third bullet of (3.2).

**b) An $\iota$-invariant, locally conformally flat metric on $\mathbb{T}_K$**

Suppose that Y is a complete, finite volume hyperbolic 3-manifold with a single cusp end. This hyperbolic metric on Y is denoted here by g. As it turns out (see e.g. [G]), a cusp end always has coordinates $(s, \tau_1, \tau_2)$ for $s \in [-1, \infty)$ and $(\tau_1, \tau_2) \in \mathbb{R}^2/(2\pi\mathbb{Z}^2)$ with which the metric can be written as

$$g = e^{-2s} \mathfrak{m} + ds^2 \tag{3.3}$$

where $\mathfrak{m} = \sum_{a,b \in \{1,2\}} \mathfrak{m}_{ab} d\tau^a \oplus d\tau^b$ is an s-independent, constant flat metric on the $(\tau_1, \tau_2)$-torus. In particular, the coefficients $\{\mathfrak{m}_{ab}\}_{a,b \in \{1,2\}}$ are fixed numbers, independent of s and of $(\tau_1, \tau_2)$. Now, fix $R \geq 0$. It follows from (3.3) that the metric $g_K = g + d\theta^2$ on $Y \times S^1$ is therefore conformal where $s \geq 0$ to the metric

$$e^{-2R} \mathfrak{m} + e^{2(s-R)}(ds^2 + d\theta^2) . \tag{3.4}$$

Introducing $\rho$ to denote $e^{s-R}$, this metric can be written as

$$e^{-2R} \mathfrak{m} + d\rho^2 + \rho^2 d\theta^2, \tag{3.5}$$

which is manifestly flat because it is the product of the flat metric on the $(\tau_1, \tau_2)$ torus and the Euclidean metric from $\mathbb{R}^2$ (written using polar coordinates).

Let $\mathbb{T}$ again denote the four dimensional torus. As before, it is identified with $\mathbb{R}^4/(2\pi\mathbb{Z}^4)$ using the $\mathbb{R}/(2\pi\mathbb{Z})$ valued coordinates $(t_1, t_2, t_3, t_4)$. The torus $\mathbb{T}$ has the metric

$$g_\mathbb{T} = e^{-2R} \sum_{a,b=1,2} \mathfrak{m}_{ab} dt_a dt_b + dt_3^2 + dt_4^2 . \tag{3.6}$$

Let N again denote the part of $\mathbb{T}$ where $(t_3 - t_*)^2 + (t_4 - t_*)^2 \leq \frac{1}{16} t_*^2$. Writing $t_3 - t_*$ as $\rho \cos \theta$ and $t_4 - t_*$ as $\rho \sin \theta$, and making the identifications $t_1 = \tau_1$ and $t_2 = \tau_2$, the metric in (3.6) is observedly the same metric as that in (3.5).

With the preceding understood, let $g_K$ denote the product metric on $(S^3 - N_K) \times S^1$ with the metric on the $S^3 - N_K$ factor being the restriction from $S^3 - K$ of the finite volume, sectional curvature -1 metric, and the metric on the $S^1$ factor being $d\theta^2$. Since (3.5) and (3.6) agree when $(t_1, t_2, t_3, t_4)$ are written as $(\rho \cos \theta, \rho \sin \theta, \tau_1, \tau_2)$, it follows from (2.1) that the metric from (3.6) on the $\mathbb{T} - (N \cup \iota(N))$ part of $\mathbb{T}_K$, the metric $g_K$ on $(S^3-N_K)\times S^1$ part of $\mathbb{T}_K$, and the metric $\iota^* g_K$ on the $\iota((S^3-N_K)\times S^1)$ part of $\mathbb{T}_K$ define a locally conformally flat structure on $\mathbb{T}_K$. This locally conformally flat structure is $\iota$-invariant because the metric in (3.6) on $\mathbb{T}$ is $\iota$-invariant. It therefore defines a locally conformally flat structure on the complement of the sixteen singular points in $\mathbb{T}_K/\iota$.

The locally conformally flat structure on $\mathbb{T}_K$ that was just defined is denoted in what follows by $[\mathfrak{g}_R]$. This same notation is also used to denote the induced locally conformally flat structure on the complement of the singular points in $\mathbb{T}_K/\iota$.

### c) Anti-self dual metrics on T*S²

As explained in Part 4 of Section 2, the manifold $X_K$ is obtained from $\mathbb{T}_K/\iota$ by a surgery that replaces $\mathbb{RP}^3$ cone neighborhoods of each of the singular points in $\mathbb{T}_K/\iota$ by small radius disk subbundles in $T^*S^2$. Metrics on $X_K$ with very small self-dual Weyl curvature are constructed in the next subsection by grafting the locally conformally flat metric on the complement in $\mathbb{T}_K/\iota$ of the singular points to a suitable metric on $T^*S^2$ with anti-self dual Weyl curvature (and zero Ricci curvature). The metrics on $T^*S^2$ are the Eguchi-Hanson metrics [EH], and the grafting of these metrics to $\mathbb{T}/\iota$ was first described by Gibbons and Pope [GP] (see also [P] and [LS]). (Their purpose was to obtain metrics on the K3 surface with very small self-dual Weyl and Ricci curvature.) The upcoming Proposition 3.1 describes the salient properties of the Eguchi-Hanson metrics on $T^*S^2$.

Proposition 3.1 refers to a particular map from $T^*S^2$ to $\mathbb{R}^4/\{\pm 1\}$ that sends the complement of the zero section in $T^*S^2$ diffeomorphically onto the complement of the origin in $\mathbb{R}^4/\{\pm 1\}$. To present the desired map, it proves useful to view the 2-dimensional sphere $S^2$ as the complex projective line $\mathbb{CP}^1$; and then to view any given point in $\mathbb{CP}^1$ as an equivalence class of pairs $(u_1, u_2) \in \mathbb{C}^2 - 0$ with the equivalence being $(u_1, u_2) \sim (v_1, v_2)$ when $(v_1, v_2)$ is a non-zero complex multiple of $(u_1, u_2)$. The bundle $T^*S^2$ appears in this guise as the underlying real bundle of a complex line bundle over $\mathbb{CP}^1$ (denoted by $E^2$) that is defined to be the space of equivalence classes of triples $(u_1, u_2; \zeta) \in (\mathbb{C}^2-0) \times \mathbb{C}$ with the equivalence relation identifying $(u_1, u_2; \zeta)$ and $(\lambda u_1, \lambda u_2; \lambda^{-2}\zeta)$ when $\lambda \in \mathbb{C}-0$.

Define the map $\mathfrak{p}: E^2 \to \mathbb{C}^2/\{\pm 1\}$ by the rule

$$(u_1, u_2; \zeta) \to (\zeta^{1/2} u_1, \zeta^{1/2} u_2)$$

(3.7)

Let $\underline{0}$ denote the zero section of the bundle $E^2$. The map $\mathfrak{p}$ restricts to $E^2 - \underline{0}$ as a diffeomorphism onto $(\mathbb{C}^2-0)/\pm 1$. Now choose an isometric, constant almost complex structure on $\mathbb{R}^4$ so as to identify $\mathbb{R}^4$ with $\mathbb{C}^2$. The map in (3.7) identifies the complement of the zero section in $E^2$ with the complement of the origin in $\mathbb{R}^4/\{\pm 1\}$. (A constant isometric almost complex is a matrix J obeying $J^2 = -1$ and $J^T = -J$.) Proposition 3.1 assumes an isometric identification of $\mathbb{R}^4$ with $\mathbb{C}^2$ to view $\mathfrak{p}$ as a map from $E^2$ to $\mathbb{R}^4/\{\pm 1\}$.

Supposing that $r \in (0, \infty)$, Proposition 3.1 uses $\mathfrak{q}_r$ to denote the map $\mathbb{R}^4$ to $\mathbb{R}^4/\{\pm 1\}$ that is obtained by composing first the rescaling map from $\mathbb{R}^4$ to $\mathbb{R}^4$ sending x to $r^{-1}x$ and then the tautological projection map to $\mathbb{R}^4/\{\pm 1\}$. If $x \in \mathbb{R}^4$ is a given point, the proposition uses |x| to denote its Euclidean norm.

**Proposition 3.1**: *There exists a metric on $E^2$ to be denoted by $\hat{g}$ having vanishing anti-self dual Weyl curvature and vanishing Ricci curvature. The normalization is such that the zero section in $E^2$ has area $4\pi$. This metric has the following properties with regards to the map $\mathfrak{p}$ and any given $r > 0$ version of the map $\mathfrak{q}_r$: Let $\eth$ denote the Euclidean metric on $\mathbb{R}^4$ and let $g_r$ denote the metric on $\mathbb{R}^4-0$ that is obtained by multiplying the $(\mathfrak{p}^{-1} \circ \mathfrak{q}_r)$ pull-back of $\hat{g}$ by $r^2$. If $x \in \mathbb{R}^4 - 0$ and $r > 0$, then*
- $|g_r - \eth| \le \kappa r^4 |x|^{-4}$,
- $|\nabla g_r| \le \kappa r^4 |x|^{-5}$,
- $|\nabla \nabla g_r| \le \kappa r^4 |x|^{-6}$,

*with $\kappa$ being independent of the parameter r and the point x.*

*Proof of Proposition 3.1*: The metric $\hat{g}$ is the Eguchi-Hanson metric [EH] with a suitable choice of [EH]'s parameter a. The properties asserted by the four bullets follow from the formulas in [EH].

### d) The Fubini-Study metric on $\overline{\mathbb{CP}}^2$

The Fubini-Study metric on $\overline{\mathbb{CP}}^2$ has anti-self dual Weyl curvature; and its traceless Ricci curvature is zero. It also has positive scalar curvature. These properties were exploited by Floer [Fl] and subsequently by Donaldson/Friedman [DF] and others to construct metrics with anti-self dual Weyl curvature on 4-manifolds that are connect sums with many $\overline{\mathbb{CP}}^2$s. Floer observed that the connect sum of $\overline{\mathbb{CP}}^2$ with a manifold having a metric with anti-self dual Weyl curvature can be done so as to obtain a metric on the connnect sum manifold with very small anti-self dual Weyl curvature. Donaldson/Friedman made a second crucial observation about which more is said momentarily. The next three paragraphs describes the version of the connect sum construction that is used here; Proposition 3.2 summarizes the story.

To start, fix a point in $\overline{\mathbb{CP}}^2$. Its complement can be viewed as the total space of the complex line bundle $E \to \mathbb{CP}^1$ whose points are equivalence classes of triples $(u_1, u_2; \zeta) \in (\mathbb{C}^2 - 0) \times \mathbb{C}$ with the equivalence relation identifying $(u_1, u_2; \zeta)$ with $(\lambda u_1, \lambda u_2; \lambda^{-1}\zeta)$ when $\lambda$ is a non-zero complex number. Let $\mathfrak{p}: E \to \mathbb{C}^2$ denote now the map that is defined by the rule whereby

$$(u_1, u_2; \zeta) \to (\zeta u_1, \zeta u_2).$$
(3.8)

This is a complex analytic diffeomorphism from the complement of the zero-section in E to $\mathbb{C}^2$–0. Now choose an isometric almost complex structure on $\mathbb{R}^4$ to identify $\mathbb{R}^4$ as $\mathbb{C}^2$. This allows $\mathfrak{p}$ to be viewed as a map from E to $\mathbb{R}^4$; and in this guise, it restricts to the complement of the zero section as an orientation preserving diffeomorphism onto $\mathbb{R}^4$–0. Given $r \in (0, \infty)$, the upcoming Proposition 3.1 use $\mathfrak{q}_r$ to denote the map from $\mathbb{R}^4$ to $\mathbb{R}^4$ that multiplies any given point by $r^{-1}$; thus $\mathfrak{q}_r(x) = r^{-1}x$.

**Proposition 3.2**: *The Fubini-Study metric on $\overline{\mathbb{CP}}^2$ restricts so as to define a metric on E with vanishing anti-self dual Weyl curvature, vanishing traceless Ricci curvature and positive scalar curvature. This metric has the following properties with regards to the map $\mathfrak{p}$ and any given $r > 0$ version of the map $\mathfrak{q}_r$: Let $\eth$ denote the Euclidean metric on $\mathbb{R}^4$ and let $g_r$ denote the metric on $\mathbb{R}^4$–0 that is obtained by multiplying the $(\mathfrak{p}^{-1} \circ \mathfrak{q}_r)$ pull-back of the Fubini-Study metric by $r^{-2}|x|^4$. If $x \in \mathbb{R}^4$–0 and $r > 0$, then*

- $|g_r - \eth| \leq \kappa r^2 |x|^{-2}$,
- $|\nabla g_r| \leq \kappa r^2 |x|^{-3}$,
- $|\nabla\nabla g_r| \leq \kappa r^2 |x|^{-4}$,

*with $\kappa$ being independent of the parameter $r$ and the point $x$.*

### e) Metrics on $X_K$ with small anti-self dual Weyl curvature

Suppose as before that K is a hyperbolic knot in $S^3$. This subsection uses the metrics described in Proposition 3.1 on $T^*S^2$ and Section 3b's locally conformally flat structure $[g_R]$ on $\mathbb{T}_K/\iota$ to construct a family of conformal structures on $X_K$ with small self-dual Weyl curvature. The construction uses embeddings from disk subundles of $E^2$ into $X_K$ to implement the surgery described in Part 4 of Section 2a. (There is one embedding for each of the sixteen singular points in $\mathbb{T}_K/\iota$). Part 1 of what follows describes these embeddings. Parts 2 and 3 construct the family of conformal structures.

*Part 1*: This part of the subsection describes the required embeddings of disk subbundles of $E^2$ into $X_K$. To start, introduce by way of notation *m* to denote the $2 \times 2$ matrix whose entries are the coefficients $\{\mathfrak{m}_{ab}\}_{a,b \in \{1,2\}}$ that appear in (3.6). Since *m* is a symmetric and positive definite matrix, there is a $2 \times 2$ matrix to be denoted by $z$ with positive determinant such that $z^T m z = \mathbb{I}$. Fix such a matrix $z$.

Supposing that $o \in \mathbb{T}_K$ denotes a given fixed point of $\iota$, write its $(t_1, t_2, t_3, t_4)$ coordinates as $(o_1, o_2, o_3, o_4)$ with each entry either 0 or $\pi$ as the case may be. For each index $a \in \{1, 2\}$, let $x_a = e^{-R} \sum_{b \in \{1,2\}} (z^{-1})_{ab} (t_a - o_a)$; and for $a \in \{3, 4\}$, let $x_a = (t_a - o_a)$ with this definition valid where $\sum_{a \in \{1,2,3,4\}} |t_a - o_a|^2$ is less than $\frac{1}{16} t_*^2$. This condition defines an open neighborhood of $o$ in $\mathbb{T}_K$ and the set of functions $\{x_1, x_2, x_3, x_4\}$ are coordinates on this neighborhood of $o$. The metric $g_T$ when written with these coordinates is the Euclidean metric $dx_1^2 + dx_2^2 + dx_3^2 + dx_4^2$. The involution $\iota$ when written with these coordinates is the map $(x_1, x_2, x_3, x_4) \to (-x_1, -x_2, -x_3, -x_4)$. The coordinate functions $\{x_1, x_2, x_3, x_4\}$ are used to define complex coordinates $(u_1, u_2)$ on this neighborhood of $o$ in $\mathbb{T}_K$. This is done by first choosing an element $O \in SO(4)$, then rotating the coordinate vector $x = (x_1, x_2, x_3, x_4)$ by $O$ to obtain $x' = O \cdot x$, and with $x' = (x'_1, x'_2, x'_3, x'_4)$ in hand, then writing the complex coordinates $(u_1, u_2)$ as $u_1 = x'_1 + i x'_2$ and $u_2 = x'_3 + i x'_4$.

If $\delta$ is a sufficiently small but positive number, then the $|x| < \delta$ ball in $\mathbb{R}^4$ is well inside the coordinate chart from the preceding paragraph for a neighborhood of any given fixed point of $\iota$ in $\mathbb{T}_K$. This is to say that if $|x| < \delta$ and $o = (o_1, o_2, o_3, o_4) \in \mathbb{T}$ is a given fixed point of $\iota$, then $\sum_{a \in \{1,2,3,4\}} |t_a(x) - o_a|^2 < \frac{1}{16} t_*^2$ with $t_a(x) = o_a + e^R \sum_{b=1,2} z_{ab} x_b$ for $a \in \{1,2\}$ and $t_a(x) = o_a + x_a$ for $a \in \{3, 4\}$. Fix $\delta > 0$ so that this is the case and let $B_\delta$ denote both this $|x| < \delta$ ball in $\mathbb{R}^4$ and the corresponding ball about $o$ in $\mathbb{T}_K$. Having chosen a matrix $O \in SO(4)$, the corresponding complex coordinates (these being $(u_1, u_2)$ from the preceding paragraph) identify the $\mathbb{R}^4$ incarnation of $B_\delta$ as a ball about 0 in $\mathbb{C}^2$. This identifies $\mathbb{T}_K$ incarnation of $B_\delta$ with the radius $\delta$ ball about the origin in $\mathbb{C}^2$. Supposing that $r > 0$ has been chosen, then the rescaling map from $\mathbb{C}^2$ to $\mathbb{C}^2$ sending any given point $u$ to $r^{-1} u$ in turn identifies the incarnation of $B_\delta$ in $\mathbb{T}_K$ with the ball of radius $\delta/r$ about the origin in $\mathbb{C}^2$. The latter ball is denoted by $\mathbb{B}_{\delta/r}$.

To continue, reintroduce the complex line bundle $E^2 \to \mathbb{CP}^1$ from the previous section. With a positive number $\mathfrak{r}$ chosen, define $E^2_{\mathfrak{r}}$ to be the disk subbundle in E where $(|u_1|^2 + |u_2|^2) |\zeta| \leq \mathfrak{r}^2$. The map $\mathfrak{p}$ from $E^2$ to $\mathbb{C}^2/\{\pm 1\}$ restricts to $E^2_{\mathfrak{r}}$ as a diffeomorphism from $E^2_{\mathfrak{r}} - \underline{0}$ onto the $(\mathcal{B}_{\mathfrak{r}} - 0)$ part of $\mathbb{C}^2/\{\pm 1\}$ that comes from the complement of 0 in the radius $\mathfrak{r}$ ball about the origin.

With all of the preceding understood, the surgery in Part 4 of Section 1 is implemented at the given singular point of $\iota$ in $\mathbb{T}_K/\iota$ as follows: Fix $\delta > 0$ and $r > 0$. The quotient $B_\delta/\{\pm 1\}$ is viewed simultaneously as a neighborhood of the image of $o$ in $\mathbb{T}_K/\iota$ using the coordinate functions $(x_1, x_2, x_3, x_4)$, and as $\mathbb{B}_{\delta/r}/\{\pm 1\}$ in $\mathbb{C}^2/\{\pm 1\}$ using first a matrix $O \in SO(4)$ to define rotated coordinate functions $x´ = O \cdot x$, then the complex coordinates $u = (x´_1 + i x´_2, x´_3 + i x´_4)$, and finally the scaling map $u \to u/r$. The map $\mathfrak{p}$ from (3.7) is then used to identify $E^2_{\delta/r} - \underline{0}$ with $(B_{\delta/r} - 0)/\{\pm\}$ and thus with the complement of the image of $o$ in a neighborhood of this image in $\mathbb{T}_K/\iota$. This is summarized schematically by the following diagram:

$$\begin{array}{ccc} & (B_\delta - 0)/\{\pm 1\} & \\ \swarrow & & \searrow \\ \mathbb{T}_K / \iota & & E^2_{\delta/r} \end{array}$$

(3.9)

where both arrows are embeddings onto open sets; the right hand arrow being the composition $(\mathfrak{p}^{-1}) \circ q_r \circ O$ and the left hand arrow being the identification of first $B_\delta$ as an open set in $\mathbb{T}_K$ and then the projection from $\mathbb{T}_K$ to $\mathbb{T}_K/\iota$.

*Part 2*: This part of the subsection uses the embeddings that are indicated by the arrows in (3.8) to construct metrics on $X_K$ with very small self-dual Weyl curvature. To start the construction, choose once and for ever a smooth, non-decreasing function on $\mathbb{R}$ that is equal to 1 on $(-\infty, \frac{1}{4}]$ and equal to 0 on $[\frac{3}{4}, \infty)$. This function is denoted in what follows by $\chi$. Given $\delta > 0$, the function $\chi$ will be used to define a function on $\mathbb{R}^4$ to be denoted by $\chi_\delta$ by the rule $\chi_\delta(x) = \chi(2|x|/\delta - 1)$. Thus, $\chi_\delta = 0$ on the complement of the radius $\delta$ ball centered at the origin in $\mathbb{R}^4$; and $\chi_\delta = 1$ on the concentric, radius $\frac{1}{2}\delta$ ball.

With $r > 0$ chosen, let $q_r$ again denote the composition of first multiplication by $r^{-1}$ on $\mathbb{R}^4$ and then the projection map from $\mathbb{R}^4$ to $\mathbb{R}^4/\{\pm 1\}$. Let $g_r$ denote the metric on $\mathbb{R}^4 - 0$ given by $r^{-2}$ times the pull-back of the metric $\hat{g}$ on $T^*S^2$ via the map $(\mathfrak{p}^{-1}) \circ q_r \circ O$. Let $\partial$ again denote the Euclidean metric on $\mathbb{R}^4 - 0$. Having fixed also $\delta > 0$, define a third metric on $\mathbb{R}^4 - 0$ as follows:

$$g_{O,r,\delta} = \partial + \chi_\delta(g_r - \partial)$$

(3.10)

The metric in (3.10) descends to define a Riemannian metric on $(\mathbb{R}^4 - 0)/\{\pm 1\}$ because the function $\chi_\delta$ and both the Euclidean metric and the metric $g_r$ are invariant with respect to the action of $\{\pm 1\}$. The metric defined by (3.9) on $\mathbb{R}^4/\{\pm 1\}$ pulls back to $E^2 - \underline{0}$ via the map $\mathfrak{p}$ where it defines a metric that is identical to the metric $\hat{g}$ on the set of points $(u_1, u_2; \zeta) \in E^2 - \underline{0}$ where $(|u_1|^2 + |u_2|^2)|\zeta| \le \frac{1}{4}\delta^2 r^{-2}$. It thus extends across the 0-section of E as the metric $\hat{g}$. The metric defined by (3.9) is the Euclidean metric on the complement of $B_\delta - 0$. It therefore smoothly extends the locally conformally flat metric on $(\mathbb{T}_K-(B_\delta-0))/\iota$ to define a smooth metric on the part of $X_K$ that comes from the surgery in Part 4 of Section 2a that replaced $(B_\delta-0)/\iota$ by $E_{\delta/r}$.

The lemma that follows summarizes the properties of the metric in (3.10) that are needed in subsequent parts of this article. Having chosen a fixed point of $\iota$ in $\mathbb{T}_K$ and positive numbers $\delta$ and $r$, the lemma implicitly uses the maps in (3.9) to view $(B_\delta-0)/\iota$ on the one hand as a subspace in $\mathbb{T}_K/\iota$ and on the other as the complement of the zero section in disk bundle $E^2_{\delta/r}$.

**Lemma 3.3**: *There exists $\kappa > 1$ depending on the knot K and the choice of the parameter R with the following significance: Let o denote a fixed point of $\iota$ on $\mathbb{T}_K$. Fix $\delta \in (0, \kappa^{-1}]$ and $r \in (0, \delta]$, and then an element in SO(4) so as to define the metric $g_{o,r,\delta}$ that is depicted in (3.10). This metric is the pull-back to $B_\delta-0$ of a metric on $(B_\delta-0)/\{\pm 1\}$ with the properties listed below.*

- *It extends the flat metric defined by $g_T$ from a neighborhood of the complement of $(B_\delta-0)/\iota$ in $\mathbb{T}_K$ to $E_{\delta/r}$.*
- *The metric volume of $E_{\delta/r}$ differs from $\frac{1}{2}\pi^2\delta^4$ by at most $\kappa r^4$.*
- *The versions of the curvatures $\mathcal{W}_+$, $\mathcal{B}$ and R as defined by (3.1) are non-zero only in the $\{\pm 1\}$ quotient of the part of $B_\delta$ where $|x| \in (\frac{1}{2}\delta, \delta)$. In any event, the metric norms of these curvatures obey $|\mathcal{W}_+| + |\mathcal{B}| + |R| \leq \kappa r^4 \delta^{-6}$.*

*Proof of Lemma 3.3*: The first bullet summarizes what is said subsequent to (3.10). The second and third bullets follow directly from the formula in (3.10) using what is said by Proposition 3.1 about the metric $g_r$.

*Part 3*: Let $\Theta$ denote the set of sixteen singular points in $\mathbb{T}_K/\iota$. Each singular point has its corresponding $(B_\delta-0)/\{\pm 1\}$ neighborhood and its corresponding copy of $E^2_{\delta/r}$ in $X_K$. With an element in the group SO(4) chosen for each singular point, and with r fixed, the formula in (3.10) defines a metric on each the sixteen copies of $E^2_{\delta/r}$ in $X_K$. This set of metrics smoothly extend the locally conformally flat metric defined previously on the $(\mathbb{T}_K-(\cup_{o\in\Theta}(B_\delta-0)))/\iota$ part of $X_K$ to define a conformal structure on all of $X_K$. The precise choice of r and $\delta$ is not so important in the subsequent applications except that $\delta$ and $r/\delta$ must be small in a suitable sense. This said, it proves convenient to fix $r = \delta^2$ in what follows. The resulting conformal structure is denoted below by $[\mathfrak{g}_{R,\delta}]$. The notation does not indicate that the conformal structure $[\mathfrak{g}_{R,\delta}]$ depends on the elements in SO(4) assigned to each of the sixteen singular point in $\mathbb{T}_K/\iota$. These group element are not noted because they do not play a significant role in what follows.

The $(\mathbb{T}_K - (\cup_{o \in \Theta} (B_\delta - 0)))/\iota$ part of $X_K$ is denoted in what follows by $X_{K,\delta}$. The subspace $X_{K,\delta}$ has the conformal structure $[\mathfrak{g}_{R,\delta}]$ from $X_K$ and it has the conformal structure $[\mathfrak{g}_R]$ coming from its incarnation as a subspace of $\mathbb{T}_K/\iota$. These two conformal structures agree. As a consequence, each point in $X_{K,\delta}$ has a neighborhood where $[\mathfrak{g}_{R,\delta}]$ is defined by a flat metric.

### f) Metrics on $X_K \#_n \overline{\mathbb{CP}}^2$ with small self-dual Weyl curvature

The construction of conformal structures on the connect sum of $X_K$ with n copies of $\overline{\mathbb{CP}}^2$ with small self-dual Weyl curvature is described in this subsection. The upcoming construction is a close kin to the construction in the preceding subsection. The construction is given in four parts.

*Part 1*: To set the stage, suppose that $\delta \in (0, \frac{1}{64} t_*)$ and the required 16 elements in SO(4) have been chosen so as to define the conformal structure $[\mathfrak{g}_{R,\delta}]$ using the instructions from Parts 2 and 3 of the preceding subsection. Reintroduce the subset $X_{K,\delta}$ of $X_K$ from Part 3 of the preceding subsection and suppose that p is a given point in $X_{K,\delta}$. As noted previously, p has a neighborhood where the conformal structure $[\mathfrak{g}_{R,\delta}]$ is defined by a flat metric. Fix such a neighborhood of p with its flat metric, and then a metric ball in this neighborhood where there are coordinates that identify this flat metric with the Euclidean metric on $\mathbb{R}^4$. If p is in the $(\mathbb{T} - (N \cup \iota(N)))/\iota$ part of $X_K$, then the flat metric can be taken to be $\mathfrak{g}_T$ in which case the radius of this ball can be chosen a priori given only the knot K, the number R and the $\mathfrak{g}_T$ distance from p to the $\mathfrak{g}_T$-radius $\delta$ balls about the fixed points of $\iota$. If p is in either the $(S^3 - N_K) \times S^1$ or the $\iota((S^3 - N_K) \times S^1)$ part of $X_{K,\delta}$, then the flat metric can be chosen so that the radius of this ball can be chosen a priori given only the numbers R and the knot K. In particular, the radius in this case does not depend on the parameter $\delta$ nor does it depend on the SO(4) parameters that were chosen for each fixed point of $\iota$ on $\mathbb{T}$. In any case, given p, let $U_p$ denote such a ball in $X_K$.

There is an SO(4)'s worth of Gaussian coordinate charts at p that are defined by the chosen flat metric that defines the equivalence structure $[\mathfrak{g}_{R,\delta}]$ on $U_p$. The coordinates in any two of these charts are related by an SO(4) rotation. In particular the set of such Gaussian coordinate charts has a canonical identification with the fiber at p of the orthonormal frame bundle for any metric from the equivalence structure $[\mathfrak{g}_{R,\delta}]$. (If g and g´ are any two metrics in the equivalence class, then a g-orthonormal frame at p will be g´-orthogonal and the four frame vectors will have the same g´-length. As a consequence they define a canonical g´-orthonormal frame at p.) All of these charts embed a ball of some radius $\mathfrak{r}_p$ about the origin in $\mathbb{R}^4$ as an open neighborhood of p in $U_p$.

*Part 2*: Let $(x_1, x_2, x_3, x_4)$ denote Gaussian coordinates for a neighborhood of p in $U_p$ as defined using the chosen flat metric that gives the equivalence structure $[\mathfrak{g}_{R,\delta}]$ near p. Fix $\delta_p > 0$ so that the $|x| < \delta_p$ ball in $\mathbb{R}^4$ is in the domain of the coordinate functions. Use $B_p$ henceforth to denote both the $|x| < \delta_p$ ball in $\mathbb{R}^4$ and its corresponding image in $X_K$ via the coordinate embedding. Writing $u_1 = x_1 + ix_2$ and $u_2 = x_3 + ix_4$ identifies the $X_K$ incarnation of $B_p$ with the radius $\delta_p$ ball about the origin in $\mathbb{C}^2$. Supposing that $r > 0$ has been specified, the subsequent multiplication of the coordinates by $r^{-1}$ identifies the incarnation of $B_p$ in $X_K$ with the radius $\delta_p/r$ ball about the origin in $\mathbb{C}^2$. Let $\mathfrak{q}_r$ now denote the latter identification.

Meanwhile, the map $\mathfrak{p}$ from the complex line bundle $E \to \mathbb{CP}^1$ to $\mathbb{C}^2$ given in (3.8) identifies the complement of the zero section in E with $\mathbb{C}^2 - 0$. Supposing that $\mathfrak{r} > 0$ has been specified, let $E_\mathfrak{r}$ denote the disk subbundle in E where the homogeneous coordinates $(u_1, u_2; \zeta)$ obey $(|u_1|^2 + |u_2|^2)^{1/2}|\zeta| < \mathfrak{r}$. The map $\mathfrak{p}$ identifes $E_\mathfrak{r} - \underline{0}$ with the complement of the origin in radius $\mathfrak{r}$ ball about the origin in $\mathbb{C}^2$. As a consequence, the composition of first $\mathfrak{p}$ and then $\mathfrak{q}_r^{-1}$ identifies complement of the zero section in the disk sub-bundle $E_{\delta_p/r}$ with the complement of p in the $X_K$ incarnation of $B_p$.

With the preceding understood, the connect sum operation to obtain $X_K \# \overline{\mathbb{CP}}^2$ from $X_K$ is implemented by the surgery on $X_K$ that removes $B_p$ and then replaces it with any given $r > 0$ version of $E_{\delta_p/r}$ with the complement of the zero section in $E_{\delta_p/r}$ being identified with $B_p - p$ using the composition of first $\mathfrak{p}$ and then $\mathfrak{q}_r^{-1}$.

Supposing that n is a given positive integer, then the connect sum operation to obtain $X_K \#_n \overline{\mathbb{CP}}^2$ from $X_K$ is implemented by first choosing n distinct points in $X_\delta$ and then choosing for each point, a version of $\delta_p$ such that the resulting n copies of $B_p$ in $X_K$ are pairwise disjoint. Once this is done, then the connect sum construction from the preceding paragraphs can be done in each copy of $B_p$ simultaneously having chosen for each of the n points a Gaussian coordinate chart and a positive number $r_p$.

*Part 3*: Fix a point p in $X_{K,\delta}$ and define its open neighborhood $U_p$ as instructed in Part 1. Fix a Gaussian coordinate chart centered at p using a metric that gives the conformal structure $[\mathfrak{g}_{R,\delta}]$ on $U_p$. Then choose $\delta_p > 0$ so that the $|x| < \delta_p$ ball in $\mathbb{R}^4$ is in the domain of this coordinate chart. Also choose a positive number to be denoted by $r_p$.

Let $g_{r_p}$ denote the $r = r_p$ version of the metric $g_r$ on $\mathbb{R}^4 - (0)$ that is described in Proposition 3.2. Let $\chi_{\delta_p}$ denote the $\delta = \delta_p$ version of the function $\chi_\delta$ that was introduced in Part 2 of the preceding subsection. Use this function with $g_{r_p}$ and the Euclidean metric $\mathfrak{d}$ to define a the metric $g_p$ on $\mathbb{R}^4 - 0$ by the formula

$$g_p = \eth + \chi_{\delta_p}(g_{r_p} - \eth).$$
(3.11)

Since the pull-back of this metric to the complement of the zero section of $E_{\delta_p/r_p}$ by the map $\mathfrak{p}$ is the Fubini-Study metric near the zero section, this pull-back extends smoothly to the whole of the disk bundle $E_{\delta_p/r_p}$. Since this metric is the Euclidean metric on a neighborhood of the complement of the $|x| \le \delta_p$ ball in $\mathbb{R}^4$, it also smoothly extends as the flat metric on $U_p$ that defines the conformal structure $[\mathfrak{g}_{R,\delta}]$ on $U_p$.

The next lemma summarizes what was just said and it says somethings about the curvature of the metric $g_p$. By way of notation, the upcoming lemma implicitly uses the chosen Gaussian coordinate chart for the flat metric in the conformal structure of $[\mathfrak{g}_{R,\delta}]$ with the maps $\mathfrak{q}_{r_p}$ and $\mathfrak{p}$ to identify the set $B_p-p$ in $X_{K,\delta}$ with the complement of the zero section in the disk bundle $E_{\delta_p/r_p}$.

**Lemma 3.4**: *Fix $p \in X_{K,\delta}$ and a flat metric $U_p$ in the conformal structure $[\mathfrak{g}_{R,\delta}]$. There exists $\kappa > 1$ which is independent of the parameter $r$ used to define $[\mathfrak{g}_{R,\delta}]$ and which has the following significance: Fix $\delta_p \in (0, \kappa^{-1}]$ and $r_p \in (0, \delta]$, and then a Gaussian coordinate chart centered at $p$ so as to define the metric $g_p$ that is depicted in (3.11).*
- *The metric $g_p$ extends the chosen flat metric in the conformal structure of $[\mathfrak{g}_{R,\delta}]$ from a neighborhood of the complement of $(B_p-p)$ in $U_p$ to $E_{\delta_p/r_p}$.*
- *The $g_p$ metric volume of $E_{\delta_p/r_p}$ differs from $\tfrac{1}{2}\pi^2 \delta_p^4$ by at most $\kappa r_p^2 \delta_p^2$.*
- *The $g_p$ versions of the curvatures $\mathcal{W}_+$ and $\mathcal{B}$ as defined by (3.1) are non-zero only in the part of $B_p$ where $|x| \in (\tfrac{1}{2}\delta_p, \delta_p)$; and their norms here obey $|\mathcal{W}_+| + |\mathcal{B}| \le \kappa r_p^2 \delta_p^{-4}$.*

*Proof of Lemma 3.4*: The assertion in the top bullet follows from what is said in the preceding paragraph and the assertion in the lower bullet follows from what is said in Proposition 3.2.

*Part 4*: Fix a positive integer n and then a set of n distinct points in $X_{K,\delta}$ to be denoted by $\vartheta$. Supposing that $p \in \vartheta$, define its open neighborhood $U_p$ as instructed in Part 1. Fix a Gaussian coordinate chart centered at $p$ using a flat metric that gives the conformal structure $[\mathfrak{g}_{R,\delta}]$ on $U_p$. Then choose $\delta_p > 0$ so that the $|x| < \delta_p$ ball in $\mathbb{R}^4$ is in the domain of this coordinate chart and subject to the constraint that resulting set of balls $\{B_p\}_{p\in\vartheta}$ when viewed in $X_K$ are pairwise disjoint. Each $p \in \vartheta$ version of $r_p$ is chosen to be proportional to $\delta_p^2$ with the proportionality constant being at most 1 and taken to be independent of p. This is to say that $r_p = \varepsilon \delta_p^2$ for each $p \in \vartheta$ with $\varepsilon \in (0, 1]$. A particular choice of $\varepsilon$ is fixed in Section 4f. Given what is said in Parts 2 and 3 and Lemma 3.4, it follows that the various $p \in \vartheta$ versions of (3.11) with the conformal structure $[\mathfrak{g}_{R,\delta}]$ on $X_K - (\cup_{p \in \vartheta} B_p)$ define a conformal structure on the whole of $X_K \#_n \overline{\mathbb{CP}}^2$.

The conformal structure just defined is denoted in what follows by $[\mathfrak{g}^n_{R,\delta}]$. Keep in mind that this notation does not indicated that the conformal structure depends on various choices that have been made. By way of a summary, these choices are as follows: The set $\vartheta$ must be chosen and then a flat metric in the $[\mathfrak{g}_{R,\delta}]$ conformal structure on $X_K$ for a neighborhood of each $p \in \vartheta$. Given the flat metric, the positive number $\delta_p$ must be chosen; and a Gaussian coordinate chart must be selected for the radius $\delta_p$ ball as defined by this metric. (Only two of these choices at each $p \in \vartheta$ are relevant for what is to come; the first is $\delta_p$ and the second is the choice of the Gaussian coordinate chart at p up to the action of the U(2) subgroup in SO(4).) Finally, a number $\varepsilon \in (0, 1]$ must also be selected so as to specify each $p \in \vartheta$ version of $r_p$ as $r_p = \varepsilon \delta_p^2$.

## 4. Deformation to metrics with zero self-dual Weyl curvature

Theorem A in Kovalev and Singer's paper [KS] will be invoked so as to obtain metrics with vanishing self-dual Weyl curvature from some of the conformal structures that are described in Section 3e.

### a) Metrics in the conformal classes

Suppose as in the preceding section that K is a knot in $S^3$ with hyperbolic complement. A number $R > 1$ should be chosen first to define the manifold $\mathbb{T}_K$ and its $\iota$-invariant, locally conformally flat structure $[\mathfrak{g}_R]$ as instructed in Section 3b. Supposing that R has been chosen, fix $\delta > 0$ subject to the constraints in Section 3e. Having chosen an element in SO(4) for each singular point of $\mathbb{T}_K/\iota$, use $\delta$ to define both $X_K$ and the conformal structure $[\mathfrak{g}_{R,\delta}]$ per the instructions in Part 3 of Section 3e. Now fix a non-negative integer and then fix this same number of points $X_{K,\delta}$; the number is denoted by n and the set in $X_{K,\delta}$ is denoted by $\vartheta$. Supposing that $p \in \vartheta$, fix a Gaussian coordinate chart based at p for a flat metric in the conformal structure of $[\mathfrak{g}_{R,\delta}]$ on a neighborhood of p. Choose next a positive number $\delta_p$ subject to the constraints in Section 3f, and then choose a Gaussian coordinate chart for the radius $\delta_p$ ball in $X_K$ centered at p. Use these choices to define the conformal structure $[\mathfrak{g}^n_{R,\delta}]$ on $X_K \#_n \overline{\mathbb{CP}}^2$ according to the instructions in Part 4 of Section 3f.

A version of Theorem A in [KS] asserts that $X_K \#_n \overline{\mathbb{CP}}^2$ has a metric with vanishing self-dual Weyl curvature if $\delta$ and $\{\delta_p\}_{p \in \vartheta}$ are all sufficiently small, and provided that a certain obstruction is zero. What is meant by 'sufficiently small' is defined momentarily in the upcoming Proposition 4.1 The heart of the matter is the obstruction. It is described at length in Sections 5 and 6.

Theorem A in [KS] also implies that the new metric with anti-self dual Weyl curvature (when it exists) is close in a suitable sense to a metric that defines the conformal structure $[\mathfrak{g}^n_{R,\delta}]$. The three steps that follow specify the relevant metric from $[\mathfrak{g}^n_{R,\delta}]$; Proposition 4.1 then defines what is meant by 'close'.

<u>Step 1</u>: This step specifies a metric on $\mathbb{T}_K$ and on the complement of the singular points in $\mathbb{T}_K/\iota$ that gives the conformal structure $[\mathfrak{g}_R]$. To start, fix a smooth function of the coordinate s that appears in (3.3) and (3.5) which is equal to 0 where $s \leq \frac{1}{2}R - 1$ and which is equal to $s - R$ where $s \geq \frac{1}{2}R + 1$. Denote this function by $f$. Then the metric

$$\mathfrak{g}_f = e^{2f}(e^{-2s}\mathfrak{m} + ds^2 + d\theta^2)$$

(4.1)

on the $s \geq 0$ part of $(S^3 - N_K) \times S^1$ agrees with the product of the hyperbolic metric on $S^3 - N_K$ and the Euclidean metric on $S^1$ where $s \leq \frac{1}{2}R - 1$. This product metric where $s \leq \frac{1}{2}R - 1$ and $\mathfrak{g}_f$ define a smooth metric on $(S^3 - N_K) \times S^1$ which is conformal to the product metric. Denote this new metric on $(S^3 - N_K) \times S^1$ by $\mathfrak{g}_*$. Since $\mathfrak{g}_f$ agrees with the flat metric in (3.5) (with $\rho = e^{s-R}$) where $s \geq \frac{1}{2}R + 1$, the metric $\mathfrak{g}_*$ on $(S^3 - N_K) \times S^1$, the metric $\iota^* \mathfrak{g}_*$ on $\iota((S^3 - N_K) \times S^1)$ and the metric $\mathfrak{g}_T$ on the rest of $\mathbb{T}_K$ define a smooth, $\iota$-invariant metric on $\mathbb{T}_K$ in the conformal structure $[\mathfrak{g}_R]$. This metric descends to the complement of the singular points in $\mathbb{T}_K/\iota$ to define a metric in the conformal class of $[\mathfrak{g}_R]$ on the smooth manifold part of $\mathbb{T}_K/\iota$. This $\iota$-invariant metric on $\mathbb{T}_K$ and the corresponding induced metric on the smooth part of $\mathbb{T}_K/\iota$ are both denoted by $\mathfrak{g}_R$.

<u>Step 2</u>: Suppose that a positive number $\delta$ has been chosen subject to the constraints in Section 3e and let $B_\delta$ denote the ball of radius $\delta$ centered at the origin in $\mathbb{R}^4$. Let $o$ denote one of the 16 fixed points of $\iota$ in $\mathbb{T}_K$. Section 3e describes a certain identification between $(B_\delta - 0)/\{\pm 1\}$ and the complement of $o$ in one of its neighborhood in the quotient space $\mathbb{T}_K/\iota$. As explained in Section 3e, the manifold $X_K$ is obtained from $\mathbb{T}_K/\iota$ by removing each of these 16 incarnations of $(B_\delta - 0)/\{\pm 1\}$ from $\mathbb{T}_K/\iota$ and replacing each by a copy of the disk subbundle $E^2_{1/\delta}$ in the complex line bundle $E^2 \to \mathbb{CP}^1$. Part 2 of Section 3e explains why the $r = \delta^2$ version of the metric in (3.10) extends the metric induced by $\mathfrak{g}_T$ on the complement of the $(B_\delta - 0)/\{\pm 1\}$ part of $\mathbb{T}_K$ to the whole of $X_K$. These sixteen extension also extend the metric $\mathfrak{g}_*$ from the complement of the 16 versions of $E^2_{1/\delta}$ in $X_K$ to the whole of $X_K$. Denote this metric by $\mathfrak{g}_{R,\delta}$. By construction, this metric defines the conformal structure $[\mathfrak{g}_{R,\delta}]$ on $X_K$.

A different, but conformal extension of $\mathfrak{g}_R$ over $E^2_{1/\delta}$ is introduced in [KS]. This alternate metric extends $\mathfrak{g}_R$ as $\phi^2 \mathfrak{g}_{R,\delta}$ with the conformal factor $\phi^2$ defined as follows: Fix $\delta_0 > 0$ subject to the constraints in Section 3e. The $r = \delta^2$ version of the metric in (3.10) is then defined on the complement of 0 in $|x| \leq \delta_0$ ball in $\mathbb{R}^4$ when $\delta < \delta_0$. Assume in what follows that $\delta < \frac{1}{16}\delta_0$. Define $u = \chi(4\delta_0^{-1}(|x| + \delta^2|x|^{-1}) - 1)$ and set $\phi^2 = (1 - u + u|x|^2)^{-1}$. This function $\phi^2$ is equal to 1 where $|x| \geq \delta_0$ and where $|x| \leq \delta^2/\delta_0$. It follows as a consequence that the metric $\phi^2 \mathfrak{g}_{R,\delta}$ also extends $\mathfrak{g}_T$ over the $E^2_{1/\delta}$ parts of $X_K$ as a member of the conformal structure $[\mathfrak{g}_{R,\delta}]$. This new extension is denoted by $\mathfrak{g}_{R,\delta*}$. The pointwise curvature norms of this $\mathfrak{g}_{R,\delta*}$ metric have $\delta$-independent upper bounds; and its injectivity radius has a $\delta$-independent, positive lower bound. There are no such bounds for $\mathfrak{g}_{R,\delta}$. (These bounds arise because $\phi^2$ is equal to $|x|^{-2}$ where $|x|$ is between $16\delta^2/\delta_0$ and $\frac{1}{16}\delta_0$.) On the other hand, the $\mathfrak{g}_{R,\delta*}$ volume of $X_K$ grows like a positive multiple of $|\ln\delta|$ as $\delta$ limits to 0. By comparison, the $\mathfrak{g}_{R,\delta}$-volume of $X_K$ has a $\delta$-independent upper bound and a positive, $\delta$-independent lower bound.

<u>Step 3</u>: Fix a positive integer n and consider the conformal structure $[\mathfrak{g}^n_{R,\delta}]$ on $X_K \#_n \overline{\mathbb{CP}}^2$. Let p denote a given point from the set $\vartheta$. Keeping in mind that $\vartheta$ is a subset of $X_{K,\delta}$, there is an open set $U_p$ containing p where the metrics $\mathfrak{g}_{R,\delta}$ and $\mathfrak{g}_{R,\delta*}$ are the same, and both are conformally flat. Write $\mathfrak{g}_{R,\delta}$ on $U_p$ as $h^2 \mathfrak{g}_F$ where h is a positive function on $U_p$ and where $\mathfrak{g}_F$ is the flat metric that is used in the constructions of Section 3f. The chosen Gaussian coordinates for p writes $\mathfrak{g}_R$ as the Euclidean metric $\mathfrak{d}$, so use these coordinates to write $\mathfrak{g}_{R,\delta}$ near p as $h^2 \mathfrak{d}$ where $h^2$ is now viewed as a function on a neighborhood of the origin in $\mathbb{R}^4$. This depiction of $\mathfrak{g}_{R,\delta}$ near p is in particular valid in some radius $\delta_{p0}$ ball centered at p. It is assumed in what follows that the parameter $\delta_p$ that is used to define the conformal structure $[\mathfrak{g}^n_{R,\delta}]$ near p obeys $\delta_p \leq \frac{1}{16}\delta_{p0}$. With the preceding understood, use the function $\chi$ from Part 2 of Section 3e to construct a function to be denoted by $h_p$ that is equal to h on the $|x| \geq \frac{1}{2}\delta_{p0}$ ball in $\mathbb{R}^4$ and is equal to 1 on the $|x| \leq \frac{1}{4}\delta_{p0}$ ball. Let $\mathfrak{g}_p$ denote the metric that is depicted in (3.11). The metric $h_p^2 \mathfrak{g}_p$ smoothly extends the metric $\mathfrak{g}_{R,\delta}$ from the complement of a neighborhood of p in $X_K$ to the disk bundle $E_{\delta_p/r_p}$ in $X_K \#_n \overline{\mathbb{CP}}^2$ because it agrees with the metric in (3.11) on the $0 < |x| \leq \delta_p$ part of $\mathbb{R}^4$. Moreover, this extension is in the conformal structure that is defined near p by $[\mathfrak{g}^n_{R,\delta}]$. Using $h^2_p \mathfrak{g}_p$ for each $p \in \vartheta$ defines a metric to be denoted by $\mathfrak{g}^n_{R,\delta}$ that gives the conformal structure $[\mathfrak{g}^n_{R,\delta}]$.

An analog of $\mathfrak{g}_{R,\delta*}$ for $X_K \# \overline{\mathbb{CP}}^2$ is defined as follows: Given $p \in \vartheta$, define the function $u_p$ on $\mathbb{R}^4$ by the rule whereby $u_p(x) = \chi(4\delta_{p0}^{-1}(|x| + \delta_p^2|x|^{-1}) - 1)$. With $u_p$ in hand, define $\phi_p^2$ to be $(1 - u_p + u_p|x|^2)^{-1}$. The metric $\phi_p^2 h_p^2 \mathfrak{g}_p$ smoothly extends the metric $\mathfrak{g}_{R,\delta*}$ from the complement of a neighborhood of p in $X_K$ to the disk bundle $E_{\delta_p/r_p}$ in $X_K \#_n \overline{\mathbb{CP}}^2$. Use these extensions to define the metric $\mathfrak{g}^n_{R,\delta*}$. This metric is also in the conformal structure $[\mathfrak{g}^n_{R,\delta}]$. The curvature norms of the extension $\mathfrak{g}^n_{R,\delta*}$ enjoy $\delta$ and $\{\delta_p\}_{p\in\vartheta}$ independent upper bounds and the injectivity radius enjoys a $\delta$ and $\{\delta_p\}_{p\in\vartheta}$ indendent, positive lower bound.

**b) Theorem A of Kovalev-Singer**

Fix a positive number $\delta$ and the sixteen SO(4) parameters that are needed to construct the conformal structure $[\mathfrak{g}_{R,\delta}]$ on $X_K$. Step 2 in Section 4a describes the metric $\mathfrak{g}_{R,\delta*}$ which gives this conformal class. Supposing that n is a given positive integer, fix a set of n distinct points in $X_{K,\delta}$ to be denoted by $\vartheta$. Having chosen a number $\varepsilon \in (0, 1]$, then chose an appropriate set of positive numbers $\{\delta_p\}_{p\in\vartheta}$ and the various SO(4) parameters to construct the conformal structure $[\mathfrak{g}^n_{R,\delta}]$ for $X_K \#_n \overline{\mathbb{CP}}^2$. Step 3 in Section 4a describes the metric $\mathfrak{g}^n_{R,\delta*}$ which gives this conformal structure. In the case when n = 0, let $\mathfrak{g}_* = \mathfrak{g}_{R,\delta*}$ and in the case n > 0, let $\mathfrak{g}_* = \mathfrak{g}^n_{R,\delta*}$.

The convention in what follows has $c_R$ denoting a number that is greater than 1 with dependence only on the choice of R, the knot K and the number n. In particular, this number does not depend on $\delta$ or the data for the set $\vartheta$. The value of $c_R$ can be assumed to increase between successive appearances.

Supposing that $\delta$ and $\{\delta_p\}_{p\in\vartheta}$ are less than $c_R^{-1}$, then [KS] construct a symmetric, traceless section of $\otimes^2 T^*(X_K \#_n \overline{\mathbb{CP}}^2)$ to be denoted here by $\mathfrak{h}_{\mathfrak{g}_*}$, with sup-norm less than $\frac{1}{100}$ as measured by the metric $\mathfrak{g}_*$, and such that $\mathcal{W}_+(\mathfrak{g}_* + \mathfrak{h}_{\mathfrak{g}_*})$ lies entirely in a certain finite dimensional vector space. This vector space is called the *obstruction* space.. The vector $\mathcal{W}_+(\mathfrak{g}_* + \mathfrak{h}_{\mathfrak{g}_*})$ in the obstruction space is the (tautological) obstruction vector since it is zero if and only if $\mathcal{W}_+(\mathfrak{g}_* + \mathfrak{h}_{\mathfrak{g}_*})$ is zero. This is summarized by Theorem A in [KS]. The following proposition asserts slightly more than Theorem A in [KS].

**Proposition 4.1**: *Fix a hyperbolic knot* K. *There is a finite dimensional vector space to be denoted by* $\mathbb{V}$, *and given also a sufficiently large positive number* R *and* $\varepsilon \in (0, 1]$, *there exists* $\kappa > 1$ *such that the following is true: Fix* $\delta \in (0, \kappa^{-1})$ *and use* R *and* $\delta$ *to define the conformal structure* $[\mathfrak{g}_{R,\delta}]$ *on* $X_K$.
- *There is a metric on* $X_K$ *with anti-self dual Weyl curvature if a certain vector in* $\mathbb{V}$ *determined by* $[\mathfrak{g}_{R,\delta}]$ *is zero.*
- *Fix an integer* n > 0 *and a set of* n *distinct points in* $X_{K,\delta}$ *to be denoted by* $\vartheta$. *Given* $\vartheta$ *and the chosen flat metrics on neighborhoods of the points in* $\vartheta$, *there exists* $\kappa_\vartheta > \kappa$ *such that if* $\delta$ *and all* $p \in \vartheta$ *versions of* $\delta_p$ *are less than* $\kappa_\vartheta$ *then there is a metric with anti-self dual Weyl curvature on* $X_K \# \overline{\mathbb{CP}}^2$ *if a certain vector in* $\mathbb{V}$ *determined by* $[\mathfrak{g}^n_{R,\delta}]$ *is zero.*

*The anti-self dual Weyl curvature metric on* $X_K$ *from the top bullet is close to a metric from* $[\mathfrak{g}_{R,\delta}]$ *when* $\delta$ *is small; and the anti-self dual Weyl curvature metric on any* n ≥ 1 *version of* $X_K \#_n \overline{\mathbb{CP}}^2$ *from the second bullet is close to a metric from* $[\mathfrak{g}^n_{R,\delta}]$ *when* $\delta$ *and* *all $\{\delta_p\}$ are small. In particular, given $\mu > 0$ and a positive integer k; and supposing that $\delta$ and, if $n \geq 1$, all $p \in \vartheta$ versions of $\delta_p$ are sufficiently small (for fixed parameter $\varepsilon$), then the assertions below are true*

i) *Let $\mathfrak{g}$ denote either $\mathfrak{g}_{R,\delta}$ or some $n \geq 1$ version of $\mathfrak{g}^n_{R,\delta}$. The anti-self dual Weyl curvature metric can be written as $\mathfrak{g} + \mathfrak{h}$ where the pointwise norm of $\mathfrak{h}$ as measured by $\mathfrak{g}$ is less than $\mu$.*

ii) *Let $\mathfrak{g}_*$ denote either $\mathfrak{g}_{R,\delta*}$ or $\mathfrak{g}^n_{R,\delta*}$ as the case may be. Write $\mathfrak{g}_*$ as $\phi_*^2 \mathfrak{g}$. The $C^k$ norm of $\phi_*^2 \mathfrak{h}$ as defined by $\mathfrak{g}_*$ is also less than $\mu$.*

*Proof of Proposition 4.1*: With three caveats to be discussed in the subsequent paragraphs, the assertions of the first two bulleted items are instances of Theorem A in [KS]. Item ii) also follows directly from this same theorem. Item i) follows from Item ii) because the pointwise $\mathfrak{g}$ norm of $\mathfrak{h}$ is identical to the pointwise $\mathfrak{g}_* = \phi_*^2 \mathfrak{g}$ norm of $\phi_*^2 \mathfrak{h}$.

The first and second caveats are with regards to the appeal to [KS] in the case of $X_K \#_n \overline{\mathbb{CP}}^2$. These caveats concern two salient differences between the metric $\mathfrak{g}^n_{R,\delta*}$ and the metric that is described in Section 2.3.5 of [KS]: First, the [KS] metrics strictly speaking describe only the metrics here where $\varepsilon = 1$. Second, the [KS] metrics are the $\varepsilon = 1$ metrics here only in the case when all $p \in \vartheta$ versions of $\delta_p$ to equal $\delta$. The fact that the parameters $\delta$ and $\{\delta_p\}_{p \in \vartheta}$ can be allowed to differ is of no significance with regards to what is done in [KS]. The fact that $\varepsilon$ can be taken less than 1 has no significance also. Only notational changes to what is said in [KS] are needed to account for the variation of $\{\delta_p\}_{p \in \vartheta}$ and a choice for $\varepsilon$ that is less than 1. (The variation of $\{\delta_p\}_{p \in \vartheta}$ would require a multi-headed giraffe in Figure 1 of [KS] with different neck lengths. Meanwhile, the $\varepsilon < 1$ giraffes would have the 'damage zone' in Figure 1 of [KS] shifted by $|\ln \varepsilon|$ along the t direction towards the body of the giraffe.) (A value of $\varepsilon$ less than 1 and the freedom to vary $\{\delta_p\}_{p \in \vartheta}$ are used in the proof of the upcoming Proposition 4.2 to prove that the other parameters defining $[\mathfrak{g}^n_{R,\delta}]$ can be chosen so that the associated $\mathbb{V}$ vector is zero.)

The third caveat concerns the assertion in Proposition 4.1 to the effect that the vector space $\mathbb{V}$ does not depend on the parameter R when R is large. Theorem A in [KS] supplies a vector space that can, in principle, depend on R. The assertion that this space is independent of R when R is large is a consequence of what is said Sections 5. In particular, the R-independence of this vector space follows directly from the upcoming Propositions 5.3 and 5.4 which are in Sections 5b and 5c.

**c) The existence of metrics with anti-self dual Weyl curvature**

Sections 5a-5c describe Proposition 4.1's vector space $\mathbb{V}$; and Section 6 describes the relevant vector in $\mathbb{V}$. What is said in these subsections can be used when invoking Proposition 4.1 to obtain an existence assertion for metrics on $X_K \#_n \overline{\mathbb{CP}}^2$ for suitable n:

**Proposition 4.2**: *Suppose that* K *is a hyperbolic knot in* $S^3$. *There exists* N *such that Proposition 4.1's obstruction vector vanishes for certain choices of the parameters that are used to define the conformal structure* $[\mathfrak{g}^n_{R,\delta}]$ *when* $n \geq N$ *and when* R *is sufficiently large (given* K *and* n*) and* $\delta$ *is sufficiently small (given* K, n *and* R *and* $\varepsilon$*). These parameter choices include those with the parameters* $\varepsilon$ *and/or* $\{\delta_p\}_{p \in \vartheta}$ *bounded away from zero for choices of* $\delta$ *that limit to zero. Moreover,*
- *An upper bound for* N *is determined by the hyperbolic volume of* $S^3 - K$.
- *There exists an infinite set of distinct, hyperbolic knots in* $\mathbb{R}^3$ *whose complement in* $S^3$ *has a priori bounded volume such that the preceding is true for* $N = 4$.
- *There exists an infinite set of distinct, hyperbolic knots in* $\mathbb{R}^3$ *with no a priori bound on the volumes of their complements in* $S^3$ *such that the preceding is true for* $N = 4$.

This proposition is proved in Section 7.

The following theorem is now a corollary to Propositions 4.1 and 4.2.

**Theorem 4.3**: *Supposing that* K *is a hyperbolic knot in* $S^3$, *let* N *denote the integer that appears in Proposition 4.2. If* $n \geq N$, *then the manifold* $X_K \#_n \overline{\mathbb{CP}^2}$ *has metrics with anti-self dual Weyl curvature.*

The set of anti-self dual Weyl conformal structures on any given compact 4-manifold has an action of the group of diffeomorphisms of the 4-manifold. As explained in [KK], the quotient can be given a reasonable topology which makes it a stratified space, and a smooth manifold in favorable circumstances. The formal dimension of this space is $-\frac{1}{2}(15e + 29\tau)$ with $e$ denoting the Euler characteristing of X and $\tau$ denoting the signature of X. For $X_K \#_n \overline{\mathbb{CP}^2}$, this number is $52 + 7n$. See [KK] for more about the structure of this space.

**5. The obstruction space**

The central issue in Proposition 4.1 is the obstruction vector space $\mathbb{V}$ and the corresponding vectors that are defined in this space by $[\mathfrak{g}_{R,\delta}]$ and, when $n \geq 1$, by $[\mathfrak{g}^n_{R,\delta}]$. This section describes the space $\mathbb{V}$.

**a) The definition of the obstruction space**

Supposing that some large value for R has been fixed, define the locally conformally flat structure $[\mathfrak{g}_R]$ as directed in Section 3b. The vector space that is obtained by a direct appeal to Theorem A in [KS] is isomorphic to the ι-invariant kernel of a certain differential operator that can be defined using any metric from the locally conformally flat structure $[\mathfrak{g}_R]$. (As explained below, if $\mathfrak{g}$ and $\mathfrak{g}'$ are two metrics that define $[\mathfrak{g}_R]$, then the kernels of the corresponding operators can be canonically identified.) The rest of this subsection describes the operator in question.

Having defined the locally conformally flat structure $[\mathfrak{g}_R]$, let $\Lambda^+$ denote the subbundle in $\wedge^2 T^*\mathbb{T}_K$ of self dual 2-forms. Let $S$ denote the vector space of symmetric, traceless sections of $\Lambda^+ \otimes \Lambda^+$. With $S$ viewed as a subbundle of the space of symmetric sections of $(\wedge^2 T^*\mathbb{T}_K) \otimes (\wedge^2 T^*\mathbb{T}_K)$ and let $\wp: (\wedge^2 T^*\mathbb{T}_K) \otimes (\wedge^2 T^*\mathbb{T}_K) \to S$ denote the orthogonal projection. (The subbundle $\Lambda^+$ and the projection $\wp$ are defined canonically by the conformal structure.)

Now let $\mathfrak{g}$ denote an $\iota$-invariant metric on $\mathbb{T}_K$ that defines the locally conformally flat structure $[\mathfrak{g}_R]$. If $\mathfrak{h}$ is a symmetric, traceless section of $T^*\mathbb{T}_K \otimes T^*\mathbb{T}_K$ and if $t \in \mathbb{R}$ has sufficiently small norm, then $\mathfrak{g} + t\mathfrak{h}$ is also a metric on $\mathbb{T}_K$. The self-dual Weyl curvature of this metric is denoted in below by $\mathcal{W}_+[\mathfrak{g}+t\mathfrak{h}]$. Since the $\mathfrak{g}$-orthogonal projection of $\mathcal{W}_+[\mathfrak{g}+t\mathfrak{h}]$ to $S$ varies smoothly as a function of t, a section of $S$ to be denoted by $\mathcal{L}_\mathfrak{g}\mathfrak{h}$ is defined by the rule

$$\mathcal{L}_\mathfrak{g}\mathfrak{h} = (\tfrac{d}{dt}(\wp \cdot \mathcal{W}^+[\mathfrak{g}+t\mathfrak{h}]))|_{t=0} \ .$$

(5.1)

The assignment of $\mathfrak{h}$ to $\mathcal{L}_\mathfrak{g}\mathfrak{h}$ defines a second order differential operator that maps symmetric, traceless sections of $T^*X \otimes T^*X$ to symmetric, traceless sections of $\Lambda^+ \otimes \Lambda^+$.

The operator whose kernel appears in Theorem A of [KS] is the formal, $L^2$ adjoint of $\mathcal{L}_\mathfrak{g}$. This formal $L^2$ adjoint is denoted by $\mathcal{L}_\mathfrak{g}^\dagger$. This definition characterizes $\mathcal{L}_\mathfrak{g}^\dagger \mathcal{X}$ for a given section, $\mathcal{X}$, of $S$ by the following requirement:

$$\int_X \langle \mathcal{L}_\mathfrak{g}^\dagger \mathcal{X}, \mathfrak{h} \rangle \, \mathrm{dvol}_\mathfrak{g} = \int_X \langle \mathcal{X}, \mathcal{L}_\mathfrak{g}\mathfrak{h} \rangle \, \mathrm{dvol}_\mathfrak{g}$$

(5.2)

when $\mathfrak{h}$ is a symmetric section of $T^*X \otimes T^*X$. In this equation and henceforth, what is denoted by $\langle \, , \, \rangle$ signifies the fiberwise inner product on the relevant vector bundle. Since the composite operator $\mathcal{L}_\mathfrak{g}\mathcal{L}_\mathfrak{g}^\dagger$ is elliptic, the kernel of $\mathcal{L}_\mathfrak{g}^\dagger$ is a finite dimensional vector space of smooth sections of $S$.

There is a canonical identification between the respective kernels of $\mathcal{L}_\mathfrak{g}$ and $\mathcal{L}_{\mathfrak{g}'}$ when $\mathfrak{g}$ and $\mathfrak{g}'$ are conformal. To give this identification, write $\mathfrak{g}'$ as $\phi^2 \mathfrak{g}$ with $\phi$ being a positive function on X. The promised identification is as follows:

*If $\mathcal{X} \in \mathrm{kernel}(\mathcal{L}_\mathfrak{g}^\dagger)$, then $\phi^2 \mathcal{X} \in \mathrm{kernel}(\mathcal{L}_{\mathfrak{g}'}^\dagger)$* .

(5.3)

As explained in the next paragraph, this is a direct consequence of (5.1) and (5.2) and the fact that $\mathcal{W}_+[\phi^2 \mathfrak{g}]$ and $\phi^2 \mathcal{W}_+[\mathfrak{g}]$ are equal as sections of $S$. It also follows from (5.3) that the $L^2$ norm of $\mathcal{X}$ as measured by the metric $\mathfrak{g}$ is the same as the $L^2$ norm of $\phi^2 \mathcal{X}$ as measured by the metric $\mathfrak{g}' = \phi^2 \mathfrak{g}$.

To explain the preceding remarks, suppose for the moment that $\{\omega^a\}_{a=1,2,3}$ is a local, $\mathfrak{g}$-orthonormal frame for $\Lambda^+$. Then $\mathcal{W}_+[\mathfrak{g}]$ has components $\{\mathcal{W}_+^{ab}\}_{a,b=1,2,3}$ with respect to the induced frame for $\Lambda^+ \otimes \Lambda^+$. Meanwhile, the set $\{\omega'^a = \phi^2 \omega^a\}_{a=1,2,3}$ is a $\mathfrak{g}' = \phi^2 \mathfrak{g}$ orthonormal frame for $\Lambda^+$; and $\{\mathcal{W}_+'^{ab} = \phi^{-2} \mathcal{W}_+^{ab}\}_{a,b=1,2,3}$ are the components of $\mathcal{W}_+[\mathfrak{g}']$ with respect to the corresponding basis $\{\omega'^a \otimes \omega'^b\}_{a,b \in \{1,2,3\}}$. Also, the volume 4-form of the $\mathfrak{g}'$ metric is $\phi^4$ times that of the $\mathfrak{g}$-metric. With the preceding understood, suppose that $\{\mathcal{X}'^{ab}\}_{a,b=1,2,3}$ are the components of an element in the kernel of $\mathcal{L}_{\mathfrak{g}'}^\dagger$ when written using the basis $\{\omega'^a \otimes \omega'^b\}_{a,b \in \{1,2,3\}}$. Then the integration by parts identity in (5.2) with what was just said about the respective $\mathfrak{g}$ and $\mathfrak{g}'$ volume forms and self dual Weyl curvatures implies that $\mathcal{X}^{ab} = \phi^2 \mathcal{X}'^{ab}$ are the components of an element in the kernel of $\mathcal{L}_\mathfrak{g}^\dagger$ when written using the basis $\{\omega^a \otimes \omega^b\}_{a,b \in \{1,2,3\}}$. This is to say that $\mathcal{X}'$ and $\phi^2 \mathcal{X}$ are equal as sections of $S$. These observations also imply that the $L^2$ norm of $\mathcal{X}$ as measured using the metric $\mathfrak{g}$ is the same as the $L^2$ norm of $\mathcal{X}'$ as measured using the metric $\mathfrak{g}'$.

The involution $\iota$ defines an action of the group $\{\pm 1\}$ on the kernel of $\mathcal{L}_\mathfrak{g}^\dagger$ by virtue of the fact that $\mathfrak{g}$ is $\iota$-invariant. Let $\mathcal{H}_{K,R} \subset \text{kernel}(\mathcal{L}_\mathfrak{g}^\dagger)$ denote the $\iota$-invariant subspace. The obstruction vector space that is supplied by Theorem A in [KS] is the vector space $\mathcal{H}_{K,R}$ with the identification in (5.3) implicit. (See also Theorem B in [KS] and Theorems 4.12 and 4.13 in [KS].)

The space $\mathcal{H}_{K,R}$ can, in principle, depend on R. The space $\mathbb{V}$ that appears in Proposition 4.1 is a finite dimensional, R independent vectors space which is the direct sum of spaces $\mathcal{H}_T$ and $\mathcal{H}_K$ that are defined in the upcoming Propositions 5.3 and 5.4. These propositions also describe an injective homomorphism from $\mathcal{H}_{K,R}$ to $\mathbb{V} = \mathcal{H}_T \oplus \mathcal{H}_K$. This homomorphism is used in the subsequent sections to view the obstruction vector from Theorem A of [KS] as a vector in the fixed vector space $\mathbb{V}$.

**b) A 'Mayer-Vietoris' decomposition of $\mathcal{H}_{K,R}$**

Part 1 identified Proposition 4.1's obstruction space with the $\iota$-invariant kernel of $\mathcal{L}_\mathfrak{g}^\dagger$ with $\mathfrak{g}$ being any chosen $\iota$-invariant metric $\mathfrak{g}$ on $\mathbb{T}_K$ that defines the conformal structure $[\mathfrak{g}_R]$. This subsection writes elements in $\mathcal{H}_{K,R}$ so as to exploit the decomposition of $\mathbb{T}_K$ as the union of a part of $\mathbb{T}$ and $(S^3 - N_K) \times S^1$ and $\iota((S^3 - N_K) \times S^1)$. To start, reintroduce the $\mathbb{R}^4/(2\pi\mathbb{Z}^4)$ coordinates $(t_1, t_2, t_3, t_4)$ for $\mathbb{T}$; then use these coordinates to write $\mathbb{T}$ as the product $T \times T'$ with $T$ having the coordinates $(t_1, t_2)$ and $T'$ having the coordinates $(t_3, t_4)$. Let $D \subset T'$ denote the disk where $(t_3-t_*)^2+(t_4-t_*)^2 < e^{-2R}$; and use $\iota(D)$ to denote the reflected disk in $T'$ centered at $(t_3 = -t_*, t_4 = -t_*)$.

With $D$ and $\iota(D)$ defined, write the conformal structure $[\mathfrak{g}_R]$ on $\mathbb{T}_K$ using the decomposition of $\mathbb{T}_K$ as the union of three open sets, the $(S^3-N_K) \times S^1$ and $\iota((S^3-N_K) \times S^1)$ parts of $\mathbb{T}_K$, and the $T \times (T'-(D \cup \iota(D)))$ part of $\mathbb{T}_K$. By way of a reminder, the metric that defines $[\mathfrak{g}_R]$ on $(S^3-N_K) \times S^1$ is the product of the hyperbolic metric with sectional curvature on $S^3-N_K$ with the Euclidean metric on $S^1$ that comes via the identification of $S^1$ with $\mathbb{R}/2\pi\mathbb{Z}$. This product metric was previously denoted by $g_K$. The metric that defines the conformal structure $[\mathfrak{g}_R]$ on $\iota((S^3-N_K) \times S^1)$ is the pull-back metric $\iota^*(g_K)$. And, the metric that defines the conformal structure $[\mathfrak{g}_R]$ on $T \times (T'-(D \cup \iota(D)))$ is the flat metric $g_T$ that is depicted in (3.6).

Use $\mathcal{L}_{g_K}^\dagger$ and $\mathcal{L}_{g_T}^\dagger$ in what follows to denote the respective $\mathfrak{g} = g_K$ and $\mathfrak{g} = g_T$ versions of the operator $\mathcal{L}_\mathfrak{g}^\dagger$. (The operator $\mathcal{L}_{g_K}^\dagger$ obeys the $(S^3-N_K) \times S^1$ version of (4.3) with $\mathfrak{g} = g_K$ when $\mathfrak{h}$ has compact support on $(S^3-N_K) \times S^1$.) Now suppose that $\mathfrak{g}$ is an $\iota$-invariant metric on $\mathbb{T}_K$ that defines the conformal structure $[\mathfrak{g}_R]$. Then there is a positive function $\phi_T$ on the $T \times (T'-(D \cup \iota(D)))$ part of $\mathbb{T}_K$ such that $\mathfrak{g} = \phi_T^2 g_T$ on this same part of $\mathbb{T}_K$. There is also a positive function $\phi_K$ on the $(S^3-N_K) \times S^1$ part of $\mathbb{T}_K$ such that $\mathfrak{g} = \phi_K^2 g_K$ on the $(S^3-N_K) \times S^1$ part of $\mathbb{T}_K$ and such that $\mathfrak{g} = \iota^*(\phi_K^2 g_K)$ on the $\iota((S^3-N_K) \times S^1)$ part of $\mathbb{T}_K$. Therefore, if $\mathcal{X}$ is an $\iota$-invariant element in the kernel of $\mathcal{L}_\mathfrak{g}^\dagger$, then $\mathcal{L}_T^\dagger(\phi_T^{-2}\mathcal{X}) = 0$ on the $T \times (T'-(D \cup \iota(D)))$ part of $\mathbb{T}_K$ and $\mathcal{L}_K^\dagger(\phi_K^{-2}\mathcal{X}) = 0$ on the $(S^3-N_K) \times S^1$ part of $\mathbb{T}_K$. There is a converse to this last statement: Suppose that $\mathcal{X}_T$ and $\mathcal{X}_K$ are respective sections of $S$ on $T \times (T'-(D \cup \iota(D)))$ and $(S^3-N_K) \times S^1$ that obey the following:

- $\mathcal{X}_T$ *is an $\iota$-invariant element in the kernel of* $\mathcal{L}_{g_T}^\dagger$.
- $\mathcal{X}_K$ *is in the kernel of* $\mathcal{L}_{g_K}^\dagger$.
- $\mathcal{X}_T = \phi_T^{-2}\phi_K^2 \mathcal{X}_K$ *where the $T \times (T'-(D \cup \iota(D)))$ and $(S^3-N_K) \times S^1$ parts of $\mathbb{T}_K$ intersect.*

(5.4)

Under these circumstances, the section $\mathcal{X}$ of $S$ defined on the whole of $\mathbb{T}_K$ by $\mathcal{X} = \phi_T^2 \mathcal{X}_T$ on $T \times (T'-(D \cup \iota(D)))$, by $\mathcal{X} = \phi_K^2 \mathcal{X}_K$ on $(S^3-N_K) \times S^1$ and by $\mathcal{X} = \iota^*(\phi_K^2 \mathcal{X}_K)$ on $\iota((S^3-N_K) \times S^1)$ is in the kernel of $\mathcal{L}_\mathfrak{g}^\dagger$.

To be slightly less abstract about the third bullet of (5.4), let $\mathbb{D}$ denote the disk in the $(t_3, t_4)$ torus $T'$ where $(t_3-t_*)^2 + (t_4-t_*)^2 < \frac{1}{16} t_*^2$. (What is denoted by N in Section 2a is $T \times \mathbb{D}$.) The $T \times (T'-(D \cup \iota(D)))$ part of $\mathbb{T}_K$ intersects the $(S^3-N_K) \times S^1$ part of $\mathbb{T}_K$ as the $T \times (\mathbb{D}-D)$ part of the former and as the $s > 0$ part of the latter. The identification of T×(𝔻−D) with the s > 0 part of $(S^3-N_K) \times S^1$ is implemented using the map depicted in (2.1) keeping in mind that the $(t_3, t_4)$ coordinates on 𝔻−D are written using the functions $(\rho, \theta_N)$ as $(t_3 = t_* + \rho \cos\theta_N, t_4 = t_* + \rho \sin\theta_N)$. This identification map, viewed as a map from the s ≥ 0 part of $(S^3-N_K) \times S^1$ to T×(𝔻−D) is denoted below by $\psi$. What is denoted by $\mathcal{X}_T$ in (5.4) restricts to a section of the $\mathfrak{g}_T$ version of S over T×(𝔻−D) and what is denoted by $\mathcal{X}_K$ in (5.4) restricts to a section of the $\mathfrak{g}_K$ version of S over the s > 0 part of $(S^3-N_K) \times S^1$. The pull-back via $\psi$ of the former version of S is canonically isomorphic to the $\mathfrak{g}_K$ version. With this identification understood, the third bullet in (5.4) requires that

$$\mathcal{X}_K = e^{-2(s-R)} \psi^* \mathcal{X}_T$$

(5.5)

hold as an equivalence between section of S over the s > 0 part of $(S^3-N_K) \times S^1$.

*Part 3*: The purpose of writing an element from the kernel of $\mathcal{L}_\mathfrak{g}^\dagger$ as a pair of elements, one from the kernel of $\mathcal{L}_{\mathfrak{g}_T}^\dagger$ on T × (T′−(D∪ι(D))) and one from the kernel of $\mathcal{L}_{\mathfrak{g}_K}^\dagger$ on $(S^3-N_K) \times S^1$ is as follows: The Fourier transform can be used to study the kernels of $\mathcal{L}_{\mathfrak{g}_T}^\dagger$ and $\mathcal{L}_{\mathfrak{g}_K}^\dagger$ because both T × (T′−(D∪ι(D))) and $(S^3-N_K) \times S^1$ have isometric Abelian lie group actions. In the former case, the action is that of the torus T on the T factor; and in the latter case, the action is that of $S^1$ on the $S^1$ factor. The appendix of this article explains how these group actions with (5.5) can be used to describe the ι-invariant part of the kernel of $\mathcal{L}_\mathfrak{g}^\dagger$ when R is sufficiently large. (Note that the T and $S^1$ actions are mutually orthogonal on the intersection of their two domains, which is 𝔻−D; and as a consequence there is no circle action on the whole of $\mathbb{T}_K$.) The upcoming Propositions 5.3 and 5.4 summarize what is needed from the appendix.

### c) The $\mathcal{X}_T$ part of an ι-invariant element in the kernel of $\mathcal{L}_\mathfrak{g}^\dagger$

Suppose that K is a hyperbolic knot and that R has been chosen to define $\mathbb{T}_K$ and its the conformal structure $[\mathfrak{g}_R]$. Let $\mathfrak{g}$ denote an ι-invariant metric from this conformal structure, and let $\mathcal{X}$ is an ι-invariant element in the kernel of $\mathcal{L}_\mathfrak{g}^\dagger$. The three parts of this subsection describe what is denoted by $\mathcal{X}_T$ in (5.4).

*Part 1*: The story begins with a description of the T action on the domain of $\mathcal{L}_{\mathfrak{g}_T}^\dagger$. To describe this action and to set the stage for the subsequent observations, let $\{e^1, e^2\}$ denote a constant, oriented orthonormal basis for the metric

$$\mathfrak{m}_T = e^{-2R} \sum_{a,b=1,2} \mathfrak{m}_{ab}\, dt_a dt_b$$

(5.6)

on the $(t_1, t_2)$ torus T. This basis with the basis $\{e^3 = dt_3, e^4 = dt_4\}$ define an oriented $g_T$-orthonormal and covariantly constant basis for $T^*(T \times T')$. Use this basis to construct a covariantly constant basis for any given tensor bundle on $T \times T'$ constructed from the tangent and cotangent spaces of $T \times T'$. In particular, this basis induces one for $\Lambda^+$, $\Lambda^-$ and for tensor bundles that are constructed from them. The basis 2-forms for $\Lambda^+$ and $\Lambda^-$ are as follows:

- *For* $\Lambda^+$: $\omega^1 = \frac{1}{\sqrt{2}}(e^1 \wedge e^4 + e^2 \wedge e^3)$   $\omega^2 = \frac{1}{\sqrt{2}}(e^2 \wedge e^4 + e^3 \wedge e^1)$   $\omega^3 = \frac{1}{\sqrt{2}}(e^3 \wedge e^4 + e^1 \wedge e^2)$.
- *For* $\Lambda^-$: $\omega^1 = \frac{1}{\sqrt{2}}(e^1 \wedge e^4 - e^2 \wedge e^3)$   $\omega^2 = \frac{1}{\sqrt{2}}(e^2 \wedge e^4 - e^3 \wedge e^1)$   $\omega^3 = \frac{1}{\sqrt{2}}(e^3 \wedge e^4 - e^1 \wedge e^2)$.

(5.7)

A section of $\Lambda^+ \otimes \Lambda^+$ can be written using the corresponding basis $\{\omega^a \otimes \omega^b\}_{a,b \in \{1,2,3\}}$; if $\mathcal{X}$ denotes a section, then it is written as $\mathcal{X} = \mathcal{X}^{ab} \omega^a \otimes \omega^b$ with $\{\mathcal{X}^{ab}\}_{a,b \in \{1,2,3\}}$ being functions on $T \times (T' - (D \cup \iota(D)))$. The section $\mathcal{X}$ is $\iota$-invariant if and only if all of the coefficient functions $\{\mathcal{X}^{ab}\}_{a,b \in \{1,2,3\}}$ are $\iota$-invariant functions.

*Part 2*: The kernel of $\mathcal{L}_{g_T}^\dagger$ on the whole 4-torus $T \times T'$ is the span of the constant sections of $S$ with respect to the basis $\{\omega^a \otimes \omega^b\}_{a,b \in \{1,2,3\}}$. This is a 5-dimensional vector space because $\mathcal{X}$ must a priori be symmetric and traceless. Note in this regard that all of the constant sections of $S$ are $\iota$-invariant. These constant sections restrict to give elements in the kernel of $\mathcal{L}_{g_T}^\dagger$ on $T \times (T' - (D \cup \iota(D)))$; and the four dimensional subspace where $\mathcal{X}^{33} = 0$ (and $\mathcal{X}^{11} = -\mathcal{X}^{22}$ because $\mathcal{X}$ is traceless) plays a role in the subsequent story.

These constant sections of $S$ are part of a larger subspace in the kernel of $\mathcal{L}_{g_T}^\dagger$ on $T \times (T' - (D \cup \iota(D)))$ that also plays a role. To say more, let $p_*$ denote the point $(t_*, t_*) \in T'$. This second subspace is the restriction to the domain $T \times (T' - (D \cup \iota(D)))$ of elements in the kernel of $\mathcal{L}_{g_T}^\dagger$ on $T \times (T' - \{p_*, \iota(p_*)\})$ whose coefficient functions $\{\mathcal{X}^{ab}\}_{a,b \in \{1,2,3\}}$ are independent of the $(t_1, t_2)$ coordinates on T and have first order poles at the points $p_*$ and $\iota(p_*)$ as functions of the coordinates $(t_3, t_4)$ on T'. This larger subspace is 7 dimensional. Proposition 5.3 in Part 3 of this subsection describes this space.

To set the stage for Proposition 5.3, let $\iota$ now denote the involution on T' that acts as $(t_3, t_4) \to (-t_3, -t_4)$. Introduce a complex structure on T' by declaring that $dt_3 + i dt_4$ span the holomorphic cotangent bundle $T^{1,0}$. Two lemmas about meromorphic functions on T' are needed in what follows.

**Lemma 5.1**: *An $\iota$-invariant, meromorphic function on $T'$ with poles of order 1 at the points $p_*$ and $-p_*$ can be written as $\alpha_0 + \alpha_1 x$ with $(\alpha_0, \alpha_1) \in \mathbb{C}^2$ and with $x$ having poles of order 1 at $p_*$ and $\iota(p_*)$ with residue 1, and a zero of order 2 at the point $(t_3 = 0, t_4 = 0)$.*

*Proof of Lemma 5.1*: This follows from the basic theory of elliptic functions (see for example [A]). (With $T'$ viewed as $\mathbb{R}^2/(2\pi\mathbb{Z}^2)$, the function $x$ can be written in terms of a Weierstrasse function $\wp$ as $x = \frac{\alpha}{\wp - \beta}$ for suitable $\alpha, \beta \in \mathbb{C}$.)

By way of a parenthetical remark for now, any $\iota$-invariant, meromorphic function on $T'$ that is analytic on the complement of $\{p_*, \iota(p_*)\}$ and has poles of order at most $n$ at $p_*$ and $\iota(p_*)$ can be written as $a_0 + a_1 x + \cdots + a_n x^n$ with $(a_0, \ldots, a_n) \in \mathbb{C}^{n+1}$. That this is so can be proved by induction: Suppose that $X$ is such a meromorphic function, with a pole of order $n$. Then one can find $a_n \in \mathbb{C} - 0$ so that $X - a_n x^n$ has a pole of order less than $n$ at $p_*$. It is also $\iota$ invariant, so it has poles of order less than $n$ at $p_*$ and $\iota(p_*)$

The next lemma concerns the derivative of the function $x$ in Lemma 5.1.

**Lemma 5.2**: *There exists a unique, $\iota$-invariant function on $T' - \{p_*, \iota(p_*)\}$ that obeys the bulleted conditions that follow. The function is denoted by $u$.*
- $\bar{\partial} u = \partial x$.
- $\int_{T'} u = 0$.
- *Let $z$ denote a local holomorphic coordinate for a neighborhood of $p_*$ obeying $z(p_*) = 0$ and $dz = dt_3 + i\, dt_4$. Then $u$ near $p_*$ obeys $u = \frac{1}{z^2} \bar{z} + \mathcal{O}(1)$.*

*Proof of Lemma 5.2*: Recall that $\chi$ is a smooth, non-increasing function on $\mathbb{R}$ that equals 1 on $(-\infty, \frac{1}{4}]$ and equals zero on $[\frac{3}{4}, \infty)$. Fix $\delta < \frac{1}{100} t_*$ and let $\chi_\delta$ now denote the function of the coordinate $z$ given by $\chi(2|z|/\delta - 1)$. This function equals 1 where $|z| \leq \frac{1}{2}\delta$ and it equals zero where $|z| \geq \delta$. Define $u_0$ to be the function on $T'$ that is equal to 0 on the complement of $\mathbb{D}$ and $\iota(\mathbb{D})$; on $\mathbb{D}$ set it equal to $\chi_\delta \bar{z}\, \partial x$ and on $\iota(\mathbb{D})$, set it equal to the $\iota$-pull-back of $\chi_\delta \bar{z}\, \partial x$. This makes $u_0$ to be $\iota$-invariant. Write $u = u_0 + u_1$. The top bullet of the lemma is obeyed if and only if $u_1$ obeys

$$\bar{\partial} u_1 = -(\bar{\partial} \chi_\delta)\, \bar{z}\, \partial x - \iota^*((\bar{\partial}\chi_\delta)\, \bar{z}\, \partial x) + (1 - \chi_\delta - \iota^*\chi_\delta)\, \partial x\,.$$

(5.8)

The function on the right hand side of (5.8) has compact support in $T' - \{p_*, \iota(p_*)\}$. This understood, it can be viewed as a smooth function on the whole of $T'$. Since the operator $\bar\partial$ on T´ is Fredholm with index 0 and since its kernel and also cokernel are the constants, there is a unique solution to (5.8) with integral zero on T´ if the integral over T´ of the function on the right hand side of (5.8) is zero. This is the case because the function on the right hand side of (5.8) changes sign when pulled back by $\iota$. With $u_1$ obeying (5.8) and having integral zero on T´, then $u = u_0 + u_1$ obeys the three bullets of the lemma.

*Part 3*: Suppose that $\mathfrak{g}$ is an $\iota$-invariant metric in $[\mathfrak{g}_R]$. Proposition 5.3 which is stated directly uses the functions $x$ and $u$ to describe the $\mathcal{X}_T$ part of any $\iota$-invariant element in the kernel on $\mathbb{T}_K$ of $\mathcal{L}_\mathfrak{g}^\dagger$. To set the stage for the proposition, suppose that $\mathfrak{X}$ is a given element in the metric $g_T$ version $\Lambda^+ \otimes \Lambda^+$ on $T \times (T´ - \{p_*, \iota(p_*)\})$. The basis in (5.7) for $\Lambda^+$ is used in the upcoming (5.9) to write the components of $\mathfrak{X}$ as $\{\mathfrak{X}^{ab}\}_{a,b \in \{1,2,3\}}$.

Let $\mathcal{H}_T$ denote the (real) seven dimensional vector space $\mathbb{R} \times \mathbb{C}^3$. A canonical incarnation of $\mathcal{H}_T$ as a subspace of symmetric, traceless sections of $\Lambda^+ \otimes \Lambda^+$ over $T \times (T´ - \{p_*, \iota(p_*)\})$ is defined as follows: This incarnation is spanned by elements with components $\{\mathfrak{X}^{ab}\}_{a,b \in \{1,2,3\}}$ that can be written using a real number $q$ and three complex numbers $c, \flat_0, \flat_1$ as

- $\mathfrak{X}^{33} = 0$ *and* $\mathfrak{X}^{11} = -\mathfrak{X}^{22}$,
- $\mathfrak{X}^{13} - i\mathfrak{X}^{23} = c + q\mathfrak{z}_K x$,
- $\mathfrak{X}^{11} - \mathfrak{X}^{22} - 2i\mathfrak{X}^{12} = \flat_0 + \flat_1 x + \overline{\flat}_1 u$.

(5.9)

The incarnation of $\mathcal{H}_T$ depicted in (5.9) as a subspace of sections of $S$ is in the $\iota$-invariant kernel of $\mathcal{L}_{g_T}^\dagger$ on $T \times (T´ - \{p_*, \iota(p_*)\})$. This is a straightforward calculation, but in any event, it is proved in the Section A2 of the appendix. The initial incarnation of $\mathcal{H}_T$ as $\mathbb{R} \times \mathbb{C}^3$ and its subsequent incarnation as the vector space depicted by (5.9) will not be notationally distinguished in what follows; the symbol $\mathcal{H}_T$ will be used for both incarnations.

There is one last item of notation for Proposition 5.3: This proposition uses $\Delta$ to denote the T invariant function on $T \times T´$ giving the distance on the T´ factor to $D \cup \iota(D)$.

**Proposition 5.3**: *Given a hyperbolic knot* K, *there exists* $\kappa > 100$ *and a norm 1 complex number denoted by* $\mathfrak{z}_K$ *whose significance is described directly. Use* $\mathfrak{z}_K$ *in (5.9) to define the vector space* $\mathcal{H}_T$. *Fix* $R > \kappa^2$ *to define* $\mathbb{T}_K$ *and its locally conformally flat structure* $[\mathfrak{g}_R]$, *and let* $\mathfrak{g}$ *denote an* $\iota$-*invariant metric in* $[\mathfrak{g}_R]$. *There is a homomorphism from the* $\iota$ *invariant kernel of* $\mathcal{L}_\mathfrak{g}^\dagger$ *on* $\mathbb{T}_K$ *to* $\mathcal{H}_T$, *this denoted by* $\mathcal{Q}_T$, *with the following property: Let* $X$ *denote an* $\iota$-*invariant element in the kernel of* $\mathcal{L}_\mathfrak{g}^\dagger$ *on* $\mathbb{T}_K$. *Write* $X$ *as a pair* $(X_T, X_K)$ *in the manner of (5.4) and (5.5). Let* $z_T$ *denote the* $L^2$ *norm of* $X_T$ *on the part of* $T \times T´$ *where the distance to both* $p_*$ *and* $\iota(p_*)$ *is greater than* $\frac{1}{4} t_* e^{-\kappa}$; *and let* $z_K$ *denote the* $L^2$ *norm of* $\mathcal{X}_K$ *where* $s \leq \kappa$ *on* $(S^3 - N_K) \times S^1$. *Assume that* $z_T + z_K = 1$. *The* $\mathcal{X}_T$ *part of* $\mathcal{X}$ *can be written as* $\mathcal{Q}_T(\mathcal{X}) + \mathfrak{x}$ *with the* $g_T$ *norm of* $\mathfrak{x}$ *obeying* $|\mathfrak{x}| \leq \kappa(\Delta^{-2} e^{-2R} + e^{-R})$ *where* $\Delta \geq \kappa e^{-R}$.

This proposition is proved in Section A5 of the appendix. Section A1 of the appendix contains an outline of the proof.

### d) The $\mathcal{X}_K$ part of an $\iota$-invariant element in the kernel of $\mathcal{L}_\mathfrak{g}^\dagger$

Proposition 5.4 (which is stated momentarily) describes the part of the $\iota$-invariant kernel of $\mathcal{L}_\mathfrak{g}^\dagger$ that is not visible in $\mathcal{H}_T$ via the homorphism $\mathcal{Q}_T$. To set the stage for this new proposition, introduce $\mathbb{R}^{3,1}$ to denote the 4-dimensional Minkowski space. This is $\mathbb{R}^4$ as a vector space, but with the signature (-, +, +, +) metric that is written using coordinates $(y_0, y_1, y_2, y_3)$ for $\mathbb{R}^4$ as

$$-dy_0 \otimes dy_0 + dy_1 \otimes dy_1 + dy_2 \otimes dy_2 + dy_3 \otimes dy_3.$$
(5.10)

The universal cover of $S^3 - N_K$ can be viewed as the $y_0 > 0$ part of the hyperbola where the coordinates obey $y_0^2 = 1 + y_1^2 + y_2^2 + y_3^2$. This component of the hyperbola is denoted below by $\mathbb{H}^+$. The fundamental group of $S^3 - N_K$ acts on $\mathbb{H}^+$ as a discrete subgroup of the Lie group SO(3,1) which is the group of isometries of $\mathbb{R}^{3,1}$ that fix the origin. (The action of the latter group on $\mathbb{R}^{3,1}$ maps $\mathbb{H}^+$ to itself.) Let $\Gamma$ denote for the moment the fundamental group of $S^3 - N_K$.

With the preceding understood, the quotient $\mathbb{H}^+ \times_\Gamma \mathbb{R}^{3,1}$ defines a flat $\mathbb{R}^{3,1}$ bundle over $S^3 - N_K$. Denote this bundle by V and let $H^1(S^3 - K; V)$ denote the first (Čech) cohomology of $S^3 - K$ with coefficients in the bundle V. The inclusion $\hat{\imath}: N_K - K \to S^3 - K$ defines by pull-back a homomorphism $\hat{\imath}^*: H^1(S^3 - K; V) \to H^1(N_K - K; V)$.

To continue the stage setting, reintroduce the metric $g_K$ on $(S^3 - K) \times S^1$ that is given by the product of the constant sectional curvature -1 metric on $S^3 - K$ and the Euclidean metric $d\theta^2$ on $S^1$. This metric $g_K$ is used in Proposition 5.4 to define the bundle $\Lambda^+ \otimes \Lambda^+$ over $(S^3 - N_K) \times S^1$ and its subbundle of symmetric, traceless elements and the operator $\mathcal{L}_{g_K}^\dagger$ on the space of sections of this subbundle.

Proposition 5.4 refers to a norm on the vector space $\mathcal{H}_T$ that is depicted in (5.9) This norm is defined as follows: The norm of any given element $\mathcal{Q} = (q, c, \mathfrak{b}_0, \mathfrak{b}_1) \in \mathcal{H}_T$ is taken to be $(q^2 + |c|^2 + |\mathfrak{b}_0|^2 + |\mathfrak{b}_1|^2)^{1/2}$. This is denoted by $|\mathcal{Q}|_T$. By way of some more notation, the Proposition 5.4 uses $\kappa_T$ to denote the version of the number $\kappa$ that appears in Proposition 5.3.

**Proposition 5.4**: *Given a hyperbolic knot K, there exists $\kappa > \kappa_T$ and a subspace of square integrable elements in the kernel of $\mathcal{L}_{g_K}^\dagger$ on $(S^3-K) \times S^1$ whose significance is described directly. Use $\mathcal{H}_K$ to denote this subspace in the kernel of $\mathcal{L}_{g_K}^\dagger$.*

- *The vector space $\mathcal{H}_K$ is naturally isomorphic to the kernel of $\hat{\imath}^*$ in $H^1(S^3-K; V) \otimes_\mathbb{R} \mathbb{C}$.*
- *Fix $R > \kappa$ to construct $\mathbb{T}_K$ and its conformal structure $[g_R]$. Let $g$ denote a given metric from this conformal structure. There is a homomorphism from the $\iota$ invariant kernel of $\mathcal{L}_g^\dagger$ on $\mathbb{T}_K$ to $\mathcal{H}_K$, this denoted by $Q_K$, with the following property: Let $X$ denote an $\iota$-invariant element in the kernel of $\mathcal{L}_g^\dagger$ on $\mathbb{T}_K$. Write $X$ as a pair $(X_T, X_K)$ in the manner of (5.4) and (5.5). Let $z_T$ denote the $L^2$ norm of $X_T$ on the part of $T \times T'$ where the distance to both $p_*$ and $\iota(p_*)$ is greater than $\frac{1}{4} t_* e^{-\kappa}$; and let $z_K$ denote the $L^2$ norm of $X_K$ where $s \leq \kappa$ on $(S^3-N_K) \times S^1$. Assume that $z_T + z_K = 1$. The $X_K$ part of $X$ can be written as $Q_K(X) + \mathfrak{x}$ with the $g_K$ norm of $\mathfrak{x}$ on the $s \leq \frac{1}{2} R$ part of $(S^3-N_K) \times S^1$ obeying the bound $|\mathfrak{x}| \leq \kappa e^{2s-R}$*

This proposition is also proved in Section A5 of the appendix. Section A1 of the appendix contains an outline of the proof.

The following is an immediate corollary to Propositions 5.3 and 5.4:

**Corollary 5.5**: *Given a hyperbolic knot K, there exists $\kappa > 100$ with the following significance: If $R \geq \kappa$, then the conclusions of Propositions 5.3 and 5.4 are valid. Moreover, if $g$ is a metric on $\mathbb{T}_K$ from the conformal class $[g_R]$, then Proposition 5.3's homomorphism $Q_T$ from the $\iota$ invariant kernel of $\mathcal{L}_g^\dagger$ to $\mathcal{H}_T$ and Proposition 5.4's homorphims $Q_K$ from the $\iota$ invariant kernel of $\mathcal{L}_g^\dagger$ to $\mathcal{H}_K$ define an injective homomorphism from the $\iota$ invariant kernel of $\mathcal{L}_g^\dagger$ to $\mathcal{H}_T \oplus \mathcal{H}_K$.*

As explained in Part 4 of Section A3 in the appendix,, the twisted cohomology with coefficients in the bundle V appears in Proposition 5.4 for much the same reason that it appears in Sections 5 and 6 of [AV]: Use of the Fourier transform on $(S^3-N_K) \times S^1$ identifies a part of the $n = \pm 1$ Fourier modes in the kernel of $\mathcal{L}_{g_K}^\dagger$ with the vector space of Codazzi tensors on $S^3-N_K$. These are symmetric, traceless sections of $\otimes^2 T^*(S^3-N_K)$ that are closed with respect to the exterior covariant derivative defined by the Levi-Civita connection when they are viewed as $T^*(S^3-N_K)$ valued sections of $T^*(S^3-N_K)$. A theorem of Lafontaine [Laf] (which is based on an observation of Ferus [Fe]) relates Codazzi tensors to elements of $H^1(S^3-N_K; V)$.

The Codazzi tensor modes in the kernel of $\mathcal{L}_{g_K}^{\dagger}$ are also elements in the kernel of the operator $\mathcal{L}_{g_K}$. The appearance of $H^1(S^3-K;V)$ in the latter guise is explained by a theorem of [Ka] (see also [Sc]): Kapovich observed that the twisted cohomology $H^1(S^3-K;V)$ classifies the first order deformations of the canonical flat $SO(3,1)$ conformal structure of $S^3-N_K$ in the larger conformal group $SO(4,1)$. (Mostow's rigidity theorem asserts in part that the canonical flat conformal structure in $SO(3;1)$ has no non-trivial deformations in $SO(3,1)$.)

The next proposition can be used to obtain an a priori bound on the dimension of Proposition 5.4's vector space $\mathcal{H}_K$.

**Proposition 5.6**: *Given* $N > 0$, *there exists* $\mathcal{V} > 0$ *with the following significance: If the hyperbolic volume of* $S^3-K$ *is less than* $\mathcal{V}$, *then the dimension of* $H^1(S^3-K; V))$ *is less than* $N$. *As a consequence, the dimension of* $\mathcal{H}_K$ *is less than* $N$.

*Proof of Proposition 5.6*: It follows from work of Thurston (see e.g [G]) that the set of finite volume, complete hyperbolic 3-manifolds with a given a priori volume bound is compact in the Gromov-Hausdorff topology. With the preceding understood, suppose that M is a given finite volume, complete hyperbolic 3-manifold and that $\{M_i\}_{i \in \{1,2,...\}}$ is a sequence of such manifolds that converges in the Gromov-Hausdorff topology to M. As explained in Gromov's Bourbaki exposition of Thurston and Jørgensen's work [G], there is a subsequence where the convergence can be described more or less explicitly as a process whereby each $M_i$ for $i \gg 1$ looks very much like what one would obtain from M by gluing some pairs of cusp ends to form a long tube with a central, very short geodesic. This picture of convergence from [G] and a corresponding Mayer-Vietoris sequence can be used to bound the dimension of the large i versions of $H^1(M_i;V)$ in terms of the dimension of $H^1(M;V)$.

The final proposition in this subsection talks about a certain infinite family of hyperbolic knots with kernel($i*$) equal to 0. This family consists of the hyperbolic 2-bridge knot. The definition is as follows: A knot $\mathbb{R}^3 = S^3$−point is said to be a 2-bridge knot when it can be isotoped so that one of the coordinate functions has only 2 maxima and 2 minima on the knot. A 2-bridge knot has hyperbolic knot complement unless it is a torus knot, in which case it must be a $(2, n)$ torus knot (see [Schb], or for those who don't read German, [Schl1].) It follows from the classification of 2-bridge knots (see [Schl2] that infinitely many are not torus knots. Let $\mathcal{K}_2$ denote the latter set of 2-bridge knots.

**Proposition 5.7**: *If* $K \in \mathcal{K}_2$, *then the kernel of* $i*$ *in* $H^1(S^3-K,V)$ *is trivial. Moreover, there exists a positive number* $\nu$ *and an infinite set* $\mathcal{K}_{2,\nu} \subset \mathcal{K}_2$ *such if* $K \subset \mathcal{K}_{2,\nu}$, *then the* *hyperbolic volume of* $S^3$–K *is less than* v. *There is also an infinite set of knots in* $\mathcal{K}_2$ *with no a priori bound on the hyperbolic volume of their complements in* $S^3$.

***Proof of Proposition 5.7***: The assertion about the kernel of $\hat{\imath}*$ in $H^1(S^3$–K; V) was proved by Kapovich in [Ka]. The existence of an a priori upper bound on the hyperbolic volume of the knot complement for an infinite set of knots in $\mathcal{K}_2$ follows from the classification of these knots (see, e.g. [Schl 2]) and the volume upper bound in [Lac]. Schultens [Schl2] attributes the observation that there is no a priori hyperbolic volume of the knot complement for knots in $\mathcal{K}_2$ to work of Hatcher and Thurston.

## 6. The obstruction vector

Proposition 4.1's vector space $\mathbb{V}$ is defined to be the direct sum $\mathcal{H}_T \oplus \mathcal{H}_K$. The upcoming Section 6c defines Proposition 4.1's obstruction vector (a vector in $\mathbb{V}$). The intervening subsections set the stage.

### a) Review of the Kovalev-Singer construction

Let K again denote the given hyperbolic knot in $S^3$. Fix R > 1 so as to construct $\mathbb{T}_K$ and its locally conformally flat structure $[\mathfrak{g}_R]$ on $\mathbb{T}_K$. Let $\mathfrak{g}_R$ denote the ι invariant metric in this conformal equivalence class that is described in Step 1 of Section 4a. Use $\mathcal{H}_{K,R}$ to denote the ι-invariant kernel of the $\mathfrak{g} = \mathfrak{g}_R$ version of $\mathcal{L}_\mathfrak{g}^\dagger$ on $\mathbb{T}_K$. This section reviews how $\mathcal{H}_{K,R}$ comes to be the obstruction space in Theorem A in [KS].

By way of preliminaries, fix a positive number δ and the sixteen SO(4) parameters that are needed to construct the conformal structure $[\mathfrak{g}_{R,\delta}]$ on $X_K$. Introduce the metric $\mathfrak{g}_{R,\delta*}$ from this conformal class; this metric is defined Step 2 Section 4a. Given a positive integer n, fix a set of n distinct points in $X_{K,\delta}$ to be the set $\vartheta$. With $\varepsilon \in (0, 1]$ chosen, fix a set of positive numbers $\{\delta_p\}_{p \in \vartheta}$ subject to the constraints in Section 3f and fix the required SO(4) parameters to construct the conformal structure $[\mathfrak{g}^n_{R,\delta}]$ for $X_K \#_n \overline{\mathbb{CP}}^2$. Let $\mathfrak{g}^n_{R,\delta*}$ denote the metric from this conformal class that is described in Step 3 in Section 4a. Set $\mathfrak{g}_* = \mathfrak{g}_{R,\delta*}$ when n = 0 and set $\mathfrak{g}_* = \mathfrak{g}^n_{R,\delta*}$ when n ≥ 1.

The notation in what follows uses $c_\diamond$ to denote a number greater than 1 that can be determined a priori given the knot K, the number R, the chosen ε and the set $\vartheta$. Its value can be assumed to increase between successive appearances. Supposing that δ and each $p \in \vartheta$ version of $\delta_p$ is less than $c_\diamond^{-1}$, then [KS] construct in Section 4.12 of their paper an isomorphic copy of $\mathcal{H}_{K,R}$ in the space of symmetric, traceless sections of the $\mathfrak{g}_*$ version of $\Lambda^+ \otimes \Lambda^+$ over $X_K \#_n \overline{\mathbb{CP}}^2$. This isomorphic copy is denoted here by $\mathcal{V}_*'$. The four steps that follow review the [KS] definition of $\mathcal{V}_*'$.

Step 1: The complement in $X_K \#_n \overline{\mathbb{CP}}^2$ of the zero sections in the 16 copies of $E^2_{\delta/r}$ and the n copies of $E_{\delta_p/r_p}$ are identified with the complement in $\mathbb{T}_K/\iota$ of the 16 singular points and the points $\vartheta$. Note that this identification implicitly uses the coordinates that appear in (3.10) and (3.11) from Sections 3e) and 3f) to specify the $E^2_{\delta/r}$ and $E_{\delta_p/r_p}$ parts of $X_K \#_n \overline{\mathbb{CP}}^2$.

Step 2: Let $X$ denote a given element in $\mathcal{H}_{K,R}$. Granted the identifications from Step 1, the element $X$ defines gives a section of $(\wedge^2 T^*(X_K \#_n \overline{\mathbb{CP}}^2)) \otimes (\wedge^2 T^*(X_K \#_n \overline{\mathbb{CP}}^2))$ on the complement in $X_K \#_n \overline{\mathbb{CP}}^2$ of the zero sections in the 16 copies of $E^2_{\delta/r}$ and the n copies of $E_{\delta_p/r_p}$.

Step 3: Let $\mathfrak{g}$ denote either the metric $\mathfrak{g}_{R,\delta}$ or $\mathfrak{g}^n_{R,\delta}$ as the case may be. Keep in mind that this is a metric on $X_K$ or $X_K \#_n \overline{\mathbb{CP}}^2$ that is conformal to $\mathfrak{g}_*$. Multiplying $X$ by a suitable cut-off function on $X_K \#_n \overline{\mathbb{CP}}^2$ that is zero near the zero sections in the 16 copies of $E^2_{\delta/r}$ and the n copies of $E_{\delta_p/r_p}$ gives a section of the $\mathfrak{g}$ version of $\Lambda^+ \otimes \Lambda^+$ on the whole of $X_K \#_n \overline{\mathbb{CP}}^2$. This section is denoted by $X_\mathfrak{g}$.

Near the zero section of any given $E^2_{\delta/r}$, the section $X_\mathfrak{g}$ when written with the coordinates that are used in (3.10) is

$$X_\mathfrak{g} = (1 - \chi_{e\delta}) X$$

(6.1)

with $e\delta$ denoting the product of Euler's number $e$ with the number $\delta$. The function $\chi_{e\delta}$ is defined by replacing $\delta$ with $e\delta$ in the first paragraph of Part 2 in Section 3e. Note in this regard that $1 - \chi_{e\delta}$ is equal to 0 where the function $\chi_\delta$ in (3.10) is positive which is why $X_\mathfrak{g}$, a prior a section of $(\wedge^2 T^*(X_K \#_n \overline{\mathbb{CP}}^2)) \otimes (\wedge^2 T^*(X_K \#_n \overline{\mathbb{CP}}^2))$, extends over $E^2_{\delta/r}$ as a section a traceless symmetric section of $\Lambda^+ \otimes \Lambda^+$ as defined by the metric $\mathfrak{g}$.

Near the zero section in any $p \in \vartheta_p$ version of $E_{\delta_p/r_p}$, the section $X_\mathfrak{g}$ when written with the coordinates that are used in (3.11) is

$$X_\mathfrak{g} = (1 - \chi_{e\delta_p}) X .$$

(6.2)

The fact that $1 - \chi_{e\delta_p}$ is zero where $\chi_{2\delta_p} \neq 0$ implies that $X_\mathfrak{g}$ also extends over $E_{\delta_p/r_p}$ as a section of the metric $\mathfrak{g}$ version of the bundle of traceless, symmetric sections of $\Lambda^+ \otimes \Lambda^+$.

Step 4: The metrics $\mathfrak{g}$ and $\mathfrak{g}_*$ are conformal, which is to say that $\mathfrak{g}_* = \phi^2 \mathfrak{g}$ with $\phi$ being a positive function on $X_K \#_n \overline{\mathbb{CP}}^2$. Define $X_{\mathfrak{g}_*}$ to be $\phi^2 X_\mathfrak{g}$. Since the $\mathfrak{g}$ and $\mathfrak{g}_*$ versions of $\Lambda^+ \otimes \Lambda^+$ are the same subbundle of $(\wedge^2 T^*(X_K \#_n \overline{\mathbb{CP}}^2)) \otimes (\wedge^2 T^*(X_K \#_n \overline{\mathbb{CP}}^2))$, this section $\mathcal{X}_{\mathfrak{g}*}$ is therefore a section of $\mathfrak{g}_*$'s version of $\Lambda^+ \otimes \Lambda^+$. The span of the various $\mathcal{X} \in \mathcal{H}_{K,R}$ versions of $\mathcal{X}_{\mathfrak{g}*}$ defines the vector bundle $\mathcal{V}_*'$ from Theorem A of [KS] when $\delta$ and $\{\delta_p\}_{p \in \vartheta}$ are less than $c_\diamond^{-1}$.

**b) Tweaking the definition of $\mathcal{V}_*'$**

Given that $\delta$ and $\{\delta_p\}_{p \in \vartheta}$ are less than $c_\diamond^{-1}$, then [KS] construct a symmetric, traceless section of $\otimes^2 T^*(X_K \#_n \overline{\mathbb{CP}}^2)$, denoted here by $\mathfrak{h}_\mathfrak{g}$, with sup-norm less than $\frac{1}{100}$ as measured by $\mathfrak{g}_*$, and such that $\mathcal{W}_+(\mathfrak{g}_* + \mathfrak{h}_{\mathfrak{g}*})$ is in the vector space $\mathcal{V}_*'$. Since the vector $\mathcal{W}_+(\mathfrak{g}_* + \mathfrak{h}_{\mathfrak{g}*})$ in $\mathcal{V}_*'$ vanishes if and only if $\mathcal{W}_+(\mathfrak{g}_* + \mathfrak{h}_{\mathfrak{g}*}) = 0$, it is tautologically an obstruction vector. This first incarnation of the obstruction vector will ultimately determine Proposition 4.1's obstruction vector in $\mathbb{V}$.

The definition of $\mathcal{W}_+(\mathfrak{g}_* + \mathfrak{h}_{\mathfrak{g}*})$ in $\mathcal{V}_*'$ as an obstruction vector becomes more than a tautology by using Taylor's theorem (with remainder) to write this vector as done in [KS] just prior to their Proposition 5.12. This appeal to Taylor's theorem writes

$$\mathcal{W}_+(\mathfrak{g}_* + \mathfrak{h}_{\mathfrak{g}*}) = \Pi \mathcal{W}_+(\mathfrak{g}_*) + \Pi \mathcal{L}_{\mathfrak{g}_*} \mathfrak{h}_{\mathfrak{g}*} + \mathcal{Q}_{\mathfrak{g}*}$$

(6.3)

where $\Pi$ denotes the $L^2$ orthogonal projection to $\mathcal{V}_*'$ and where $\mathcal{Q}_{\mathfrak{g}*}$ is the remainder term in the Taylor's approximation. The norm of this term is bounded by a quadratic function of $\mathfrak{h}_{\mathfrak{g}*}$ the following sense: Use the $L^2$ inner product for the metric $\mathfrak{g}_*$ to define an inner product and norm on $\mathcal{V}_*'$. This norm is denoted by $\|\cdot\|_2$. Meanwhile, let $\|\cdot\|_6$ denote the $L^6$ norm on the space of symmetric, traceless sections of the $\mathfrak{g}_*$ version of $\Lambda^+ \otimes \Lambda^+$. Then

$$\|\mathcal{Q}_{\mathfrak{g}*}\|_2 \leq c_\diamond \|\mathcal{W}_+(\mathfrak{g}_*)\|_6^2.$$

(6.4)

As it turns out, the definition of $\mathcal{V}_*'$ in [KS] is such that the $\Pi \mathcal{L}_{\mathfrak{g}_*} \mathfrak{h}_{\mathfrak{g}*}$ term in (6.3) has potentially the largest norm of the three terms. (If $\varepsilon = 1$ is used in the definition of $[\mathfrak{g}^n{}_{R,\delta}]$, then the $\mathcal{Q}_{\mathfrak{g}*}$ term can be just as large.) This is unfortunate because $\mathfrak{h}_{\mathfrak{g}*}$ is obtained by solving a non-linear differential equation and so little is known about this section but for bounds on its $L^6{}_2$ Sobolev norm. (This is the sum of the $L^6$ norms of $\mathfrak{h}_{\mathfrak{g}*}$ and its covariant derivatives to second order.) The $L^6{}_2$ norm of $\mathfrak{h}_{\mathfrak{g}*}$ is bounded by $c_\diamond \|\mathcal{W}_+(\mathfrak{g}_*)\|_6$.

To make something of equation in (6.3), it is necessary to tweak the definition of $\mathcal{V}_*{'}$ so as to make the $\Pi W_+(\mathfrak{g}_*)$ by far the largest term in (6.3). The tweaking of the definition of $\mathcal{V}_*{'}$ has two steps. The first step uses the parameter $\varepsilon$ to replace (6.2) with

$$\mathcal{X}_\mathfrak{g} = (1 - \chi_{\varepsilon\delta_p})\mathcal{X}.$$

(6.5)

Where as (6.2) has $\mathcal{X}_\mathfrak{g} = 0$ where the coordinate x obey $|x| \leq \tfrac{1}{2} e \, \delta_p$, the formula in (6.5) has $\mathcal{X}_\mathfrak{g} = 0$ where $|x| \leq \tfrac{1}{2}\varepsilon\delta_p$. Meanwhile, (6.5) has $\mathcal{X}_\mathfrak{g} = \mathcal{X}$ where $|x| \geq \varepsilon\delta_p$. Having replaced (6.2) by (6.5) at each $p \in \vartheta$, the new version of $\mathcal{X}_\mathfrak{g}$ still extends smoothly across the zero sections of each $p \in \vartheta$ version of $E_{\delta_p/r_p}$.

The extension $\mathcal{X}_\mathfrak{g}$ just defined does not lie entirely in the metric $\mathfrak{g}$ version of the bundle $\Lambda^+ \otimes \Lambda^+$. (Even so, its projection in $\otimes^2(\wedge^2 T^*(X_K \#_n \overline{\mathbb{CP}}^2))$ orthogonal to this subbundle is non-zero only near each $p \in \vartheta$ where the relevant version of the coordinate x obeys $\tfrac{1}{2}\varepsilon\delta_p \leq |x| \leq \delta_p$.) The second step in the definition corrects this problem: Introduce $\underline{\mathcal{X}}_\mathfrak{g}$ to denote the $\mathfrak{g}$-orthogonal projection in $\otimes^2(\wedge^2 T^*(X_K \#_n \overline{\mathbb{CP}}^2))$ of $\mathcal{X}_\mathfrak{g}$ onto the subbundle of symmetric, traceless elements in $\Lambda^+ \otimes \Lambda^+$. Again write $\mathfrak{g}_*$ as $\phi^2 \mathfrak{g}$. The span of the various $\mathcal{X} \in \mathcal{H}_{K,R}$ versions of $\phi^2 \underline{\mathcal{X}}_\mathfrak{g}$ is the tweaked version of $\mathcal{V}_*{'}$. This new vector space is denoted in what follows by $\underline{\mathcal{V}}_*{'}$.

The lemma that follows directly asserts in effect that the constructions in [KS] that led to (6.3) can be repeated with $\underline{\mathcal{V}}_*{'}$ used in lieu of $\mathcal{V}_*{'}$. This lemma uses $\underline{\Pi}$ to denote the $L^2$-orthogonal projection onto $\underline{\mathcal{V}}_*{'}$ as defined by the metric $\mathfrak{g}_*$. The norm on $\underline{\mathcal{V}}_*{'}$ in this lemma is the $L^2$ norm as defined by $\mathfrak{g}_*$.

**Lemma 6.1**: *Given the knot* K, *a sufficiently large number* $R > 1$, *a non-negative integer* n *and a set* $\vartheta \subset X_K$ *of* n *points, there exists* $\kappa > 100$; *and given in addition,* $\varepsilon \in (0, 1]$, *there exists* $\kappa_\varepsilon > \kappa$ *with the following significance: Set the parameters* $\delta$ *and* $\{\delta_p\}_{p \in \vartheta}$ *to be less than* $\kappa_\varepsilon^{-1}$; *and then choose the remaining SO(4) parameters to construct the metric* $\mathfrak{g}_*$ *on* $X_K \#_n \overline{\mathbb{CP}}^2$. *There is a symmetric, traceless section of* $\otimes^2(T^*(X_K \#_n \overline{\mathbb{CP}}^2))$, *denoted by* $\underline{\mathfrak{h}}_{\mathfrak{g}_*}$, *with sup norm less than* $\tfrac{1}{100}$ *as measured by* $\mathfrak{g}_*$ *and with the properties listed below*.

- *The* $L^6_2$ *norm of* $\underline{\mathfrak{h}}_{\mathfrak{g}_*}$ *is less than* $\kappa(\delta^2 + \varepsilon^2 \sum_{p \in \vartheta} \delta_p^{\,2})$.
- *The tensor* $W_+(\mathfrak{g}_* + \underline{\mathfrak{h}}_{\mathfrak{g}_*})$ *is in the subspace* $\underline{\mathcal{V}}_{\mathfrak{g}*}{'}$.
- $\|W_+(\mathfrak{g}_* + \underline{\mathfrak{h}}_{\mathfrak{g}_*}) - \underline{\Pi} W_+(\mathfrak{g}_*)\|_2 \leq \kappa(\delta^2 + \varepsilon^3 \sum_{p \in \vartheta} \delta_p^{\,2})$.

*Proof of Lemma 6.1*: The proof has three steps.

Step 1: The assertion in Proposition 4.2 of [KS] still holds. Indeed, this follows because the $L^6_2$ norm of $\mathcal{X}_{\mathfrak{g}*} - \underline{\mathcal{X}}_{\mathfrak{g}*}$ is bounded by $c_\diamond \varepsilon^2 \sum_{p\in\vartheta} \delta_p^2$. The latter bound follows directly from three facts: First, if $\mathcal{X} \in \mathcal{H}_{K,R}$, then $\underline{\mathcal{X}}_{\mathfrak{g}*} - \mathcal{X}_{\mathfrak{g}*}$ is not zero only where each $p \in \vartheta$ version of the coordinate x that is used in (6.5) obeys $\frac{1}{2}\varepsilon\delta_p < |x| < \delta_p$. The second fact is that there are a priori $C^2$ bounds near p in $X_K$ for $\mathcal{X}$. The third fact is as follows: Write $\mathfrak{g}_*$ as $\phi^2 \mathfrak{g}$. Now suppose that $\mathcal{X}$ when written using a $\mathfrak{g}_R$ orthonormal frame for $\Lambda^+$ near p has components $\{\mathcal{X}^{ab}\}_{a,b \in \{1,2,3\}}$. Fix $p \in \vartheta$. Then $\phi^2 \mathcal{X}$ where p's version of the coordinate x in (6.5) obeys $\frac{1}{2}\varepsilon\delta_p \le |x| < \delta_p$ when written using a $\phi^2 \mathfrak{g}_R$ orthonormal frame for $\Lambda^+$ has components $\{\phi^{-2}\mathcal{X}^{ab}\}_{a,b \in \vartheta}$.

Step 2: Granted that Proposition 4.2 in [KS] holds, then what is said in Proposition 5.8 of [KS] holds with $\underline{\mathcal{V}}_*'$ used in lieu of $\mathcal{V}_*'$ when $\delta$ and $\{\delta_p\}_{p \in \vartheta}$ are less than $c_\diamond^{-1}$. Proposition 5.8 in [KS] implies in turn Proposition 5.12 in [KS] which, when $\delta$ and $\{\delta_p\}_{p \in \vartheta}$ are less than $c_\diamond^{-1}$, gives the desired section $\underline{\mathfrak{h}}_{\mathfrak{g}*}$ that obeys the first two bullets of Lemma 6.1.

Step 3: Equation (5.9) in [KS] leads to an analog of the Taylor expansion (6.3),

$$\mathcal{W}_+(\mathfrak{g}_* + \underline{\mathfrak{h}}_{\mathfrak{g}*}) = \Pi \mathcal{W}_+(\mathfrak{g}_*) + \Pi \mathcal{L}_{\mathfrak{g}*} \underline{\mathfrak{h}}_{\mathfrak{g}*} + \underline{\mathcal{Q}}_{\mathfrak{g}*}$$

(6.6)

with the $L^2$ norm of $\underline{\mathcal{Q}}_{\mathfrak{g}*}$ also bounded by what is written on the right hand side of (6.4). In particular, Lemmas 3.3 and 3.4 with the bound from the right hand side of (6.4) imply

$$\|\underline{\mathcal{Q}}_{\mathfrak{g}*}\|_2 \le c_\diamond(\delta^2 + \varepsilon^3 \sum_{p\in\vartheta} \delta_p^2) \ .$$

(6.7)

(This inequality with $\varepsilon$ having exponent 3 (as opposed to 2) is the reason for introducing $\varepsilon$ in Section 3f.). Meanwhile, an integration by parts can be used with the three facts listed in Step 1 and the fact that the $\mathfrak{g} = \mathfrak{g}_R$ version of $\mathcal{X}$ obeys $\mathcal{L}_{\mathfrak{g}}^\dagger \mathcal{X} = 0$ on $\mathbb{T}_K$ to see that

$$\|\Pi \mathcal{L}_{\mathfrak{g}*} \underline{\mathfrak{h}}_{\mathfrak{g}*}\|_2 \le c_\diamond(\delta^2 + \varepsilon^3 \sum_{p\in\vartheta} \delta_p^2) \ .$$

(6.8)

(This inequality with $\varepsilon$ having exponent 3 (as opposed to exponent 2) is the reason for introducing $\varepsilon$ in (6.5).) Inequalities in (6.7) and (6.8) give the lemma's third bullet

### c) Definition of the obstruction vector

As noted previously, the obstruction space $\mathbb{V}$ from Proposition 4.1 is the direct sum of the spaces $\mathcal{H}_T$ and $\mathcal{H}_K$ that are described in Propositions 5.3 and 5.4. To define Proposition 4.1's obstruction vector in $\mathbb{V}$, fix $R > 1$, a non-negative integer n and a number $\varepsilon \in (0, 1)$. The number R is chosen in particular with an appeal to Lemma 6.1 in mind. For this appeal, fix parameters $\delta$ and $\{\delta_p\}_{p \in \vartheta}$, each less than $\kappa_\varepsilon^{-1}$ with $\kappa_\varepsilon$ from Lemma 6.1. Use R with the parameters in $\vartheta$ to to define the metric $\mathfrak{g}_*$ from the conformal class $[\mathfrak{g}^n_{R,\delta}]$, and use this data to define the vector space $\underline{\mathcal{V}}_*'$. Lemma 6.1 supplies the tensor $\underline{\mathfrak{h}}_{\mathfrak{g}_*}$ which obeys Lemma 6.1's three bullets. The corresponding vector $\mathcal{W}_+(\mathfrak{g}_* + \underline{\mathfrak{h}}_{\mathfrak{g}_*})$, which is in $\underline{\mathcal{V}}_{\mathfrak{g}*}'$, defines a vector in $\mathbb{V}$ via the composition of first the isomorphism between $\underline{\mathcal{V}}_{\mathfrak{g}*}'$ and $\mathcal{H}_{K,R}$ and then the monomorphism $\mathcal{Q}_T \oplus \mathcal{Q}_K$ from $\mathcal{H}_{K,R}$ to $\mathbb{V}$ that is described in Propositions 5.3 and 5.4. The image of $\mathcal{W}_+(\mathfrak{g}_* + \underline{\mathfrak{h}}_{\mathfrak{g}_*})$ in $\mathbb{V}$ via the composition of these to homomorphisms is Proposition 4.1's obstruction vector.

### 7. Proof of Proposition 4.2

The proof of Proposition 4.2 in the case when the kernel of the homomorphism $\hat{\imath}*$ is trivial is given below in Section 7c and the proof when the kernel of $\hat{\imath}*$ is non-trivial is given in Section 7d. The basic idea for the proof comes from Section 6 in the paper by Donaldson and Friedman [DF]. Sections 7a and 7b explain what is needed from Section 6 of this Donaldson-Friedman paper.

### a) The $\underline{\Pi}\mathcal{W}_+(\mathfrak{g}_*)$ term in (6.6)

As explained in Section 6 of Donaldson and Friedman's paper [DF], the term $\underline{\Pi}\mathcal{W}_+(\mathfrak{g}_*)$ in (6.6) can be written more or less explicitly because the support of $\mathcal{W}_+(\mathfrak{g}_*)$ lies entirely in the sixteen copies of $E^2_{\delta/r}$ and the n copies of $E_{\delta_p/r_p}$. The steps that follow directly outline how this is done and the upcoming Lemma 7.1 summarizes the resulting formula for $\underline{\Pi}\mathcal{W}_+[\mathfrak{g}_*]$.

<u>Step 1</u>: Write $\mathfrak{g}_*$ as $\phi^2\mathfrak{g}$ again. Supposing that $\mathcal{X} \in \mathcal{H}_{K,R}$, then the corresponding element in $\underline{\mathcal{V}}_*'$ is $\phi^2\underline{\mathcal{X}}_\mathfrak{g}$. The projection of $\mathcal{W}_+[\mathfrak{g}_*]$ along the unit vector in $\underline{\mathcal{V}}_*'$ in the $\phi^2\underline{\mathcal{X}}_\mathfrak{g}$ direction is the ratio whose numerator is the inner product between $\phi^2\underline{\mathcal{X}}_\mathfrak{g}$ and $\mathcal{W}_+(\mathfrak{g}_*)$ as defined by the metric $\phi^2\mathfrak{g}$, and whose denominator is the $L^2$ norm of $\phi^2\underline{\mathcal{X}}_\mathfrak{g}$ as defined by $\phi^2\mathfrak{g}$. The numerator of this ratio is also equal to the inner product between $\underline{\mathcal{X}}_\mathfrak{g}$ and $\mathcal{W}_+(\mathfrak{g})$ as defined by the metric $\mathfrak{g}$ because $\mathcal{W}_+(\phi^2\mathfrak{g}) = \phi^2\mathcal{W}_+(\mathfrak{g})$ and because of how the volume forms for $\phi^2 \mathfrak{g}$ and $\mathfrak{g}$ are related. Likewise, the denominator is the $L^2$ norm of $\underline{\mathcal{X}}_\mathfrak{g}$ as defined by $\mathfrak{g}$.

Step 2: Let $o \in \mathbb{T}_K$ denote a fixed point of the involution $\iota$. Since $\mathcal{X}$ is smooth, it has a Taylor expansion near $o$ when $\mathcal{X}$ is written with the Euclidean coordinates that are introduced in Part 1 of Section 3e. This Taylors expansion with the formula in (3.10) for the metric and (4.13) for $\mathcal{X}_\mathfrak{g}$ can be used to see that the contribution from $o$'s version of $E^2_{\delta/r}$ to the $L^2$ inner product between $\underline{\mathcal{X}}_\mathfrak{g}$ and $\mathcal{W}_+(\mathfrak{g})$ is bounded by $c_0 \delta^3$ with $c_0$ denoting a number that is independent of $\delta$ and of the other parameters that are used to define the metric $\mathfrak{g}$. (One might a priori expect this contribution to be $\mathcal{O}(\delta^2)$, but an integration by parts can be used to see that the terms that could make an $\mathcal{O}(\delta^2)$ contribution are zero.)

Step 3: Taylor's theorem with remainder is also used in the case of $X_K \#_n \overline{\mathbb{CP}}^2$ to write the leading order contribution to the inner product between $\underline{\mathcal{X}}_\mathfrak{g}$ and $\mathcal{W}_+(\mathfrak{g})$ from the points in $\vartheta$. To say more about this leading order contribution, remember that the construction of the conformal structure requires a choice of a suitable Gaussian coordinate chart centered at each point in $\vartheta$. These are parametrized by the elements in the SO(4) principal orthonormal frame bundle for the metric $\mathfrak{g}_R$ at the given point in $\vartheta$. An identification of the fiber of this orthonormal frame bundle at a given point $p \in \vartheta$ also identifies the vector space of traceless, symmetric $3 \times 3$ matrices with the fiber at p of the bundle of traceless, symmetric elements in $\Lambda^+ \oplus \Lambda^+$. Let $\mathcal{I}$ denote the diagonal matrix with upper and middle diagonal entries being 1 and with lower diagonal entry being -2. The choice of a Gaussian coordinate chart for a given $p \in \vartheta$ writes $\mathcal{I}$ as an element in $(\Lambda^+ \oplus \Lambda^+)|_p$. This element is denoted in what follows by $\mathcal{I}_p$.

Step 4: Suppose first that p is a point from $\vartheta$ in the $\mathbb{T}-(N\cup\iota(N))$ part of $\mathbb{T}_K$. A calculation using (3.11) and the Taylor's expansion from Step 3 can be used (see Equations (6.9) and (6.14) in [DF]) to write the contribution from p's version of $E_{\delta_p/r_p}$ to the $L^2$ inner product between $\underline{\mathcal{X}}_\mathfrak{g}$ and $\mathcal{W}_+(\mathfrak{g})$ as

$$c_\# \varepsilon^2 \delta_p^2 \langle \mathcal{X}|_p, \mathcal{I}_p \rangle + r_p ,$$
(7.1)

where $\langle \, , \, \rangle$ indicates the metric inner product on $\Lambda^+ \oplus \Lambda^+$; where $c_\#$ is non-zero and independent of all parameter choices; and where $r_p$ is a number with norm bound

$$|r_p| \leq c_{\diamond\diamond} \varepsilon^3 \delta_p^2 \|\mathcal{X}\|_2 .$$
(7.2)

What is denoted here by $c_{\diamond\diamond}$ depends on K and R, but not on the other parameters; it is independent of $\varepsilon$ and $\delta_p$ in particular.

If p is from either the $(S^3 - N_K) \times S^1$ or $\iota((S^3 - N_K) \times S^1)$ part of $\mathbb{T}_K$, then the contribution from p's version of $E_{\delta_p / r_p}$ to the inner product of $\underline{X}_{\mathfrak{g}}$ and $\mathcal{W}_+(\mathfrak{g})$ has a form that is analogous to (7.1):

$$z_p c_\# \varepsilon^2 \delta_p^2 \langle \mathcal{X}|_p, \mathcal{I}_p \rangle + r_p ,$$

(7.3)

with $z_p$ being a positive number that depends only on p, K and R; and with $r_p$ obeying (7.2) with the understanding that $c_{\diamond\diamond}$ can now depend on p, K and R and $\vartheta$, but not on the other parameters.

The following lemma summarizes the observations in the preceding paragraphs.

**Lemma 7.1**: *Fix the knot K, a sufficiently large number R > 1, a non-negative integer n and a set $\vartheta \subset X_K$ of n points. Given this data, there exists $\kappa > 100$; and given in addition, $\varepsilon \in (0, 1]$, there exists $\kappa_\varepsilon > \kappa$ with the following significance: Set the parameters $\delta$ and $\{\delta_p\}_{p \in \vartheta}$ to be less than $\kappa_\varepsilon^{-1}$; and then choose the sixteen SO(4) parameters labeled by the fixed points of $\iota$ on $\mathbb{T}_K$ and then the n Gaussian coordinate chart parameters labeled by the points in $\vartheta$ to construct the metric $\mathfrak{g}_*$ on $X_K \#_n \overline{\mathbb{CP}^2}$. Let $\mathcal{X}$ denote a given element in $\mathcal{H}_{K,R}$. Then the $L^2$ inner product between the corresponding element in $\underline{\mathcal{V}}_*'$ and $\mathcal{W}_+(\mathfrak{g}_*)$ can be written as*

$$\sum_{p \in \vartheta} (z_p c_\# \varepsilon^2 \delta_p^2 \langle \mathcal{X}|_p, \mathcal{I}_p \rangle + r_p) + r_\diamond$$

*where $z_p$ is independent of $\varepsilon$, $\{\delta_{p'}\}_{p' \in \vartheta}$ and the Gaussian coordinate charts; and it obeys $z_p > \kappa^{-1}$. It is also independent of $p \in \vartheta$ if p is defined by a point from the $\mathbb{T} - (N \cup \iota(N))$ part of $\mathbb{T}_K$. Meanwhile, the number $r_p$ obeys $|r_p| \leq \kappa \varepsilon^3 \delta_p^2 \|\mathcal{X}\|_2$. Finally, what is denoted by $r_\diamond$ obeys $|r_\diamond| \leq \kappa \delta^2 \|\mathcal{X}\|_2$.*

**b) Calculations from Donaldson/Friedman**

This subsection is a digression to review some computations from [DF] concerning the $\{\langle \mathcal{X}|_p, \mathcal{I}_p \rangle\}_{p \in \Theta}$ terms that appear in the formula from Lemma 7.1. The digression starts by introducing by way of notation $\mathbb{S}$ to denote the vector space of $3 \times 3$, symmetric traceless matrices. View $S^2$ as the unit sphere in $\mathbb{R}^3$. The sphere $S^2$ has an SO(3) equivariant embedding in $\mathbb{S}$ as the subspace of matrices with the three eigenvalues $\{1, 1, -2\}$. This embedding is denoted by $\Phi_1$; it is defined by the rule where by any given

unit length vector $\mathfrak{n} \in \mathbb{R}^3$ is sent to the matrix $\mathbb{I} - 3\mathfrak{n} \otimes \mathfrak{n}^T$ with $\mathbb{I}$ denoting here the $3 \times 3$ identity matrix. Note that the differential of $\Phi_1$ at a given unit vector $\mathfrak{n}$ sends a vector $\mathfrak{v}$ orthogonal to $\mathfrak{n}$ to the matrix $-3(\mathfrak{v} \otimes \mathfrak{n}^T + \mathfrak{n} \otimes \mathfrak{v}^T)$.

A map from $(0,1) \times S^2 \times S^2$ to $\mathbb{S}$ to be denoted by $\Phi_2$ is defined by the rule that sends $(t, \mathfrak{n}_1, \mathfrak{n}_2)$ to $(1-t)\Phi_1(\mathfrak{n}_1) + t\Phi_1(\mathfrak{n}_2)$. This is to say that

$$\Phi_2(t, \mathfrak{n}_1, \mathfrak{n}_2) = \mathbb{I} - 3(1-t)\mathfrak{n}_1 \otimes \mathfrak{n}_1^T - 3t\mathfrak{n}_2 \otimes \mathfrak{n}_2^T .$$

(7.4)

Supposing that $\mathfrak{n}_1$ and $\mathfrak{n}_2$ are orthogonal, let $\mathfrak{n}_3$ denote their vector cross product. The following observations about $\Phi_2$ are valid when $\mathfrak{n}_1$ and $\mathfrak{n}_2$ are orthogonal:

- $\Phi_2(\frac{1}{2}, \mathfrak{n}_1, \mathfrak{n}_2) = -\frac{1}{2}\Phi_1(\mathfrak{n}_3)$.
- *The differential of $\Phi_2$ defines a surjective map from the tangent space of $(0,1) \times S^2 \times S^2$ at $(t = \frac{1}{2}, \mathfrak{n}_1, \mathfrak{n}_2)$ to the orthogonal complement in $\mathbb{S}$ of the span of $\Phi_1(\mathfrak{n}_3)$.*

(7.5)

The first bullet's assertion follows by inspection. The assertion in the second bullet follows from the fact that the differential of $\Phi_1$ at a given point $\mathfrak{n} \in S^2$ sends a vector $\mathfrak{v}$, orthogonal to $\mathfrak{n}$, to the matrix $-3(\mathfrak{v} \otimes \mathfrak{n} + \mathfrak{n} \otimes \mathfrak{v})$.

A third map also plays a role. This one is denoted by $\Phi_3$. To define this map, introduce by way of notation $\Delta^2$ to denote the simplex in $\mathbb{R}^3$ where the coordinates $(t_1, t_2)$ obey $t_1, t_2 > 0$ and $t_1 + t_2 < 1$. The map $\Phi_3$ sends $\Delta^2 \times (\times^3 S^2)$ into $\mathbb{S}$ according to the rule

$$((t_1, t_2), \mathfrak{n}_1, \mathfrak{n}_2, \mathfrak{n}_3) \to t_1 \Phi_1(\mathfrak{n}_1) + t_2 \Phi_2(\mathfrak{n}_2) + (1 - t_1 - t_2)\Phi_3(\mathfrak{n}_3) .$$

(7.6)

The following assertions about $\Phi_3$ are true when $\mathfrak{n}_1$, $\mathfrak{n}_2$ and $\mathfrak{n}_3$ are orthogonal.

- $\Phi_3((\frac{1}{3}, \frac{1}{3}), \mathfrak{n}_1, \mathfrak{n}_2, \mathfrak{n}_3) = 0$.
- *The differential of $\Phi_3$ defines a surjective map from the tangent space of $\Delta^2 \times^3 S^2$ at $((\frac{1}{3}, \frac{1}{3}), \mathfrak{n}_1, \mathfrak{n}_2, \mathfrak{n}_3)$ to $\mathbb{S}$.*

(7.7)

The verification of these assertions is left to the reader, or see Section 6 in [DF].

Let $\mathfrak{n}^0_1$, $\mathfrak{n}^0_2$ and $\mathfrak{n}^0_3$ denote the respective unit vectors in $\mathbb{R}^3$ along the positive $x_1$, $x_2$ and $x_3$ axis in $\mathbb{R}^3$. For each $a \in \{1, 2, 3\}$, let $\mathfrak{f}_a$ to denote the map from $SO(4)$ to $\mathbb{S}$ that is obtained by composing three maps: The first map sends points in $SO(4)$ to $SO(3)$ via the self dual representation (denoted above by $O \to \hat{o}$); and the second map is the map from $SO(3)$ to $S^2$ given by the rule that sends $\hat{o} \in SO(3)$ to the vector in $\mathbb{R}^3$ that is obtained from acting on $\mathfrak{n}^0_a$ by $\hat{o}$. This is denoted by $\hat{o}\mathfrak{n}^0_a$. The third map is $\Phi_1$. Thus, $\mathfrak{f}_a(\mathcal{O}) = \Phi_1(\hat{o} \mathfrak{n}^0{}_a)$. Now define respective maps $\mathbb{O}_2$ and $\mathbb{O}_3$ from $(0, 1) \times (\times^2 SO(4))$ to $\mathbb{S}$ and from $\Delta^2 \times (\times^3 SO(4))$ to $\mathbb{S}$ by setting

- $\mathbb{O}_2(t, O_1, O_2) = t \mathfrak{f}_1(O_1) + (1 - t) \mathfrak{f}_2(O_2)$.
- $\mathbb{O}_3((t_1, t_1), O_1, O_2, O_3) = t_1 \mathfrak{f}_1(O_1) + t_2 \mathfrak{f}_2(O_2) + (1 - t_1 - t_2) \mathfrak{f}_3(O_3)$ .

(7.8)

The following lemma summarizes the salient features of these maps. This lemma uses $\hat{I}$ to denote the identity element in SO(4).

**Lemma 7.2**: *The maps $\mathbb{O}_2$ and $\mathbb{O}_3$ have the following properties:*
- $\mathbb{O}_2(\frac{1}{2}, \hat{I}, \hat{I}) = 0$ *and the differential of $\mathbb{O}_2$ at $(\frac{1}{2}, \hat{I}, \hat{I})$ is surjective onto the orthogonal complement in $\mathbb{S}$ of $\mathfrak{f}_3(\hat{I})$.*
- $\mathbb{O}_3((\frac{1}{3}, \frac{1}{3}), \hat{I}, \hat{I}, \hat{I}) = 0$ *and the differential of $\mathbb{O}_3$ at $((\frac{1}{3}, \frac{1}{3}), \hat{I}, \hat{I}, \hat{I})$ is surjective onto $\mathbb{S}$.*

*Proof of Lemma 7.2*: The assertions in the first and second bullets follow directly from (7.5) and (7.7) respectively.

### c) The proof when the kernel of $i^*$ is zero

The proof of Proposition 4.2 when the kernel of $i^*$ is zero has seven parts.

*Part 1*: Define a four (real dimensional) subspace in the vector space $\mathcal{H}_T$ from (5.9) using two complex parameters $(a_1, a_2)$ as follows: The components $\{\mathcal{X}^{ab}\}_{a,b \in \{1,2,3\}}$ of its elements have the form

- $\mathcal{X}^{33} = 0$ *and* $\mathcal{X}^{11} = -\mathcal{X}^{22}$.
- $\mathcal{X}^{13} - i\mathcal{X}^{23} = a_1$ .
- $\mathcal{X}^{11} - \mathcal{X}^{22} - 2i \mathcal{X}^{12} = a_2$ .

(7.9)

Denote this subspace by $\mathcal{H}_a$.

The involution $\iota$ induces an involution on T´ (to be called $\iota$ also) that sends the $\mathbb{R}/(2\pi\mathbb{Z})$ coordinates $(t_3, t_4)$ to $(-t_3, -t_4)$. Reintroduce $p_*$ to denote the point $(t_3 = t_*, t_4 = t_*)$ in the T´ torus. Let $x$ denote the function on T´$-\{p_*, \iota(p_*)\}$ that is described in Lemma 5.1 and let $u$ denote the function on T´$-\{p_*, \iota(p_*)\}$ that is described in Lemma 5.2. Fix a point $q \in$ T´$-\{p_*, \iota(p_*)\}$ and having chosen q, define $\alpha_1 = -x(q)$ and $\alpha_2 = -u(q)$. These are defined so that the complex functions $\alpha_1 + x$ and $\alpha_2 + u$ vanish at q. Define a three dimensional subspace in $\mathcal{H}_T$ using a real number $s$ and a complex number $\flat$ as follows: The components of its elements $\{\mathcal{X}^{ab}\}_{a,b \in \{1,2,3\}}$ have the form

- $\mathcal{X}^{33} = 0$ and $\mathcal{X}^{11} = -\mathcal{X}^{22}$.
- $\mathcal{X}^{13} - i\mathcal{X}^{23} = s\mathfrak{z}_K(\alpha_1 + x)$ .
- $\mathcal{X}^{11} - \mathcal{X}^{22} - 2i\mathcal{X}^{12} = \flat(\alpha_1 + x) + \overline{\flat}(\alpha_2 + u)$ .

(7.10)

Denote this suspace by $\mathcal{H}_\flat$. The vector space $\mathcal{H}_T$ from Proposition 5.3 is the direct sum of its subspace of $\mathcal{H}_a$ and $\mathcal{H}_\flat$. A norm on $\mathcal{H}_T$ is needed later and the Euclidean norm that takes the length of $((a_1, a_2), (s, \flat))$ to be $(|a_1|^2 + |a_2|^2 + s^2 + |\flat|^2)^{1/2}$ serves this purpose.

*Part 2*: Let K denote a hyperbolic knot in $S^3$ with kernel($i^*$) $\subset H^1(S^3-K; V)$ equal to 0. Fix R to be very much greater than 1 and use it to construct $\mathbb{T}_K$, its conformal structure $[\mathfrak{g}_R]$ and the metric $\mathfrak{g}_R$ in this conformal structure from Step 1 in Section 4a. Let $\mathcal{H}_{K,R}$ again denote the $\iota$-invariant kernel of the $\mathfrak{g}=\mathfrak{g}_R$ version of $\mathcal{L}_\mathfrak{g}^\dagger$. Propositions 5.3 and 5.4 imply the following: Let D denote the disk T´ centered at $p_*$ with radius $e^{-2R}$. This is the disk where the coordinates $(t_3, t_4)$ obey $(t_3-t_*)^2 + (t_4-t_*)^2 < e^{-2R}$. A given element $\mathcal{X} \in \mathcal{H}_{K,R}$ can be written on the domain $T \times (T´-(D \cup \iota(D)))$ as $\mathcal{Q}_T(\mathcal{X}) + \mathfrak{r}_T$ with $\mathfrak{r}_T$ and $\mathcal{Q}_T(\mathcal{X})$ uniformly large in the following sense:

- *Write $\mathcal{H}_T$ as $\mathcal{H}_a \oplus \mathcal{H}_\flat$. The $(a_1, a_2)$ and $(s, \flat)$ coordinates of $\mathcal{Q}_T(\mathcal{X})$ with respect to the $\mathcal{H}_a$ and $\mathcal{H}_\flat$ basis in (7.1) and (7.2) obey*

$$c_K^2 \|\mathcal{X}\|_2 \geq e^{-2R}(|a_1|^2 + |a_2|^2) + \ln(R) e^{-2R}(t^2 + |\flat|^2) \geq c_K^{-2} \|\mathcal{X}\|_2.$$

- *The element $\mathfrak{r}$ obeys the bound $|\mathfrak{r}_T| \leq \kappa(\Delta^{-2} e^{-2R} + e^{-R})\|\mathcal{X}\|_2$ where $\Delta \geq c_K e^{-R}$ .*

(7.11)

Here and in what follows, $c_K$ is a number greater than 1 that depends only on the knot K. In particular, it is independent of a given element in $\mathcal{H}_{K,R}$ and it is independent of R. (The $L^2$ norm of $\mathcal{X}$ will be less than 1 if the lower bound in the first bullet of (7.11) is violated, and greater than 1 if the upper bound is violated. The factor $e^{-2R}$ appears here because the area of the T factor is proportional to $e^{-2R}$. The factor of $\ln(R)$ appears before $t^2 + |\flat|^2$ because $|x|$ and $|u|$ are both commensurate with $|z|^{-1}$ on $\mathbb{D}-p_*$.)

*Part 3*: A point $q \in T´-\{p_*, \iota(p_*)\}$ was specified in Part 1. Choose a second point in $T´-\{p_*, \iota(p_*)\}$ which is distinct from $q$, $\iota(q)$ and the fixed points of the involution ι on T´. This second point in T´ is denoted by q´. With regards to (q, q´) and the choice of R: If R is sufficiently large, then the points q and q´ will lie in T´−(D∪ι(D)). It is assumed in what follows that R is large enough so that this is the case. With R fixed, take δ > 0 but very small and then choose the required SO(4) parameters for each of the sixteen fixed points of ι on $\mathbb{T}_K$ to construct the conformal structure $[\mathfrak{g}_{R,\delta}]$ on $X_K$.

Let $q_1$ and $q_2$ denote distinct points in the torus T that are not fixed points of the involution on T induced by ι. This is the involution that sends $(t_1, t_2)$ to $(-t_1, -t_2)$. Define $\vartheta_q$ to be the following two element set $\{(q_1,q),(q_2,q)\}$ in $T \times (T´-(D\cup\iota(D)))$. Define a second two element set $\vartheta_{q´}$ to be $\{(q_1,q´),(q_2,q´)\}$. The sets $\vartheta_q$ and $\vartheta_{q´}$ when viewed in $\mathbb{T}_K$ project to distinct two element sets in $\mathbb{T}_K/\iota$. If the number δ that is used to define the conformal structure $[\mathfrak{g}_{R,\delta}]$ on $X_K$ is sufficiently small, then the incarnations of both $\vartheta_q$ and $\vartheta_{q´}$ in $\mathbb{T}_K/\iota$ can be viewed as sets in $X_{K,\delta}$. It is assumed in what follows that δ is small enough so that this is the case. Let ϑ denote the union of the points from the $X_K$ versions of $\vartheta_q$ and $\vartheta_{q´}$.

The set ϑ will be used to construct the n = 4 version of the conformal structure $[\mathfrak{g}^n_{R,\delta}]$ on $X_K \#_4 \overline{\mathbb{CP}}^2$. The construction requires the choice of a number $\varepsilon \in (0, 1]$. Having chosen ε, the construction also requires the choice of the set $\{\delta_p\}_{p\in\vartheta} \subset (0, 1]$. To choose this set, first fix $\delta_\diamond > 0$ but small; the set $\{\delta_p\}_{p\in\vartheta}$ is determined below $\delta_\diamond$ and by parameters $t, t´ \in (0, 1)$ according to the following rule:

- *For p = $(q_1, q)$, take $\delta_p = \delta_\diamond \sqrt{t}$; and for p = $(q_2,q)$, take $\delta_p = \delta_\diamond \sqrt{1-t}$.*
- *For p = $(q_1,q´)$, take $\delta_p = \delta_\diamond \sqrt{t´}$; and for p = $(q_2,q´)$, take $\delta_p = \delta_\diamond \sqrt{1-t´}$.*

(7.12)

The definition of $[\mathfrak{g}^4_{R,\delta}]$ also requires a choice of a Gaussian coordinate charts at each point in ϑ. This is done in the next part of the proof.

*Part 4*: There is an almost canonical choice of Gaussian coordinate chart for a point in the $\mathbb{T}-(N\cup\iota(N))$ part of $\mathbb{T}_K$. To say more, write the $(t_1, t_2, t_3, t_4)$ coordinates of p as $(p_1, p_2, p_3, p_4)$. Let z denote the 2 × 2 matrix from Part 1 of Section 3e. Local Euclidean coordinates near p are $(x_1, x_2, x_3, x_4)$ with $x_a = e^{-R}\sum_{b\in\{1,2\}}(z^{-1})_{ab}(t_a-p_a)$ for $a \in \{1,2\}$ and $x_a = (t_a - p_a)$ for $a \in \{3,4\}$. These coordinates write $\mathfrak{g}_T$ as the standard Euclidean metric; they are Gaussian coordinates and this is the canonical Gaussian coordinate choice. With this choice of coordinates, what was defined by $\mathcal{I}_p$ in the preceding step is denoted simply by $\mathcal{I}$.

A different choice of Gaussian coordinates is obtained from the canonical choice through the action of an element in SO(4) on the $(x_1, x_2, x_3, x_4)$ coordinates. Supposing that O is the element in question, then the corresponding version of $\mathcal{I}_p$ is obtained from $\mathcal{I}$ via the adjoint action of an SO(3) matrix defined by $\mathcal{O}$. This is to say that $\mathcal{I}_p$ can be written as $\hat{o}\mathcal{I}\hat{o}^{-1}$ with $\hat{o}$ denoting the matrix in SO(3) that is defined by $\mathcal{O}$ via the homomorphism to SO(3) that gives the action of SO(4) on the self-dual summand in the vector space $\wedge^2\mathbb{R}^4$. Henceforth, the choice of a Gaussian coordinate chart for a point in $\vartheta$ from the $\mathbb{T}-(N\cup\iota(N))$ part of $\mathbb{T}_K$ will be viewed as a choice of an element in SO(4).

The two SO(4) elements for the points in $\vartheta$ from $\vartheta_q$ and the $(0,1)$ parameter $t$ define a point in a copy of $(0,1)\times SO(4)\times SO(4)$. The version of this space labeled by q is denoted by $Z_q$. By the same token, the two SO(4) elements for the points in $\vartheta$ from $\vartheta_{q'}$ and the $(0,1)$ parameter $t'$ define a point in another copy of space $Z_{q'}$.

With $\varepsilon$ chosen and supposing that $\delta$ and $\delta_\Diamond$ are sufficiently small, then the data from the set $Z_q\times Z_{q'}$ determines a conformal structure $[\mathfrak{g}^{n=4}{}_{R,\delta}]$ of the sort that is described in Propositions 4.1 and 4.2. Supposing that $(z_q, z_{q'})\in Z_q\times Z_{q'}$ has been chosen, let $\mathfrak{g}_*$ denote the corresponding version of the metric $\mathfrak{g}^n{}_{R,\delta*}$ that appears in Proposition 4.1.

*Part* 5: If $\delta$ and $\delta_\Diamond$ are sufficiently small, and supposing a point $(z_q, z_{q'})\in Z_q\times Z_{q'}$ have been chosen, then the writing of $\mathcal{X}\in\mathcal{H}_{K,R}$ as $\mathcal{Q}_T(\mathcal{X})+\mathfrak{r}$ and what is said by Lemmas 6.1 and 7.1 lead to a particularly useful description of the obstruction vector $\mathcal{W}_+(\mathfrak{g}_*+\underline{\mathfrak{h}}_{\mathfrak{g}_*})$ in the $(z_q, z_{q'})$ version of $\underline{\mathcal{V}}_*'$. The next lemma supplies this description. The lemma supposes implicitly that R, $\varepsilon$ and the various parameter choices for the points in $\vartheta$ have been made so that Lemmas 6.1 and 7.1 can be invoked.

**Lemma 7.3**: *Fix $\mathcal{X}\in\mathcal{H}_{K,R}$ and write it near the points in $\vartheta$ as $\mathcal{X}=\mathcal{Q}_T(\mathcal{X})+\mathfrak{r}$ with $\mathfrak{r}$ given in (7.11). The $L^2$ inner product (as defined by $\mathfrak{g}_*$) between $\mathcal{W}_+(\mathfrak{g}_*+\underline{\mathfrak{h}}_{\mathfrak{g}_*})$ and the element in $\underline{\mathcal{V}}_*'$ defined by $\mathcal{X}$ can be written as*

$$c\varepsilon^2\delta_\Diamond^2(\langle\mathcal{Q}_T(\mathcal{X})|_q, \mathbb{O}_2(z_q)\rangle + \langle\mathcal{Q}_T(\mathcal{X})|_{q'}, \mathbb{O}_2(z_{q'})\rangle) + \varepsilon^2\delta_\Diamond^2\,\wp(\mathcal{Q}_T(\mathcal{X})),$$

*where c is positive and independent of all choices and parameters, and where $\wp$ is a linear functional on $\mathcal{H}_a\oplus\mathcal{H}_b$ with the following properties:*
- *This functional depends continuously on the points in $Z_q\times Z_{q'}$.*
- *The norm of $\wp$ obeys the bound $|\wp|\leq c_\Delta e^{-R}+c_R(\varepsilon+\varepsilon^{-2}\delta_\Diamond^{-2}\delta^2)$ with $c_\Delta$ being independent of the parameters R, $\delta$, $\varepsilon$, $\delta_\Diamond$ and the parameters in $Z_q\times Z_{q'}$; and with $c_R$ being independent of $\delta, \varepsilon, \delta_\Diamond$ but not necessarily on R or on the parameters in $Z_q\times Z_{q'}$.*

**Proof of Lemma 7.3**: The linear functional $\wp$ on $\mathcal{H}_a\oplus\mathcal{H}_b$ is defined as follows: Let $\mathcal{H}\subset\mathcal{H}_a\oplus\mathcal{H}_b$ denote the image via $\mathcal{Q}_T$ of the vector space $\mathcal{H}_{K,R}$. Write $\mathcal{H}_a\oplus\mathcal{H}_b$ as $\mathcal{H} \oplus \mathcal{H}^\perp$ with $\mathcal{H}^\perp$ being the orthogonal complement as defined using the metric from the identification of $\mathcal{H}_a \oplus \mathcal{H}_b$ with $\mathbb{C}^2 \oplus (\mathbb{R} \oplus \mathbb{C})$ that comes by using the coordinates $(a_1, a_2)$ on $\mathcal{H}_a$ and the coordinates $(s, b)$ on $\mathcal{H}_b$. The linear functional $\wp$ is equal zero on $\mathcal{H}^\perp$. The fact that $\mathcal{Q}_T$ is a monomorphism implies that there is an inverse homormorphism from $\mathcal{H}$ to $\mathcal{H}_{K.R}$. This is denoted by $\mathcal{Q}_T^{-1}$. With this notation in hand, then $\wp(\mathfrak{X})$ for $\mathfrak{X} \in \mathcal{H}$ can be written as $A(\mathfrak{X}) + B(\mathfrak{X}) + C(\mathfrak{X})$ where A, B and C are described in the subsequent paragraph.

What is denoted by $A(\mathfrak{X})$ is the $L^2$ inner product as defined by the metric $\mathfrak{g}_*$ between the element in $\underline{\mathcal{V}}_*'$ defined by $\mathcal{Q}_T^{-1}(\mathfrak{X})$ and $W_+(\mathfrak{g}_* + \underline{\mathfrak{h}}_{-\mathfrak{g}_*}) - \Pi W_+(\mathfrak{g}_*)$. It follows from Lemma 6.1 and (7.11) that the norm of A is bounded by $\gamma_R (\delta^2 + \varepsilon^3 \delta_\diamond^2)$ with $\gamma_R$ being independent of $\delta$, $\varepsilon$ and $\delta_\diamond$. What is denoted by $B(\mathfrak{X})$ is the $\mathcal{X} = \mathcal{Q}_T^{-1}(\mathfrak{X})$ version of the contributions to Lemma 7.1's formula by $\sum_{p \in \Theta} r_p + r_\diamond$. It follows from Lemma 7.1 and (7.11) that B also has norm at most $\gamma_R (\delta^2 + \varepsilon^3 \delta_\diamond^2)$. To define $C(\mathfrak{X})$, first write $\mathcal{X} = \mathcal{Q}_T^{-1}(\mathfrak{X})$ as $\mathfrak{X} + \mathfrak{x}_T$ in the manner of (7.1). What is denoted by $C(\mathfrak{X})$ is the contribution to Lemma 7.1's formula from the term $\sum_{p \in \vartheta} z_p c_\# \varepsilon^2 \delta_p^2 \langle \mathfrak{x}_T|_p, \mathcal{I}_p \rangle$. It follows from (7.11) that the norm of C is bounded by $c_\Delta e^{-R}$.

Note that the contributions to Lemma 7.3's formula from $\langle \mathcal{Q}_T(\mathcal{X})|_q, \mathbb{O}_2(z_q) \rangle$ and $\langle \mathcal{Q}_T(\mathcal{X})|_{q'}, \mathbb{O}_2(z_{q'}) \rangle$ account respectively for the contributions from the $\mathcal{Q}_T(\mathcal{X})$ part of $\mathcal{X}$ to the respective $\vartheta_q$ and $\vartheta_{q'}$ terms in Lemma 7.1's sum $\sum_{p \in \vartheta} z_p c_\# \varepsilon^2 \delta_p^2 \langle \mathcal{X}|_p, \mathcal{I}_p \rangle$.

*Part 6*: The upcoming Lemma 7.4 will be used with Lemma 7.3 to find points in $(z_q, z_{q'})$ from $Z_q \times Z_{q'}$ where $W_+(\mathfrak{g}_* + \underline{\mathfrak{h}}_{-\mathfrak{g}_*})$ is zero. To set the notation for this lemma, introduce $\hat{I}$ to denote the identity matrix in SO(4) and let $z_0$ to denote the point $(\frac{1}{2}, \hat{I}, \hat{I})$ in $(0,1) \times SO(4) \times SO(4)$.

**Lemma 7.4**: *Given $\mu \in (0, \infty)$ and the points q and q´, there exists $\kappa_\mu > 1$ with the following significance: Let $\wp$ denote a linear functional on $\mathcal{H}_a \oplus \mathcal{H}_b$ that depends continuously on the points in $Z_q \times Z_{q'}$. If $|\wp| \leq \kappa_\mu^{-1}$, there exist $(z_q, z_{q'}) \in Z_q \times Z_{q'}$ where*

$$\langle \mathfrak{X}|_q, \mathbb{O}_2(z_q) \rangle + \mu \langle \mathfrak{X}|_{q'}, \mathbb{O}_2(z_{q'}) \rangle + \wp(\mathfrak{X}) = 0$$

*for all $\mathfrak{X} \in \mathcal{H}_a \oplus \mathcal{H}_b$. Moreover, this equation holds with $(z_q, z_{q'})$ obeying*

$$\text{dist}(z_q, z_0) + \text{dist}(z_{q'}, z_0) \leq \kappa_\mu |\wp|.$$

This lemma is proved momentarily.

With the formula in Lemma 7.3 in hand, invoke Lemma 7.4 when R is large, $\varepsilon$ and $\delta_\diamond$ are small and $\delta \leq \varepsilon^{3/2}\delta_\diamond$ to see that there are parameter choices in $Z_q \times Z_{q'}$ where the corresponding $\mathcal{W}_+(\mathfrak{g}_* + \mathfrak{h}_{-\mathfrak{g}_*})$ is zero. Thus, if the kernel of $i^*$ in $H^1(S^3-K; V)$ is zero, then $X_K \#_4 \overline{\mathbb{CP}}^2$ has metrics with anti-self dual Weyl curvature. The second and third bullets of Proposition 4.2 follow from this last observation and what is said by Proposition 5.7.

*Part* 7: Metrics on $n > 4$ versions of $X_K \#_n \overline{\mathbb{CP}}^2$ with anti-self dual curvature are obtained in the following way: Let $\vartheta_0$ now denote the set of 4 points in the small $\delta$ versions of $X_{K,\delta}$ that was denoted in Parts 3-6 by $\vartheta$. Fix a set of n-4 distinct points in $X_{K,\delta}$ (to be denoted by $\vartheta_\#$) that is disjoint from $\vartheta_0$; then set $\vartheta = \vartheta_0 \cup \vartheta_\#$. Write the versions of $\delta_p$ for $p \in \vartheta_0$ in terms of $\delta_\diamond$ and $t$ and $t'$ as instructed by (7.12). As was done before, identify the Gaussian coordinate chart choices for the points in $\vartheta_0$ with the parameter space $(0, 1) \times (0, 1)$ for the pair $(t, t')$ with the space $Z_q \times Z_{q'}$. Fix small values for the $p \in \vartheta_\#$ versions of $\delta_p$ and any favorite Gaussian coordinate charts for these $\vartheta_\#$ points. With the parameters for $\vartheta_\#$ fixed, then the choice of the parameters $(z_q, z_{q'})$ in $Z_q \times Z_{q'}$ defines a family of $\mathfrak{g}^n_{R,\delta*}$ metrics parametrized by $Z_q \times Z_{q'}$. Supposing that $\delta$ and $\delta_\diamond$ and $\{\delta_p\}_{p \in \vartheta_\#}$ are sufficiently small, then Lemmas 6.1 and 7.1 can be invoked to see that there is an analog of Lemma 7.3's formula for any pair $(z_q, z_{q'}) \in Z_q \times Z_{q'}$ and any $\mathcal{X} \in \mathcal{H}_{K,R}$ but with $\wp$ depending on the chosen points in $Z_q \times Z_{q'}$ and, now, also on data for the points in $\vartheta_\#$; and with $\wp$'s norm now obeying

$$|\wp| \leq c_\Delta e^{-R} + c_R(\varepsilon + \varepsilon^{-2}\delta_\diamond^{-2}\delta^2 + \delta_\diamond^{-2}\sum_{p \in \vartheta_\#} z_p \delta_p^2) \, .$$
(7.13)

In this equation, $c_\Delta$ and $c_R$ are just like their namesakes in Lemma 7.3 except that $c_R$ can now depends also on the data that defines $\vartheta_\#$. Meanwhile, $z_p$ for each $p \in \vartheta_\#$ is positive and independent of the parameters $\varepsilon$, $\delta$, $\delta_\diamond$ and $\{\delta_p\}_{p \in \vartheta_\#}$. It is also independent of the chosen Gaussian coordinate charts for the points in $\vartheta_\#$. If R is large, $\varepsilon$ and $\delta_\diamond$ are small, $\delta < \varepsilon^{3/2}\delta_\diamond$ and $\delta_p < \varepsilon\delta_\diamond z_p^{-1}$, then Lemma 7.4 can again be invoked to see that there are parameter choices in $Z_q \times Z_{q'}$ where the $\mathfrak{g}_* = \mathfrak{g}^n_{R,\delta*}$ version of $\mathcal{W}_+(\mathfrak{g}_* + \mathfrak{h}_{-\mathfrak{g}_*})$ is zero.

*Proof of Lemma 7.4*: It proves useful to write the left hand side of the identity asserted by the lemma as the projection along a given $\mathcal{X} \in \mathcal{H}_a \oplus \mathcal{H}_b$ of a map from $Z_q \times Z_{q'}$ to $\mathcal{H}_a \oplus \mathcal{H}_b$. This map is denoted by $\Psi$. The map $\Psi$ can be written as $\Psi_0 + \psi$ where $\Psi_0$ is defined so that the projection of $\Psi_0(z_q, z_{q'})$ on any given vector $\mathcal{X}$ in $\mathcal{H}_a \oplus \mathcal{H}_b$ is

$$\langle \mathfrak{X}|_q, \mathbb{O}_2(z_q)\rangle + \mu \langle \mathfrak{X}|_{q'}, \mathbb{O}_2(z_{q'})\rangle \, .$$

(7.14)

The plan in what follows is to use the properties of the map $\Psi_0$ and some topological arguments to prove the lemma. To this end, the first observation is that $\Psi_0(z_0, z_0) = 0$. The second argument concerns the differential of $\Psi_0$ at $(z_0, z_0)$. Denote this differential by $d\Psi_0$. Suppposing that $\mathfrak{k} = (\mathfrak{z}_q, \mathfrak{z}_{q'})$ is a given tangent vector to $Z_q \times Z_{q'}$ at $(z_0, z_0)$, then $d\Psi_0(\mathfrak{k})$ is a vector in $\mathcal{H}_a \oplus \mathcal{H}_b$. Moreover, if $\mathfrak{X}_a \in \mathcal{H}_a$ and $\mathfrak{X}_b \in \mathcal{H}_b$, then the projection of $d\Psi_0(\mathfrak{k})$ to the vector $\mathfrak{X}_a + \mathfrak{X}_b$ can be written as

$$d\Psi_0(\mathfrak{k}) = \langle \mathfrak{X}_a, B\mathfrak{z}_{q'} + \mu A\mathfrak{z}_q\rangle + \langle \mathfrak{X}_b, C\mathfrak{z}_q\rangle \, ,$$

(7.15)

with A and B being respective homomorphisms from $TZ_q|_{z_0}$ and $TZ_{q'}|_{z_0}$ to $\mathcal{H}_a$, and with C being a homomorphism from $TZ_{q'}|_{z_0}$ to $\mathcal{H}_b$. There is no component to $d\Psi_0$ mapping $\mathfrak{z}_q$ to $\mathcal{H}_b$ because all of the elements in $\mathcal{H}_b$ vanish at the point q.

It follows from the second part of the first bullet in Lemma 7.2 that A, B and C are all isomorphisms. This implies that $d\Psi_0$ is also an isomorphism. Granted that $d\Psi_0$ is invertible, then Lemma 7.4 would follow from the inverse function theorem if the $C^1$ norm of $\psi = \Psi - \Psi_0$ is small. Since the assumptions of the lemma talk only about the sup-norm of $\psi$, this inverse function theorem argument can not be used. A homological argument is used instead. To begin the argument, let $\mathbb{V} \subset T(Z_q \times Z_{q'})|_{(z_0, z_0)}$ denote the orthogonal complement to the kernel of $d\Psi_0$, this being a 8 dimensional subspace. Use the exponential map for the product metric on $Z_q \times Z_{q'}$ to identify a ball about the origin in $\mathbb{V}$ with an embedded, 8 dimensional submanifold in $Z_q \times Z_{q'}$. Let $B_0$ denote this submanifold and let $\Psi_{00}$ denote the restriction of $\Psi_0$ to $B_0$. Since $d\Psi_{00}$ is invertible at $(z_0, z_0)$, there exists an open ball (to be denoted by $B_1$) in $B_0$ with the following properties: It has compact closure in $B_0$, its center is $(z_0, z_0)$, and it is mapped diffeomorphically by $\Psi_{00}$ onto an open neighborhood of the origin in $\mathcal{H}_a \oplus \mathcal{H}_b$. Let $B_2$ denote the closed ball in $Z_q \times Z_{q'}$ that is concentric to $B_1$ with half the radius of $B_1$. Having defined $B_2$, introduce $d_*$ to denote the smaller of the following two numbers:

- *The distance between $\Psi_{00}(\partial B_2)$ and the complement of $\Psi_{00}(B_1)$ in $\mathcal{H}_a \oplus \mathcal{H}_b$.*
- *The distance between $\Psi_{00}(\partial B_2)$ and the origin in $\mathcal{H}_a \oplus \mathcal{H}_b$.*

(7.16)

Suppose now that $\mu|\wp_q| + |\wp_{q'}| < \frac{1}{100} d_*$. Under this assumption, $\Psi_{00} + \lambda\psi$ maps $B_2$ into $\Psi_{00}(B_1)$ for all $\lambda \in [0, 1]$ because of the top bullet in (7.14). Moreover, when $\lambda \in [0, 1]$, then the whole boundary of $B_2$ is mapped by $\Psi_{00} + \lambda\psi$ to points in the complement of the origin in $\mathcal{H}_a \oplus \mathcal{H}_b$ because of the second bullet in (7.14). The following is a consequence of these two observations: If the origin in $\mathcal{H}_a \oplus \mathcal{H}_b$ were not in the image of the map $(\Psi_{00} + \lambda\psi)$ on $B_2$ for all values of $\lambda$ in $[0,1]$, then the family of maps $\{(\Psi_{00}|_{B_1})^{-1} \circ (\Psi_{00} + \lambda\psi)\}$ from $B_2$ into $B_1$ could be used to construct a homotopy of the identity map from $B_2$ to itself (rel $\partial B_2$) whose end member was a map from $B_2$ to $\partial B_2$. There is no such map because balls don't deformation retract onto their boundaries.

### d) The proof when the kernel of $\imath^*$ in $H^1(S^3-K; V)$ is not zero

The proof of Proposition 4.2 when kernel of $\imath^*$ is not zero has four parts.

*Part 1*: Fix $R > 1$ and large to define $\mathbb{T}_K$, its conformal structure $[\mathfrak{g}_R]$ and the $\iota$-invariant metric $\mathfrak{g}_R$ in the conformal structure $[\mathfrak{g}_R]$. Let $\mathcal{H}_{K,R}$ again denote the $\iota$-invariant kernel of the $\mathfrak{g} = \mathfrak{g}_R$ version of $\mathcal{L}_\mathfrak{g}^\dagger$. If R is sufficiently large (as will be assumed henceforth), then what is said by Propositions 5.3 and 5.4 about the vector space $\mathcal{H}_{K,R}$ can be invoked. In particular, Propositions 5.3 and 5.4 describes respective R independent vector spaces vector space $\mathcal{H}_T$ and $\mathcal{H}_K$, and homomorphisms $\mathcal{Q}_T$ and $\mathcal{Q}_K$ from $\mathcal{H}_{K,R}$ to $\mathcal{H}_T$ and $\mathcal{H}_K$ such that the direct product homomorphism $\mathcal{Q}_T \otimes \mathcal{Q}_K: \mathcal{H}_{K,R} \to \mathcal{H}_T \oplus \mathcal{H}_K$ is injective. The vector space $\mathcal{H}_K$ is isomorphic to the kernel of $\imath^*$ in $H^1(S^3-K; V)$. Let $N_1$ denote the dimension of this kernel (it is 1 less than $\dim(H^1(S^3-K; V))$).

*Part 2*: Let $\vartheta_0$ denote the set of 4 points in the small $\delta$ versions of $X_{K,\delta}$ that Part 3 of the previous subsection denoted by $\vartheta$. A second set of points is needed from the $(S^3-N_K) \times S^1$ part of $X_K$. This set is chosen using the following methodology: Fix for the moment a point from $(S^3-N_K) \times S^1$ where the function s is less than 1. Denoting this point by p, let $\mathbb{S}_p \subset (\Lambda^+ \otimes \Lambda^+)|_p$ denote the subspace of symmetric, traceless elements and let $\mathfrak{e}_p: \mathcal{H}_K \to \mathbb{S}_p$ denote the restriction map. Fix a finite set in the $s \leq 1$ part of $(S^3-N_K) \times S^1$ to be denoted by $\vartheta_*$ so that the homomorphism

$$\mathfrak{e} = \oplus_{p \in \vartheta_*} \mathfrak{e}_p: \mathcal{H}_K \to \oplus_{p \in \vartheta_*} \mathbb{S}_p$$

(7.17)

is injective. Sets with this property exist with $N_1$ or fewer points and one of the latter should be chosen for $\vartheta_*$. Define from $\vartheta_*$ a larger set (to be denoted by $\vartheta_1$) that has three points for each point in $\vartheta_*$. If $p \in \vartheta_*$, then the corresponding three points are distinct and they are much closer to p then they are to any other point in $\vartheta_*$. Having specified a priori a positive number to be denoted by $\rho$, then the three points associated to p should have distance $\rho$ or less from p. An upper bound for $\rho$ will be specified momentarily.

The set $\vartheta = \vartheta_0 \cup \vartheta_1$ will be used to construct a metric on $X_K \#_{4+3N_1} \overline{\mathbb{CP}}^2$ with anti-self dual Weyl curvature.

*Part 3*:  The construction of an $n = 4 + 3N_1$ version of the conformal class $[\mathfrak{g}^n{}_{R,\delta}]$ using the set $\vartheta$ requires the choice of the positive numbers $\{\delta_p\}_{p \in \vartheta}$ and a suitable Gaussian coordinate chart at each point from $\vartheta$. The $p \in \vartheta_0$ versions of $\delta_p$ are again written in terms of $\delta_\lozenge$ and $t$ and $t'$ as in (7.12). Also as before, the space of Gaussian coordinate charts for the points in $\vartheta_0$ and the parameter space $(0,1) \times (0,1)$ for the pair $(t, t')$ are identified with the points in the space $Z_q \times Z_{q'}$.

Now suppose that $p \in \vartheta_*$. The point p corresponds to three points in $\vartheta_1$ that will be labeled as $\{p_1, p_2, p_3\}$. The parameters $\{\delta_{p_a}\}_{a=1,2,3}$ for these points will be written using the number $\delta_\lozenge$ and a pair of numbers $(t_1, t_2)$ in the simplex $\Delta^2$ (from Section 7b). The upcoming formula in (7.18) uses $t_3$ to denote $1 - t_1 - t_2$. The formula that follows also involves the numbers $\{z_{p_a}\}_{a=1,2,3}$ that appear in (7.3) and Lemma 7.1:

$$\delta_{p_a} = (z_{p_a}^{-1} t_a)^{1/2} \delta_\lozenge.$$

(7.18)

A Gaussian coordinate chart centered at each point in the set $\{p_a\}_{a=1,2,3}$ must also be chosen to define the desired version of $[\mathfrak{g}^n{}_{R,\delta}]$. This is a chart for a flat metric that defines the conformal class $[\mathfrak{g}_R]$ near the relevant point. Since the points are close together and very close to p, an oriented, orthonormal frame for the tangent space at p defines respective orthonormal frames for the tangent spaces at each of the three nearby points $\{p_a\}_{a=1,2,3}$. This is done by parallel transporting the frame at p using $\mathfrak{g}_R$'s Levi-Civita connection to the nearby point along the short $\mathfrak{g}_R$ geodesic between them. Choose once and for all an orthonormal frame for the tangent space at p and use the induced frame at each of point from the set $\{p_a\}_{a=1,2,3}$ to identify the orthonormal frame bundle at the point with a copy of the group SO(4). This in turn identifies the set of choices of Gaussian coordinate chart at any given point in $\{p_a\}_{a=1,2,3}$ with SO(4). The three SO(4) parameters and the simplex parameters $(t_1, t_2)$ define a point in the space $\Delta^2 \times (\times_3 SO(4))$. The version of this space that is associated to a given point p from $\vartheta_*$ will be denoted by $W_p$.

Suppose again that $p \in \vartheta_*$. The a priori choice of a point in the orthonormal frame bundle for the tangent bundle at p identifies the space of symmetric, traceless elements in $(\Lambda^+ \otimes \Lambda^+)|_p$ with the vector space $\mathbb{S}$. Meanwhile, the parallel transport alluded to in the previous paragraph identifies the subspace of symmetric, traceless elements in each $a \in \{1, 2, 3\}$ version of $(\Lambda^+ \otimes \Lambda^+)|_{p_a}$ with the corresponding subspace of $(\Lambda^+ \otimes \Lambda^+)|_p$, and thus also with the vector space $\mathbb{S}$. These identifications are used implicitly in the upcoming Equation (7.20).

*Part 4*: Fix R to be large and $\varepsilon \in (0, 1]$. Supposing that $\delta$ and $\delta_\diamond$ are sufficiently small, and supposing that points $(z_q, z_{q'}) \in Z_q \times Z_{q'}$ have been chosen and likewise a point (to be denoted by $w_p$) in each $p \in \vartheta_*$ version of $W_p$, then the constructions in Section 3 supply a conformal structure $[\mathfrak{g}^n_{R,\delta}]$ for $n = 4 + 3N_1$ on the space $X_K \#_{4+3N_1} \overline{\mathbb{CP}^2}$ and the corresponding version of the metric $\mathfrak{g}^n_{R,\delta*}$ in this conformal structure. Lemmas 6.1 and 7.1 can again be used to study the resulting obstruction vector (the $\mathfrak{g}_* = \mathfrak{g}^n_{R,\delta*}$ version of $\mathcal{W}_+(\mathfrak{g}_* + \underline{\mathfrak{h}}_{-\mathfrak{g}_*})$) when $\delta$ and $\delta_\diamond$ are small.

Lemmas 6.1 and 7.1 imply in particular that if $X \in \mathcal{H}_{K,R}$, then the $L^2$ orthogonal projection of the obstruction vector $\mathcal{W}_+(\mathfrak{g}_* + \underline{\mathfrak{h}}_{-\mathfrak{g}_*})$ to the element that is defined by $X$ in the $\mathfrak{g}^n_{R,\delta*}$ version of $\underline{\mathcal{V}}_*'$ can be written as a sum of three terms which are denoted here by $\mathcal{W}_1$, $\mathcal{W}_2$ and $\mathcal{W}_3$. These terms are described in the subsequent paragraphs.

What was denoted by $\mathcal{W}_1$ in the preceding paragraph is written below in (7.20). This upcoming formula for $\mathcal{W}_1$ writes a given element $X \in \mathcal{H}_{K,R}$ with $L^2$ norm equal to 1 (as defined by $\mathfrak{g}_R$) near the points from $\vartheta_\diamond$ as $\mathcal{Q}_T(X) + \mathfrak{r}_T$ such that (7.11) holds. Meanwhile, it writes $X$ near the points of $\vartheta_1$ as $\mathcal{Q}_K(X) + \mathfrak{r}_K$ with the $g_K$ norm of $\mathfrak{r}_K$ obeying

$$|\mathfrak{r}_K| \leq \ \leq \kappa e^{2s-R} \tag{7.19}$$

on the $s \leq \frac{1}{2} R$ part of $(S^3 - N_K) \times S^1$.

The promised formula for $\mathcal{W}_1$ is

$$\mathcal{W}_1 = c\varepsilon^2 \delta_\diamond^2 \left( \langle \mathcal{Q}_T(X)|_q, \mathbb{O}_2(z_q) \rangle + \langle \mathcal{Q}_T(X)|_{q'}, \mathbb{O}_2(z_{q'}) \rangle \right) \tag{7.20}$$

with $c$ being the positive number that appears in Lemma 7.3's formula. What was denoted by $\mathcal{W}_2$ in the preceding paragraph is

$$\mathcal{W}_2 = \varepsilon^2 \delta_\diamond^2 \sum_{p \in \vartheta_*} \langle \mathcal{Q}_K(X)|_p, \mathbb{O}_3(w_p) \rangle \tag{7.21}$$

with $\mathbb{O}_3$ being the function on $\Delta^2 \times (\times_3 SO(4))$ that is defined in (7.8). What was denoted above by $\mathcal{W}_3$ is the value on $X$ of a linear functional on $\mathcal{H}_{K,R}$ (to denoted by $\wp$) whose dual norm obeys the bound

$$|\wp| \leq c_\Delta \varepsilon^2 \delta_\diamond^2 e^{-R} + c_R \varepsilon^2 \delta_\diamond^2 (\rho + \varepsilon + \varepsilon^{-2} \delta_\diamond^{-2} \delta^2) \tag{7.22}$$

with $\rho$ denoting the minimum of the distances between any given $p \in \vartheta_*$ version of the three element set $\{p_a\}_{a=1,2,3}$. What is denoted in (7.22) by $c_\Delta$ is a positive number that is independent of the parameters R, ρ, δ, ε, $δ_◊$ and the chosen points in $Z_q \times Z_{q'} \times (\times_{p \in \vartheta_*} W_p)$. The number $c_R$ can depend on R, but not on ρ, δ, ε or $δ_◊$. One last point of note: Although the norm of $\wp$ has an upper bound that is independent of the chosen point in $Z_q \times Z_{q'} \times (\times_{p \in \vartheta_*} W_p)$, the linear functional $\wp$ can depend on the chosen point in this space. In any event, $\wp$ necessarily varies continuously as a Hom($\mathcal{H}_{K,R}$; ℝ) valued function on the space $Z_q \times Z_{q'} \times (\times_{p \in \vartheta_*} W_p)$. The definition of $\wp$ mimicks the definition of its namesake in Lemma 7.3.

*Part 5*: Let $\hat{I}$ again denote the identity matrix in SO(4) and let $z_0$ again denote the point ($\frac{1}{2}$, $\hat{I}$, $\hat{I}$) in $(0,1) \times SO(4) \times SO(4)$. Introduce $w_0$ to denote the point $((\frac{1}{3}, \frac{1}{3}), \hat{I}, \hat{I}, \hat{I})$ in $\Delta^2 \times (\times_3 SO(4))$. Let $\Psi_0$ now denote the map from $Z_q \times Z_{q'} \times (\times_{p \in \vartheta_*} W_p)$ to $\mathcal{H}_T \oplus \mathcal{H}_K$ that is defined as follows: View $\mathcal{H}_T$ and $\mathcal{H}_K$ as respective subspaces in the kernels of $\mathcal{L}^\dagger_{g_T}$ and $\mathcal{L}^\dagger_{g_K}$. Let $(\mathfrak{X}_T, \mathfrak{X}_K)$ denote a given vector in $\mathcal{H}_T \oplus \mathcal{H}_K$. Then the $L^2$-orthogonal projection along $(\mathfrak{X}_T, \mathfrak{X}_K)$ (as defined by the metrics $g_T$ and $g_K$) of the $\Psi_0$ image of any given point in its domain is given by

$$c\varepsilon^2 δ_◊^2 \left( \langle \mathfrak{X}_T|_q, \mathbb{O}_2(z_q) \rangle + \langle \mathfrak{X}_T|_{q'}, \mathbb{O}_2(z_{q'}) \rangle + \varepsilon^2 δ_◊^2 \sum_{p \in \vartheta_*} \langle \mathfrak{X}_K|_p, \mathbb{O}_3(w_p) \rangle \right).$$
(7.23)

Since the map in (7.17) is surjective, it follows from Lemma 7.2 that the differential of $\Psi_0$ is surjective onto $\mathcal{H}_T \oplus \mathcal{H}_K$ at the point in the domain where $z_q$ and $z_{q'}$ are equal to $z_0$ and each $p \in \vartheta_*$ version of $w_p$ is equal to $w_0$. Granted this observation, granted (7.22) and granted that $\wp$ varies continuously as a function on the space $Z_q \times Z_{q'} \times (\times_{p \in \vartheta_*} W_p)$, then much the same argument that is used to prove Lemma 7.4 proves the following: If R is sufficiently large and ε is sufficiently small (with an upper bound that is independent of R), and if $δ_◊$, $(\varepsilon δ_◊)^{-1}δ$ and ρ are sufficiently small, then there are points in $Z_q \times Z_{q'} \times (\times_{p \in \vartheta_*} W_p)$ where $\mathcal{W}_1 + \mathcal{W}_2 + \mathcal{W}_3$ is zero. By definition, these points parametrize metrics on $X_K \#_{4+3N_1} \overline{\mathbb{CP}}^2$ with anti-self dual Weyl curvature tensor.

*Part 6*: Suppose now that $n > 4+3N_1$. Metrics on $X_K \#_n \overline{\mathbb{CP}}^2$ with anti-self dual Weyl curvature tensor can be constructed by mimicking the arguments from Part 7 of the previous subsection. By way of a summary, a set $\vartheta_\#$ of n - $(4+3N_1)$ additional points must be chosen in $X_K$ that are distinct from the points in the set $\vartheta_0 \cup \vartheta_1$. The set $\vartheta = \vartheta_0 \cup \vartheta_1 \cup \vartheta_\#$ is then used to construct the desired metric. The $p \in \vartheta_0 \cup \vartheta_1$ versions of $δ_p$ and the Gaussian coordinate charts at these points are parametrized as before using $δ_◊$ and the points in the space $Z_q \times Z_{q'} \times (\times_{p \in \vartheta_*} W_p)$. Meanwhile, if $p \in \vartheta_\#$, then the corresponding $δ_p$ parameter should be chosen to be very small. Any favorite Gaussian coordinate chart can be chosen for these points from $\vartheta_*$. With all of the $\vartheta_*$ parameters fixed, then the corresponding obstruction vector for the associated metric $\mathfrak{g}^n_{R,\delta*}$ can again be written as $\mathcal{W}_1 + \mathcal{W}_2 + \mathcal{W}_3$ with $\mathcal{W}_1$ and $\mathcal{W}_2$ as before (see (7.20) and (7.21); and with $\mathcal{W}_3$ defined by a linear function $\wp$ on $\mathcal{H}_T \oplus \mathcal{H}_K$ that varies continuously with respect to the parameters in $Z_q \times Z_{q'} \times (\times_{p \in \vartheta_*} W_p)$ and obeys

$$|\wp| \leq c_\Delta \varepsilon^2 \delta_\Diamond^2 e^{-R} + c_R(\rho + \varepsilon + \varepsilon^{-2} \delta_\Diamond^{-2} \delta^2 + \delta_\Diamond^{-2} \sum_{p \in \vartheta_\#} z_p \, \delta_p^{\,2}) \, .$$

(7.24)

In this equation, $z_p$ for $p \in \vartheta_\#$ comes from Lemma 7.1.

Granted the bound in (7.24), then the argument used in Part 5 can be repeated when R is large, $\varepsilon$ and $\delta_\Diamond$ are small, $\delta < \varepsilon^{3/2} \delta_\Diamond$ and $\delta_p < \varepsilon \delta_\Diamond z_p^{-1}$ to see that there are parameters in $Z_q \times Z_{q'} \times (\times_{p \in \vartheta_*} W_p)$ where the $\mathfrak{g}_* = \mathfrak{g}^n_{R,\delta*}$ version of $\mathcal{W}_+(\mathfrak{g}_* + \underline{\mathfrak{h}}_{\mathfrak{g}_*})$ is zero.

## 8. Gromov-Hausdorff limits

The introduction to this paper (Section 1) remarked about the Gromov-Hausdorff limit of anti-self dual Weyl curvature metrics on $X_K$. These remarks refer to the limit metric spaces that are obtained from Proposition 4.2's metrics by taking $\delta \to 0$ and all of the collection $\{\delta_p\}_{p \in \vartheta}$ limit to zero also. The proposition that follows makes a precise statement to this effect.

**Proposition 8.1**: *Let K denote a hyperbolic knot in $S^3$ and let N denote the integer that appears in K's version of Proposition 4.2. If $n \geq N$, then $X_K \#_n \overline{\mathbb{CP}}^2$ has a sequence of anti-self dual Weyl curvature metrics with the following properties:*
- *All elements in the sequence have the same volume and there is an priori $L^2$ bound for the Riemann curvature tensor.*
- *The sequence converges in the Gromov-Hausdorff topology to the complement of n-4 points in $(S^3-K) \times S^1$ with the metric space structure coming from the product of the hyperbolic metric on $S^3 - K$ and the length $2\pi$ metric on $S^1$.*

*Proof of Proposition 8.1*: The metrics in the desired sequence have the form $\mathfrak{g}^n_{R,\delta} + \mathfrak{h}$ as described in Proposition 4.1 with the parameters that define $\mathfrak{g}^n_{R,\delta}$ chosen according to the rules that are listed in the next paragraph.

Choose a sequence $\{R_k\}_{k \in \{1,2...\}}$ that is increasing and unbounded. Supposing that k is a positive integer, then the k'th metric in the desired sequence will have $R = R_k$. The parameters $(\varepsilon, \delta, \delta_\Diamond, \{\delta_p\}_{p \in \vartheta_\#})$ for the k'th metric are chosen subject to the constraints listed below.

- $\varepsilon = \varepsilon_0$ *is independent of* k.
- $\delta_\diamond = \delta_{\diamond(k)}$ *is sufficiently small given* $R_k$*, and* $\lim_{k\to\infty} \delta_{\diamond(k)} = 0$ .
- $\delta \le \frac{1}{1000}\varepsilon^{3/2}\delta_\diamond$.
- $\delta_p \le \frac{1}{1000}\varepsilon\delta_\diamond$ *for all* $p \in \vartheta_\#$

(8.1)

If k is sufficiently large, then what is said in Lemmas 7.3 and 7.4 can be invoke when kernel $\hat{\imath}*$ is zero, and what is said in Parts 4-6 of Section 7e when the kernel of $\hat{\imath}*$ is not zero to find SO(4) parameters and the $p \in \vartheta_0 \cup \vartheta_1$ versions of $\delta_p$ so as to obtain a metric of the form $\mathfrak{g}^n_{R,\delta}+\mathfrak{h}$ as described in Item i) of Proposition 4.1 with anti-self dual Weyl curvature. Note in particular that the construction of this metric $\mathfrak{g}^n_{R,\delta} + \mathfrak{h}$ uses values of $\delta_p$ for $p \in \vartheta_0 \cup \vartheta_1$ that are bounded by $c_\diamond \delta_\diamond$ with $c_\diamond$ here (and below) indicating a number that is greater than 1 and independent of $R_n$ and $\delta_\diamond$. Therefore, all versions of $p \in \vartheta$ versions of $\delta_p$ are bounded by $c_\diamond \delta_\diamond$.

Fix a positive integer k and write $R = R_k$ and $\delta_\diamond = \delta_{\diamond(k)}$. It follows from the definition of $\mathfrak{g}^n_{R,\delta*}$ in Section 4a) and from what is said in Lemmas 3.3 and 3.4 that the $L^6$ norm of the self-dual part of the Weyl curvature of $\mathfrak{g}^n_{R,\delta*}$ is bounded by $c_\diamond \delta_\diamond^2$. This leads via the analysis in [KS] to a corresponding $c_\diamond \delta_\diamond^2$ bound for the $L^6_2$ norm (as measured using $\mathfrak{g}^n_{R,\delta*}$) of what is denoted by $\mathfrak{h}_*$ in Item ii) of Proposition 4.2. There are also to $c_\diamond \delta_\diamond^2$ pointwise bounds for $\mathfrak{h}_*$ as is noted by Equation (5.6) in [KS].

Define $\mathfrak{h}$ as in Item i) of Proposition 4.2. The Gromov-Hausdorff convergence with the asserted limit of the sequence whose k'th element is the integer k version of $\mathfrak{g}^n_{R,\delta}+\mathfrak{h}$ follows from what is said by Items i) and ii) of Proposition 4.2. This sequence has volume bounded from above and below because this is the case for the corresponding $\mathfrak{g}^n_{R,\delta}$ sequence.

To see about the $L^2$ norm bound for the Riemann curvature tensor, first use the fact that the volume of $X_K$ as measured by $\mathfrak{g}^n_{R,\delta*}$ is bounded by $c_\diamond |\ln \delta_\diamond|$ and the $L^6_2$ norm of $\mathfrak{h}_*$ is bounded by $c_\diamond \delta_\diamond^2$ to obtain a $c_\diamond \delta_\diamond^2 |\ln \delta_\diamond|^{1/3}$ bound for the $L^2_2$ norm of $\mathfrak{h}_*$. This $L^2_2$ bound leads to a corresponding $c_\diamond \delta_\diamond^2 |\ln(\delta_\diamond)|^{1/3}$ bound for the $L^2_2$ norm (as measured by $\mathfrak{g}^n_{R,\delta}$) of $\mathfrak{h}$. (The $L^2_2$ norm of sections of $\Lambda^+ \otimes \Lambda^-$ is invariant with respect to constant conformal changes of the metric; but it is not invariant with respect to non-constant conformal changes. As a consequence, the asserted bound for $L^2_2$ norm of $\mathfrak{h}$ requires a calculation to verify. This task is left to the reader.) Since the $L^2$ norm of the Riemann curvature of $\mathfrak{g}^n_{R,\delta}$ is bounded independent of the index k, the afore-mentioned $L^2_2$ norm bound for $\mathfrak{h}$ implies in turn that there is a k independent bound for the $L^2$ norm of the Riemann curvature tensor of $\mathfrak{g}^n_{R,\delta}+\mathfrak{h}$.

**APPENDIX: The obstruction vector space**

This appendix has five sections, with the first containing an outline of the proof of Propositions 5.3 and 5.4 and the last containing the details of the proof. The intervening three sections supply the background that is needed. Here is the formal table of contents for the appendix:

SECTION A1: Outline of the proofs of Propositions 5.3 and 5.4.

SECTION A2: The $\mathcal{X}_T$ part of a pair $(\mathcal{X}_T, \mathcal{X}_K)$.

SECTION A3: The $\mathcal{X}_K$ part of a pair $(\mathcal{X}_T, \mathcal{X}_K)$.

SECTION A4: Approximation by $(S^3-K) \times S^1$ solutions with finitely many modes.

SECTION A5: Proofs of Propositions 5.3 and 5.4.

Note that Sections A3 and A4 may have independent interest, especially Section A3 which gives a detailed description of the kernel of the operator $\mathcal{L}_{g_K}^\dagger$ on $(S^3-K) \times S^1$.

**A1. Outline of the proofs of Propositions 5.3 and 5.4**

What follows is a four part outline of the arguments for these propositions.

*Part 1*: The metric $g_T$ on $T \times T'$ is invariant under the action of T by translations the T factor in $T \times T'$. It follows as a consequence that if $U \subset T'$ is any given open set, then any element in the kernel of the operator $\mathcal{L}_{g_T}^\dagger$ on $T \times U$ can be written as a sum of Fourier modes with respect to this T action with each mode in the kernel of this $\mathcal{L}_{g_T}^\dagger$. The T-action Fourier decomposition of a given element $\mathcal{X}$ has the form

$$\mathcal{X} = \sum_{(k_1, k_2) \in \mathbb{Z}} \mathfrak{X}_{(k_1, k_2)}(t_3, t_4)\, e^{i(k_1 t_1 + k_2 t_2)} \ .$$

(A1.1)

The proofs of both Propositions 5.3 and 5.4 use the fact that if R is large, then an element in the kernel of $\mathcal{L}_{g_T}^\dagger$ on $T \times (\mathbb{D}-D)$ with no $(k_1 = 0, k_2 = 0)$ Fourier mode is necessarily very small away from the boundaries of this domain relative to its size near these boundaries. To say this precisely, let $\Delta$ denote the distance to the boundary of the annulus $\mathbb{D}-D$. Then let $A \subset \mathbb{D}-D$ denote the set where $\Delta \leq c e^{-R}$ with $c$ being somewhat greater than 1, but independent of the value of R that is used to define the metric $g_T$. The following is proved momentarily in Section A2b: Let $X$ denote an element in the kernel of $\mathcal{L}_{g_T}^\dagger$ on $T\times(\mathbb{D}-D)$ with no $(k_1=0, k_2=0)$ Fourier mode. Then

$$|X| \leq c_* \exp(-e^R \eth/c_*) \, \Big( \int_{T\times A} |X|^2 \Big)^{1/2}$$

(A1.2)

with $c_*$ being greater than 1 and independent of both $X$ and R. This bound is (ultimately) a consequence of the fact that the norm (as measured using $g_T$) of the matrix of second derivatives of the function $(t_1, t_2) \to e^{i(k_1 t_1 + k_2 t_2)}$ is greater than an R and $X$ independent positive multiple of $e^{2R}$ when at least one of $k_1$ or $k_2$ is not zero. A Bochner-Weitzenboch formula for $\mathcal{L}_{g_T}^\dagger$ is also needed for the proof of (A1.2).

The bound in (A1.2) leads to a similar bound for the norm of an $\iota$ invariant element in the kernel of $\mathcal{L}_{g_T}^\dagger$ on the whole $T\times(T'-(D\cup\iota(D)))$. The latter bound is given by (A1.2) with $\eth$ being the distance in $T'-(D\cup\iota(D)))$ to $\partial D\cup\iota(\partial D)$ and with A being the subset of $\mathbb{D}-D$ where the distance to $\partial D$ is less than $ce^{-R}$. The $T\times(T'-(D\cup\iota(D)))$ version of (A1.2) leads to the following observation: Supposing that $X$ is an $\iota$-invariant element in the kernel of $\mathcal{L}_{g_T}^\dagger$ on $T\times(T'-(D\cup\iota(D)))$, write it as a pair $(X_T, X_K)$ in the manner of (5.4) and (5.5). Then the $X_T$ part of $X$ is mostly the $(k_1=0, k_2=0)$ Fourier mode for the T action on the domain $T\times(T'-(D\cup\iota(D)))$.

*Part 2*: The observations from the preceding paragraph is important by virtue of the fact that the $(k_1=0, k_2=0)$ Fourier mode are independent of the coordinates $(t_1, t_2)$ on T; and the equation $\mathcal{L}_{g_T}^\dagger(\cdot) = 0$ for such a mode is thus an equation on $T'-(D\cup\iota(D)))$. As it turns out, the solutions to the latter equation can be written explicitly in terms of the functions $x$ and $u$ from Lemmas 5.1 and 5.2: Any given solution is a pointwise limit of solutions whose components have the form

- $X^{33} = -(X^{11} + X^{22}) = s$
- $X^{13} - iX^{23} = c_0 + c_1 x + c_2 x^2 + \cdots$
- $X^{11} - iX^{12} = (a_0 + a_1 x + a_2 x^2 + \cdots) u + (\flat_0 + \flat_1 x + \flat_2 x^2 + \cdots)$

(A1.3)

with $s$ being a real number; and with $\{c_0, c_1, \ldots\}$ and $\{a_0, a_1, \ldots\}$ and $\{\flat_0, \flat_1, \ldots\}$ being complex numbers $\mathbb{C}$. This is proved in the upcoming Section A2c.

*Part 3*: There is no global T action on the $(S^3-K)\times S^1$ part of $X_K$, but there is the action of the $S^1$ factor on itself; and this is an isometric action for the metric $g_K$. As a consequence, if $U \subset S^3 - K$ is a given open set, then each element in the kernel of $\mathcal{L}_{g_K}^\dagger$ on $U \times S^1$ can be written as a sum of $S^1$ action Fourier modes from the kernel of $\mathcal{L}_{g_K}^\dagger$. The Fourier modes in this case are labeled by the integers; and supposing that n is an integer, then the corresponding Fourier mode has the form

$$\mathfrak{X} e^{in\theta} \tag{A1.4}$$

with $\theta$ denoting here the $\mathbb{R}/(2\pi\mathbb{Z})$ coordinate for the $S^1$ factor of $U \times S^1$ and with $\mathfrak{X}$ denoting a tensor that depends only on the points in the U factor. Of particular interest here are the cases where $U = S^3 - N_K$ (and, as it turns out, where U is the whole of $S^3 - K$.)

Now, the $s \geq 0$ part of $S^3 - K$ (and hence of $S^3 - N_K$) does have an isometric action of the torus T because the metric $g_K$ has the form $g + d\theta^2$ with g as depicted in (3.3). The resulting T action on the $s \geq 0$ part of $(S^3 - N_K) \times S^1$ is the same as that on $T \times (\mathbb{D} - D)$ when these spaces are identified in the manner of Section 3b. A key fact is that this T action commutes with the $S^1$ action. Therefore,, any given $S^1$ action Fourier mode of an element in the kernel of $\mathcal{L}_{g_K}^\dagger$ on the $s \geq 0$ part of $(S^3 - N_K) \times S^1$ or on the $s \geq 0$ part of $(S^3 - K) \times S^1$ can be written as a sum of Fourier modes for the T action. The inequality in (A1.2) suggests that the $(k_1 = 0, k_2 = 0)$ mode for the T action will dominate the other T action modes where s is large, and this turns out to be the case except for the $n = \pm 1$ Fourier modes when the kernel of $i^*$ is not zero. More is said about this momentarily.

In any event, (A1.2) has the following consequence (proved in Section A4): If R is sufficiently large, then any $S^1$ action Fourier mode on $(S^3 - N_K) \times S^1$ that is consistent with (A1.2) is well approximated (except near the $\partial N_K \times S^1$) by an element in the kernel of $\mathcal{L}_{g_K}^\dagger$ on the whole of $(S^3 - K) \times S^1$ whose non-trivial T action modes are bounded and square integrable on the $s \geq 0$ part of $(S^3 - K) \times S^1$. This last observation is useful by virtue of the fact that the elements in the kernel of $\mathcal{L}_{g_K}^\dagger$ on $(S^3 - K) \times S^1$ with the non-trivial T action modes being bounded (or square integrable) can be completely determined for each integer n. With regards to $n = \pm 1$: It is only in this case that there are elements in the kernel of $\mathcal{L}_{g_K}^\dagger$ on $(S^3 - K) \times S^1$ of the sort just described that are square integrable on the whole of $(S^3 - K) \times S^1$. These elements account for the appearance of the kernel of $i^*$ in the Proposition 5.4. Section A3 proves the preceding assertions about the kernel of the $\mathcal{L}_{g_K}^\dagger$ on $(S^3 - K) \times S^1$. The formalism from [AV] helps to do this.

*Part 4*: As noted, the elements in the kernel of $\mathcal{L}_{g_K}^\dagger$ on $(S^3 - N_K) \times S^1$ that are consistent with (A1.2) are well approximated by elements in the kernel of $\mathcal{L}_{g_K}^\dagger$ on $(S^3 - K) \times S^1$ with the non-trivial T action modes being bounded where $s \geq 0$. As also noted, the latter set can be completely determined. Moreover, the components of the relevant version of the tensor $\mathfrak{X}$ in (A1.4) for any given mode number n can be well approximated for s large by a specific n-dependent functions of s (it is a sum of exponentials). This behavior of $\mathfrak{X}$ on the $s \geq 0$ part of $(S^3 - N_K) \times S^1$ can be compared with the behavior of (A1.3) on $\mathbb{D} - D$ using the formulas in (5.4) and (5.5). The comparison finds that most of the terms in (A1.3) do not have extensions over $(S^3 - N_K) \times S^1$; and the comparison finds that those that do extend are accounted for by Propositions 5.3 and 5.4.

## A2. The $\mathcal{X}_T$ part of a pair $(\mathcal{X}_T, \mathcal{X}_K)$

Fix an element (to be denoted by $\mathcal{X}$) in the $\iota$-invariant kernel of the $\mathfrak{g} = \mathfrak{g}_R$ version of $\mathcal{L}^\dagger_\mathfrak{g}$. This element is written in the manner of (5.4) and (5.5) as a pair $(\mathcal{X}_T, \mathcal{X}_K)$. This section helps set the stage for the proofs of Propositions 5.3 and 5.4 by analyzing the $\mathcal{X}_T$ part of the pair $(\mathcal{X}_T, \mathcal{X}_K)$. Propositions A2.1, A2.6 and A2.7 are the key results in this section of the appendix.

### a) The operator $\mathcal{L}^\dagger_\mathfrak{g}$

Supposing that X is an oriented, Riemannian 4-manifold, denote its metric by $\mathfrak{g}$. The formal definition in (5.1) and (5.2) of the operator $\mathcal{L}^\dagger_\mathfrak{g}$ is mostly useless for the task of finding its cokernel. A useful definition writes $\mathcal{L}^\dagger_\mathfrak{g}$ as a differential operator using a chosen (local) orthonormal frame for $T^*X$ and a corresponding orthonormal frame for $\Lambda^+$. This is what is done in this subsection.

Fix an open set in X where there is an oriented, orthonormal frame for $T^*X$. The four basis 1-forms are denoted by $\{e^1, e^2, e^3, e^4\}$. These can be used to define respective orthonormal frames $\{\omega^a\}_{a=1,2,3}$ and $\{\underline{\omega}^a\}_{a=1,2,3}$ for $\Lambda^+$ and $\Lambda^-$ using the formulas in (5.7). These frames then define corresponding orthonormal frames for the tensor bundles $\Lambda^+ \otimes T^*X$ and $\Lambda^+ \otimes \Lambda^+$ and $\Lambda^+ \otimes \Lambda^-$.

Differential operators taking sections of one tensor bundle to another appear as matrix valued operators acting on maps to Euclidean vector spaces when the sections of the bundles are depicted using orthonormal frames. Two matrices play a central role in the upcoming formula for $\mathcal{L}^\dagger_\mathfrak{g}$. These are the two matrices that map $\wedge^2 \mathbb{R}^4$ to $\mathbb{R}^3$ giving the respective projections to the vector subspaces of self dual and anti-self dual elements in $\wedge^2 \mathbb{R}^4$. (These are known to physicists as the self dual and anti-self dual 't Hooft symbols.) The self dual version has components $\{\eta^a_{ik}\}_{a=1,2,3; i,k=1,2,3,4}$ and are defined by writing the self dual basis in (5.7) as

$$\{\omega^a = \tfrac{1}{2\sqrt{2}} \, \eta^a_{ij} e^i \wedge e^j \, \}_{a=1,2,3}$$

(A2.1)

with it understood that the repeated Latin indices are summed over their 4 values. This summation convention for repeated indices is used subsequently with no further comment. The anti-self dual version has components $\{\eta^a_{ik}\}_{a=1,2,3; i,k \in \{1,2,3,4\}}$ and it is defined by writing the anti-self dual basis in (5.7) as

$$\{\omega^a = \tfrac{1}{2\sqrt{2}} \, \eta^a_{ij} e^i \wedge e^j \}_{a=1,2,3} \, .$$

(A2.2)

The bundle of symmetric, traceless elements in $T^*X \otimes T^*X$ is isometric to the bundle $\Lambda^+ \otimes \Lambda^-$; and the 't Hooft symbols implement this isometry as follows: Supposing that $h_{ik} e^i \otimes e^k$ is a symmetric, traceless section of $T^*X \otimes T^*X$, then the corresponding section of $\Lambda^+ \otimes \Lambda^-$ is $\tfrac{1}{2}(\eta^a_{im} \eta^c_{km} h_{ik}) \omega^a \otimes \omega^c$. The inverse of this identification writes any given section $t^{ac} \omega^a \otimes \omega^c$ of $\Lambda^+ \otimes \Lambda^-$ as $\tfrac{1}{2}(\eta^a_{im} \eta^c_{km} t^{ac}) e^i \otimes e^k$. This identification of bundles is used implicitly in what follows. By way of a relevant example, the Ricci tensor for the metric $\mathfrak{g}$ can be viewed as a symmetric section of $T^*X \otimes T^*X$ and so its traceless part can be viewed as a section of $\Lambda^+ \otimes \Lambda^-$. Half of the latter section is denoted by $\mathcal{B}$ in (3.1). The components of $\mathcal{B}$ with respect to the given frames for $\Lambda^+$ and $\Lambda^-$ are $\{\mathcal{B}^{ac}\}_{a,c=1,2,3}$.

The directional covariant derivatives along the frame vectors for $TX$ that give the dual basis to $\{e^i\}_{i=1,2,3,4}$ are written as $\{\nabla_i\}_{i=1,2,3,4}$. The notation in what follows writes covariant derivatives of sections as follows: Supposing that $\mathfrak{F}$ is a tensor bundle of some rank N with a local frame $\{\mu^A\}_{A=1,2,\ldots,N}$, and supposing that a section $\phi$ is written with respect to this basis as $\phi^A \mu^A$, then the components of the covariant derivative $\nabla \phi$ are written as $\nabla_k \phi^A$ with respect to the basis $\{\mu^A \otimes e^k\}_{A=1,\ldots,N; k=1,2,3,4}$ for $\mathfrak{F} \otimes T^*X$. The components $\nabla \nabla \phi$ are then written as $\nabla_i \nabla_k \phi^A$, and so on.

The Riemann curvature tensor for the metric $\mathfrak{g}$ has components $\{R_{nmik}\}_{n,m,i,k=1,2,3,4}$. These appear in the commutator identity

$$[\nabla_i, \nabla_k] e^n = - R_{nmik} e^m$$

(A2.3)

for the covariant derivatives of the $T^*X$ frame vectors. The self-dual Weyl tensor is the traceless part of the section of $\Lambda^+ \otimes \Lambda^+$ that is defined by the symmetric, $3 \times 3$ matrix with the components $\{\tfrac{1}{8} \eta^a_{nm} \eta^b_{ik} R_{nmik}\}_{a,b=1,2,3}$. The trace of this matrix is $\tfrac{1}{4}$ times the scalar curvature of the metric $\mathfrak{g}$. The $\Lambda^+ \otimes \Lambda^-$ components of the tensor $\mathcal{B}$ that appears in (3.1) can also be written as $\{\mathcal{B}^{ac} = \tfrac{1}{8} \eta^a_{nm} \eta^c_{ik} R_{nmik}\}_{a,c=1,2,3}$.

Let $S$ denote the subbundle in $\Lambda^+ \otimes \Lambda^+$ of traceless, symmetric elements. On the open set where $T^*X$ has the given oriented orthonormal frame, a give section of $S$ (to be denoted by $\mathcal{X}$) can be written as $\mathcal{X} = \mathcal{X}^{ab} \omega^a \otimes \omega^b$ with $\{\mathcal{X}^{ab}\}_{a,b \in \{1,2,3\}}$ being the components of a traceless, symmetric matrix valued function on the set. The components of the section $\mathcal{L}_T^\dagger \mathcal{X}$ of $\Lambda^+ \otimes \Lambda^-$ with respect to the basis $\{\omega^a \otimes \omega^c\}_{a,c \in \{1,2,3\}}$ are given by the rule

$$-(\mathcal{L}_{\mathfrak{g}}^{\dagger}\mathcal{X})^{ac} = \eta^{b}{}_{im}\eta^{c}{}_{km}\nabla_{i}\nabla_{k}\mathcal{X}^{ab} + \mathcal{B}^{bc}\mathcal{X}^{ab}.$$

(A2.4)

Of interest in what follows are $\mathcal{L}_{\mathfrak{g}}^{\dagger}$ when X is the $T\times(T'-(D\cup\iota(D)))$ part of $\mathbb{T}^k$ with the metric $g_T$ from (3.6) and when $X = (S^3-N_K)\times S^1$ with the product metric $g_K$.

**b) The Fourier decomposition of $\mathcal{X}_T$**

Reintroduce the metric $g_T$ on the torus $\mathbb{T} = T\times T'$ from (3.6). The $\mathcal{X}_T$ part of the pair $(\mathcal{X}_T, \mathcal{X}_K)$ from (5.4) and (5.5) is an $\iota$-invariant element in the kernel of the $\mathfrak{g} = g_T$ version of $\mathcal{L}_{\mathfrak{g}}^{\dagger}$ on $T\times(T'-(D\cup\iota(D)))$. By way of a reminder, D is the disk in T' with radius $e^{-R}$ whose center is the point $p_*$ with coordinates $(t_3 = t_*, t_4 = t_*)$. The operator $\mathcal{L}_{g_T}^{\dagger}$ has constant coefficients when written with the basis $\{e^1, e^2, e^3, e^4\}$ from Part 1 of Section 5c; it is given by (A2.4) with $\mathcal{B} = 0$ and with the covariant derivatives being ordinary derivative. Because $\mathcal{L}_{g_T}^{\dagger}$ has constant coefficients and because $T\times(T'-(D\cup\iota(D)))$ is invariant with respect to the T action on T by constant translations, any given element in the kernel of $\mathcal{L}_{g_T}^{\dagger}$ on this domain can be written as a Fourier sum whose terms are indexed by pairs of integers with each term of the sum in the kernel of $\mathcal{L}_{g_T}^{\dagger}$. In this regard, a given $(k_1, k_2) \in \mathbb{Z}\times\mathbb{Z}$ term has the form

$$\mathcal{X}_{(k_1,k_2)}(t_3, t_4)\, e^{i(k_1 t_1 + k_2 t_2)}$$

(A2.5)

with $\mathcal{X}_{(k_1,k_2)}$ being a function on $T'-(D\cup\iota(D))$ with values in the complexification of the vector space $\mathbb{S}$ of $3\times 3$, traceless symmetric matrices. The following proposition is used to justify the subsequent focus on solely the $(k_1 = 0, k_2 = 0)$ Fourier mode. This proposition uses $\Delta$ to denote the T invariant function on $T\times T'$ that gives the distance along the T' factor to the subspace $D\cup\iota(D)$ in T'.

**Proposition A2.1**: *There exists $\kappa > 1$ that is independent of R and has the following significance: Having fixed $R > \kappa^2$, let D denote the radius $e^{-R}$ disk in T' centered at $p_*$ and let A denote the annulus in T' where the distance to $p_*$ is between $e^{-R}$ and $\kappa e^{-R}$. Suppose that X is an $\iota$-invariant element in the kernel of $\mathcal{L}_{g_T}^{\dagger}$ on $T\times(T'-(D\cup\iota(D)))$ with zero $(k_1 = 0, k_2 = 0)$ Fourier mode and normalized so that*

$$\int_{T\times(A\cup\iota(A))} |X|^2 = 1.$$

*Then $|\mathcal{X}| \leq \kappa e^{2R} \exp(-e^R \Delta/\kappa^2)$ where $\Delta > \kappa e^{-R}$.*

The upcoming proof of this proposition (and subsequent proofs) invokes the upcoming Lemma A2.2. By way of notation, this lemma has X denoting a smooth, oriented 4-manifold and $\mathfrak{g}$ being a metic on X. Given an open set $V \subset X$ with compact closure, let $\mathfrak{f}$ denote a non-negative, $C^2$ function on the closure of V that vanishes on $\partial V$. Let $\varepsilon > 0$ be such that $|d\mathfrak{f}| \leq \varepsilon^{-1}$ and $|d^\dagger d\mathfrak{f}| \leq \varepsilon^{-2}$. Set $V_\mathfrak{f} \subset V$ to be the set of points in V where $\mathfrak{f} \geq 1$. Let $|\mathcal{R}_\mathfrak{g}|_0$ denote the norm of the Riemann curvature tensor of the metric $\mathfrak{g}$ and let $|\nabla \mathcal{R}_\mathfrak{g}|_0$ denote the norm of its covariant derivative.

**Lemma A2.2**: *There exist $\kappa > 1$ with the following significance: Let X denote an oriented 4 dimensional manifold and let $\mathfrak{g}$ denote a given Riemannian metric on X. Fix an open set V, a number $\varepsilon > 0$, and then fix a function $\mathfrak{f}$ as described above. Suppose that $\mathcal{X}$ is a symmetric, traceless section of $\otimes^2 \Lambda^+$ on V with $\mathcal{X}$ and $\mathcal{L}_\mathfrak{g}^\dagger \mathcal{X}$ being square integrable. Then*
$$\int_{V_\mathfrak{f}} |\nabla \nabla \mathcal{X}|^2 \leq \kappa \left( \int_V |\mathcal{L}_\mathfrak{g}^\dagger \mathcal{X}|^2 + (|\mathcal{R}_\mathfrak{g}|_0^2 + |\nabla \mathcal{R}_\mathfrak{g}|_0^{4/3}) \int_V |\mathcal{X}|^2 + \kappa \varepsilon^{-4} \int_{V-V_\mathfrak{f}} |\mathcal{X}|^2 \right).$$

The proof of this lemma is given in Section A2f. The key point to note is that the number $\kappa$ does not depend on X, its metric, the set V, the function $\mathfrak{f}$, and also not on $\mathcal{X}$.

*Proof of Proposition A2.1*: Suppose that U is an open set in T´. If $\mathcal{X}$ is a tensor on $T \times U$ with no $(k_1 = 0, k_2 = 0)$ term in its T-action Fourier expansion, then

$$\int_{T \times U} |\nabla \nabla \mathcal{X}|^2 \geq c_*^{-1} e^{4R} \int_{T \times U} |\mathcal{X}|^2 ,$$

(A2.6)

with $c_*$ denoting here and subsequently in this proof a number that is greater than 1 and independent of $\mathcal{X}$, R and U. The precise value of $c_*$ can be assumed to increase between subsequent appearances. (The number $c_*$ in most instances will depend on the knot K via the matrix $\mathfrak{m}$ that from (3.6). and on the tensor bundle in question) The inequality in (A2.6) is used below with the tensor bundle being the bundle of symmetric, traceless sections of $\otimes^2 \Lambda^+$.

Let $\mathbb{D}$ denote as before the disk in T´ centered at the point $p_*$ with radius equal to $\frac{1}{4} t_*$ and let $\rho$ denote the radial coordinate in $\mathbb{D}$. Fix $L > 10$ and let N denote the largest integer less $\frac{1}{8} e^R L^{-1} t_*$. Supposing that n is an integer from the set $\{1, 2, \ldots, N\}$, let $A_n$ denote the annulus in $\mathbb{D} - D$ where the radial function $\rho$ obeys $\rho \in (nLe^{-R}, (n+1)Le^{-R})$. With $\mathcal{X}$ now denoting a symmetric, traceless section of $\otimes^2 \Lambda^+$ over $\mathbb{D} - D$, define the number $f_n$ to be

$$f_n = \int_{T \times A_n} |\mathcal{X}|^2 .$$

(A2.7)

Now assume in addition that $\mathcal{X}$ is annilated by $\mathcal{L}_{g_T}^\dagger$. Invoke Lemma A2.2 using for V the set $A_{n-1} \cup A_n \cup A_{n+1}$ and $\mathfrak{f} = (\rho(Le^{-R})^{-1} - (n-1))(n+2 - \rho(Le^{-R})^{-1})$. Noting that the number $\varepsilon$ in this instance of Lemma A2.2 can be taken to be $c_*(Le^{-R})^{-1}$, the inequality asserted by Lemma A2.2 and the $U = A_n$ version of (A2.6) together lead to the inequality

$$c_*^{-1} e^{4R} f_n \leq c_* L^{-4} e^{4R} (f_{n-1} + f_n + f_{n+1}) .$$

(A2.8)

It follows as a consequence that if $L > 4c_*$, then (A2.8) implies the following:

$$0 > c_*^{-2} L^4 f_n - (f_{n+1} + f_{n-1} - 2f_n) .$$

(A2.9)

The maximum/comparison principle can be invoked using the latter inequality to see that

$$f_n \leq f_N e^{-(N-n)L^2/c_*} + f_1 e^{-(n-1)L^2/c_*}$$

(A2.10)

when $n \in \{1, \ldots, N\}$. By way of an explanation, the right hand side of (A2.9) has the form of a discretized version of the Helmholtz operator $-\frac{d^2}{dx^2} + m^2$ acting on a function, $f$, of a variable $x$. A maximum principle principle argument can be invoked with this discrete Helmoltz operator in (A2.9) using the function on the set $\{1, \ldots, N\}$ whose value at any given n is $f_n - g_n$ with $g_n$ being the function that appears on the right hand side of (A2.10).

Let $A_\Diamond$ now denote the annulus in $\mathbb{D} - D$ where $\rho \in (\frac{1}{32} t_*, \frac{1}{16} t_*)$. The inequality in (A2.10) leads to an $L^2$ bound for $|\mathcal{X}|$ on $T \times A_\Diamond$ that reads

$$\int_{T \times A_\Diamond} |\mathcal{X}|^2 \leq (f_N + f_1) \exp(-c_*^{-1} e^R L) .$$

(A2.11)

Let $\mathbb{D}_* \subset \mathbb{D}$ denote the disk that is concentric to $\mathbb{D}$ but has radius $\frac{1}{12} t_*$. Suppose that $\mathcal{X}$ is an $\iota$-invariant element in the kernel of $\mathcal{L}_{g_T}^\dagger$ on $T \times (D \cup \iota(D))$. Given $R \geq c_*$, invoke Lemma A2.2 with V being $T \times (T' - (\mathbb{D}_* \cup \iota(\mathbb{D}_*)))$ and with $\mathfrak{f}$ being a smooth, $\iota$-invariant function that is equal to the $1000 t_*^{-1}(\rho - \frac{1}{12} t_*)$ on $\mathbb{D} - \mathbb{D}_*$. With this input, Lemma A2.2 and (A2.6) bound the square of the $L^2$ norm of $\mathcal{X}$ on V by $c_*$ times the integral on the left hand side of (A2.11), and thus by $c_*$ times what is written on the right hand side (A2.11).  This has the following immediate consequence when L and R are larger than $c_*$: Since the annulus $A_N$ is in $T'-(\mathbb{D}_* \cup \iota(\mathbb{D}_*))$, the number $f_N$ that appears in (A2.10) and (A2.11) is a priori bounded by $c_* f_1 \exp(-c_*^{-1} e^R L)$.  This last bound with (A2.10) and (A2.11) lead directly to the pointwise norm asserted by Proposition A2.1 because standard elliptic regularity techniques can be used to bound the sup norm of the function $|\mathcal{X}|$ on any radius $e^{-R}$ ball in $T \times (T'-(D \cup \iota(D)))$ with distance no less than $4e^{-R}$ from D and $\iota(D)$ by $c_0 e^{2R}$ times the $L^2$ norm of this function on the concentric, radius $2e^{-R}$ ball.

### c) The ($k_1 = 0, k_2 = 0$) Fourier mode of $\mathcal{X}_T$.

Assume in this subsection that $\mathcal{X}$ is a T-independent element in the kernel of $\mathcal{L}_{g_T}^\dagger$ on $T \times (T'-(D \cup \iota(D)))$.  To see what $\mathcal{X}$ looks like, it is useful to first look at the section $\mathcal{A}$ of $\Lambda^+ \otimes T^*X$ whose components when written using the basis $\{e_1, e_2, e_3, e_4\}$ of $T^*X$ from Part 1 of Section 5c and the basis for $\Lambda^+$ in (5.7) are

$$\mathcal{A}^a_i = \eta^b_{ik} \nabla_k \mathcal{X}^{ab} \ .$$

(A2.12)

The superscript a and the subscript i take respective values in $\{1, 2, 3\}$ and in $\{1,2,3,4\}$. The twelve components of $\mathcal{A}$ obey the linear constraint

$$\eta^a_{ni} \mathcal{A}^a_i = 0$$

(A2.13)

because the components of $\mathcal{X}$ define a symmetric and traceless matrix.  Meanwhile, each of the three 1-forms $\{\mathcal{A}^a_i e^i\}_{i=1,2,3}$ obeys the first order system of differential equations

$$\eta^c_{ik} \nabla_i \mathcal{A}^a_k = 0 \quad \text{and} \quad \nabla_i \mathcal{A}^a_i = 0 \ .$$

(A2.14)

The three versions of the left hand equation express the fact that $\mathcal{L}_{g_T}^\dagger \mathcal{X} = 0$, and the right hand equation follows directly from (A2.12).  With the preceding understood, there are two parts in what follows:  Part 1 writes the general solution to (A2.14); and Part 2 puts the solutions to (A2.14) from Part 1 on the left hand side of (A2.12) and then solves the resulting equation for the desired $\mathcal{X}$.

*Part 1*:  Since derivatives with respect to the coordinates $t_1$ and $t_2$ are zero, the equations in (A2.14) for the components of $\mathcal{A}$ decouple to give a pair of Cauchy-Riemann equations on T′:

- $\nabla_3 \mathcal{A}^a_4 - \nabla_4 \mathcal{A}^a_3 = 0$ *and* $\nabla_3 \mathcal{A}^a_3 + \nabla_4 \mathcal{A}^a_4 = 0$.
- $\nabla_4 \mathcal{A}^a_1 - \nabla_3 \mathcal{A}^a_2 = 0$ *and* $\nabla_4 \mathcal{A}^a_2 + \nabla_3 \mathcal{A}^a_1 = 0$.

(A2.15)

The equations in the top bullet restate the respective c = 3 version of the left hand equation in (A2.14) and the right hand equation in (A2.14). The equations in the second bullet restate the respective c = 1 and c = 2 versions of the left hand equation in (A2.14). Let $\mathcal{K}^a = \mathcal{A}^a_3 - i\mathcal{A}^a_4$, this to be viewed as a function on T´−(D∪ι(D)). The top bullet equations in (A2.15) says that $\mathcal{K}^a$ is a holomorphic function of the complex coordinate $t_3 + it_4$ on any given disk in T´−(D∪ι(D)). The second bullet equations in (A2.16) makes the same assertion for the function $\mathcal{J}^a = \mathcal{A}^a_1 - i\mathcal{A}^a_2$ with the latter also viewed as a function on the domain T´−(D∪ι(D)).

Let $\mathbb{D}$ again denote the radius $\frac{1}{4} t_*$ disk in T´ centered at the point $p_*$. Fix an index a ∈ {1, 2, 3}. Since $\mathcal{K}^a$ and $\mathcal{J}^a$ are holomorphic on $\mathbb{D}$−D, they have a Laurent expansion on this domain. This expansion can be written using the coordinate $z = t_3 - t_* + i(t_4 - t_*)$ as follows: Let $\mathcal{I}^a$ denote either $\mathcal{K}^a$ or $\mathcal{J}^a$. Then

$$\mathcal{I}^a = \ldots a_1 z + a_0 + a_{-1} \frac{1}{z} + a_{-2} \frac{1}{z^2} + \cdots$$

(A2.16)

with $\{a_k\}_{k \in \mathbb{Z}} \subset \mathbb{C}$. There is a similar expansion near ι($\mathbb{D}$−D).

The next lemma makes the formal assertion that the coefficient $a_{-1}$ that appears in (A2.16) is necessarily equal to 0 if $\mathcal{I}^a$.

**Lemma A2.3**: *Let $\mathcal{I}^a$ denote either $\mathcal{K}^a$ or $\mathcal{J}^a$. The Laurent expansion of $\mathcal{I}^a$ on $\mathbb{D}$−D that is depicted in (A2.16) has $a_{-1} = 0$. The analogous term in the Laurent expansion on ι($\mathbb{D}$−D) is also zero.*

*Proof of Lemma A2.3*: Supposing that $a_{-1} \neq 0$, then $w = \mathcal{I}^a + a_{-1} x$ will have a non-zero residue at $p_*$ and zero residue at ι($p_*$) because $x$ is ι-invariant and $\mathcal{I}^a$ is minus its pull-back by ι. This is impossible for the following reason: If the expansion in (A2.16) has but a finite number of non-zero inverse powers of z, then $w$ is a meromorphic function on the torus if $\mathcal{I}^a$, and the sum of the residues of a meromorphic function must vanish (see Theorem 4 in Chapter 7 of Ahlfor's book [A]). But $w$ has residue $2a_{-1}$ at $p_*$ and residue 0 at ι($p_*$). In the general case, the integration by parts argument for Theorem 4 in this same chapter of Ahlfor's book can be used to see that the sum of the path integrals of $w dz$ on ∂D and ι(∂D) are zero. This sum has zero contribution from the ι(∂D) integral and contribution $4\pi i a_{-1}$ from the ∂D integral.

The next lemma depicts $\mathcal{I}^a$ on the whole of T´−(D∪ι(D)) when the expansion in (A2.16) has a finite order pole. This lemma uses $\partial$ to denote $\frac{1}{2}(\frac{\partial}{\partial t_3} - i\frac{\partial}{\partial t_4})$ which is the holomorphic derivative on T´.

**Lemma A2.4**: *Suppose that $X$ is a T and $\iota$ invariant element in the kernel of $\mathcal{L}_{g_T}^\dagger$ on the domain T×(T´−(D∪ι(D))). Define $\mathcal{A}$ by (A2.12) and then define $\{\mathcal{K}^a\}_{a=1,2,3}$ and $\{\mathcal{J}^a\}_{a=1,2,3}$ from $\mathcal{A}$ as done above. If the expansion in (A2.14) for a given function $\mathcal{I}^a$ from $\{\mathcal{K}^a\}_{a=1,2,3}$ or $\{\mathcal{J}^a\}_{a=1,2,3}$ has a finite order pole, then $\mathcal{I}^a$ can be written on the whole of T´−(D∪ι(D)) in terms of the function $x$ from Lemma 5.1 as*

$$\mathcal{I}^a = (a_0 + a_1 x + a_2 x^2 + \cdots + a_n x^n)\partial x$$

*for some positive integer $n$ and for some set $\{a_0, \ldots, a_n\} \subset \mathbb{C}$.*

**Proof of Lemma A2.4**: The function $\mathcal{I}^a$ obeys $\iota^*\mathcal{I}^a = -\mathcal{I}^a$ because $\mathcal{I}^a$ is a linear combination of first derivatives of the components of $X$ and $\iota^*X = X$. It follows as a consequence that $\mathcal{I}^a = 0$ if there are no poles at $p_*$ (and thus at $\iota(p_*)$) because the constant functions are $\iota$-invariant. Now suppose that some version of $\mathcal{I}^a$ has a non-trivial Laurent expansion as in (A2.16). As noted by Lemma A2.3, this implies that some $k \geq 2$ version of $a_{-k}$ must be non-zero. Supposing that only finitely many inverse powers of $z$ appear in (A2.14), let $n \in \{0, 1, \ldots\}$ be such that the $a_{-(n+2)} \neq 0$ in (A2.14) but $a_{-k} = 0$ for $k > n+2$. Since $\partial x$ has a pole of order 2 at $p_*$, there exists $a_n \in \mathbb{C}$ such that $\mathcal{I}^a - a_n x^n \partial x$ has a pole of order at most $n+1$ at $p_*$; and thus at $\iota(p_*)$ also because this function changes sign under the action of $\iota$. Because $\mathcal{I}^a - a_n x^n \partial x$ has pole at $p_*$ and $\iota(p_*)$ of order less than that of $\mathcal{I}^a$ and because it changes sign under the action of $\iota$, an induction argument on the integer $n$ proves that $\mathcal{I}^a$ can be written as claimed on T´−(D∪ι(D)).

*Part 2*: Suppose now that $\mathcal{A}$ is as described in Part 1 and that it can be written as in (A2.12) with $X$ being independent of the coordinates on T. Granted this assumption, then the various $a \in \{1,2,3\}$ and $k \in \{1, 2, 3, 4\}$ versions of $\mathcal{A}^a_k$ and the various $X^{ab}$ will be viewed as functions on T´−(D∪ι(D)). Viewed in this light, the equation in (A2.12) asserts the following:

- $\mathcal{A}^a_3 = \nabla_4 X^{a3}$ and $\mathcal{A}^a_4 = -\nabla_3 X^{a3}$.
- $\mathcal{A}^a_1 = \nabla_4 X^{a1} - \nabla_3 X^{a2}$ and $\mathcal{A}^a_2 = \nabla_4 X^{a2} + \nabla_3 X^{a1}$.

(A2.17)

When written using $\mathcal{K}$ and $\mathcal{J}$, these equations say that

- $\mathcal{K}^a = 2i\partial \mathcal{X}^{a3}$,
- $\mathcal{J}^a = -2i\bar{\partial}(\mathcal{X}^{a1} - i\mathcal{X}^{a2})$,

(A2.18)

with $2\partial = \frac{\partial}{\partial t_3} - i\frac{\partial}{\partial t_4}$ and $2\bar{\partial} = \frac{\partial}{\partial t_3} + i\frac{\partial}{\partial t_4}$.

The constraint in (A2.13) (which follows from the fact that $\mathcal{X}$ is traceless and symmetric) can be used to write $\mathcal{K}^3$ and $\mathcal{J}^3$ as follows:

$$\mathcal{K}^3 = -\bar{\mathcal{J}}^1 + i\bar{\mathcal{J}}^2 \quad \text{and} \quad \mathcal{J}^3 = \bar{\mathcal{K}}^1 - i\bar{\mathcal{K}}^2.$$

(A2.19)

This linear constraint requires that $\mathcal{K}^3$ and $\mathcal{J}^3$ be constant since $\{\mathcal{J}^a\}_{a=1,2,3}$ and $\{\mathcal{K}^a\}_{a=1,2,3}$ are holomorphic and only constant functions are both holomorphic and anti-holomorphic. Moreover, these constants must be zero when $\iota^*\mathcal{X} = \mathcal{X}$ because the latter identity implies that $\iota^*\mathcal{A}^a_i = -\mathcal{A}^a_i$ for all index pairs (a,i) since each such $\mathcal{A}^a_i$ is a linear combination of derivatives of components of $\mathcal{X}$. The vanishing of $\mathcal{K}^3$ and $\mathcal{J}^3$ require that both $\mathcal{J}^1 + i\mathcal{J}^2$ and $\mathcal{K}^1 + i\mathcal{K}^2$ are also zero.

Since $\mathcal{K}^3 = 0$ and $\mathcal{X}^{33}$ is real, the top bullet in (A2.18) requires that $\mathcal{X}^{33}$ be constant. Since $\mathcal{X}$ is traceless, it follows that $\mathcal{X}^{11} + \mathcal{X}^{22}$ is -1 times this same constant. Meanwhile, the vanishing of $\mathcal{K}^1 + i\mathcal{K}^2$ requires that $\mathcal{X}^{13} - i\mathcal{X}^{23}$ be holomorphic on $T'-(D\cup\iota(D))$. The lower bullet of (A2.18) implies this because $\mathcal{X}^{13} = \mathcal{X}^{31}$ and $\mathcal{X}^{23} = \mathcal{X}^{32}$. Since $\mathcal{X}^{13} - i\mathcal{X}^{23}$ is holomorphic, it can be depicted on $\mathbb{D}-D$ as convergent Laurent expansion

$$\mathcal{X}^{13} - i\mathcal{X}^{23} = \cdots c_1 z + c_0 + c_{-1}\frac{1}{z} + c_{-2}\frac{1}{z^2} + \cdots.$$

(A2.20)

The pull-back of (A2.20) by $\iota$ depicts $\mathcal{X}^{13} - i\mathcal{X}^{23}$ on $T \times \iota(\mathbb{D}-D)$ because $\mathcal{X}^{13} - i\mathcal{X}^{23}$ is $\iota$ invariant. If (A2.20) has but finitely many nonzero powers of $\frac{1}{z}$, then $\mathcal{X}^{13} - i\mathcal{X}^{23}$ on $T'-(D\cup\iota(D))$ can be depicted as a finite polynomial in the function $x$:

$$\mathcal{X}^{13} - i\mathcal{X}^{23} = c_0 + c_1 x + c_2 x^2 + \cdots$$

(A2.21)

with $\{c_0, c_1, \ldots\}$ being a finite set of complex numbers. This depiction of $\mathcal{X}^{13} - i\mathcal{X}^{23}$ follows because $x$ generates the ring of $\iota$-invariant, meromorphic functions on $T'$ with poles at $p_*$ and $\iota(p_*)$. The detailed argument differs little from the proof of Lemma A2.4.

Since $\mathcal{J}^3 = 0$ and $\mathcal{J}^1 + i\mathcal{J}^2 = 0$, it is enough to consider only the a = 1 case in the second bullet of (A2.18). The function $\mathcal{J}^1$ on $\mathbb{D}-D$ has the convergent Laurent expansion that is depicted in (A2.16). This understood, then the second bullet in (A2.18) implies that the function $\mathcal{X}^{11} - i\mathcal{X}^{12}$ can be written on $\mathbb{D}$–D as the sum of two convergent series:

$$\mathcal{X}^{11} - i\mathcal{X}^{12} = \tfrac{i}{2}((\cdots a_1 z + a_0 + a_{-2}\tfrac{1}{z^2} + \cdots)\bar{z} + (\cdots b_1 z + b_0 + b_{-1}\tfrac{1}{z} + b_{-2}\tfrac{1}{z^2} + \cdots)) \ .$$

(A2.22)

The pull-back by $\iota$ of (A2.22) depicts $\mathcal{X}^{11} - i\mathcal{X}^{12}$ on $\iota(\mathbb{D}$–D) because $\mathcal{X}^{11} - i\mathcal{X}^{12}$ is $\iota$ invariant.

Suppose now that $\mathcal{J}^1$ is described by Lemma A2.4 for some integer n and coefficient set $\{a_1, \ldots, a_n\} \subset \mathbb{C}$ (This is the case if and only if the left most series in (A2.20) has but finitely many non-zero powers of $\tfrac{1}{z}$). Let $u$ denote the function that is described by Lemma 5.2. Then the function

$$w = (a_0 + a_1 x + a_2 x^2 + \cdots + a_n x^n)\, u$$

(A2.23)

obeys $\bar{\partial} w = \mathcal{J}^1$. It follows as a consequence that if both sums in (A2.22) have but finitely many powers of $\tfrac{1}{z}$, then $\mathcal{X}^{11} - i\mathcal{X}^{12}$ can be written on the whole of T´−(D∪$\iota$(D)) as

$$\mathcal{X}^{11} - i\mathcal{X}^{12} = \tfrac{i}{2}((a_0 + a_1 x + a_2 x^2 + \cdots + a_n x^n)\, u + (b_0 + b_1 x + b_2 x^2 + \cdots + b_m x^m)$$

(A2.24)

for some non-negative integers n, m and complex numbers $\{a_0, \ldots, a_n\}$ and $\{b_0, \ldots, b_m\}$.

An $\iota$-invariant and T-invariant element in the kernel of $\mathcal{L}_{g_T}^\dagger$ on T×(T´−(D∪$\iota$(D))) is said in what follows to be *regular* if the expansions in (A2.20) and (A2.22) have but a finite number of powers of $\tfrac{1}{z}$. The following proposition summarizes what was said in the preceding paragraphs about these elements.

**Proposition A2.5**: *The vector space of* T *and* $\iota$ *invariant regular elements in the kernel of* $\mathcal{L}_{g_T}^\dagger$ *on the domain* T×(T´−($p_*$∪$\iota(p_*)$)) *is the linear span of solutions having the form*
- $\mathcal{X}^{33} = -(\mathcal{X}^{11} + \mathcal{X}^{22}) = s$ ,
- $\mathcal{X}^{13} - i\mathcal{X}^{23} = c_0 + c_1 x + c_2 x^2 + \cdots$ ,
- $\mathcal{X}^{11} - i\mathcal{X}^{12} = (a_0 + a_1 x + a_2 x^2 + \cdots)\, u + (b_0 + b_1 x + b_2 x^2 + \cdots)$ ,

*with* $s \in \mathbb{R}$ *being constant; and with* $\{c_0, c_1, \ldots\} \subset \mathbb{C}$ *and* $\{a_0, a_1, \ldots\} \subset \mathbb{C}$ *and* $\{b_0, b_1, \ldots\} \subset \mathbb{C}$ *all being constant.*

### d) An R-invariant reformulation of Proposition A2.5

It is convenient for the subsequent applications to rewrite Proposition A2.5 so as to remove all references to the choice of the parameter R. (The parameter R is used to define the metric $g_T$, and the metric $g_T$ in turn is needed first to define the bundles $\Lambda^\pm$ and their tensor powers, and then to define the operator $\mathcal{L}_{g_T}^\dagger$.) What is said below is admittedly pedantic, but even so, it is important to keep in mind for later.

Having fix $R > c_*$, let $\{e^1, e^2\}$ denote a constant, oriented orthonormal frame for the metric on T that is depicted in (5.6) and let $e^3 = dt_3$ and $e^4 = dt_4$. This basis identifies the bundle $T^*(T \times T')$ with the product vector bundle with fiber $\mathbb{R}^4$. Meanwhile, the corresponding basis for the $g_T$ versions of $\Lambda^+$ and $\Lambda^-$ that are depicted in (5.7) identify these bundles with the product vector bundle whose fiber is $\mathbb{R}^3$. A section over a domain in $T \times T'$ of the $g_T$ version of $\Lambda^+$ or $\Lambda^-$ with $g_T$ defined by any given $R > c_*$ can be viewed (when convenient) as a map from the domain to the R independent space $\mathbb{R}^3$. Letting $\mathbb{M}_\diamond \subset \mathbb{R}^3 \otimes \mathbb{R}^3$ denote the space of symmetric, traceless matrices, then a traceless, symmetric section of any given $g_T$ version of $\Lambda^+ \otimes \Lambda^+$ over a domain in $T \times T'$ can likewise be viewed as a map from the domain to the R-independent vector spaces $\mathbb{M}_\diamond$. Moreover, if the domain has the form $T \times U$ with $U \subset T'$ being an open set, and if the original section of $\Lambda^+ \otimes \Lambda^+$ is T invariant, then the corresponding map from $T \times U$ to $\mathbb{M}_\diamond$ can be viewed as a map from the domain U to $\mathbb{M}_\diamond$.

Having specified R so as to define $g_T$, suppose that $U \subset T'$ is a given open set and that $\mathcal{X}$ is a T invariant, traceless and symmetric section of $\Lambda^+ \otimes \Lambda^+$ on the domain $T \times U$; which is to say that $\mathcal{X}$ is a map from U to the vector space $\mathbb{M}_\diamond$. Viewing $\mathcal{X}$ as a map from U to $\mathbb{M}_\diamond$ writes the action of the operator $\mathcal{L}_{g_T}^\dagger$ as that of a certain second order, R-independent operator; it is the operator on maps from U to $\mathbb{M}_\diamond$ that is obtained from $\mathcal{L}_{g_T}^\dagger$ by setting all derivatives along T equal to 0. (This operator is the composition of the operator that defines $\mathcal{A}$ from $\mathcal{X}$ using (A2.17) and then takes the linear combinations of the derivatives of $\mathcal{A}$ that are written in the second bullet of (A2.15) and in the left hand equation of the first bullet in (A2.15). The right hand equation in this bullet is automatically obeyed when $\mathcal{A}$ comes from $\mathcal{X}$ via (A2.17).) This operator on maps from domains in T' to $\mathbb{M}_\diamond$ is denoted in what follows by $\mathcal{L}^\dagger$.

What follows are three important points to keep in mind. The first point is that there is a canonical bijection between the T-invariant kernel of any given $g_T$ version of $\mathcal{L}_{g_T}^\dagger$ on a domain $T \times U$ in $T \times T'$ and the kernel of $\mathcal{L}^\dagger$ on U. The second point concerns the involution $\iota$: If the involution $\iota$ is defined to act on T' via the rule $(t_3, t_4) \to (-t_3, -t_4)$ and if $U \subset T'$ is an $\iota$ invariant domain, then the aforementioned bijection between the T invariant maps from $T \times U$ to $\mathbb{M}_\diamond$ and the maps from U to $\mathbb{M}_\diamond$ commutes with the corresponding actions of $\iota$ on $T \times U$ and on U. The final point concerns $L^2$ inner products: This canonical bijection is not isometric with respect to the $L^2$ inner products on T×U and on U; but it is conformal in the following sense: Let $\mathfrak{m}$ denote the matrix that is used to define $g_T$ in (3.6). The $L^2$ inner product on T×U of two T invariant maps to $\mathbb{M}_\diamond$ is $4\pi^2 \det(\mathfrak{m}) e^{-2R}$ times their inner product on U as maps from U to $\mathbb{M}_\diamond$.

Supposing now that $\mathcal{X}$ is an $\iota$-invariant element in the kernel of $\mathcal{L}^\dagger$ on the domain T´−($p_*\cup\iota(p_*)$), then $\mathcal{X}$ can be depicted as in (A2.20) and (A2.22) on $\mathbb{D}$−$p_*$ because the corresponding $\iota$-invariant element in the kernel of any $R \geq c_*$ version of $\mathcal{L}^\dagger_{g_T}$ has such a depiction. This understood, the element $\mathcal{X}$ is said to be *regular* when the expansions in (A2.20) and (A2.22) have but a finite number of powers of $\frac{1}{z}$.

Here is the promised rewording of Proposition A2.5:

**Proposition A2.6**: *The vector space of $\iota$ invariant regular elements in the kernel of $\mathcal{L}^\dagger$ on the domain T´−($p_*\cup\iota(p_*)$) is the linear span of solutions having the form*
- $\mathcal{X}^{33} = -(\mathcal{X}^{11} + \mathcal{X}^{22}) = s$,
- $\mathcal{X}^{13} - i\mathcal{X}^{23} = c_0 + c_1 x + c_2 x^2 + \cdots$,
- $\mathcal{X}^{11} - i\mathcal{X}^{12} = (a_0 + a_1 x + a_2 x^2 + \cdots)u + (\mathfrak{b}_0 + \mathfrak{b}_1 x + \mathfrak{b}_2 x^2 + \cdots)$,

*with $s \in \mathbb{R}$ being constant; and with $\{c_0, c_1, \ldots\} \subset \mathbb{C}$ and $\{a_0, a_1, \ldots\} \subset \mathbb{C}$ and $\{\mathfrak{b}_0, \mathfrak{b}_1, \ldots\} \subset \mathbb{C}$ all being constant.*

### e) The approximation by regular elements

The upcoming Proposition A2.7 (with Proposition A2.1) implies that any given $\iota$ invariant element in the kernel of $\mathcal{L}^\dagger_{g_T}$ on T×(T´−(D$\cup\iota$(D))) can be well approximated on a slightly smaller domain by a regular element. The statement of the proposition uses the following notation: Given a non-negative integer N and a data set consisting of a real number $r$, and sets of N complex numbers $\{c_0, \ldots, c_N\}$ and $\{a_0, \ldots, a_{N-1}\}$, $\{\mathfrak{b}_0, \ldots, \mathfrak{b}_N\}$, the proposition uses $\mathfrak{X}_N$ to denote the $\iota$ invariant element in the kernel of the operator $\mathcal{L}^\dagger$ on T´−($p_*\cup\iota(p_*)$) whose components are defined using this data by the three bullets of Proposition A2.6.

To set more of the stage for Proposition A2.7, suppose for the moment that $\mathcal{X}$ is some $\iota$ invariant element in the kernel of $\mathcal{L}^\dagger$ on T´−(D$\cup\iota$(D)). As noted subsequent to (A2.19), the element $\mathcal{X}$ has $\mathcal{X}^{33} = -(\mathcal{X}^{11} + \mathcal{X}^{22})$ being constant. The combination $\mathcal{X}^{13} - i\mathcal{X}^{23}$ has the depiction on $\mathbb{D}$−D given by (A2.20) and the combination $\mathcal{X}^{11} - i\mathcal{X}^{12}$ has the depiction on $\mathbb{D}$−D given in (A2.22). Let $\mathcal{X}_-$ denote the following part of $\mathcal{X}$ on $\mathbb{D}$−D:

- $\mathcal{X}_-^{33} = -(\mathcal{X}_-^{11} + \mathcal{X}_-^{22}) = s$,

- $\mathcal{X}_-^{13} - i\mathcal{X}_-^{23} = c_0 + c_{-1} z^{-1} + c_{-2} z^{-2} + \cdots$ ,
- $\mathcal{X}_-^{11} - i\mathcal{X}_-^{12} = (a_{-2} z^{-2} + a_{-3} z^{-3} \cdots) \bar{z} + (b_0 + b_{-1} z^{-1} + b_{-2} z^{-2} + \cdots)$ ,

(A2.25)

To be sure, the second bullet in (A2.25) consists of $\sum_{k \geq 0} c_{-k} z^{-k}$, and the the third bullet consists of $b_0 + b_{-1} z^{-1} + b_{-2} z^{-2} + \sum_{k \geq 3} (a_{-k+1} |z|^2 + b_{-k}) z^{-k}$.

By way of a final bit of notation, suppose that $R \geq c_*$ has been fixed and then a number $r \in (e^{-R}, \frac{1}{4} t_*)$. The proposition uses $D_r$ to denote the disk in $T'$ with radius $r$ and center at the point $p_*$. For example, the disk $D$ is the $r = e^{-R}$ version of $D_r$.

**Proposition A2.7**: *There exists $\kappa > 1$ with the following significance: Having fixed $R \geq \kappa$, let $D$ denote the disk in $T'$ with center $p_*$ and radius $e^{-R}$ and let $\mathcal{X}$ denote an $\iota$ invariant element in the kernel of $\mathcal{L}^\dagger$ on the domain $T' - (D \cup \iota(D))$. Use $z_{T,K}$ to denote the $L^2$ norm of $\mathcal{X}$ on the annulus in $T'$ where the distance to $p_*$ is between $e^{-R}$ and $2 e^{-R}$. Given a positive integer $N$ greater than $1$, there exist a data set consisting of a real number $s$ and three sets of complex numbers $\{c_0, \ldots, c_N\}$ and $\{a_0, \ldots, a_{N-1}\}$, $\{b_0, \ldots, b_N\}$ such that the resulting version of $\mathfrak{X}_N$ differs from $\mathcal{X}$ by an $\iota$ invariant element in the kernel of $\mathcal{L}^\dagger$ having the properties in the next two bullets.*
- *The version of (A2.25) for $(\mathcal{X} - \mathfrak{X}_N)_-$ has $s$, $(c_0, \ldots, c_N)$, $(a_{-2}, \ldots, a_{-N-1})$ and $(b_0, \ldots, b_{-N})$ replaced by $0$; but the other coefficients are the same as those in the original $\mathcal{X}$ version of (A2.25).*
- *Supposing that $r \in (4 e^{-R}, \frac{1}{4} t_*)$, then the $L^2$ norm of $\mathcal{X} - \mathfrak{X}_N$ on $T' - (D_r \cup \iota(D_r))$ is at most $\kappa r^{-N} e^{-NR} z_{T,K}$ and the pointwise norm on $T' - (D_r \cup \iota(D_r))$ is at most $\kappa r^{-N-1} e^{-(N-1)R} z_{T,K}$.*

The proof of this proposition occupies the remainder of this subsection.

*Proof of Proposition A2.7*: The proof of the proposition has three parts. Since no generality is lost by assuming that $z_{T,K} = 1$, this condition is assumed in what follows.

*Part 1*: The lemma that follows momentarily will be used in Part 3 to bound the size of the coefficients of the positive powers of $z$ in (A2.20) and (A2.22) by the size of the coefficients of the non-positive powers of $z$. By way of notation, the lemma uses $D'$ to denote the disk in $T'$ with center $p_*$ and radius $2 e^{-R}$. This disk contains $D$ and it is contained in $\mathbb{D}$.

**Lemma A2.8**: *There exists $\kappa > 1$ with the following significance: Fix $R \geq \kappa$ to define the disks $D$ and $D'$ in the torus $T'$. Let $X$ denote an $\iota$ invariant element in the kernel of $\mathcal{L}^\dagger$ on $T' - (D \cup \iota(D))$. Define $X_-$ from $X$ as in (A2.25). Then* $\int_{T'-(D'\cup\iota(D'))} |X|^2 \leq \kappa \int_{\mathbb{D}-D} |X_-|^2$.

*Proof of Lemma A2.8*: Suppose that no such $\kappa$ exists so as to derive nonsense. In this instance, there would be sequences $\{R_n\}_{n=1,2,\ldots}$ and $\{X^{(n)}\}_{n=1,2,\ldots}$ of the following sort: The first sequence consists of positive real numbers, with the n'th member $R_n > n$. To describe the second sequence, fix a positive integer n; and let $D^{(n)}$ denote the disk in $\mathbb{D}$ with center $p_*$ and radius $e^{-R_n}$. What is denoted by $X^{(n)}$ is an $\iota$ invariant in the kernel of $\mathcal{L}^\dagger$ on the domain $T' - (D^{(n)} \cup \iota(D^{(n)}))$ with $L^2$ norm equal to 1. In addition, the norm of the corresponding $X^{(n)}_-$ on the domain $\mathbb{D} - D_n$ is less than $\frac{1}{n}$. Here and in what follows, $X^{(n)}_-$ is defined by taking $X$ in (A2.25) to be $X^{(n)}$. The derivation of nonsense from these sequences has four steps.

<u>Step 1</u>: Given $r \in (0, \frac{1}{4} t_*)$, let $A_r \subset \mathbb{D}$ denote the concentric annulus with inner radius $\frac{1}{2} r$ and outer radius $r$. Keep in mind that $\rho$ when written using the coordinate z on $\mathbb{D}$ and is $|z|$. Note that the annulus $A_r$ is contained in $\mathbb{D} - D^{(n)}$ when $n > |\ln r|$.

Terms from the right hand side of (A2.20) with different powers of z are orthogonal with respect to the $L^2$ inner product on any constant $\rho$ slice of $\mathbb{D}$. They are therefore orthogonal with respect to the $L^2$ inner product on the whole $A_r$. Supposing that $k \geq 1$, then the square of the $L^2$ norm on $A_r$ of the term $c_k z^k$ from (A2.20) is no smaller than

$$c_*^{-1} \frac{1}{k+1} r^{2(k+1)} |c_k|^2 \ .$$

(A2.26)

Supposing that $n \geq |\ln r|$, then $r > 2e^{-R_n}$ and thus $A_r \subset \mathbb{D} - D^{(n)}$. Fix n so this bound is obeyed. Then, the $X^{(n)}$ version of (A2.26) can not be greater than 1 because the $L^2$ norm of $X^{(n)}$ is equal to 1. This is the case in particular for $r = \frac{1}{4} t_*$. Using $r = \frac{1}{4} t_*$ in (A2.26) leads to the bound $|c_k| \leq c_*^{k+1}$ for the $X^{(n)}$ version of $c_k$.

The terms $a_0 \bar{z}$ and $a_1 |z|^2$, and those from the set $\{(a_{k+1}|z|^2 + b_k) z^k\}_{k=1,2,\ldots}$ that appear in (A2.22) are pairwise orthogonal for the $L^2$ norm on any given constant $\rho$ disk in $\mathbb{D}$. Meanwhile, the square of the $L^2$ norms of $a_0 \bar{z}$ and $a_1 |z|^2$ over $A_r$ are no smaller than

$$c_*^{-1} |a_0|^2 r^4 \quad and \quad c_*^{-1} |a_1|^2 r^6 \ ;$$

(A2.27)

and the square of the $L^2$ norm on $A_r$ of any given $k \geq 1$ version of $(a_{k+1}|z|^2 + b_k)z^k$ is at least

$$c_*^{-1} \frac{1}{(k+1)^2} r^{2(k+1)} (r^2|a_{k+1}|^2 + |b_k|^2) .$$

(A2.28)

If $n > |\ln r|$ and if these coefficients come from the $\mathcal{X}_{T0} = \mathcal{X}^{(n)}$ version of (A2.22), then none of the expressions in (A2.27) and (A2.28) can be greater than 1. In particular, taking $r = \frac{1}{4}t_*$ leads to the requirement that the $k \geq 0$ versions of $a_k$ and the $k \geq 1$ versions of $b_k$ obey the norm bound $|a_k| + |b_k| \leq c_*^{k+1}$.

Step 2: Given $r \in (0, \frac{1}{4}t_*)$, let $D_r$ denote the radius $r$ disk in $\mathbb{D}$ with center $p_*$. Supposing that $k \geq 1$, then the square of the $L^2$ norm on $D_r - D$ of the $c_k z^k$ term in (A2.20) is no greater than what is written in (A2.26) with the $c_*^{-1}$ factor replaced by $c_*$. Likewise, the square of the $L^2$ norms over $D_r$ of the terms $a_0 \bar{z}$ and $a_1|z|^2$ from (A2.22) are no greater than what is written in (A2.27) with the $c_*^{-1}$ factor replaced by $c_*$. Meanwhile, the square of the $L^2$ norms over $D_r$ of the $k \geq 1$ version of $(a_{k+1}|z|^2 + b_k)z^k$ is no greater than what is written in (A2.28) with the factor of $c_*^{-1}$ replaced by $c_*$.

Fix $n \geq |\ln r|$ so as to apply the preceding observations to the case when $\mathcal{X}$ is $\mathcal{X}^{(n)}$. In this case, the observations in the preceding paragraph with the norm bounds from Step 1 on the coefficients $\{c_k\}_{k\geq 1}$ and $\{a_k\}_{k\geq 0}$ and $\{b_k\}_{k\geq 1}$, plus the $\frac{1}{n}$ upper bound for the $L^2$ norm $\mathcal{X}^{(n)}_-$ on $\mathbb{D}-D^{(n)}$ leads to the $L^2$ norm bound

$$\int_{D_r - D^{(n)}} |\mathcal{X}^{(n)}|^2 \leq c_*(r^4 + \Sigma_{k\geq 2} c_*^k r^{2k+2} + \frac{1}{n^2}) .$$

(A2.29)

Note in particular that the right hand side of this is less than $c_*^{-1}$ if $r < c_*^{-1}$ and $n \geq c_*$. The key point here is that the left hand side is going to be much less than 1 if $r$ is uniformly small (independent of $n$) and $n$ is large (independent of $r$).

Fix $r \in (0, \frac{1}{4}t_*)$ again. If $n \geq |\ln r|$, then (A2.29) implies in particular that

$$\int_{A_r} |\mathcal{X}^{(n)}|^2 \leq c_*(r^4 + \Sigma_{k\geq 2} c_*^k r^{2k+2} + \frac{1}{n^2}) .$$

(A2.30)

This is because $A_r \subset D_r - D^{(n)}$ when $n \geq |\ln r|$.

Step 3: Use what is said in the previous section (but in reverse) to view $\mathcal{X}^{(n)}$ for the moment as an $\iota$ invariant and $T$ invariant element in the kernel of the $R = R^{(n)}$ version of the operator $\mathcal{L}_{g_T}^\dagger$ on the domain $T \times (T' - (D^{(n)} \cup \iota(D^{(n)})))$. Invoke Lemma A2.2 taking X to be $T \times T'$ with the metric being the $R = R_n$ version of $g_T$, taking $V$ to be the subset $T \times (T' - (D_{r/2} \cup \iota(D_{r/2})))$, and taking for the function $\mathfrak{f}$ a smooth function that is greater than 1 on $T \times (T' - (\mathbb{D} \cup \iota(\mathbb{D})))$ and is equal to $8 r^{-1}(\rho - \frac{1}{2} r)$ on $\mathbb{D}$. Since the lemma's $\varepsilon$ in this case can be taken greater than $c_*^{-1} r$, the use of $\mathcal{X} = \mathcal{X}^{(n)}$ in this instance of Lemma A2.2 leads from (A2.30) to the following bound for the $T' - (D^{(n)} \cup \iota(D^{(n)}))$ incarnation of $\mathcal{X}^{(n)}$:

$$\int_{T' - (D_r \cup \iota(D_r))} |\nabla \nabla \mathcal{X}^{(n)}|^2 \le c_* (1 + c_*(r^2 + r^{-4} \tfrac{1}{n^2})).$$

(A2.31)

To continue, note that the average of any component of $\nabla \mathcal{X}^{(n)}$ over the domain $T' - (D_\rho \cup \iota(D_\rho))$ is equal to zero because $\mathcal{X}^{(n)}$ is $\iota$ invariant. It follows as a consequence that there is version of $c_*$ such that

$$\int_{T' - (D_r \cup \iota(D_r))} |\nabla \mathcal{X}^{(n)}|^2 \le c_* \int_{T' - (D_r \cup \iota(D_r))} |\nabla \nabla \mathcal{X}^{(n)}|^2 .$$

(A2.32)

Thus, the square of the $L^2$ norm of $\nabla \mathcal{X}^{(n)}$ on $T' - (D_r \cup \iota(D_r))$ is also bounded by the product of $c_*$ and what is written on the right hand side of (A2.31).

<u>Step 4</u>: With $r \in (0, \tfrac{1}{4} t_*)$ chosen, fix $n \ge |\ln r|$. The following bullets rewrite the salient conclusions of Steps 1-3.

- $\mathcal{L}^\dagger \mathcal{X}^{(n)} = 0$.
- $\displaystyle\int_{T' - (D^{(n)} \cup \iota(D^{(n)}))} |\mathcal{X}^{(n)}|^2 = 1.$
- $\displaystyle\int_{(D_r - D^{(n)})} |\mathcal{X}^{(n)}|^2 \le c_* (r^4 + \tfrac{1}{n^2}) .$
- $\displaystyle\int_{T' - (D_r \cup \iota(D_r))} (|\nabla \nabla \mathcal{X}^{(n)}|^2 + |\nabla \mathcal{X}^{(n)}|^2) \le c_* (1 + r^{-4} \tfrac{1}{n^2}).$
- $\displaystyle\int_{(\mathbb{D} - D_t)} |\mathcal{X}_-^{(n)}|^2 \le c_* \tfrac{1}{n^2} .$

(A2.33)

The bounds in the first four bullets of (A2.33) with some standard elliptic regularity theorems imply that the sequence $\{\mathcal{X}^{(n)}\}_{n=1,2,\ldots}$ has a subsequence that converges weakly in the $L^2_2$ topology on the domain $T' - (p_* \cup \iota(p_*))$ and strongly in the $C^\infty$ topology on compact subsets of in $T' - (p_* \cup \iota(p_*))$ to a symmetric, traceless, $\iota$-invariant $3 \times 3$ matrix with $L^2$ norm equal to 1 that is annihilated by $\mathcal{L}^\dagger$. The limit matrix for such a subsequence is denoted by $\mathcal{X}$. The fifth bullet in (A2.33) has the following implication: The corresponding $\mathcal{X}_-$ is zero. This is the desired nonsense because $\mathcal{X}$ is described by Proposition A2.6 and all of the elements described by this proposition have singularities at $p_*$ and $\iota(p_*)$.

*Part 2*: Fix $R > c_*$ to define the disk $D \subset T'$ and suppose that $\mathcal{X}$ is an $\iota$-invariant element in the kernel of $\mathcal{L}^\dagger$ on $T' - (D \cup \iota(D))$ with $z_{T,K} = 1$. Supposing that $r \in (2e^{-R}, \frac{1}{4} t_*)$, let $D_r \subset \mathbb{D}$ again denote the disk with center $p_*$ and radius $r$. If $k \geq 2$, then the square of the $L^2$ norm on $D_{2r} - D_r$ of the term $c_{-k} z^{-k}$ from (A2.25) is no less than

$$c_*^{-1} \frac{1}{k-1} r^{-2(k-1)} |c_{-k}|^2 .$$

(A2.34)

This can be at most 1 for $r = e^{-R}$. If (A2.34) in this case is no bigger than 1 then $|c_{-k}|$ is no bigger than $c_* \sqrt{k} e^{-(k-1)R}$.

The terms of the form $(a_{-k+1} |z|^2 + b_{-k}) z^{-k}$ for $k \geq 2$ from the third bullet of (A2.25) with different integer values of $k$ are mutually orthogonal with respect to the $L^2$ inner product on any constant $\rho$ slice of $\mathbb{D} - D$. Meanwhile, if $k \geq 2$, then the square of the $L^2$ norm of $(a_{-k+1} |z|^2 + b_{-k}) z^{-k}$ on $D_{2r} - D_r$ is no less than

$$c_*^{-1} \frac{1}{k^2} (r^2 |a_{-k+1}|^2 + |b_{-k}|^2) r^{-2(k-1)} .$$

(A2.35)

Since this can be at most 1 for $r = e^{-R}$, it follows that $e^{-R} |a_{-k+1}| + |b_{-k}| \leq c_* k e^{-(k-1)R}$.

Fix an integer $N \geq 2$. The sum of the individual terms in (A2.25) with norm greater than a constant multiple of $|z|^{-N}$ are the $\sum_{k>N} c_{-k} z^{-k}$ part of $\mathcal{X}^{13} - i\mathcal{X}^{23}$ and the part of $\mathcal{X}^{11} - i\mathcal{X}^{12}$ given by $\sum_{k>N} a_{-k-1} \bar{z} z^{-k-1}$ and $\sum_{k>N} b_{-k} z^{-k}$. It follows from the coefficient norm bounds of the preceding two paragraphs that the square of the $L^2$ norm on $\mathbb{D} - D_r$ of the $\sum_{k>N} c_{-k} z^{-k}$ part of $\mathcal{X}^{13} - i\mathcal{X}^{23}$ and the $\sum_{k>N} a_{-k-1} \bar{z} z^{-k-1}$ and $\sum_{k>N} b_{-k} z^{-k}$ parts of $\mathcal{X}^{11} - i\mathcal{X}^{12}$ are at most

$$c_* \sum_{k>N} r^{2(1-k)} e^{-2(k-1)R} .$$

(A2.36)

This in turn is no greater than $c_* r^{-N} e^{-NR}$.

*Part 3*: Choose the real number $r$ for defining $\mathfrak{X}_N$ to be the number that appears in the $\mathcal{X}$ version of the top bullet in (A2.25). Meanwhile, given $N \geq 2$, there are sets of complex numbers $\{c_0, \ldots, c_N\}$ and $\{a_0, \ldots, a_{N-1}\}$ and $\{b_0, \ldots, b_N\}$ so that the corresponding $\mathcal{X}$-$\mathfrak{X}_N$ version of (A2.25) has coefficients $c_{-k}$, $a_{-k+1}$ and $b_{-k}$ equal to 0 for all values of k from the set $\{0, \ldots, N\}$. These sets exist because the function $x$ that appears in Proposition A2.6 has a pole of order 1 at $p_*$ and no critical points near the disk $\mathbb{D}$.

Fix $r \in [2e^{-R}, \frac{1}{4} t_*)$  Given the preceding definition of $\mathfrak{X}_N$, it then follows from what is said in Part 2 that the square of the $L^2$ norm of $(\mathcal{X}-\mathfrak{X}_N)_-$ on $\mathbb{D}-D_r$ is no greater than $c_* r^{-N} e^{-NR}$. This bound and Lemma A2.8 give a $c_* r^{-N} e^{-NR}$ bound for the square of the $L^2$ norm of $\mathcal{X}-\mathfrak{X}_N$ on the domain $T'-(D_r \cup \iota(D_r))$. This is the $L^2$ norm bound that is asserted by Proposition A2.7 when $\mathcal{X}$ is normalized so that $z_{T,K} = 1$.

To obtain the asserted pointwise bound for $\mathcal{X}-\mathfrak{X}_N$ change the point of view for the moment and consider $\mathcal{X}-\mathfrak{X}_N$ as an $\iota$ invariant and T invariant element in the kernel of $\mathcal{L}_{g_T}^\dagger$ on $T \times (T'-(D \cup \iota(D)))$. The purpose is to invoke Lemma A2.2. To this end, take X for Lemma A2.2 to be $T \times T'$ with its metric being $g_T$ and take V to be the set $T \times (T'-(D_r \cup \iota(D_r)))$. Take the function $\mathfrak{f}$ to be a smooth function, $\iota$ invariant function that is equal to $r^{-1}(\rho - r)$ on $\mathbb{D}-D$. Invoke Lemma A2.2 with this data and with $\mathcal{X}-\mathfrak{X}_N$ used in lieu of $\mathcal{X}$. Because the corresponding version of the lemma's $\varepsilon$ can be taken larger than $c_*^{-1} r$, this instance of Lemma A2.2 with the $L^2$ bound in the preceding paragraph leads to a second derivative bound for the $T'-(D \cup \iota(D))$ incarnation of $\mathcal{X}-\mathfrak{X}_D$ that reads

$$\int_{T'-(D_{2r} \cup \iota(D_{2r}))} |\nabla\nabla(\mathcal{X}-\mathfrak{X}_N)|^2 \le c_* \, r^{-2N-4} \, e^{-2NR} \; .$$

(A2.37)

A standard Sobolev inequality in dimension 2 with the bound in (A2.37) and the $c_* c_* r^{-N} e^{-NR}$ bound on the $L^2$ norm of $\mathcal{X}-\mathfrak{X}_N$ on $T'-(D_r \cup \iota(D_r))$ lead to the pointwise norm bound that is asserted by the proposition.

**f) Proof of Lemma A2.2**

By way of notation, $c_\diamond$ is used in this proof to denote a number that is greater than 1 and is independent of X, its Riemannian metric, the set V, the function $\mathfrak{f}$ and $\mathcal{X}$. This number can be assumed to increase between successive appearances.

The proof requires a compactly supported function on V with certain specific properties. This function is denoted by $\sigma$ and the required properties are listed below.

- *$\sigma$ is smooth ($C^{1,1}$ is sufficient) and non-negative with compact support in V.*
- *$\sigma = 1$ on $V_\varepsilon$.*
- *$\varepsilon |d\sigma| + \varepsilon^2 |d^\dagger d\sigma| \le c_\diamond$.*

- $|d\sigma|^2 \le c_\diamond \varepsilon^{-2} \sigma(1-\sigma)$.

(A2.38)

The proof that such a function exists is given momentarily. Assume there is one for now.

Let $\mathfrak{g}$ denote the given metric on X. Suppose that $\mathcal{X}$ is a symmetric, traceless section of $\otimes^2 \Lambda^+$ on V that is annihilated by the operator $\mathcal{L}_\mathfrak{g}^\dagger$. The Bochner-Weitzenbock formula for $\mathcal{L}_\mathfrak{g}^\dagger$ writes $\mathcal{L}_\mathfrak{g}\mathcal{L}_\mathfrak{g}^\dagger \mathcal{X}$ using a (local) orthonormal frame for T*X can be written schematically as

$$\mathcal{L}_\mathfrak{g}\mathcal{L}_\mathfrak{g}^\dagger \mathcal{X} = \nabla_k \nabla_i \nabla_i \nabla_k \mathcal{X} + \mathcal{R}_2(\nabla\nabla\mathcal{X}) + \mathcal{R}_1(\nabla\mathcal{X})$$

(A2.39)

where $\mathcal{R}_0$ and $\mathcal{R}_1$ are tensors whose norms are bounded by $c_*$ times the respective norms of the Riemann curvature tensor and its covariant derivative. Take the inner product of both sides of (A2.16) with $\sigma^2 \mathcal{X}$ and then integrate both sides of the resulting equation over V. Having done this, then judicious (and multiple) integration by parts leads directly to an inequality of the following sort:

$$\int_V \sigma^2 |\nabla\nabla\mathcal{X}|^2 \le \int_V \sigma^2 |\mathcal{L}_\mathfrak{g}^\dagger \mathcal{X}|^2 + \mathcal{T}_1' + \mathcal{T}_2' + \mathcal{T}_3' + \mathcal{T}_4',$$

(A2.40)

where the terms $\{\mathcal{T}_i'\}_{i=1,2,3}$ are as follows:

- $\mathcal{T}_1' = c_\diamond \int_V |d\sigma|^2 |\mathcal{X}||\nabla\nabla\mathcal{X}|$.
- $\mathcal{T}_2' = c_\diamond \int_V (\sigma|d^\dagger d\sigma| + \sigma^2 |\mathcal{R}_\mathfrak{g}|)|\mathcal{X}||\nabla\nabla\mathcal{X}|$.
- $\mathcal{T}_3' = c_\diamond \int_V \sigma|d\sigma||\nabla\mathcal{X}||\nabla\nabla\mathcal{X}|$.
- $\mathcal{T}_4' = c_\diamond \int_V \sigma^2 (|\mathcal{R}_\mathfrak{g}||\nabla\mathcal{X}|^2 + |\nabla\mathcal{R}_\mathfrak{g}||\mathcal{X}||\nabla\mathcal{X}|)$.

(A2.41)

This paragraph talks about the $\mathcal{T}_1$ and $\mathcal{T}_2$ terms in (A2.41), and the next two paragraphs say more about other terms. With regards to $\mathcal{T}_1$: Since the support of $d\sigma$ is in $V-V_\varepsilon$, the triangle inequality with the fourth bullet in (A2.38) bounds $\mathcal{T}_1'$ by the sum

$$\tfrac{1}{100} \int_V \sigma^2 |\nabla\nabla\mathcal{X}|^2 + c_\diamond \varepsilon^{-4} \int_{V-V_\varepsilon} |\mathcal{X}|^2.$$

(A2.42)

The triangle inequality with the third bullet in (A2.38) bounds $\mathcal{T}_2'$ by the sum of what is written in (A2.42) (with a larger version of $c_\diamond$) and

$$c_\Diamond |\mathcal{R}_g|_0^2 \int_V |\mathcal{X}|^2 \ .$$

(A2.43)

Meanwhile, the $\mathcal{T}_3'$ term in (A2.41) is observedly no larger than

$$\tfrac{1}{100} \int_V \sigma^2 |\nabla\nabla\mathcal{X}|^2 + c_\Diamond \int_{V-V_\varepsilon} |d\sigma|^2 |\nabla\mathcal{X}|^2 \ .$$

(A2.44)

An appeal to the fourth bullet of (A2.38) bounds the integrand in the right most term of (A2.44) by $c_\Diamond \varepsilon^{-2} \sigma(1-\sigma)|\nabla\mathcal{X}|^2$. After writing $|\nabla\mathcal{X}|^2$ as $-\tfrac{1}{2} d^\dagger d |\mathcal{X}|^2 + \langle \mathcal{X}, \nabla^\dagger \nabla \mathcal{X}\rangle$, the triangle inequality, integration by parts and the third bullet of (A2.38) lead to the bound on the left most term in (A2.44) by the expression in (A2.42) (with a larger version of $c_\Diamond$).

The $\mathcal{T}_4'$ term in (A2.41) is observedly no larger than

$$c_\Diamond |\nabla\mathcal{R}_g|_0^{4/3} \int_V |\mathcal{X}|^2 + |\nabla\mathcal{R}_g|_0^{2/3} \int_V |\sigma|^2 |\nabla\mathcal{X}|^2 \ .$$

(A2.45)

Writing the $|\nabla\mathcal{X}|^2$ in the right most term of (A2.45) as $-\tfrac{1}{2} d^\dagger d |\mathcal{X}|^2 + \langle \mathcal{X}, \nabla^\dagger \nabla \mathcal{X}\rangle$, integration by parts and the third bullet of (A2.38) bounds the right most term in (A2.45) by

$$\tfrac{1}{100} \int_V \sigma^2 |\nabla\nabla\mathcal{X}|^2 + c_\Diamond |\nabla\mathcal{R}_g|_1^{4/3} \int_V |\mathcal{X}|^2 + c_\Diamond \varepsilon^{-4} \int_{V-V_\varepsilon} |\mathcal{X}|^2 \ .$$

(A2.46)

Add up the preceding bounds for $\mathcal{T}_1'$, $\mathcal{T}_2'$, $\mathcal{T}_3'$ and $\mathcal{T}_4'$ to obtain the bound in Lemma A2.2.

The construction of a function $\sigma$ that obeys (A2.38) completes the proof of Lemma A2.2. To start the construction, introduce the standard 'cut-off' function on $\mathbb{R}$, denoted here by $\eta$, which is defined to be zero on $(-\infty, 0]$ and defined by the rule $t \to \eta(t) = e^{-1/t}$ on $(0, \infty)$. This function has the property that

$$|d\eta|^2 \le (\tfrac{4}{e})^4 \eta.$$

(A2.47)

With $\eta$ in hand, let $\eta_\Diamond$ denote the function on $\mathbb{R}$ that is defined by the rules whereby $\eta_\Diamond(t) = 0$ for $t \le \tfrac{1}{3}$ and

$$\eta_\Diamond(t) = \exp\left(-\tfrac{\eta(2/3 - t)}{t - 1/3}\right)$$

(A2.48)

for $t \geq \frac{1}{3}$. Since $\eta(\frac{2}{3} - t) = 0$ for $t \geq \frac{2}{3}$, this function $\eta_\diamond$ is equal to 1 where $t \geq \frac{2}{3}$. It follows from (A2.47) that $|d\eta_\diamond|^2 \leq c_\diamond \eta_\diamond$. Now define $\sigma$ by the rule $x \to \eta_\diamond(\mathfrak{f}(x))$.

## A3. The $\mathcal{X}_K$ part of a pair $(\mathcal{X}_T, \mathcal{X}_K)$

Fix an element (to be denoted by $\mathcal{X}$) in the $\iota$-invariant kernel of the $\mathfrak{g} = \mathfrak{g}_R$ version of $\mathcal{L}^\dagger_\mathfrak{g}$. This element is again written in the manner of (5.4) and (5.5) as a pair $(\mathcal{X}_T, \mathcal{X}_K)$. This section helps set the stage for the proofs of Propositions 5.3 and 5.4 by analyzing the $\mathcal{X}_K$ part of the pair $(\mathcal{X}_T, \mathcal{X}_K)$. Propositions A3.1 and A3.4 and A3.5 and A3.7 and A3.13 summarize the main results in this section.

### a) Fourier modes on $(S^3 - N_K) \times S^1$

Let $g_K$ denote the product metric on $(S^3 - K) \times S^1$ with the metric on $S^3 - K$ being the finite volume, constant sectional curvature -1 metric; and with the metric on $S^1$ being the Euclidean metric comes from writing $S^1$ as $\mathbb{R}/(2\pi\mathbb{Z})$. The group $S^1$ acts on $(S^3 - K) \times S^1$ via the constant rotations of the $S^1$ factor; and the metric $g_K$ is invariant with respect to this $S^1$ action. Since $g_K$ is $S^1$ invariant, the canonically lift of the $S^1$ action to the tangent bundle of $(S^3 - K) \times S^1$ defines corresponding lifts to the cotangent bundle and to other tensor bundles, in particular to $\Lambda^+ \otimes \Lambda^+$ and $\Lambda^+ \otimes \Lambda^-$. The resulting $S^1$ action on the spaces of sections of these bundles will be used to simplify the analysis of the kernel of the operator $\mathcal{L}^\dagger_{g_K}$ by first writing this kernel as a direct sum of $S^1$-character subspaces; and then writing the operator $\mathcal{L}^\dagger_{g_K}$ on each character subspace as an operator between tensor bundles on the 3-manifold $S^3 - K$. This is a five part task that occupies the remainder of this subsection. Note in this regard that what is done below is completely analogous to what is done in [AV] except with regards to notation. (There is no practical benefit in using [AV]'s abstract notation because the subsequent calculations must be done using *some* chosen frame for the tangent bundle and associated tensor bundles, which is what [AV] do in their proofs anyway.)

*Part 1*: The equation $\mathcal{L}^\dagger_{g_K} \mathcal{X} = 0$ is a second order equation that can be written as a coupled system of first order equations for a pair $(\mathcal{X}, \mathcal{A})$ with $\mathcal{A}$ being a section of the bundle $\Lambda^+ \otimes T^*((S^3 - N_K) \times S^1)$. These equations assert that

$$\eta^b_{nk} \nabla_k \mathcal{X}^{ab} - \mathcal{A}^a_n = 0 \quad and \quad \eta^c_{ni} \nabla_i \mathcal{A}^a_n + \mathcal{B}^{b,c} \mathcal{X}^{ab} = 0 \,.$$
(A3.1)

The left most equation is the $g_K$ analog of (A2.12) and, given the left most equation, then the right most equation says that the expression on the right hand side of (A2.4) is zero.

*Part 2*: Let s again denote coordinate function that appears in (3.3). The $s \geq 0$ part of $S^3-K$ is diffeomorphic to $[0, \infty) \times T$ with T being a torus. As noted in Section 3b, there is a diffeomorphism that writes the constant sectional curvature -1 metric on $S^3-K$ as in (3.3). Since this metric is invariant under the translation action of the torus T on itself, the metric $g_K$ on the $s \geq 0$ part of $(S^3-N_K) \times S^1$ is invariant under the T action. It follows as a consequence that any element in the kernel of $\mathcal{L}_{g_K}^\dagger$ on the $s \geq 0$ part of $(S^3-N_K) \times S^1$ can be written as a sum of Fourier modes with respect to this T action, and that each Fourier mode is in the kernel of $\mathcal{L}_{g_K}^\dagger$. Such a Fourier sum is indexed by integer pairs; the term in the sum that is labeled by an integer pair $(k_1, k_2)$ has the form

$$\mathfrak{X}_{(k_1,k_2)}(s,\theta) e^{i(k_1\tau_1 + k_2\tau_2)}$$

(A3.2)

with $(\tau_1, \tau_2)$ being the $\mathbb{R}^2/(2\pi\mathbb{Z}^2)$ coordinates for T. The following proposition will be used to justify a subsequent focus on the behavior of the $(k_1 = 0, k_2 = 0)$ Fourier mode.

**Proposition A3.1**: *There exists $\kappa > 10$ with the following significance: Fix $r > \kappa^2$.*

- *Suppose that $\mathcal{X}$ is in the kernel of $\mathcal{L}_{g_K}^\dagger$ on the $s \in [0, r]$ part of $(S^3-K) \times S^1$. Assume that the sum of the $L^2$ norms of $\mathcal{X}$ on the $s \in [0, \kappa]$ and the $s \in [r - \kappa, r]$ parts of $(S^3-N_K) \times S^1$ is equal to 1; and assume that the $(k_1 = 0, k_2 = 0)$ Fourier mode from $\mathcal{X}$ for the action of the torus T is absent. Then $|\mathcal{X}| \leq \kappa(\exp(-\kappa^{-1}e^s) + \exp(-\kappa^{-1}e^{(r-s)}))$ at the points where $s \in [\kappa, r - \kappa]$.*
- *Let $\mathcal{X}_K$ denote the $(S^3-N_K) \times S^1$ part of a pair $(\mathcal{X}_T, \mathcal{X}_K)$ that is described by (5.4) and (5.5) with the integral of $|\mathcal{X}_K|^2$ on the $s \leq \kappa$ part of $(S^3-N_K) \times S^1$ equal to 1. Define $\mathcal{X}_K^\perp$ on the $s \geq 0$ part of $(S^3-N_K) \times S^1$ to be the sum of the $(k_1, k_2) \neq (0,0)$ Fourier modes from $\mathcal{X}_K$ for the T action. Then $|\mathcal{X}_K^\perp|^2 \leq \kappa \exp(-\kappa^{-1}e^s)$ where $\kappa \leq s \leq R + \ln(\frac{1}{4}t_*) - \kappa$ in $(S^3-N_K) \times S^1$.*
- *Suppose that $\mathcal{X}$ is in the kernel of $\mathcal{L}_{g_K}^\dagger$ on the whole of the $s \geq 0$ part of $(S^3-K) \times S^1$. Assume that the integral of $|\mathcal{X}|^2$ on the $s \geq 0$ part of $(S^3-K) \times S^1$ is equal to 1 and that the $(k_1 = 0, k_2 = 0)$ Fourier mode from $\mathcal{X}$ for the action of the torus T is absent. Then $|\mathcal{X}| \leq \kappa \exp(-\kappa^{-1}e^s)$ on the $s \geq \kappa$ part of $(S^3-K) \times S^1$.*

***Proof of Proposition A3.1***: Define $R = r - \ln(\frac{1}{4}t_*)$. Use this value of R to define the disk $D \subset T'$ with center $p_*$ and radius $e^{-R}$. Let $\mathbb{D}$ again denote the disk in the torus $T'$ centered at the point $p_*$ with radius $\frac{1}{4}t_*$; and use this same R to define the metric $g_T$ on $T \times \mathbb{D}$.

The claim in the first bullet follows from (A2.10) with some some standard elliptic regularity theorems because the $s \in [0, r]$ part of $(S^3-K) \times S^1$ with the metric $g_K$ is conformal to $T \times (\mathbb{D}-D)$ with the metric $g_T$. The claim in the second bullet of the Proposition follows from the first using the conclusions of Proposition A2.1. The third bullet's claim follows from what is said by the first bullet for values of $r$ in an increasing, unbounded set.

*Part 3*: Let $\theta$ denote the Euclidean $\mathbb{R}/(2\pi\mathbb{Z})$ valued coordinate on $S^1$. The projection from $(S^3-K) \times S^1$ to $S^1$ allows $\theta$ to be viewed as an $\mathbb{R}/(2\pi\mathbb{Z})$ coordinate on the manifold $(S^3-K) \times S^1$. The convention in what follows is to identify the 3-manifold $S^3-K$ with the $\theta = 0$ slice in $(S^3-K) \times S^1$. The lifted $S^1$ action to the various tensor bundles on $(S^3-K) \times S^1$ is used implicitly to identify these bundles with their restriction to the $\theta = 0$ slice, which is to say $S^3-K$. Granted these identifications, any section of a tensor bundle over $(S^3-K) \times S^1$ is viewed as an $S^1$ dependent section of the bundle's restriction to $S^3-K$. The derivative with respect to $\theta$ acting on such a section is denoted by $\partial_\theta$ in what follows.

A related convention is to use an oriented orthonormal frame for $T^*(S^3-K)$ to define corresponding frames for $T^*((S^3-K) \times S^1)$ and for the bundles $\Lambda^+$ and $\Lambda^-$. To say more, suppose that $\{e^1, e^2, e^3\}$ is an oriented orthonormal frame for $T^*(S^3-K)$. Then $\{e^1, e^2, e^3, e^4 = d\theta\}$ defines the corresponding frame for $T^*((S^3-K) \times S^1)$; and with this understood, then the formulae in (5.7) define respective frames for $\Lambda^+$ and $\Lambda^-$. Of particular note is that the formulae in (5.7) are SO(3) equivariant and so they identify $\Lambda^+$ and $\Lambda^-$ on $S^3-K$ with $T^*(S^3-K)$. Meanwhile the frame $\{e^1, e^2, e^3, e^4 = d\theta\}$ identifies $T^*((S^3-K) \times S^1)$ on $S^3-K$ with the direct sum of $T^*(S^3-K)$ and the product $\mathbb{R}$ bundle (the span of $d\theta$). All of these identifications are used implicitly below to rewrite (A3.1).

*Part 4*: Let U again denote an open set in $S^3-K$ and let $C = (\mathcal{X}, \mathcal{A})$ denote a pair consisting of a section over $U \times S^1$ of the bundle $\Lambda^+ \otimes \Lambda^+$ and a section of the bundle $\Lambda^+ \otimes T^*(U \times S^1)$. This pair is assumed to obey (A2.15). What is said by the first two paragraphs of Part 3 can be used to view $\mathcal{A}$ as a pair $(\mathfrak{s}, \mathfrak{j})$ with $\mathfrak{s}$ being an $S^1$ dependent, symmetric, traceless section of $T^*U \otimes T^*U$ and with $\mathfrak{j}$ being an $S^1$ dependent section of $T^*U$: The pair $(\mathfrak{s}, \mathfrak{j})$ determines $\mathcal{A}$ (and vice-versa) by the rule depicted below in (A3.3). The convention in (A3.3) and subsequently has Latin indices taking values from $\{1,2,3\}$. As done previously, repeated indices are summed. Equation (A3.3) has $\{\varepsilon^{abc}\}_{a,b,c=1,2,3}$ denoting the components of the completely anti-symmetric 3 tensor with $\varepsilon^{123}$ equal to 1.

- $\mathcal{A}^a{}_4 = \mathfrak{j}^a$.
- $\mathcal{A}^a{}_b = \mathfrak{s}^{ab} + \frac{1}{2} \varepsilon^{abc} \mathfrak{j}^c$.

(A3.3)

By way of an explanation for the second bullet, the definition of $\mathcal{A}$ given in (A3.1) implies that it obeys $\eta^a_{ik}\mathcal{A}^a_k = 0$. Take $i = 4$ in this identity to see that the $3 \times 3$ matrix with components $\{\mathcal{A}^a_b\}_{a,b=1,2,3}$ is traceless, and thus so is its symmetric part (which is by definition $\mathfrak{s}$). Take $i \in \{1, 2, 3\}$ to see that the anti-symmetric part of this matrix can be written in terms of $j$ as $\frac{1}{2}\varepsilon^{abc}j^c$.

With $\mathcal{X}$ understood to be an $S^1$ dependent section of $T^*U \otimes T^*U$ and with $j$ and $\mathfrak{s}$ as just described, then the left most equation in (A3.1) can be written as follows:

- $\nabla_b \mathcal{X}^{ab} + j^a = 0$.
- $\partial_\theta \mathcal{X}^{ab} + \frac{1}{2}(\varepsilon^{bcd}\nabla_c \mathcal{X}^{ad} + \varepsilon^{acd}\nabla_c \mathcal{X}^{bd}) - \mathfrak{s}^{ab} = 0$.

(A3.4)

Given what is said in Part 2, the right most equation in (A3.1) asserts the vanishing of an $S^1$ dependent section of $T^*U \otimes T^*U$, this being the section given by the expression on the equation's left hand side. Because the matrix $\mathcal{B}^{ac}$ in this expression is $\delta^{ac}$, the respective symmetric and anti-symmetric parts of the right most equation in (A3.1) say

- $\partial_\theta \mathfrak{s}^{ab} + \frac{1}{2}(\varepsilon^{bcd}\nabla_c \mathfrak{s}^{ad} + \varepsilon^{acd}\nabla_c \mathfrak{s}^{bd}) - \frac{1}{4}(\nabla_b j^a + \nabla_a j^b) - \frac{1}{2}\delta^{ab}(\nabla_d j^d) + \mathcal{X}^{ab} = 0$.
- $\partial_\theta j^a + \frac{3}{2}\varepsilon^{abc}\nabla_b j^c - \nabla_d \mathfrak{s}^{da} = 0$.

(A3.5)

The trace of the top equation in (A3.5) leads to the identity $\nabla_a j^a = 0$; and so the top equation in (A3.5) is equivalent to the equation

$$\partial_\theta \mathfrak{s}^{ab} + \frac{1}{2}(\varepsilon^{bcd}\nabla_c \mathfrak{s}^{ad} + \varepsilon^{acd}\nabla_c \mathfrak{s}^{bd}) - \frac{1}{4}(\nabla_b j^a + \nabla_a j^b) + \mathcal{X}^{ab} = 0 .$$

(A3.6)

Meanwhile, the identity $\nabla_k \mathcal{A}^a_k = 0$ (due to the left most equation in (A3.1)) says that

$$\partial_\theta j^a + \nabla_b \mathfrak{s}^{ab} + \frac{1}{2}\varepsilon^{abc}\nabla_b j^c = 0 ;$$

(A3.7)

and this with the lower equation in (A3.5) are equivalent to the identities

$$\partial_\theta j^a + \varepsilon^{abc}\nabla_b j^c = 0 \quad and \quad \nabla_b \mathfrak{s}^{ab} = \frac{1}{2}\varepsilon^{abc}\nabla_b j^c .$$

(A3.8)

By way of a summary, the equations in (A3.4) and (A3.5) are equivalent to the autonomous equations for $j$ that follow

- $\partial_\theta j^a + \varepsilon^{abc}\nabla_b j^c = 0$ ,
- $\nabla_d j^d = 0$ ;

(A3.9)

and equations for $\mathfrak{s}$ and $\mathcal{X}$ that have j appearing as a 'source' term:

- $\nabla_d \mathfrak{s}^{ad} = \frac{1}{2} \varepsilon^{abc} \nabla_b j^c$,
- $\nabla_b \mathcal{X}^{ab} = -j^a$,
- $\partial_\theta \mathfrak{s}^{ab} + \frac{1}{2} (\varepsilon^{bcd} \nabla_c \mathfrak{s}^{ad} + \varepsilon^{acd} \nabla_c \mathfrak{s}^{bd}) + \mathcal{X}^{ab} = \frac{1}{4} (\nabla_b j^a + \nabla_a j^b)$.
- $\partial_\theta \mathcal{X}^{ab} + \frac{1}{2} (\varepsilon^{bcd} \nabla_c \mathcal{X}^{ad} + \varepsilon^{acd} \nabla_c \mathcal{X}^{bd}) - \mathfrak{s}^{ab} = 0$.

(A3.10)

*Part 5*: Let U again denote an open set in $S^3 - K$ and let $C = (\mathcal{X}, \mathfrak{s}, j)$ denote an $S^1$-dependent solution to (A3.9) and (A3.10) on U. Since the equations in (A3.9) and (A3.10) are $S^1$ invariant, the solution $\mathfrak{C}$ can be written as a sum of Fourier modes with respect to the $S^1$ action with each mode obeying (A3.9) and (A3.10). The modes are indexed by $\mathbb{Z}$ with a given $n \in \mathbb{Z}$ mode having the form $\mathfrak{C} e^{in\theta}$ with $\partial_\theta \mathfrak{C} = 0$.

So as to keep the notation in check, the following convention will be used when dealing with a single Fourier mode solution to (A3.9) and (A3.10): Supposing that $n \in \mathbb{Z}$ and that $C$ is a ($\mathbb{C}$-valued) solution to (A3.9) and (A3.10) with single Fourier mode $S^1$ dependence given by $e^{in\theta}$, then $C$ will be written as $(\mathcal{X}, \mathfrak{s}, j) e^{in\theta}$ with $\mathcal{X}, \mathfrak{s}$ and $j$ being independent of the $S^1$ parameter $\theta$. The triple $(\mathcal{X}, \mathfrak{s}, j)$ obey the following equations in lieu of (A3.9) and (A3.10):

EQUATIONS FOR j:

- $in j^a + \varepsilon^{abc} \nabla_b j^c = 0$,
- $\nabla_d j^d = 0$.

(A3.11)

EQUATIONS FOR $\mathcal{X}$ AND $\mathfrak{s}$:

- $\nabla_d \mathfrak{s}^{ad} = \frac{1}{2} \varepsilon^{abc} \nabla_b j^c$,
- $\nabla_b \mathcal{X}^{ab} = -j^a$,
- $in \mathfrak{s}^{ab} + \frac{1}{2} (\varepsilon^{bcd} \nabla_c \mathfrak{s}^{ad} + \varepsilon^{acd} \nabla_c \mathfrak{s}^{bd}) + \mathcal{X}^{ab} = \frac{1}{4} (\nabla_b j^a + \nabla_a j^b)$.
- $in \mathcal{X}^{ab} + \frac{1}{2} (\varepsilon^{bcd} \nabla_c \mathcal{X}^{ad} + \varepsilon^{acd} \nabla_c \mathcal{X}^{bd}) - \mathfrak{s}^{ab} = 0$.

(A3.12)

This single mode version of (A3.9) and (A3.10) are used below.

### e) Solutions $(\mathcal{X}, \mathfrak{s})$ to (A3.12) with j = 0

This subsection talks about the j = 0 version of (A3.12). To start, suppose U is an open set in $S^3 - K$, that $n \in \mathbb{Z}$, and that $(\mathcal{X}, \mathfrak{s})$ is a pair of symmetric, traceless sections of $\otimes^2 T^*U$ over U that obey the $j = 0$ version of (A3.12). Let $\mathfrak{T}_+$ and $\mathfrak{T}_-$ denote the symmetric, traceless sections $\mathcal{X}+i\mathfrak{s}$ and $\mathcal{X}-i\mathfrak{s}$ of $\otimes^2 T^*U$. Keep in mind that these need not be complex conjugate pairs. In any event, the complex conjugates of $\mathfrak{T}_+$ and $\mathfrak{T}_-$ are denoted respectively by $\overline{\mathfrak{T}}_+$ and $\overline{\mathfrak{T}}_-$. The $j = 0$ version of the equations in (A3.12) are equivalent to the two equations

- $\varepsilon^{bcd} \nabla_c \mathfrak{T}_+{}^{ad} + i(n+1) \mathfrak{T}_+{}^{ab} = 0$.
- $\varepsilon^{bcd} \nabla_c \mathfrak{T}_-{}^{ad} + i(n-1) \mathfrak{T}_-{}^{ab} = 0$.

(A3.13)

By way of an explanation, the anti-symmetric part of the two equations (A3.13) gives the $j = 0$ version of the equations in the first and second bullets of (A3.12) because $\mathfrak{T}_+{}^{ab}$ and $\mathfrak{T}_-{}^{ab}$ are traceless. The symmetric part of the equations in (A3.13) gives the $j = 0$ version of the equations in the third and fourth bullets of (A3.12).

The equations in (A3.13) have discrete symmetries that identify versions with differing integer n: If $\mathfrak{T}_+$ solves the top bullet's equation for a given integer n, then it solves the lower bullet's equation using $n+2$ in lieu of n. Conversely, if $\mathfrak{T}_-$ solves the lower bullet's equation for a given n, then it solves the upper bullet's equation with $n-2$ used in lieu of n. This is why it is sufficient when discussing (A3.13) to consider only the top equation. By way of a second symmetry, if $\mathfrak{T}_+$ solves the equation in the top bullet of (A3.13) for a given integer n, then its complex conjugate solves this same equation using the integer $-(n+2)$ in lieu of n; and thus it solves the lower equation for the integer $-n$. Likewise, if $\mathfrak{T}_-$ solves the lower equation for a given n, then its complex conjugate solves the upper equation for $-n$.

The rest of this subsection has five parts; they focus on the solutions to the equations in the top bullet of (A3.13).

*Part 1*: The lemma that follows directly talks about s-dependence of the $L^2$ norm of solutions to (A3.13) on $S^3 - N_K$.

**Lemma A3.2**: *Fix $r \geq 2$. Supposing that* n *is a given integer, let $\mathfrak{T}_+$ denote a solution on the $s \leq r$ part of $S^3 - K$ to the corresponding equation in the top bullet of (A3.13). If $\mathfrak{r} \in (1, r]$, then the following are true:*

- $\frac{i}{2} \int\limits_{s=\mathfrak{r}} \varepsilon^{3db} \overline{\mathfrak{T}}_+{}^{ab} \mathfrak{T}_+{}^{ad} = (n+1) \int\limits_{s \leq \mathfrak{r}} |\mathfrak{T}_+|^2$.
- $\int\limits_{s \leq r} |\mathfrak{T}_+|^2 \geq e^{2|n+1|(r-\mathfrak{r})} \int\limits_{s \leq \mathfrak{r}} |\mathfrak{T}_+|^2$.

*Proof of Lemma A3.2*: To prove the lemma's top bullet, take the inner product of the equation in the top bullet of (A3.13) with $\frac{i}{2} \overline{\mathfrak{T}}_+$ and then integrate the resulting identity over the s ≤ r part of $S^3-N_K$. Add the complex conjugate of the identity to get an $\mathbb{R}$ valued integration identity whose right hand side is the integral of $|\mathfrak{T}_+|^2$ over the s ≤ r part of $S^3-N_K$. Meanwhile, the left hand side of this integral identity is the integral of $\frac{i}{2}\varepsilon^{bcd}(\overline{\mathfrak{T}}_+^{ab}\nabla_c\mathfrak{T}_+^{ad} - \mathfrak{T}_+^{ad}\nabla_c\overline{\mathfrak{T}}_+^{ab})$ over the same domain. An integration by parts identifies the latter integral with the term on the left hand side of the lemma's top bullet.

To prove the second bullet of the lemma, use the top bullet to first see that

$$\int_{s=r}|\mathfrak{T}_+|^2 \geq 2|n+1|\int_{s\leq r}|\mathfrak{T}_+|^2 .$$

(A3.14)

Let $f$ denote the function on $[1, R+\ln(\frac{1}{4}t_*)]$ whose value at r is the integral on the right hand side of this inequality. The left integral in (A3.14) is the derivative of $f$. Therefore, (A3.14) says that $\frac{d}{dr}f \geq 2|n+1|f$ which implies that $f(r) \geq e^{2|n+1|(r-r)}f(r)$.

*Part 2*: With Proposition A3.1 in mind, the next lemma describes the solutions to (A3.13) on the s ≥ 0 part of $S^3-N_K$ that are invariant under the action of the torus T on this part of $S^3-N_K$. (These are the ($k_1 = 0, k_2 = 0$) versions of (A3.2).) The upcoming lemma writes the solutions using an oriented orthonormal frame for $T^*(S^3-N_K)$ on this s ≥ 0 part of $S^3-N_K$ that has the third frame 1-form $e^3$ being ds and the other frame 1-forms $\{e^1, e^2\}$ having the form $e^1 = e^{-s}\hat{e}^1$ and $e^2 = e^{-s}\hat{e}^2$ with $\{\hat{e}^1, \hat{e}^2\}$ being a T and s independent, oriented orthonormal frame for the flat metric m on T that appears in (3.3).

**Lemma A3.3**: *Fix* $n \in \mathbb{Z}$ *and let* $\mathfrak{T}$ *denote a* T-*invariant, symmetric, trace zero solution on the* s ≥ 0 *part of* $S^3-N_K$ *to the equation in the top bullet of (A3.13).*
- $\mathfrak{T}$ *extends to the whole of the* s ≥ 0 *part of* $S^3-K$ *as a solution to the equation in the top bullet of (A3.13).*
- $\mathfrak{T}$ *has the following form on the* s ≥ 0 *part of* $S^3-K$:
  a) $\mathfrak{T}^{33} = -(\mathfrak{T}^{11} + \mathfrak{T}^{22}) = 0$ *unless* n = -1 *when* $\mathfrak{T}^{33} = -(\mathfrak{T}^{11} + \mathfrak{T}^{22}) = e^{3s}t_a$ *with* $t_a \in \mathbb{R}$.
  b) $\mathfrak{T}^{13} = \mathfrak{T}^{23} = 0$ *unless* $n \in \{0, -2\}$.
     i) *If* n = 0, *then* $\mathfrak{T}^{13} = -i\mathfrak{T}^{23} = e^{3s}t_{b+}$ *with* $t_{b+} \in \mathbb{C}$.
     ii) *If* n = -2, *then* $\mathfrak{T}^{13} = i\mathfrak{T}^{23} = e^{3s}t_{b-}$ *with* $t_{b-} \in \mathbb{C}$
  $c_+$) $\mathfrak{T}^{11} - \mathfrak{T}^{22} + 2i\mathfrak{T}^{12} = e^{-ns}t_{c+}$ *with* $t_{c+} \in \mathbb{C}$.
  $c_-$) $\mathfrak{T}^{11} - \mathfrak{T}^{22} - 2i\mathfrak{T}^{12} = e^{(n+2)s}t_{c-}$ *with* $t_{c-} \in \mathbb{C}$.

*Proof of Lemma A3.3*: The equations in the top bullet of (A3.13) for $b \in \{1,2,3\}$ are as follows:

- $(\nabla_2 \mathfrak{T})^{a3} - (\nabla_3 \mathfrak{T})^{a2} + i(n+1)\mathfrak{T}^{a1} = 0$.
- $(\nabla_1 \mathfrak{T})^{a3} - (\nabla_3 \mathfrak{T})^{a1} - i(n+1)\mathfrak{T}^{a2} = 0$.
- $(\nabla_1 \mathfrak{T})^{a2} - (\nabla_2 \mathfrak{T})^{a1} + i(n+1)\mathfrak{T}^{a3} = 0$.

(A3.15)

Since $\frac{\partial}{\partial \tau_1}\mathfrak{T} = 0$ and $\frac{\partial}{\partial \tau_2}\mathfrak{T} = 0$ and since $\mathfrak{T}^{ab} = \mathfrak{T}^{ba}$, the a = 1, 2, 3 versions of the third equation in (A3.15) are algebraic equations that are equivalent to the following:

- $-\mathfrak{T}^{23} + i(n+1)\mathfrak{T}^{13} = 0$.
- $\mathfrak{T}^{13} + i(n+1)\mathfrak{T}^{23} = 0$.
- $i(n+1)\mathfrak{T}^{33} = 0$.

(A3.16)

The first two of these say that $\mathfrak{T}^{13} = \mathfrak{T}^{23} = 0$ if n ∉ {0, -2}. In the case n = 0, they say that $\mathfrak{T}^{13} = -i\mathfrak{T}^{23}$; and in the case n = -2, they say that $\mathfrak{T}^{13} = i\mathfrak{T}^{23}$. The third equation says that $\mathfrak{T}^{33} = 0$ unless n = -1, which implies the same for $\mathfrak{T}^{11} + \mathfrak{T}^{22}$ because $\mathfrak{T}$ is traceless. With the preceding understood, note that the a = 3 versions of the first and second bullet equations in (A3.15) assert that

- $2\mathfrak{T}^{23} - \frac{d}{ds}\mathfrak{T}^{23} + i(n+1)\mathfrak{T}^{13} = 0$.
- $2\mathfrak{T}^{13} - \frac{d}{ds}\mathfrak{T}^{13} - i(n+1)\mathfrak{T}^{23} = 0$.

(A3.17)

Keeping in mind that $\mathfrak{T}^{13} = -i\mathfrak{T}^{23}$ when n = 0 and that $\mathfrak{T}^{13} = i\mathfrak{T}^{23}$ when n = -2, these equations in either case assert that $\frac{d}{ds}\mathfrak{T}^{23} - 3\mathfrak{T}^{23} = 0$ whose solution is described by the two parts of Item b) of the lemma.

The a = 1 and a = 2 versions of the equations in the first two bullets of (A3.15) are jointly equivalent to equations for $\mathfrak{w}_+ = \mathfrak{T}^{11} - \mathfrak{T}^{22} + 2i\mathfrak{T}^{12}$ and $\mathfrak{w}_- = \mathfrak{T}^{11} - \mathfrak{T}^{22} - 2i\mathfrak{T}^{12}$ and the function $\mathfrak{T}^{33}$ (which is $-\mathfrak{T}^{11} - \mathfrak{T}^{22}$) that assert the following:

- $\frac{d}{ds}\mathfrak{w}_+ + n\mathfrak{w}_+ = 0$.
- $\frac{d}{ds}\mathfrak{w}_- - (n+2)\mathfrak{w}_- = 0$.
- $\frac{d}{ds}\mathfrak{T}^{33} - 3\mathfrak{T}^{33} = 0$ *when* n = -1 *and* $\mathfrak{T}^{33} = 0$ *when* n ≠ -1.

(A3.18)

The solutions to these equations are described by Items a) and c_+) and c_-) of the lemma.

*Part 3*: The following proposition says more about the significance with regards to $S^3$−K of the solutions from Lemma A3.3 when n ∉ {-2, -1, 0}. The cases when n is -2 or -1 or 0 are special and are discussed in Parts 4 and 5 of this subsection.

**Proposition A3.4:** *There exists* κ ≥ 1 *with the following significance. Fix* n ≥ 1 *or* n ≤ -3.

- *There is a unique symmetric, trace zero solution to the equation in the top bullet of (A3.13) on $S^3-K$ with the following properties: Denote this solution by $\mathfrak{T}_+^{(n)}$.*
  
  i) *If $n \geq 1$, then $\mathfrak{T}_+^{(n)}$ on the $s \geq 0$ part of $S^3-K$ can be written as*

$$\mathfrak{T}_+^{(n)} = e^{(n+2)s} \begin{pmatrix} 1 & i & 0 \\ i & -1 & 0 \\ 0 & 0 & 0 \end{pmatrix} + \alpha_n e^{-ns} \begin{pmatrix} 1 & -i & 0 \\ -i & -1 & 0 \\ 0 & 0 & 0 \end{pmatrix} + \mathfrak{r}_n .$$

  ii) *If $n \leq -3$, then $\mathfrak{T}_+^{(n)}$ on the $s \geq 0$ part of $S^3-K$ can be written as*

$$\mathfrak{T}_+^{(n)} = e^{|n|s} \begin{pmatrix} 1 & -i & 0 \\ -i & -1 & 0 \\ 0 & 0 & 0 \end{pmatrix} + \alpha_n e^{-|n+2|s} \begin{pmatrix} 1 & i & 0 \\ i & -1 & 0 \\ 0 & 0 & 0 \end{pmatrix} + \mathfrak{r}_n .$$

  *In both cases, $\alpha_n \in \mathbb{C}$ and the function $\mathfrak{r}_n$ is bounded and has no $(k_1=0, k_2=0)$ Fourier mode for the T action on the $s \geq 0$ part of $S^3$-K.*
- *In both cases above, the number $\alpha_n$ that appears obeys $|\alpha_n| \leq \kappa e^n$, and the function $\mathfrak{r}_n$ that appears obeys $|\mathfrak{r}_n| \leq \kappa \exp(-\kappa^{-1} e^s)$.*

*Proof of Proposition A3.4*: The proof will be given only for the case $n \geq 1$ because the argument for the case $n \leq -3$ is identical but for changing n to $-(n+2)$ and changing some signs. The argument for the $n \geq 1$ case has four steps.

<u>Step 1</u>: Fix a smooth, nondecreasing function of s to be denoted by $\beta$ with $\beta = 0$ for $s \leq 0$ and with $\beta = 1$ for $s \geq 1$. View $\beta$ as a function on $S^3-K$. The derivative of $\beta$ is denoted in what follows by $\beta'$. Let $\mathfrak{T}_\diamond$ denote the symmetric, traceless section of $\otimes^2 T^*(S^3-K)$ on the $s \geq 0$ part of $S^3-K$ given by

$$\mathfrak{T}_\diamond = e^{(n+2)s} \begin{pmatrix} 1 & i & 0 \\ i & -1 & 0 \\ 0 & 0 & 0 \end{pmatrix} .$$

(A3.19)

The plan is to find a symmetric, traceless section (to be denoted by q) of $\otimes^2 T^*(S^3-K)$ with a suitable $s \to \infty$ limit that solves the equation

$$\varepsilon^{bcd} \nabla_c q^{ad} + i(n+1) q = -\varepsilon^{b3d} \beta' e^{(n+2)s} \mathfrak{T}_\diamond^{ad}$$

(A3.20)

If q obeys (A3.20), then $\mathfrak{T}_+^{(n)} = \beta \mathfrak{T}_\diamond + q$ obeys the equation in the top bullet of (A3.13). Assume for the moment that q obeys (A3.20) and that

$$\int_{s \leq r} |q|^2 \leq c_* e^{2n}.$$
(A3.21)

(This condition rules out the stupid solution $-\beta e^{(n+2)s} \mathfrak{T}_\diamond$ to (A3.21).) Granted (A3.21), then what is asserted by the lemma about the behavior of q on the $s \geq 0$ part of $S^3 - K$ follows by invoking Lemma A3.2 to bound the T-invariant part of q where $s \geq 0$ on $S^3 - K$; and by invoking the the third bullet of Proposition A3.1 for the rest of q.

Step 2: Let $\mathfrak{L}_n$ denote the operator that maps $C^\infty(S^3 - K; \otimes^2 T^*(S^3 - K))$ to itself by the rule whereby if $\mathfrak{b}$ is a section of $\otimes^2 T^*(S^3 - K)$, then $\mathfrak{L}_n \mathfrak{b}$ is the section whose components are

$$(\mathfrak{L}_n \mathfrak{b})^{ab} = \varepsilon^{bcd} \nabla_c \mathfrak{b}^{ad} + i(n+1) \mathfrak{b}^{ab}.$$
(A3.22)

The equation in (A3.20) when written using this notation says that $(\mathfrak{L}_n q)^{ab} = -\varepsilon^{b3d} \beta' \mathfrak{T}_\diamond^{ad}$.

Let $S_0 \subset C^\infty(S^3 - K; \otimes^2 T^*(S^3 - K))$ denote the subspace of symmetric, traceless sections. A section $\mathfrak{b} \in S_0$ is said to have zero divergence when $\nabla_b \mathfrak{b}^{ab} = 0$. Let $S_{0\perp} \subset S_0$ denote the subspace of zero divergence sections. As it turns out, the operator $\mathfrak{L}_n$ maps $S_{0\perp}$ to itself. (This is is because the metric on $S^3 - K$ has zero traceless Ricci tensor.)

The formal $L^2$ adjoint of $\mathfrak{L}_n$ is $\mathfrak{L}_{-(n+2)}$. If $\mathfrak{u}$ is in $S_{0\perp}$ and solves the equation

$$(\mathfrak{L}_n \mathfrak{L}_{-(n+2)} \mathfrak{u})^{ab} = -\varepsilon^{b3d} \beta' e^{(n+2)s} \mathfrak{T}_\diamond^{ad},$$
(A3.23)

then $q = \mathfrak{L}_{-(n+2)} \mathfrak{u}$ is a section of $S_{0\perp}$ that obeys (A3.20). The equation in (A3.23) will be solved by first using a Bochner-Weitzenbock formula to write it (when $\mathfrak{u} \in S_{0\perp}$) as

$$(\nabla^\dagger \nabla \mathfrak{u})^{ab} + ((n+1)^2 - 3) \mathfrak{u}^{ab} = -\varepsilon^{b3d} \beta' e^{(n+2)s} \mathfrak{T}_\diamond^{ad}.$$
(A3.24)

The next step finds a smooth, square integrable solution to (A3.24) with square integrable first derivatives in the space $S_0$ of symmetric traceless sections of $\otimes^2 T^*(S^3 - K)$. The subsequent step proves that this $S_0$ solution has zero divergence (and thus it is in $S_{0\perp}$).

Step 3: If $n \geq 1$ or $n \leq -3$, then coefficient $(n+1)^2 - 3$ that multiplies $\mathfrak{u}$ on the left hand side of (A3.24) is positive. Granted this positivity, then a square integrable solution to (A3.24) can be found in the space of $L^2_1$ symmetric, traceless sections of $\otimes^2 T^*(S^3 - K)$ as the limit of a minimizing sequence for the functional

$$\mathfrak{u} \to \mathfrak{E}(\mathfrak{u}) = \tfrac{1}{2} \int_{S^3-K} (|\nabla \mathfrak{u}|^2 + ((n+1)^2 - 3)|\mathfrak{u}|^2) + \int_{S^3-K} \mathfrak{u}^{ab} \varepsilon^{b3d} \beta' e^{(n+2)s} \mathfrak{T}_{\lozenge}{}^{ad} \, .$$

(A3.25)

Since $\beta' e^{(n+2)s} \mathfrak{T}_{\lozenge}$ has compact support and since $(n+1)^2 - 3$ is positive, this functional is bounded below on the $L^2_1$ completion of the compactly supported elements in $S_0$; and it is also coercive on this same $L^2_1$ Hilbert space. It follows as a consequence that the function $\mathfrak{E}$ has a unique minimizer in this $L^2_1$ Hilbert space. Let $\mathfrak{u}$ denote the $L^2_1$ minimizer. Standard elliptic regularity proves that $\mathfrak{u}$ is smooth and that it solves (A3.24).

Step 4: This step completes the proof by showing that the solution $\mathfrak{u}$ from the previous step has zero divergence. To this end, let $v$ denote the divergence of $\mathfrak{u}$, this being the 1-form with components $v^a = \nabla_b \mathfrak{u}^{ab}$. To prove that $v = 0$, take the divergence of both sides of (A3.24). The divergence of the right hand side is zero and the divergence of the left hand side is

$$\nabla_b (\nabla^\dagger \nabla \mathfrak{u})^{ab} + ((n+1)^2 - 3) v^a \, .$$

(A3.26)

Commuting derivative writes the term $\nabla_b (\nabla^\dagger \nabla \mathfrak{u})^{ab}$ in (A2.39) as $\nabla^\dagger \nabla v^a + 4 v^a$. As a consequence, the expression in (A3.26) is zero if and only if

$$\nabla^\dagger \nabla v^a + ((n+1)^2 + 1) v^a = 0.$$

(A3.27)

This equation implies that $v = 0$ because, being a linear combination of derivatives of $\mathfrak{u}$, the 1-form $v$ is square integrable on $S^3-K$ and, as explained next, there are no non-trivial square integrable solutions to (A3.27).

To prove that (A3.27) has no non-trivial square integrable solutions, suppose for the sake of argument that $v$ is a square integrable solution. Fix a large integer N and then a function $\chi_N: S^3-K \to [0,1]$ that is equal to 1 where $s \leq N$, equal to 0 where $s \geq N+1$ and obeys $|d\chi_N| \leq 4$. Multiply both sides of (A3.27) by $\chi_N^2$ and then take the $L^2$ norm of both sides. An integration by parts and an appeal to the triangle inequality leads to the bound

$$\int_{s \leq N} (|\nabla v|^2 + |v|^2) \leq 128 \int_{N \leq s \leq N+1} |v|^2 \, .$$

(A3.28)

Since the $N \to \infty$ limit of the right hand side of (A3.28) is zero, it follows that $v$ is identically zero.

*Part 4*: The proposition that follows ties up one loose end by describing the $S^3$–K significance of the solutions from Lemma A2.7 for the cases when n = 0 and n = -2.

**Proposition A3.5**: *Suppose that* n = 0 *or that* n = -2.
- *The* n = 0 *and* n = -2 *versions of the equation in the top bullet of (A3.13) has a unique symmetric, traceless solution on the whole $S^3$–K with the following property: Denote the* n = 0 *solution by* $\mathfrak{T}_+^{(b)}$. *This solution can be written on the* s ≥ 0 *part of $S^3$–K as*

$$\mathfrak{T}_+^{(b)} = e^{3s} \begin{pmatrix} 0 & 0 & 1 \\ 0 & 0 & i \\ 1 & i & 0 \end{pmatrix} + \mathfrak{r}_b$$

with $|\mathfrak{r}_b|$ *bounded by* $\kappa \exp(-e^s/\kappa)$ *for some* $\kappa > 1$. *Meanwhile, the* n = -2 *solution is the complex conjugate of* $\mathfrak{T}_+^{(b)}$.
- *The* n = 0 *and* n = -2 *versions of the equation in the top bullet of (A3.13) have a unique symmetric, traceless solution on the whole $S^3$–K with the following property: Denote the* n = 0 *solution by* $\mathfrak{T}_+^{(c)}$. *This solution can be written on the* s ≥ 0 *part of $S^3$–K as*

$$\mathfrak{T}_+^{(c)} = e^{2s} \begin{pmatrix} 1 & i & 0 \\ i & -1 & 0 \\ 0 & 0 & 0 \end{pmatrix} + \alpha_0 \begin{pmatrix} 1 & -i & 0 \\ -i & -1 & 0 \\ 0 & 0 & 0 \end{pmatrix} + \mathfrak{r}_c$$

*with* $\alpha_0 \in \mathbb{C}$ *and with* $\mathfrak{r}_c$ *obeying* $|\mathfrak{r}_c| \leq \kappa \exp(-e^s/\kappa)$ *for some* $\kappa > 1$. *Meanwhile, the* n = -2 *solution is the complex conjugate of* $\mathfrak{T}_+^{(c)}$.

*Proof of Proposition A3.5*: The proof has six steps.

<u>Step 1</u>: In this step, the number n in (A3.13) can be any given integer. Let E denote for the moment the complexification of the vector bundle $T^*(S^3–K)$. An $SO(3; \mathbb{C})$ connection on E is defined by its covariant derivative; and of direct interest is the connection whose covariant derivative takes a section $\phi$ (with components $\phi^a$ with respect to an oriented orthonormal frame) to the section $\nabla^\mathbb{C} \phi$ with components

$$(\nabla^\mathbb{C}_b \phi)^a = \nabla_b \phi^a - i(n+1) \varepsilon^{abc} \phi^c.$$

(A3.29)

This curvature of this $SO(3;\mathbb{C})$ connection vanishes precisely when n = 0 and n = -2. This connection in the n = 0 (and n = -2) case is the canonical $SO(3:\mathbb{C})$ flat connection that comes from the representation (and its complex conjugate) of $\pi_1(S^3-K)$ in $PSL(2;\mathbb{C})$. The bundle E in this case is isomorphic to the vector bundle with fiber the Lie algebra of $SL(2;\mathbb{C})$ that is associated to the universal cover of M via this representation of $\pi_1(SO(3))$.

The covariant derivative in (A3.29) can be used to define an exterior covariant derivative that maps E valued 1-forms to E-valued 2-forms. This exterior covariant derivative is denoted in what follows by $d^\mathbb{C}$. The components of the metric Hodge dual of $d^\mathbb{C}\flat$ with respect to any given oriented, orthonormal frame for $T^*(S^3-K)$ are as follows:

$$(*d^\mathbb{C}\flat)^{ab} = \varepsilon^{bcd}\nabla_c\flat^{ad} + i(n+1)(\flat^{ba} - \delta^{ab}\flat^{cc}).$$

(A3.30)

The preceding paragraphs are relevant to (A3.13) because a $\mathbb{C}$-valued section of $\otimes^2 T^*(S^3-K)$ (such as $\mathfrak{T}_+$) can be viewed as an E-valued 1-form by writing $\otimes^2 T^*(S^3-K)$ which is $T^*(S^3-K) \otimes T^*(S^3-K)$ as $E \otimes T^*(S^3-K)$. Moreover, if $\mathfrak{T}_+$ is symmetric and traceless as a section of $\otimes^2 T^*(S^3-K)$, applying the operator $\mathfrak{L}_n$ to $\mathfrak{T}_+$ (with $\mathfrak{L}_n$ defined by (A3.22)) and viewing the result as a section of $E \otimes T^*(S^3-K)$ is identical to what is obtained by first viewing $\mathfrak{T}_+$ as a section of $E \otimes T^*(S^3-K)$ and then acting by $*d^\mathbb{C}$.

Step 2: Focus now on the n = 0 case. (The case n = -2 is identical but for notation.) In this cases, the exterior derivatives obey $d^\mathbb{C}\nabla^\mathbb{C} = 0$ because the $SO(3;\mathbb{C})$ connection is flat. Let $H^1(S^3-K;E)$ denote the kernel of $d^\mathbb{C}$ modulo the image of $\nabla^\mathbb{C}$. The identifications of the preceding step define a homomorphism from the kernel of $\mathfrak{L}_{n=0}$ to $H^1(S^3-K;E)$. This homomorphism is denoted by $\Phi$ in what follows. It is used in the next step to analyze the kernel of $\mathfrak{L}_{n=0}$ (and $\mathfrak{L}_{n=-2}$). The rest of this step proves that $\dim_\mathbb{C}(H^1(S^3-K;E)) = 1$.

Lefschetz duality for 3-manifolds with boundary can be used to derive the formula below for for $\dim(H^1(S^3-K;E))$ (see, e.g. Proposition 2.3 in [Sc]).

$$\dim(H^1(S^3-K;E) = \dim\ker(\hat{\imath}^*) + \dim(H^0(T;E)).$$

(A3.31)

This formula uses $\hat{\imath}: T \to S^3-K$ to denote the inclusion map of T into any favorite constant s slice of the $s \geq 0$ part of $S^3-K$; and it uses $H^0(T;E)$ to denote the vector space of sections of E along T that are annihilated by the $\nabla^\mathbb{C}$ directional covariant derivatives along T. It follows from Mostow rigidity that $\ker(\hat{\imath}^*) = 0$. (Mostow rigidity can be invoked because $H^1(S^3-K;E)$ can be viewed as the space of first order formal deformations of the flat $SL(2;\mathbb{C})$ connection from the canonical representation of $\pi_1(S^3-K)$ into $SL(2;\mathbb{C})$.) A calculation finds that $\dim_{\mathbb{C}}(H^0(T;E)) = 1$. Thus, $H^1(S^3-K;E)$ has dimension 1 over $\mathbb{C}$ also.

Step 3: This step analyzes the kernel of the map $\Phi$; and it then uses the results to prove the first bullet of Proposition A3.5.

To begin, suppose that $\mathfrak{T}$ is a symmetric, traceless section of $\otimes^2 T^*(S^3-K)$ that is annihilated by the n = 0 version of $\mathfrak{L}_n$. Assume in addition that $\Phi(\mathfrak{T}) = 0$. This means that there is a section of E to be denoted by $\phi$ such that the incarnation of $\mathfrak{T}$ as a section of $E \otimes T^*(S^3-K)$ obeys $\mathfrak{T} = \nabla^{\mathbb{C}}\phi$. This one equation is equivalent to the following three assertions:

- $\mathfrak{T}^{ab} = \frac{1}{2}(\nabla_b \phi^a + \nabla_a \phi^b)$.
- $\varepsilon^{abc}\nabla_b \phi^c + 2i\phi^e = 0$.

(A3.32)

An argument much like that used to prove Lemma A3.2 leads from the second bullet in (A3.32) to the inequality

$$\int_{s\leq r} |\phi|^2 \geq e^{4r} \int_{s\leq 0} |\phi|^2 \ .$$

(A3.33)

This inequality has the following implication with regards to what is said by Lemma A3.3: Let $\mathfrak{T}_\diamond$ denote the solution on the $s \geq 0$ part of $S^3-K$ to the equation $\mathfrak{L}_0 \mathfrak{T} = 0$ that comes from Lemma A3.3 by taking only $t_{c-}$ from Item c ) not zero. Let $\beta$ denote the function from Step 1 of the proof of Proposition A3.4, and suppose that there exists a bounded, traceless, symmetric section, q, of $\otimes^2 T^*(S^3-K)$ such that $\beta \mathfrak{T}_\diamond + q$ is annihilated by $\mathfrak{L}_0$ on the whole of $S^3-K$. Then $\Phi(\beta\mathfrak{T}_\diamond + q) \neq 0$.

By way of a contrast, let $\mathfrak{T}_\diamond$ now denote the solution on the $s \geq 0$ part of $S^3-K$ to the equation $\mathfrak{L}_0 \mathfrak{T} = 0$ that comes from Lemma A3.3 by taking only $t_{b+}$ from Item bii) not zero. Using the notation from Lemma A3.3, let $\phi_\diamond$ denote the 1-form on the $s \geq 0$ part of $S^3-K$ with components $\phi_\diamond^1 = -i\phi_\diamond^2 = \frac{1}{2} e^{3s} t_{b+}$ and with $\phi_\diamond^3 = 0$. This $\phi_\diamond$ obeys (A3.32) with $\mathfrak{T} = \mathfrak{T}_\diamond$. This last fact has the following consequence: The upcoming Lemma A3.6 supplies a section $\psi$ of $T^*(S^3-K)$ such that $\phi = \beta\phi_\diamond + \psi$ obeys the second bullet of (A3.32) and such that $|\nabla\psi|^2 + e^{-s}|\nabla\nabla\psi|^2$ is integrable on $S^3-K$. Construct $\mathfrak{T}$ from $\phi$ using the formula in the top bullet of (A3.32). Then $\mathfrak{T}$ obeys the equation $\mathfrak{L}_0 \mathfrak{T} = 0$ on all of $S^3-K$ and $\mathfrak{T} - \mathfrak{T}_\diamond$ and $e^{-s}\nabla(\mathfrak{T}-\mathfrak{T}_\diamond)$ are square integrable on the $s \geq 0$ part of $S^3-K$. By construction, $\Phi(\mathfrak{T}) = 0$ in $H^1(S^3-K;E)$.

It follows as a consequence of the third bullet of Proposition A3.1 and Lemma A3.3 that $\mathfrak{T}$ is the $\mathfrak{T}_+^{(b)}$ extension of $\mathfrak{T}_\diamond$ to the whole of $S^3-K$ that is sought by the top bullet of Proposition A3.5. The fact that $\mathfrak{T}$ is the unique extension follows from Lemma

A3.2 because the large s asymptotics of the difference between any two extensions (if not zero) would run afoul of this lemma.

    Step 4: The lemma below supplies the 1-form $\psi$ for Step 3. It is also used subsequently for a different purpose.

**Lemma A3.6**: *Fix a non-zero, real number to be denoted by* $m$. *Let* $\eta$ *denote a smooth 1-form on* $S^3-K$ *with* $\eta$ *and the function* $*d*\eta$ *being square integrable.*
- *There exists a unique section* $\psi$ *of* $T^*(S^3-K)$ *obeying the following two conditions:*
  i) $\varepsilon^{abc}\nabla_b\psi_c + mi\psi^a = \eta^a$.
  ii) $\int_{S^3-K} (|\nabla\psi|^2 + |\psi|^2) < \infty$.
- *There exists* $\kappa \geq 1$ *that can be taken to depend only on a positive lower bound for* $m$ *such that the integral in Item ii) of the first bullet is bounded by* $\kappa(\int_{S^3-K} (|\eta|^2 + |d*\eta|^2))$.

*Proof of Lemma A3.6*: Supposing that $\upsilon$ is s smooth 1-form on $S^3-K$ with compact support, then

$$\int_{S^3-K} |\nabla\upsilon|^2 \geq 2 \int_{S^3-K} |\upsilon|^2$$

(A3.34)

because of the Bochner-Weitzenboch formula $\nabla^\dagger\nabla\upsilon = (d^\dagger d + dd^\dagger)\upsilon + 2\upsilon$. Now suppose that $\eta$ is a square integrable, smooth 1-form on $S^3-K$. Let $\mathfrak{E}$ denote the functional on the space of smooth square integrable 1-forms on $S^3-K$ with square integrable covariant derivative that is defined by the rule

$$\upsilon \to \mathfrak{E}_\eta(\upsilon) = \tfrac{1}{2} \int_{S^3-K} (|\nabla\upsilon|^2 + (m^2-2)|\upsilon|^2) + \int_{S^3-K} \langle\upsilon,\eta\rangle .$$

(A3.35)

Here and elsewhere, $\langle\,,\,\rangle$ denotes the metric inner product. It follows from (A3.34) that

$$\mathfrak{E}_\eta(\upsilon) \geq c_m^{-1}(\int_{S^3-K} (|\nabla\upsilon|^2 + |\upsilon|^2)) - c_m \int_{S^3-K} |\eta|^2$$

(A3.36)

with $c_m \geq 1$ being a number that depends only on $m$. It follows as a consequence that $\mathfrak{E}_\eta$ is bounded from below. This being the case, then standard variational arguments and standard elliptic regularity arguments (much like the ones used to prove Proposition

A3.4) can be used to prove the following: The functional $\mathfrak{E}_\eta$ has a unique minimizer; and this minimizer (denoted by $\upsilon$) is the unique, square integrable 1-form on $S^3-K$ that obeys

$$\nabla^\dagger \nabla \upsilon + (m^2 - 2)\upsilon = \eta \ .$$

(A3.37)

The sequence of cut-off functions $\{\chi_N\}_{N=1,2,\ldots}$ from Step 4 of the proof of Proposition A3.4 can be used to prove that $\nabla \upsilon$ and $\nabla \nabla \upsilon$ are square integrable; and that their $L^2$ norms and the $L^2$ norm of $\upsilon$ are bounded by a $\eta$-independent multiple of the $L^2$ norm of $\eta$.

Let $\mathfrak{D}$ denote the operator mapping smooth 1-forms to smooth 1-forms by the rule

$$\psi \to \mathfrak{D}\psi = *d\psi + mi\psi.$$

(A3.38)

Use $\mathfrak{D}^\dagger$ to denote the formal adjoint of $\mathfrak{D}$. The formula for $\mathfrak{D}^\dagger$ is obtained from (A3.38) by changing i to -i. Lemma A2.10 is asking for a 1-form $\psi$ that solves the equation $\mathfrak{D}\psi = \eta$. If $\upsilon$ obeys (A3.37), then a 1-form that does this is $\psi = \mathfrak{D}^\dagger \upsilon - \frac{1}{2i} d(*d*\upsilon)$. Since $\upsilon$, $\nabla \upsilon$ and $\nabla \nabla \upsilon$ are all square integrable, this 1-form $\psi$ is square integrable. Since $\nabla \nabla \upsilon$ is square integrable, an a priori bound for the $L^2$ norm of $\nabla \psi$ follows from an a priori bound on the $L^2$ norm of $\nabla \nabla (*d*\upsilon)$. There is such a bound, it is a $\eta$-independent multiple of the $L^2$ norms of $\eta$ and $*d*\eta$. Such a bound is derived by first taking the divergence of both sides of (A3.37) to see that the function $*d*\upsilon$ obeys

$$\nabla^\dagger \nabla(*d*\upsilon) + m^2(*d*\upsilon) = *d*\eta \ .$$

(A3.39)

The desired bounds on the $L^2$ norms of the second derivatives of $(*d*\upsilon)$ can be obtained from (A3.39) by using the cut-off functions $\{\chi_N\}_{N=1,2,\ldots}$ as follows: Multiply both sides of (A3.39) by ever larger N versions of $\chi_N$; and then take the $L^2$ norm of the resulting identity. Two instances of integration by parts on the left hand side will equate the square of the $L^2$ norm of $\chi_N \nabla^\dagger \nabla(*d*\upsilon)$ with the square of the $L^2$ norm of $\chi_N \nabla \nabla(*d*\upsilon)$ plus terms that involve derivatives of $\chi_N$ but fewer derivatives of $*d*\upsilon$. Since $\nabla \upsilon$ and $\nabla \nabla \upsilon$ are square integrable these terms with derivatives of $\chi_N$ already have N-independent upper bounds. This being the case, then the various $N = 1, 2, \ldots$ versions of the preceding exercise leads to a bound on the $L^2$ norm of $\nabla \nabla(*d*\upsilon)$ by a $\eta$-independent multiple of the sum of the $L^2$ norms of $\eta$ and $*d*\eta$.

Step 5: This step and Step 6 analyze the cokernel of the homomorphism $\Phi$. Step 6 uses the results to prove the second bullet of Proposition A3.5.

To start this task, choose $s \geq 0$ and let $i_s : T \to S^3 - K$ denote the inclusion as the constant s slice of the $T \times [0, \infty)$ part of $S^3 - K$. As noted in Step 2, the space $H^1(S^3 - K; E)$ is 1-dimensional over $\mathbb{C}$; and because of Mostow rigidity, the homomorphism

$$i_s^*: H^1(S^3-K; E) \to H^1(T; E)$$

(A3.40)

is injective. Thus, the image of $i_s^*$ is a 1-dimensional subspace of $H^1(T;E)$. The latter space has dimension 2 over $\mathbb{C}$ (see, for example Lemma 2.2 in [Sc]). Writing the $s \geq 0$ part of $S^3-K$ as $T \times [0, \infty)$, it follows as a consequence that $H^1(T \times [0, \infty); E)$ has dimension 2 over $\mathbb{C}$. Moreover, a calculation finds that $H^1(T \times [0, \infty); E)$ is generated by two elements that can be written as follows: Use the oriented, orthonormal frame $\{e^1, e^2, e^3 = ds\}$ to identify E and $T^*(S^3-K)$ on the $T \times [0, \infty)$ part of $S^3-K$. This writes the bundle $E \otimes T^*(S^3-K)$ over the $s \geq 0$ part of $S^3-K$ as $\otimes^2 T^*(S^3-K)$ and it writes a section of the latter bundle over this same part of $S^3-K$ as a $3 \times 3$ matrix valued function. Using these identifications, two generators of $H^1(T \times [0, \infty); E)$ are

$$e^{2s} \begin{pmatrix} 1 & i & 0 \\ i & -1 & 0 \\ 0 & 0 & 0 \end{pmatrix} \quad \text{and} \quad \begin{pmatrix} 1 & -i & 0 \\ -i & -1 & 0 \\ 0 & 0 & 0 \end{pmatrix}.$$

(A3.41)

Note that it follows from the $n = 0$ version of Items $d_+$) and $d_-$) of Lemma A2.7 that these element are annihilated by $d^{\mathbb{C}}$. A pair $(\alpha_-, \alpha_+) \in \mathbb{C}^2 - 0$ are defined by writing the image of $H^1(S^3-K; E)$ in $H^1(T \times [0, \infty); E)$ using the basis in (A3.41) as the span of the element

$$\alpha_- e^{2s} \begin{pmatrix} 1 & i & 0 \\ i & -1 & 0 \\ 0 & 0 & 0 \end{pmatrix} + \alpha_+ \begin{pmatrix} 1 & -i & 0 \\ -i & -1 & 0 \\ 0 & 0 & 0 \end{pmatrix}.$$

(A3.42)

Step 6: It follows from (A3.42) that the solution in the $n = 0$ version of Lemma A3.3 with $t_a = t_{b+} = t_{b-} = 0$ and $t_{c+} = 4\alpha_-$ and $t_{c-} = 4\alpha_+$ extends over the whole of $S^3-K$ as a closed but not exact section of $E \otimes T^*(S^3-K)$. Denote this E-valued 1-form as $\mathfrak{T}_\lozenge$. This $\mathfrak{T}_\lozenge$ can also be depicted as a section of $\otimes^2 T^*(S^3-K)$; and in the latter guise, it may or may not be symmetric and traceless. If it is symmetric and traceless, then $\mathfrak{T}_\lozenge$ is the extension that proves the second assertion of Proposition A3.5. If it is not both traceless and symmetric, then it is necessary to find a section $\psi$ of E over $S^3-K$ so that $\mathfrak{T}_\lozenge + \nabla^{\mathbb{C}}\psi$ is symmetric and traceless as a section of $\otimes^2 T^*(S^3-K)$. This will be the case if $\psi$ obeys

- $\varepsilon^{abc} \mathfrak{T}_\lozenge{}^{bc} + \varepsilon^{abc} \nabla_b \psi^c + 2i \psi^a = 0$,
- $\mathfrak{T}_\lozenge{}^{cc} + \nabla_c \psi^c = 0$.

(A3.43)

Note that if $\psi$ obeys the top bullet's equation in (A3.43), then it also obeys the bottom bullet's equation because the identity $\varepsilon^{abc}\nabla_a \mathfrak{T}_\diamond{}^{bc} = 2i\mathfrak{T}_\diamond{}^{cc}$ follows by virtue of the fact that $d^{\mathbb{C}}\mathfrak{T}_\diamond = 0$. This understood, then the version of Lemma A3.6 with $\eta^a = -\varepsilon^{abc}\mathfrak{T}_\diamond{}^{bc}$ can be invoked to find the desired 1-form $\psi$. Since $\psi$ and $\nabla\psi$ are square integrable, it follows as a consequence of Proposition A3.1 and Lemma A3.3 that $\mathfrak{T} = \mathfrak{T}_\diamond + \nabla^{\mathbb{C}}\psi$ has the properties that are required of $\mathfrak{T}_+^{(c)}$ by the Proposition A3.5's second bullet. Meanwhile, the uniqueness of $\mathfrak{T}_+^{(c)}$ follows from Lemma A3.2 because the difference of 2 solutions (if not zero) will violate the bounds asserted by Lemma A3.2.

One more thing needs to be said to complete the proof of Proposition A3.5: The number $\alpha_-$ that appears in (A3.42) can not be zero because the event that $\alpha_- = 0$ would run afoul of Lemma A3.2.

*Part 5*: The upcoming Proposition A3.7 describes the $S^3$–K significance of the solutions from Lemma A3.3 when $n = -1$. What follows directly supplies some background for Proposition A3.7.

The $n = -1$ version of the top bullet of (A3.13) asserts that

$$\varepsilon^{bcd}\nabla_c \mathfrak{T}_+{}^{ad} = 0 .$$

(A3.44)

A solution to this equation is said to be a *Codazzi tensor*. Ferus observed [Fe] that the solutions to (A3.44) on a 3-manifold with constant sectional curvature -1 are locally determined by a function in the following sense: Let Y denote the given hyperbolic manifold. Fix a priori a locally finite, open cover of Y by balls that each sit well inside some Gaussian coordinate chart. Let $\mathfrak{U}$ denote this cover. Now suppose that $\mathfrak{T}_+$ is a symmetric section of $\otimes^2 T^*Y$ that obeys (A3.44). For each $B \in \mathfrak{U}$, there is a function on B, to be denoted by $f_B$, such that

$$\mathfrak{T}_+{}^{ab}|_B = \nabla_a \nabla_b f_B - \delta_{ab} f_B .$$

(A3.45)

Lafontaine [Laf] subsequently observed that the data $\{(B, f_B)\}_{B \in \mathfrak{U}}$ obeying (A3.45) for a given $\mathfrak{T}_+$ determines and is determined (up to a natural equivalence) by an element in the Čech cohomology group $H^1(Y; V)$ that is described in Section 5d for the case when $Y = S^3$–K. Since Lafontaine's paper can be difficult to obtain, an account of Lafontaine's observations are given in the subsequent two paragraphs.

Lafontaine's observed that two solutions to (A3.45) for a fixed $\mathfrak{T}_+$ differ by a function that obeys the equation

$$\nabla_a \nabla_b h - \delta_{ab} h = 0.$$

(A3.46)

Let $\mathcal{V}$ denote the presheaf that assigns to any given open set $U \subset Y$ the vector space $\mathcal{V}(U)$ of solutions to (A3.46) on U. If $B \subset Y$ is a sufficiently small radius ball, then (A3.46) has a 4-dimensional space of solutions on B and the choice of an orthonormal frame for T*Y at B's center point gives this space a canonical basis. To explain why this is, reintroduce $\mathbb{R}^{3,1}$ from Section 5d to denote the vector space $\mathbb{R}^4$ with the indefinite metric that is depicted in (5.10). Let $\mathbb{H}^+ \subset \mathbb{R}^{3,1}$ again denote the locus where the four Euclidean coordinates $(y_0, y_1, y_2, y_3)$ obey $y_0 > 0$ and $y_0^2 = 1 + y_1^2 + y_2^2 + y_3^2$. The Minkowski metric in (5.10) induces a metric on $\mathbb{H}^+$ that is positive definite with constant sectional curvature -1. The equation in (A3.46) is defined on $\mathbb{H}^+$ with this metric. As it turns out, its solution space on any domain in $\mathbb{H}^+$ is spanned by the 4 coordinate functions $\{y_0, y_1, y_2, y_3\}$. This is relevant to the task at hand because $\mathbb{H}^+$ (with its hyperbolic metric) is an isometric model for the universal covering space of Y. As a consequence, if $p \in Y$ is a given point and $B \subset Y$ is a small radius ball centered at p, then the choice of an oriented, orthonormal frame for $TY|_p$ canonically defines an isometric map from B to $\mathbb{H}^+$ taking p to the point $(1, 0, 0, 0)$. This fact explains why there is a 4-dimensional solution to (A3.46) over B and it identifies this solution space with the coordinate functions on $\mathbb{R}^{3,1}$.

With the preceding in mind, let $\mathfrak{T}_+$ denote as before a symmetric section of $\otimes^2 T^*Y$ obeying (A3.44). If $B \subset \mathfrak{U}$, let $f_B$ denote a function on B that obeys (A3.45). If B and B´ are intersecting balls from $\mathfrak{U}$, use $h_{BB'}$ to denote $f_{B'} - f_B$ on $B \cap B'$. This function obeys (A3.46) on $B \cap B'$ so it is in the vector space $\mathcal{V}(B \cap B')$. The collection $\{h_{BB'}\}_{B,B' \in \mathfrak{U}}$ obeys the Čech cocycle condition $h_{BB'} + h_{B'B''} + h_{B''B} = 0$ on the intersection of any three balls B, B´ and B´´ from $\mathfrak{U}$. Since $\mathfrak{T}_+|_B$ is of primary interest (not $f_B$), the only relevant data is the Čech cocycle $\{h_{BB'}\}_{B,B' \in \mathfrak{U}}$ modulo the equivalence relation that is obtained by adding coboundaries; this being the equivalence $\{h_{BB'}\}_{B,B' \in \mathfrak{U}} \sim \{h_{BB'} + h_B - h_{B'}\}_{B,B' \in \mathfrak{U}}$ with $\{h_B\}_{B \in \mathfrak{U}} \in \times_{B \in \mathfrak{U}} \mathcal{V}(B)$. This equivalence class of the data $\{h_{BB'}\}_{B,B'}$ is (by definition) a class in the 1-dimensional Čech cohomology on Y for the sheaf defined by $\mathcal{V}$, which is $H^1(Y; V)$. The class determined by $\mathfrak{T}_+$ in $H^1(Y; V)$ is said in what follows to be the *Ferus//Lafontaine class* of $\mathfrak{T}_+$.

Since compact manifolds have finite covers, the vector space $H^1(Y; V)$ is finite dimensional over $\mathbb{R}$ when Y is compact. It is also finite dimensional when $Y = S^3 - K$; this is by virtue of the fact that Y deformation retracts onto its $s \leq 1$ part.

Proposition A3.7 refers to a norm on $H^1(S^3 - K; V)$; and it denotes the norm of a given element $\mathfrak{h}$ by $|\mathfrak{h}|$. Since $H^1(S^3 - K; V)$ is finite dimensional, any two norms are commensurate. This being the case, any norm that is chosen a priori will serve for the purposes of the proposition.

**Proposition A3.7**: *The significance of the solutions that are described by the* $n = -1$ *version of Lemma A3.3 with regards to* $S^3 - K$ *is described in the two bullets below.*

- *There is a unique real valued function $f$ on $S^3-K$ obeying the following two items:*
  i) $\Delta f - 3f = 0$
  ii) *$f$ on the $s \geq 0$ part of $S^3-K$ obeys $f = e^{3s} + \mathfrak{r}_f$ with $\mathfrak{r}_f$ being square integrable.*
  *Let $\mathfrak{T}_+$ denote the tensor on $S^3-K$ with local frame components $\nabla_a\nabla_b f - f\delta_{ab}$. This is a traceless Codazzi tensor that can be written on the $s \geq 0$ part of $S^3-K$ as*

$$\mathfrak{T}_+ = 3e^{3s}\begin{pmatrix} -1 & 0 & 0 \\ 0 & -1 & 0 \\ 0 & 0 & 2 \end{pmatrix} + \mathfrak{r}$$

  *where $\mathfrak{r}$ obeys $|\mathfrak{r}| \leq \kappa\exp(-e^s/\kappa)$ with $\kappa \geq 1$ being independent of $s$.*

- *There exists $\mathfrak{t} \in \mathbb{C}-\{0\}$, a linear functional $L: H^1(S^3-K; V) \to \mathbb{R}$ and an injective, linear map from $H^1(S^3-K; V)$ to the space of ($\mathbb{R}$ valued) traceless, Codazzi tensors on $S^3-K$ obeying the following: Given $\mathfrak{h} \in H^1(S^3-K; V)$, let $\mathfrak{T}_C(\mathfrak{h})$ denote the corresponding traceless Codazzi tensor. This is the unique traceless Codazzi tensor on $S^3-K$ having Ferus/Lafontaine class $\mathfrak{h}$ and, on the $s \geq 0$ part of $S^3-K$, obeying*

$$\mathfrak{T}_C(\mathfrak{h}) = L(\mathfrak{h})(e^s\mathfrak{t}\begin{pmatrix} 1 & i & 0 \\ i & -1 & 0 \\ 0 & 0 & 0 \end{pmatrix} + e^s\overline{\mathfrak{t}}\begin{pmatrix} 1 & -i & 0 \\ -i & -1 & 0 \\ 0 & 0 & 0 \end{pmatrix}) + \mathfrak{r}_\mathfrak{h}$$

  *where $\mathfrak{r}_\mathfrak{h}$ is bounded. In fact, this $\mathfrak{r}_\mathfrak{h}$ obeys $|\mathfrak{r}_\mathfrak{h}| \leq \kappa\exp(-e^s/\kappa)|\mathfrak{h}|$ with $\kappa \geq 1$ being independent of $s$.*

**Proof of Proposition A3.7**: The proof that follows has six steps. (See [OS] for more analysis with regards to Codazzi tensors.)

    Step 1: Suppose that $\mathfrak{T}_+$ is a symmetric, traceless solution to (A3.44) on the $s \geq 0$ part $S^3-K$ with its Ferus/Lafontaine class in $H^1((0,\infty) \times T; V)$ equal to 0. This is to say that there is function $f$ on $S^3-K$ such that $\mathfrak{T}_+{}^{ab} = \nabla_a\nabla_b f - f\delta_{ab}$ on $[0, \infty) \times T$. Since $\mathfrak{T}_+$ is traceless, the function $f$ also obeys $\Delta f - 3f = 0$. This equation on the $s \geq 0$ part of $S^3-K$ is invariant with respect to the T action and so each $f$'s Fourier modes with respect to this T action is also a solution to the equation $\Delta f - 3f = 0$. Keeping this in mind, let $f_0$ denote the ($k_1 = 0, k_2 = 0$) Fourier mode. The equation $\Delta f - 3f = 0$ for $f_0$ is the ordinary differential equation

$$\tfrac{d^2}{ds^2}f_0 - 2\tfrac{d}{ds}f_0 - 3f_0 = 0$$

    (A3.47)

whose general solution is $f_0 = a_+ e^{3s} + a_- e^{-s}$. It follows as a consequence that the n = -1 solutions from Items c$_+$) and c$_-$) of the second bullet in Lemma A3.3 define non-zero Ferus/Lafontaine classes in $H^1((0, \infty); V)$. Meanwhile, the solution from Item a) of Lemma A3.3 can be written as $\nabla_a \nabla_b f_\Diamond - f_\Diamond \delta_{ab}$ with $f_\Diamond = \frac{1}{3} t_a e^{3s}$. Note that the T independent function on $(0, \infty) \times T$ given by the rule $s \to e^{-s}$ obeys the equation in (A3.46); and as a consequence, adding any multiple of $e^{-s}$ to $\frac{1}{3} t_a e^{3s}$ will not change the corresponding Codazzi tensor.

Step 2: This step proves the assertion in the first bullet of Proposition A3.7. The proof finds a unique function (to be denoted by $f$) on $S^3 - K$ that obeys Items a) and b) of the proposition's first bullet. Granted such a function, $f$, define $\mathfrak{T}_+$ in any given orthonormal frame for $T^*(S^3 - K)$ by the rule $\mathfrak{T}_+{}^{ab} = \nabla_a \nabla_b f - f \delta_{ab}$. This is a symmetric, traceless section of $\otimes^2 T^*(S^3 - K)$ that obeys the equation in (A3.44) because $f$ obeys $\Delta f - 3f = 0$. Let $\mathfrak{T}$ for the moment denote the solution to the equation in (A3.44) on the $s \geq 0$ part of $S^3 - K$ given by Item a) of Lemma A3.3. Then $\mathfrak{T}_+ - \mathfrak{T}$ is square integrable on the $s \geq 0$ part of $S^3 - K$. It then follows from what is said in Proposition A3.1 and Lemma A3.3 that $\mathfrak{T}_+$ obeys the assertions in Proposition A3.7's first bullet.

To find a function $f$ that obeys Items a) and b) of the Proposition A3.7's first bullet, let $f_\Diamond$ denote the function $s \to f_\Diamond(s) = \frac{1}{3} t_a e^{3s}$ on $[0, \infty) \times T$. Reintroduce the cut-off function β from Step 1 of the proof of Proposition A3.4. (This function on $\mathbb{R}$ is non-decreasing, zero where $s \leq 0$ and 1 where $s \geq 0$.) The function $\beta f_\Diamond$ extends the definition of $f_\Diamond$ to the whole of $S^3 - K$, but this extension is not in the kernel of $\Delta - 3$. As explained in the next paragraph, there exists a unique $L^2$ function on $S^3 - K$ (to be denoted by û) that obeys $\Delta(\beta f_\Diamond + \hat{u}) - 3(\beta f_\Diamond + \hat{u}) = 0$. Moreover, this function û has $L^2$ first and second derivatives. It follows from this that $f = \beta f_\Diamond + \hat{u}$ obeys the conditions set forth by the first bullet of Proposition A3.7.

To find û, let $w$ denote for the moment any given smooth, $L^2$ function on $S^3 - K$ so as to consider solving for an $L^2$ function û obeying $-\Delta \hat{u} + 3\hat{u} = w$. In the case at hand $w$ is the function $\Delta \beta f_\Diamond + 2 \langle d\beta, df_\Diamond \rangle$ with $\langle , \rangle$.. Introduce by way of notation $\mathfrak{C}_0$ to denote the space of $L^2_1$ functions on $S^3 - K$. Define the function $\mathfrak{E}_w$ on $\mathfrak{C}_0$ by the rule

$$u \to \mathfrak{E}_w(u) = \frac{1}{2} \int_{S^3 - K} (|\nabla u|^2 + 3|u|^2) - \int_{S^3 - K} uw .$$

(A3.48)

This function $\mathfrak{E}_w$ has a unique minimizer in $\mathfrak{C}_0$. The argument for this differs only in notation from arguments in the proof of Lemma A3.6 and also from arguments in Step 3 of the proof of Proposition A3.4 for the existence of a minimizer of the functional that is depicted in (A3.25). This function, û, obeys the desired equation $-\Delta \hat{u} + 3\hat{u} = w$. Being that $\hat{u} \in \mathfrak{C}_0$, it is square integrable with square integrable first derivatives, an argument using the cut-off functions $\{\chi_N\}_{N=1,2,\ldots}$ from Step 4 of the proof of Proposition A3.4 can be used to prove that $\nabla\nabla\hat{u}$ is also an $L^2$ function. The details are straightforward and thus omitted.

<u>Step 3</u>: This step and Steps 4-6 prove the second bullet of Proposition A3.7. To start, write the $s > 0$ part of $S^3 - K$ as $(0, \infty) \times T$ as done before. Given $\mathfrak{t} \in \mathbb{C}$, define an $\mathbb{R}$ valued, traceless and symmetric section, $\mathfrak{T}_\diamond$, of $\otimes^2 T^*((0, \infty) \times T)$ using the formula

$$\mathfrak{T}_\diamond = e^s\, \mathfrak{t} \begin{pmatrix} 1 & i & 0 \\ i & -1 & 0 \\ 0 & 0 & 0 \end{pmatrix} + e^s\, \overline{\mathfrak{t}} \begin{pmatrix} 1 & -i & 0 \\ -i & -1 & 0 \\ 0 & 0 & 0 \end{pmatrix}.$$

(A3.49)

This is a Codazzi tensor since it comes from the $n = -1$ version of Items $c_+$) and $c_-$) in Lemma A3.3. Since $H^1((0, \infty) \times T; V)$ is 2-dimensional over $\mathbb{R}$, (see Lemma 2.2 in [Sc]), it follows from what is said in Step 1 that the Ferus/Lafontaine classes of the various $\mathfrak{t} \in \mathbb{C}$ versions of (A3.49) generate $H^1((0, \infty) \times T; V)$. By Proposition 2.3 in [Sc] and Lemma 4.1 in [Sc], the image of $H^1(S^3 - K; V)$ in $H^1((0, \infty) \times T; V)$ via the restriction homomorphism is 1 dimensional (over $\mathbb{R}$). It follows as a consequence that the image is the $\mathbb{R}$-linear span of the Ferus/Lafontaine classes of a fixed $\mathfrak{t} \in \mathbb{C} - 0$ version of (A3.49). Assume henceforth that $\mathfrak{T}_\diamond$ is defined using this particular value of $\mathfrak{t}$. Said differently, there exists a linear functional on $H^1(S^3 - K; V)$, this being the functional $L$, such that if $\mathfrak{h} \in H^1(S^3 - K; V)$, then its restriction to $H^1((0, \infty) \times T; V)$ is the Ferus/Lafontaine class of $L(\mathfrak{h})\mathfrak{T}_\diamond$.

<u>Step 4</u>: As argued directly, every class in $H^1(S^3 - K; V)$ is the Ferus/Lafontaine class of a Codazzi tensor. (The argument that follows works on any locally compact hyperbolic 3-manifold.) To see this, let $\mathfrak{h}$ denote a given class. Fix a locally finite cover of $S^3 - K$ by small radius balls and denote this cover by $\mathfrak{U}$. The cover should be chosen so that $\mathfrak{h}$ is represented by the cocycle data $\{h_{BB'}\}_{B,B' \in \mathfrak{U}}$ with each $h_{BB'}$ defined on the corresponding intersection $B \cap B'$ where it obeys (A2.61). Let $\{\varphi_B\}_{B \in \mathfrak{U}}$ denote a partition of unity that is subordinate to the cover $\mathfrak{U}$. For $B \in \mathfrak{U}$, define $f_B = \sum_{B'} h_{BB'} \varphi_{B'}$ and define the symmetric section of $\otimes^2 T^*(S^3 - K)$ on $B$ (to be denoted by $\mathfrak{T}_\mathfrak{h}|_B$) by taking its components with respect to a local orthonormal frame to be $(\mathfrak{T}_\mathfrak{h}|_B)_{ab} = \nabla_a\nabla_b f_B - f_B \delta_{ab}$. If $B'$ is another ball from $\mathfrak{U}$ that intersects $B$, then $f_B - f_{B'} = h_{BB'}$ on $B \cap B'$; and since $h_{BB'}$ obeys (A2.61), it follows that $\mathfrak{T}_\mathfrak{h}|_B = \mathfrak{T}_\mathfrak{h}|_{B'}$ on $B \cap B'$. This implies that the various $B \in \mathfrak{U}$ versions of $\mathfrak{T}_\mathfrak{h}|_B$ agree on intersections and so define a Codazzi tensor on all of $S^3 - K$.

Suppose that $\mathfrak{h} \in H^1(S^3 - K; V)$ and that $\mathfrak{T}_\mathfrak{h}$ is a Codazzi tensor with Ferus/Lafontaine class equal to $\mathfrak{h}$. It follows from what is said in Step 3 that $\mathfrak{T}_\mathfrak{h}$ differs from $L(\mathfrak{h})\mathfrak{T}_\diamond$ on $(0,\infty)\times T$ by a Codazzi tensor with components $\nabla_a\nabla_b f - f\delta_{ab}$ for some function $f$ on $(0,\infty)\times T$. Let $\beta$ denote the cut-off function from Step 1 of the proof of Lemma A2.8. Define $\mathfrak{T}_{\mathfrak{h}\diamond}$ to be $\mathfrak{T}_\mathfrak{h}$ on the $s \le 0$ part of $S^3-K$ and define it on the $s \ge 0$ part by taking its components to be $(\mathfrak{T}_{\mathfrak{h}\diamond})_{ab} = (\mathfrak{T}_\mathfrak{h})_{ab} - (\nabla_a\nabla_b(\beta f) - \beta f \delta_{ab})$. This new Codazzi tensor has the same Ferus/Lafontaine class $\mathfrak{h}$ and it is equal to $L(\mathfrak{h})\mathfrak{T}_\diamond$ on the $s \ge 1$ part of $S^3-K$.

Step 5: Were the Codazzi tensor $\mathfrak{T}_{\mathfrak{h}\diamond}$ traceless, then it would be perfect with regards to the second bullet of Proposition A3.7. Since it is not necessarily traceless, the task at hand is to find a function $\hat{u}$ on $S^3-K$ with suitable behavior where $s \ge 0$ such that the Codazzi tensor with the components in a local frame given by

$$(\mathfrak{T}_C)_{ab} = (\mathfrak{T}_{\mathfrak{h}\diamond})_{ab} + (\nabla_a\nabla_b\hat{u} - \hat{u}\delta_{ab})$$
(A3.50)

has zero trace everywhere. (No matter what $\hat{u}$ is, the tensor on the right hand side of (A3.50) has Ferus/Lafontaine class $\mathfrak{h}$.) The right hand side of (A3.50) is traceless when $\hat{u}$ obeys

$$-\Delta\hat{u} + 3\hat{u} = \text{trace}(\mathfrak{T}_{\mathfrak{h}\diamond}) \ .$$
(A3.51)

Since $\mathfrak{T}_{\mathfrak{h}\diamond} = L(\mathfrak{h})T_\diamond$ on the $s \ge 1$ part of $S^3-K$, its trace is an $L^2$ function. This understood, then what is said in the final paragraph of Step 2 can be invoked using $w = \text{trace}(\mathfrak{T}_{\mathfrak{h}\diamond})$ to find a unique $L^2$ function $\hat{u}$ that obeys (A3.51). This function is smooth and has $L^2$ first and second derivatives. It follows in particular from the latter observation that $\mathfrak{T}_+$ differs from $\mathfrak{T}_\diamond$ on the $s \ge 0$ part of $S^3-K$ by a square integrable Codazzi tensor.

Step 6: The assignment of $\mathfrak{h}$ to the Codazzi tensor $\mathfrak{T}_C$ from Step 5 will serve for the second bullet Proposition A3.7 if the norm of $\mathfrak{T}_C - L(\mathfrak{h})\mathfrak{T}_{+\diamond}$ on the $s \ge c_*$ part of $S^3-K$ is a priori bounded by $c_*\exp(-e^{-s}/c_*)|\mathfrak{h}|$.

To see such a bound, let N denote for the moment the dimension of $H^1(S^3-K;V)$. Fix once and for all a basis $\{\mathfrak{h}_i\}_{i=1,\ldots,N}$ for this vector space. (In particular, this basis should be fixed a priori with no reference to $\mathfrak{h}$.) It proves convenient to use the chosen basis to define a norm for $H^1(S^3-K;V)$ as follows: Write any given element $\mathfrak{h}_*$ as $\sum_{i=1,\ldots,N} a_i\mathfrak{h}_i$ with $\{a_i\}_{i=1,\ldots,N} \subset \mathbb{R}$. The norm of $\mathfrak{h}_*$ is defined to be $\sum_{i=1,\ldots,N}|a_i|$. Since any two norms on $H^1(S^3-K;V)$ are commensurate, this norm will be used for $|\cdot|$ in what follows.

For each $i \in \{1, \ldots, N\}$, let $\mathfrak{T}_{Ci}$ denote the $\mathfrak{h}_i$ version of what is denoted by $\mathfrak{T}_C$ in Step 5. Proposition A3.1 implies that $|\mathfrak{T}_{Ci} - L(\mathfrak{h}_i)\mathfrak{T}_\diamond|$ is bounded on the $s \ge c_*$ part of $S^3-K$ by $c_i\exp(-e^{-s}/c_*)$ with $c_i$ being independent of $s$.

With the preceding understood, write the given element $\mathfrak{h}$ as $\mathfrak{h} = \sum_{i=1,\ldots,N} a_i \mathfrak{h}_i$ with $\{a_i\}_{i \in \{1,\ldots,N\}} \subset \mathbb{R}$. Since the corresponding $\mathfrak{T}_C$ as depicted in (A3.50) can be written as $\mathfrak{T}_+ = \sum_{i=1,\ldots,N} a_i \mathfrak{T}_{+i}$ and since $L(\mathfrak{h}) = \sum_{i=1,\ldots,N} a_i L(\mathfrak{h}_i)$, it follows that $|\mathfrak{T}_C - L(\mathfrak{h})\mathfrak{T}_0|$ on the $s \geq c_*$ part of $S^3-K$ is bounded by $(\sum_{i=1,\ldots,N} |a_i|) c_\diamond \exp(-e^{-s}/c_*)$ with $c_\diamond = \sup_{i \in \{1,\ldots,N\}} |c_i|$. This is the required norm bound.

### f) Solutions $(\mathcal{X}, \mathfrak{s})$ to (A3.12) with $\mathfrak{j} \neq 0$

This subsection considers the $\mathfrak{j} \neq 0$ solutions to (A3.11) and (A3.12). The first four parts of the subsection talk abut the autonomous equations for $\mathfrak{j}$ given in (A3.11). Supposing that $\mathfrak{j}$ obeys (A3.11), then the last part of the subsection states and then proves a lemma that says when a pair $(\mathcal{X}, \mathfrak{s})$ obeys the corresponding version of (A3.12).

*Part 1*: The two lemmas in this part of the subsection concern solutions to (A3.11) on the $s \geq 0$ part of $S^3-K$. Keep in mind that these equations are invariant with respect to the T action on this part of $S^3-K$ that acts by translation on each level set of $s$. As a consequence, each Fourier mode with respect to this T action of any given solution to (A3.11) also obeys (A3.11). The lemma that follows is the $\mathfrak{j}$ analog of Proposition A3.1

**Lemma A3.8**: *There exists $\kappa > 10$ with the following significance: Fix an integer to be denoted by $\mathrm{n}$ and fix $r > \kappa^2$.*
- *Suppose that $\mathfrak{j}$ obeys (A3.11) on the $s \in [0, r]$ part of $S^3-K$. Assume that its $L^2$ norm is equal to 1 and that $\mathfrak{j}$ lacks the $(k_1 = 0, k_2 = 0)$ Fourier mode for the T action. Then $|\mathfrak{j}| + |\nabla \mathfrak{j}| \leq \kappa (\exp(-e^s/\kappa) + \exp(-e^{-(r-s)}/\kappa))$ at the points where $s \in [\kappa, r - \kappa]$.*
- *Suppose that $\mathfrak{j}$ obeys (A3.11) on the $s \in [0, \infty)$ part of $S^3-K$. Assume that its $L^2$ norm is equal to 1 and that $\mathfrak{j}$ lacks the $(k_1 = 0, k_2 = 0)$ Fourier mode for the T action. Then $|\mathfrak{j}| + |\nabla \mathfrak{j}| \leq \kappa \exp(-e^s/\kappa)$ at the points where $s \in [\kappa, \infty)$.*

This lemma is proved momentarily. The next lemma is the $\mathfrak{j}$ analog of Lemma A3.2.

**Lemma A3.9**: *Fix a non-zero integer to be denoted by $\mathrm{n}$. Suppose that $r \geq 1$ and that $\mathfrak{j}$ obeys (A3.11) on the $s \leq r$ part of $S^3-K$. If $\mathrm{r} \in [1, r]$, then $\int_{s \leq \mathrm{r}} |\mathfrak{j}|^2 \geq e^{2|\mathrm{n}|(\mathrm{r}-r)} \int_{s \leq r} |\mathfrak{j}|^2$.*

The proof of this last lemma differs only cosmetically from the proof of Lemma A3.2 so it is left to the reader.

*Proof of Lemma A3.8*: The 1-form j obeys the $m = n$ and $\eta = 0$ version of (A3.37) on its domain of definition. Taking the inner product of both sides of this version of (A3.37) with j leads to the following identity

$$\tfrac{1}{2} d^\dagger d |j|^2 + |\nabla j|^2 + (n^2 - 2)|j|^2 = 0 .$$

(A3.52)

Let $\mathfrak{f}$ denote the function on $[0, r]$ whose value at any given $s \in [0, r]$ is the square of the $L^2$ norm of j on $\{s\} \times T$. Integrate (A3.52) on this constant s slice $\{s\} \times T$ to obtain

$$-\tfrac{1}{2} \tfrac{d}{ds}(e^{-2s} \tfrac{d}{ds}(e^{2s} \mathfrak{f})) + \int_{\{s\} \times T} |\nabla j|^2 + (n^2 - 2)\mathfrak{f} = 0.$$

(A3.53)

If j lacks the ($k_1 = 0, k_2 = 0$) Fourier mode, then the integral of $|\nabla j|^2$ that appears in this last equation is no smaller than $c_0^{-1} e^{2s}$ times the integral of $|j|^2$ on $\{s\} \times T$ with $c_0$ being greater than 1 and independent of s and j. Therefore, if j lacks this mode, then (A3.53) leads to the inequality

$$-\tfrac{1}{2} \tfrac{d}{ds}(e^{-2s} \tfrac{d}{ds}(e^{2s} \mathfrak{f})) + c_1(e^{2s} + (n^2 - 2))\mathfrak{f} \leq 0$$

(A3.54)

with $c_1$ being greater than 1 and independent of s and j. This last inequality leads directly to the bounds for |j| in the first bullet of the lemma. The asserted bounds for $|\nabla j|$ follow from the bounds for |j| using standard elliptic regularity results.

The bounds in the second bullet follow by invoking the lemma's first bullet using an unbounded sequence of values for the parameter $r$.

*Part 2*: The lemma that follows directly discusses the T-invariant solutions to (A3.11) on the $s \geq 0$ part of $S^3 - K$. This lemma is the j analog of Lemma A3.3. This lemma writes the components of j using the frame for $T^*(S^3 - K)$ on the $s \geq 0$ part of $S^3 - K$ that has the third frame 1-form $e^3$ being ds and the other frame 1-forms $\{e^1, e^2\}$ having the form $e^1 = e^{-s} \hat{e}^1$ and $e^2 = e^{-s} \hat{e}^2$ with $\{\hat{e}^1, \hat{e}^2\}$ being a T and s independent, oriented orthonormal frame for the flat metric $\mathfrak{m}$ on T that appears in (3.3).

**Lemma A3.10**: *Supposing that j is a T-invariant solution to (A3.11) on the $s \geq 0$ part of $S^3 - K$, then j has the following form*:

- *If $n \neq 0$, then* $j = a_+ e^{(1+n)s} \begin{pmatrix} 1 \\ i \\ 0 \end{pmatrix} + a_- e^{(1-n)s} \begin{pmatrix} 1 \\ -i \\ 0 \end{pmatrix}$ *with* $a_+, a_- \in \mathbb{C}$.

- *If* $n = 0$, *then* $j = a_1 e^s \begin{pmatrix} 1 \\ 0 \\ 0 \end{pmatrix} + a_2 e^s \begin{pmatrix} 0 \\ 1 \\ 0 \end{pmatrix} + a_3 e^{2s} \begin{pmatrix} 0 \\ 0 \\ 1 \end{pmatrix}$ *with* $a_1, a_2, a_3 \in \mathbb{R}$.

*Proof of Lemma A3.10*: This lemma is proved by setting the directional derivatives (but not covariant derivatives) along the T-directions of the components of j equal to zero in (A3.11) and then solving the resulting ordinary differential equations. This is a straightforward task that is left to the reader.

*Part 3*: The next lemma says which solutions from the top bullet of Lemma A3.10 look like the large s part of some solution to (A3.11) on the whole of $S^3 - K$.

**Lemma A3.11**: *There exists* $\kappa \geq 1$ *with the following significance: Fix* $n > 0$.
- *There exists a unique solution to (A3.11) on* $S^3 - K$ *(it is denoted by* j*) that can be written on the* $s \geq 0$ *part of* $S^3 - K$ *as*

$$j = e^{(1+n)s} \begin{pmatrix} 1 \\ i \\ 0 \end{pmatrix} + \hat{a}_n e^{(1-n)s} \begin{pmatrix} 1 \\ -i \\ 0 \end{pmatrix} + \mathfrak{r} \;,$$

*where* $\hat{a}_n \in \mathbb{C}$ *and where the function* $\mathfrak{r}$ *is bounded and has zero* $(k_1 = 0, k_2 = 0)$ *Fourier mode for the* T *action on the* $s \geq 0$ *part of* $S^3$-K.
- *The number* $\hat{a}_n$ *that appears in preceding formula obeys* $|\hat{a}_n| \leq \kappa n e^{ns}$ *and the function* $\mathfrak{r}$ *that appears in the preceding formula obeys* $|\mathfrak{r}| + |\nabla \mathfrak{r}| \leq \kappa \exp(-e^s/\kappa)$

*Proof of Lemma A3.11*: If j exists, then it is unique. This follows from Lemma A3.9 because the difference between any two solutions would be square integrable on $S^3 - K$. To see about the existence of j, reintroduce the bump function $\beta$, this being a function with support on the $s \geq 0$ part of $S^3 - K$ that equals 1 on the $s \geq 1$ part of $S^3 - K$. Let $j_n$ denote the $a_+ = 1$ and $a_- = 0$ version of what is written in n's version of the top bullet of Lemma A3.10. Lemma A3.6 supplies a smooth, square integrable 1-form on $S^3 - K$ (to be denoted by $\psi$) that obeys the equation

$$*d\psi + in\psi = -*(d\beta \wedge j_n) \;.$$
(A3.55)

Granted (A3.55), then $j = \beta j_n + \psi$ obeys (A3.11) on $S^3-K$ and, by virtue of Lemma A3.8, it has the required asymptotic behavior on the large s part of $S^3-K$ to serve for Lemma A3.10.

*Part 4*: This part of the subsection talks about the significance of the solutions from the second bullet of Lemma A3.10 with regards to $n = 0$ solutions to (A3.11) on the whole of $S^3-K$. The important point in this regard is that any solution to the $n = 0$ version of (A3.11) is a harmonic 1-form on its domain of definition. In particular, it is a closed 1-form and thus defines a class in $H^1(S^3-K; \mathbb{R})$. This $\mathbb{R}$-module is 1-dimensional. Moreover, the homomorphism $\hat{\imath}^*: H^1(S^3-K; \mathbb{R}) \to H^1((0,\infty) \times T; \mathbb{R})$ that is induced by the inclusion map is injective. Meanwhile, $H^1((0,\infty) \times T; \mathbb{R})$ is the 2-dimensional vector space over $\mathbb{R}$ that is represented by the 1-forms in the second bullet of Lemma A3.10 that have $a_3 = 0$. Thus, the $\hat{\imath}^*$ image of $H^1(S^3-K; \mathbb{R})$ in $H^1(0,\infty) \times T; \mathbb{R})$ is represented by the $\mathbb{R}$ span of a 1-form on $(0, \infty)$ that can be written as

$$\hat{j} = e^s (\mathfrak{z} \begin{pmatrix} 1 \\ i \\ 0 \end{pmatrix} + \bar{\mathfrak{z}} \begin{pmatrix} 1 \\ -i \\ 0 \end{pmatrix})$$

(A3.56)

with $\mathfrak{z}$ being a unit length vector in $\mathbb{C}$.

**Lemma A3.12**: *Let $j_\diamond$ denote one of the 1-forms that are described by the second bullet of Lemma A3.10. Let $j$ denote a harmonic 1-form on $S^3-K$ that can be written on the $s > 0$ part of $S^3-K$ as $j = j_\diamond + \mathfrak{r}$ with $\mathfrak{r}$ being square integrable. Then $j_\diamond$ is a real multiple of $\hat{j}$. Moreover, in the case when $j_\diamond = \hat{j}$, there exists a unique harmonic 1-form on $S^3-K$ that can be written on the $s > 0$ part of $S^3-K$ as $\hat{j} + \mathfrak{r}$ with $\mathfrak{r}$ being square integrable; and $\mathfrak{r}$ obeys $|\mathfrak{r}| + |\nabla \mathfrak{r}| \le \kappa \exp(-e^s/\kappa)$ with $\kappa$ being independent of s.*

*Proof of Lemma A3.12*: The proof of this lemma has three steps.

<u>Step 1</u>: It follows from what was said initially about $H^1(S^3-K; \mathbb{R})$ that the only versions of $j_\diamond$ that extend as closed 1-forms are linear combinations of what is written in (A3.56) and the $a_1 = a_2 = 0$ versions of what is written in the second bullet of Lemma A3.10. Consider first the question of finding $j$ obeying (A3.11) on $S^3-K$ with $j$ on the $s \ge 0$ part $S^3-K$ differing from $j_0$ by an $L^2$ differential form. To this end, extend $\hat{j}$ from the $s \ge 0$ part of $S^3-K$ to the whole of $S^3-K$ as a closed 1-form to be denoted by $j_*$.

Supposing that $\hat{u}$ is an $L^2_1$ function on $S^3-K$, then $j = j_* + d\hat{u}$ is the desired harmonic extension of $j_0$ if $\hat{u}$ obeys the equation

$$d^\dagger d\hat{u} + d^\dagger j_* = 0 .$$

(A3.57)

Meanwhile, a smooth $L^2_1$ solution to (A3.57) can be found that minimizes the functional

$$u \to \mathfrak{E}(u) = \int_{S^3-K} (\tfrac{1}{2}|du|^2 + \langle u, d^\dagger j_* \rangle)$$

(A3.58)

on the space of smooth, $L^2_1$ functions on $S^3-K$ obeying $\int_{S^3-K} u = 0$. This is because critical points of (A3.58) on this function space are solutions to (A3.57) and vice-versa. With regards to the choice of function space, note that the integral of $d^\dagger j_*$ on $S^3-K$ is equal to zero so it is $L^2$ orthogonal to the constant functions. To see this, note first that the domain of integration can be restricted to the $s \leq 2$ part of $S^3-K$; and then integration by parts can be used to equate the integral of $d^\dagger j_*$ on this part of $S^3-K$ with the integral of the ds component of $j_*$ on the $s=2$ torus. Since $j_* = \hat{j}$ on this torus, and $\hat{j}$ has no ds component, the latter integral is zero.

Step 2: A minimizer of $\mathfrak{E}$ by invoking the fact that it obeys the coercive bound

$$\mathfrak{E}(u) \geq c_0^{-1} \int_{S^3-K} (|du|^2 + |u|^2) - c_0$$

(A3.59)

when u is an $L^2_1$ function that is $L^2$ orthogonal to 1. Here, $c_0 \geq 1$ is independent of u. This bound follows from the fact that the 0 is an isolated point in the $L^2$ spectrum of the Laplace operator on a finite volume hyperbolic 3-manifold with finitely many cusp like ends. (See, e.g. Chapter 4 in [EGM].) Given (A3.59), then standard variational arguments and elliptic regularity arguments can be used to prove that $\mathfrak{E}$ has a unique $L^2_1$ minimizer (which will be $\hat{u}$); and that this minimizer is smooth and that it obeys (A3.57).

Step 3: The $a_1 = a_2 = 0$ version of the 1-form in the second bullet of Lemma A3.10 is exact, it being a constant multiple of $d(e^{2s})$. This being the case, then Lemma A3.12 follows from what is said in Steps 1 and 2 if there are no harmonic functions on $S^3-K$ that differ on the $s \geq 0$ part of $S^3-K$ from the function $e^{2s}$ by a function that is $o(e^{2s})$ at large s. Suppose to the contrary, that there is such a function so as to derive nonsense. Call this function $f$. Supposing that $r \geq 1$, integrate the equation $d^\dagger df = 0$ on the $s \leq r$ part of $S^3-K$. Integrate by parts to derive the identity:

$$\int_{s=r} \frac{\partial}{\partial s} f = 0 .$$

(A3.60)

Now write $f = e^{2s} + \mathfrak{r}$ on the $s \geq 0$ part of $S^3-K$ with $\mathfrak{r}$ being harmonic. If $|\mathfrak{r}| = o(e^{2s})$ at large s, then it follows using the Fourier mode decomposition for the T action that $|\mathfrak{r}| = \mathcal{O}(1)$ and that $|\frac{\partial}{\partial s}\mathfrak{r}| = o(1)$. These bounds for $\mathfrak{r}$ are not compatible with (A3.60) since the integral of $\frac{\partial}{\partial s} f$ on the $s = r$ torus is positive and independent of s.

*Part 5*: The lemmas in the previous parts of this subsection describe the relevant solutions to (A3.11). The proposition that follows talks about the corresponding versions of (A3.12).

**Proposition A3.13**: *Supposing that* $n \neq \pm 2$, *let* $j_n$ *denote a solution to* n's *version of (A3.11) on the $s < r$ part of $S^3-K$ for a given positive number r. Define a pair of symmetric traceless sections of $\otimes^2 T^*(S^3-K)$ by the rule whereby their coefficients with respect to any local orthonormal frame for $T^*(S^3-K)$ are*

$$\mathcal{X}^{ab} = \frac{1}{4-n^2} (\nabla_a j^b + \nabla_b j^a) \quad and \quad \mathfrak{s}^{ab} = \frac{i}{2} \frac{n}{4-n^2} (\nabla_a j^b + \nabla_b j^a) .$$

*This pair with* $\mathfrak{j}$ *obeys the integer* n *version of the equations in (A3.12) on the $s < r$ part of $S^3-K$. If* $n = \pm 2$, *then there is no pair* $(\mathcal{X}, \mathfrak{s})$ *that obeys (A3.12) on the $s < r$ part of $S^3-K$ when there is a non-zero coefficient* $a_+$ *in Lemma A3.10's depiction of the T-invariant part of* $\mathfrak{j}$ *on the $s \in [0, r)$ subset in $S^3-K$. In particular, there is no $(\mathcal{X},\mathfrak{s})$ that obeys (A3.12) on $S^3-K$ when $\mathfrak{j}$ is described by the $n = \pm 2$ version of Lemma A3.11.*

*Proof of Proposition A3.13*: It is left as an exercise for the reader to verify that $(\mathcal{X}, \mathfrak{s})$ with the specified $\mathfrak{j}$ obey the integer n version of (A3.12) if $n \neq \pm 2$. The subsequent paragraphs prove the $n = 2$ assertion of the proposition. The case $n = -2$ is obtained from the $n = 2$ case by complex conjugation.

Supposing now that $n = 2$, let $\mathfrak{Q}$ denote $\mathcal{X} - i\mathfrak{s}$. By virtue of (A3.12), this symmetric, traceless section of $\otimes^2 T^*(S^3-K)$ obeys equations that can be written locally using an oriented orthonormal frame for $T^*(S^3-K)$ as

- $\nabla_a \mathfrak{Q}^{ab} = -2 j^b$ .
- $\frac{1}{2}(\varepsilon^{bcd} \nabla_c \mathfrak{Q}^{ad} + \varepsilon^{acd} \nabla_c \mathfrak{Q}^{bd}) + i\mathfrak{Q}^{ab} = -\frac{1}{4} i(\nabla_a j^b + \nabla_b j^a) .$

(A3.61)

These equations are the respective anti-symmetric and symmetric parts of the equation

$$\varepsilon^{acd} \nabla_c \mathfrak{Q}^{bd} + i \mathfrak{Q}^{ab} = -\tfrac{3i}{4} \nabla_a j^b + \tfrac{i}{4} \nabla_b j^a .$$

(A3.62)

As was done in the proof of Proposition A3.5, let $E \to S^3 - K$ denote the complexification of $T^*(S^3 - K)$ and let $\nabla^{\mathbb{C}}$ denote the covariant derivative of the flat $SO(3;\mathbb{C})$ connection on E that is given by the $n = 0$ version of (A3.29). Reintroduce the exterior covariant derivative on E-valued 1-forms from (A3.30). With $\mathfrak{Q}$ viewed as an E-valued 1-form on $S^3 - K$, then the equation in (A3.62) can be written schematically as

$$d^{\mathbb{C}} \mathfrak{Q} = \mathfrak{J} ,$$

(A3.63)

with $\mathfrak{J}$ denoting the E-valued 2-form whose Hodge dual is defined in a local oriented orthonormal frame for $T^*(S^3 - K)$ (and hence E also) by the right hand side of (A3.62). (So that there is no misunderstanding: The index a in (A3.62) labels 1-form components and the index b labels E component.)

The equation in (A3.63) says, first of all, that $d^{\mathbb{C}} \mathfrak{J} = 0$; and granted that this is true, it says that the $d^{\mathbb{C}}$-cohomology class of $\mathfrak{J}$ is zero. The assertion that $d^{\mathbb{C}} \mathfrak{J} = 0$ can be checked directly from the formula for $*\mathfrak{J}$ in (A3.62). This is left as an exercise for the reader (but note that the $n = 2$ version of (A3.11) must be invoked.) Granted that $\mathfrak{J}$ is a $d^{\mathbb{C}}$- closed 2-form, then $\mathfrak{J}$ is not $d^{\mathbb{C}}$ exact on $S^3 - K$ if its pull-back to any constant s torus in the $s \in [0, r]$ part of $S^3 - K$ is not $d^{\mathbb{C}}$ exact. The pull-back of (A3.62) to such a torus is the index $a = 3$ version of (A3.62) if ds is the third basis vector for the orthonormal frame for $T^*(S^3 - K)$ on the $s \geq 0$ part of $S^3 - K$.

With the preceding in mind, write j on the $s \in [0, r)$ part of $S^3 - K$ as $j = j_+ + j_- + j_\perp$ with the notation as follows: The term $j_\perp$ denote the sum of the non-zero Fourier modes of j with respect to the action of the group T; and $j_+$ and $j_-$ are the respective $a_+$ and $a_-$ terms in Lemma A3.10's depiction of the T-invariant part of j. There is a corresponding decomposition of what is denoted by $\mathfrak{J}$ in (A3.63) as $\mathfrak{J} = \mathfrak{J}_+ + \mathfrak{J}_- + \mathfrak{J}_\perp$ with each defined by using the corresponding $j_+$, $j_-$ or $j_\perp$ in lieu of j on the right hand side of (A3.62). Each of $\mathfrak{J}_+, \mathfrak{J}_-$ and $\mathfrak{J}_\perp$ is $d^{\mathbb{C}}$ closed. The pull-back of the term $\mathfrak{J}_\perp$ to any constant s torus can be written as $d^{\mathbb{C}} \mathfrak{Q}_\perp$ by using its decomposition as a sum of non-constant T action Fourier modes. This task is left to the reader. The pull-back of the $\mathfrak{J}_-$ term to any constant s torus can also be written as $d^{\mathbb{C}} \mathfrak{Q}_-$. To explain, note first that the $b \in \{1, 2, 3\}$ versions of these pull-back equations in this case are

- $(\nabla_1 \mathfrak{Q})^{12} - (\nabla_2 \mathfrak{Q})^{11} + i \mathfrak{Q}^{31} = i j^1 = i a_- e^{-s}$ .
- $(\nabla_1 \mathfrak{Q})^{22} - (\nabla_2 \mathfrak{Q})^{21} + i \mathfrak{Q}^{32} = i j^2 = a_- e^{-s}$ .
- $(\nabla_1 \mathfrak{Q})^{32} - (\nabla_2 \mathfrak{Q})^{31} + i \mathfrak{Q}^{33} = 0$ .

(A3.64)

Taking $\mathfrak{Q}$ to be T-invariant with only components $\mathfrak{Q}^{31}$ and $\mathfrak{Q}^{32}$ not zero, then the third bullet is obeyed automatically and the top two bullet are obeyed if and only if

$$-\mathfrak{Q}^{32} + i\,\mathfrak{Q}^{31} = i\,a_-\,e^{-s} \quad and \quad \mathfrak{Q}^{31} + i\mathfrak{Q}^{32} = a_-\,e^{-s}\ .$$
(A3.65)

These equation are simultaneously solved using $\mathfrak{Q}^{31} = a_-\,e^{-s}$ and $\mathfrak{Q}^{32} = 0$.

By way of a contrast, there is no solution to any constant s pull-back of the equation $d^{\mathbb{C}}\mathfrak{Q} = \mathfrak{J}_+$ if $j_+$ is not zero. To see why this is, note first that these pull-back equations in this case read

- $(\nabla_1\mathfrak{Q})^{12} - (\nabla_2\mathfrak{Q})^{11} + i\mathfrak{Q}^{31} = -2i j^1 = -2i\,a_+\,e^{3s}$.
- $(\nabla_1\mathfrak{Q})^{22} - (\nabla_2\mathfrak{Q})^{21} + i\mathfrak{Q}^{32} = -2i j^2 = 2a_+\,e^{3s}$.
- $(\nabla_1\mathfrak{Q})^{32} - (\nabla_2\mathfrak{Q})^{31} + i\mathfrak{Q}^{33} = 0$.

(A3.66)

Since the right hand side of these equations is T invariant, there is nothing lost by proving that there are no T invariant versions of $\mathfrak{Q}$ that obey these equations. In this case, the third equation in (A3.66) says that $\mathfrak{Q}^{33} = 0$. Meanwhile, the T invariant version of the first two equations require that $\mathfrak{Q}^{31}$ and $\mathfrak{Q}^{32}$ obey

- $-\mathfrak{Q}^{32} + i\mathfrak{Q}^{31} = -2i\,a_+\,e^{3s}$.
- $\mathfrak{Q}^{31} + i\mathfrak{Q}^{32} = 2a_+\,e^{3s}$.

(A3.67)

These two equations can not both be solved; from which it follows that $\mathfrak{J}_+$ is not $d^{\mathbb{C}}$-exact if it is not zero.

### A4. Approximation by $(S^3-K) \times S^1$ solutions with finitely many modes

Fix $r > c_*$. The first proposition in this section makes an assertion to the effect that any element in the kernel of $\mathcal{L}^\dagger_{g_K}$ on the $s \leq r$ part of $(S^3-K) \times S^1$ can be uniformly approximated on a somewhat smaller domain by an element from the kernel of $\mathcal{L}^\dagger_{g_K}$ on all of $(S^3-K) \times S^1$.

**Proposition A4.1**: *There exists $\kappa > 1$ such that if $r > \kappa^2$ then the following is true: Fix $n \in \mathbb{Z}$ and suppose that $\mathcal{X}^{(n)}$ is an element in the kernel of $\mathcal{L}^\dagger_{g_K}$ on the $s < r$ part of $(S^3-K) \times S^1$ with Fourier mode number n. Assume that the sum of the $L^2$ norms of $\mathcal{X}^{(n)}$ on the $s \leq \kappa$ and the $s \in [r - \kappa, r]$ parts of $(S^3-K) \times S^1$ is equal to 1. There is an element*

*(denoted by $\mathcal{X}_\infty^{(n)}$) from the kernel of $\mathcal{L}_{g_K}^\dagger$ on the whole of $(S^3-K)\times S^1$ with the same Fourier mode number n such that*

$$|\mathcal{X}^{(n)} - \mathcal{X}_\infty^{(n)}| \le \kappa\, e^{-|n|\,s}\, e^{-|n|r/2}\, \exp(-\kappa^{-1} e^{r/4}) + \kappa \exp(-\kappa^{-1} e^{(r-s)})$$

*on the $s \le r - \kappa$ part of $(S^3-K)\times S^1$. In particular, this element $\mathcal{X}_\infty^{(n)}$ is chosen so that the norm of the $(k_1=0, k_2=0)$ Fourier mode from $\mathcal{X}^{(n)} - \mathcal{X}_\infty^{(n)}$ for the T action on the $s \in [0, r]$ part of $S^3-K$ is a non-increasing function of s.*

The next proposition quantifies the sense in which an element in the kernel of $\mathcal{L}_{g_K}^\dagger$ on the $s \le r$ part of $(S^3-K)\times S^1$ can be uniformly approximated on a somewhat smaller domain by an element from the kernel of $\mathcal{L}_{g_K}^\dagger$ with an a priori bound on the number of modes in its $S^1$ action Fourier decomposition.

**Proposition A4.2**: *There exists $\kappa > 1$ such that if $r > \kappa^2$ then the following is true: Suppose that $\mathcal{X}$ is an element in the kernel of $\mathcal{L}_{g_K}^\dagger$ on the $s < r$ part of $(S^3-K)\times S^1$ with sum of the $L^2$ norms of $\mathcal{X}^{(n)}$ on the $s \le \kappa$ and the $s \in [r - \kappa, r]$ parts of $(S^3-K)\times S^1$ is equal to 1. Given $N > \kappa$, let $\mathcal{X}_N$ denote the sum of the $S^1$ action Fourier modes of $\mathcal{X}$ with mode numbers obeying $|n| \le N$. Then*

$$|\mathcal{X} - \mathcal{X}_N| \le \kappa(e^{\kappa|N|}\, e^{-|N|r/4} + \exp(-\kappa^{-1} e^{r/4}))$$

*on the $s \le \tfrac{3}{4} r$ part of $(S^3-K)\times S^1$.*

Sections A4a-c contain the proof of Proposition A4.1 and Section A4d proves Proposition A4.2.

**a) Outline of the proof of Proposition A4.1**

The outline is this: The element $\mathcal{X}^{(n)}$ is first written as $\mathcal{X} e^{in\theta}$ with $\mathcal{X}$ coming from a triple $C = (\mathcal{X}, \mathfrak{s}, \mathfrak{j})$ that obeys the equations in (A3.11) and (A3.12) on the $s \le r$ part of $S^3-K$. This triple is then written on the $s \in [0, r]$ part of $S^3-K$ as $C = C_0 + C_\perp$ where $C_0$ is the $(k_1=0, k_2=0)$ Fourier mode from $C$ for the T action on this part of $S^3-K$. What is denoted by $C_\perp$ obeys

$$|C_\perp| \le c_*(\exp(-c_*^{-1} e^s) + \exp(-c_*^{-1} e^{(r-s)}))$$

(A4.1)

when s ∈ $[c_*, r - c_*]$, this being a consequence of Proposition A3.1 and Lemma A3.8 with Proposition A3.13. Meanwhile, $C_0$, whose components are functions only of the variable s, can be written as $C_0 = C_+ + C_-$ with $|C_-|$ being bounded by $c_* n^3 e^n$ and with $|C_+|$ growing with s. That this is so follows from Lemma A3.3 and Lemma A3.10 with Proposition A3.13.

The crucial observation now is that there exists a data set $C_\infty = (X_\infty, \mathfrak{s}_\infty, \mathfrak{j}_\infty)$ with the following two properties:

- *This data set obeys (A3.11) and (A3.12) on the whole of $S^3 - K$.*
- *$C_\infty$ can be written as $C_+ + C_{\infty-} + C_{\infty\perp}$ on the $s \geq 0$ part of $S^3 - K$ with $C_+$ coming from the given solution C, with $|C_{\infty\perp}|$ also obeying (A3.68) and with $|C_{\infty-}|$ also bounded by $c_* n^3 e^n$.*

(A4.2)

The important point here is that the $(k_1 = 0, k_2 = 0)$ Fourier mode in $C_\infty$ for the T action has the same $C_+$ part as that in C. The existence of this $C_\infty$ follows from the various propositions in Section A3. With $C_\infty$ as just described, it then follows as a consequence that $|C - C_\infty| \leq c_* n^3 e^n$ on the s ∈ $[c_*, r - c_*]$ part of $S^3 - K$. This norm bound runs afoul of what is said by Lemmas A3.2 and A3.8 unless

$$|C_- - C_{\infty-}| \leq c_* e^{-|n|s} e^{-|n|r/2} \exp(-c_*^{-1} e^{r/4})$$

(A4.3)

except in three instances. These instances are when Lemma A3.2 is silent (the n = -1 case and the n = +1 case for the repective upper and lower equations in (A3.13)) and one instance when Lemma A3.8 is silent (the n = 0 case when j is present). However, those parts of $C_- - C_{\infty-}$ that are not covered by Lemmas A3.2 and A3.8 in these special cases are zero on the s ∈ $[c_*, r - c_*]$ part of $S^3 - K$ because of what is said about the Codazzi tensors and $H^1(S^3 - K; V)$ in Part 5 of Section A3e and what is said about the harmonic 1-forms and $H^1(S^3 - K; \mathbb{R})$ in Part 4 of Section A3f.

Now suppose that (A4.3) holds. Except for the special cases just mentioned, (A4.3) and (A4.1) with Lemmas A3.2 and A3.8 lead to a $c_* e^{-|n|r/2} \exp(-c_*^{-1} e^{r/4})$ bound on the $L^2$ norm of $C - C_\infty$ on the whole of the $s \leq \frac{5}{8} r$ part of $S^3 - K$. The latter bound with an appeal to a suitable version of Lemma A2.2 leads in turn to the same bound (with larger $c_*$) for the $L^2$ norm of $|\nabla\nabla(C - C_*)|^2$ on the $s \leq \frac{9}{16} r$ part of $S^3 - K$. This second derivative bound can be used to start a standard elliptic boostrapping argument to obtain the pointwise bound that is asserted by the proposition. The proof that this conclusion also holds for the special cases when n ∈ {0, ±1} (after a suitably modification of $C_\infty$) requires an additional input from what is said in Part 5 of Section A3e about Codazzi tensors on $S^3-K$ and what is said in Part 4 of Section A3f about harmonic 1-forms $S^3-K$.

**b) Proof of Proposition A4.1 when j = 0**

The three parts of the proof that follow directly give the details of how this works in the cases when j is identically 0 and when $n \geq 0$. The cases when $n \leq 0$ become the $|n| \geq 0$ cases after complex conjugation. By way of notation, the proof uses $\kappa_*$ to denote the version of $\kappa$ that appears in Proposition A3.1. The proof also uses $c_*$ to denote a number that is greater than 1 and independent of n, r and any elements under consideration from the kernel of $\mathcal{L}_{g_K}^\dagger$.

*Part 1*: Suppose in what follows that $C = (\mathcal{X}, \mathfrak{s}, j = 0)$. Introduce $\mathfrak{T}_+ = \mathcal{X} + i\mathfrak{s}$ and $\mathfrak{T}_- = \mathcal{X} - i\mathfrak{s}$ as done in Section A3e. These obey the equations in (A3.13). Assume first that $n \geq 3$. The four steps that follow prove the proposition's assertion for this case.

<u>Step 1</u>: Assume here that $r \geq \kappa_*^2$ and that the sum of the $L^2$ norms of $\mathcal{X}$ on the $s \leq \kappa_*$ and the $s \in [r-\kappa_*, r]$ parts of $S^3-K$ are equal to 1. It follows from Proposition A3.1 and Lemma A3.3 that $\mathfrak{T}_+$ and $\mathfrak{T}_-$ can be written where $s \in [\kappa_*, r-\kappa_*]$ on $S^3-K$ as

- $\mathfrak{T}_+ = e^{(n+2)s} t_1 \begin{pmatrix} 1 & i & 0 \\ i & -1 & 0 \\ 0 & 0 & 0 \end{pmatrix} + e^{-ns} t_2 \begin{pmatrix} 1 & -i & 0 \\ -i & -1 & 0 \\ 0 & 0 & 0 \end{pmatrix} + \mathfrak{r}_1$,

- $\mathfrak{T}_- = e^{ns} t_3 \begin{pmatrix} 1 & i & 0 \\ i & -1 & 0 \\ 0 & 0 & 0 \end{pmatrix} + e^{-(n-2)s} t_4 \begin{pmatrix} 1 & -i & 0 \\ -i & -1 & 0 \\ 0 & 0 & 0 \end{pmatrix} + \mathfrak{r}_2$,

(A4.4)

where the notation is as follows: What is denoted by $\{t_k\}_{k=1,2,3,4}$ are complex number; and what are denoted by $\mathfrak{r}_1$ and $\mathfrak{r}_2$ are tensors with no $(k_1 = 0, k_2 = 0)$ Fourier mode for the T action and with norms that are bounded by

$$c_*(\exp(-c_*^{-1} e^s) + \exp(-c_*^{-1} e^{(r-s)})).$$

(A4.5)

Because $\mathcal{X}$ has $L^2$ norm at most 1 where $s \in [r-\kappa_*, \kappa_*]$, the norms of the coefficients $t_1$ and $t_3$ must be bounded respectively by $c_* n^{3/2} e^{-(n+2)r}$ and $c_* n^{3/2} e^{-nr}$. The norms of $t_2$ and $t_4$ are bounded by $c_* n^{3/2} e^n$ because the $L^2$ norm of $\mathcal{X}$ where $s \leq \kappa_*$ is at most 1.

Step 2: Let $\mathfrak{T}_+^{(n)}$ denote the solution to (A3.13) given by Item i) of the top bullet in Proposition A3.4; and let $\mathfrak{T}_+^{(n-2)}$ denote the analogous solution with n replaced by n - 2 in the case when n ≥ 3 and given by the lower bullet in Proposition A3.5 when n = 2. Then $\mathfrak{P}_+ = \mathfrak{T}_+ - \mathfrak{t}_1 \mathfrak{T}_+^{(n)}$ and $\mathfrak{P}_- = \mathfrak{T}_- - \mathfrak{t}_3 \mathfrak{T}_+^{(n-2)}$ can be written where $s \in [\kappa_*, r - \kappa_*]$ on $S^3 - K$ as

- $\mathfrak{P}_+ = e^{-ns} (\mathfrak{t}_2 - \mathfrak{t}_1 \alpha_n) \begin{pmatrix} 1 & -i & 0 \\ -i & -1 & 0 \\ 0 & 0 & 0 \end{pmatrix} + \mathfrak{r}_1'$,

- $\mathfrak{P}_- = e^{-(n-2)s} (\mathfrak{t}_4 - \mathfrak{t}_2 \alpha_{n-2}) \begin{pmatrix} 1 & -i & 0 \\ -i & -1 & 0 \\ 0 & 0 & 0 \end{pmatrix} + \mathfrak{r}_2'$,

(A4.6)

where $\alpha_n$ and $\alpha_{n-2}$ have norm less than $c_* e^n$, and where $\mathfrak{r}_1'$ and $\mathfrak{r}_2'$ have the same properties as do the $\mathfrak{r}_1$ and $\mathfrak{r}_2$ terms in (A4.4). In particular, their norms are bounded by what is written in (A4.5).

Step 3: Supposing that r and r´ are both in the interval $[\kappa_*, r - \kappa_*]$ and that $r < r'$, then the integral of $|\mathfrak{P}_+|^2$ on the $s \in [r, r']$ part of $S^3 - K$ obeys

$$\int_{r \le s \le r'} |\mathfrak{P}_+|^2 \le c_* e^{-2nr} |\mathfrak{t}_2 - \mathfrak{t}_1 \alpha_n|^2 + c_* \exp(-c_*^{-1} e^r) + c_* \exp(-c_*^{-1} e^{(r-r')}).$$

(A4.7)

Meanwhile, the bound in the second bullet of Lemma A3.2 with r´ replacing r and $\mathfrak{P}_+$ replacing $\mathfrak{T}_+$ impies that

$$e^{2n(r'-r)} \int_{s \le r} |\mathfrak{P}_+|^2 \le \int_{r \le s \le r'} |\mathfrak{P}_+|^2 .$$

(A4.8)

Since the integral of $|\mathfrak{P}_+|^2$ on the $s \le r$ part of $S^3 - K$ is not less than $\frac{1}{n+1} c_*^{-1} |\mathfrak{t}_2 - \mathfrak{t}_1 \alpha_n|^2$, the inequalities in (A4.7) and (A4.8) are not mutually compatible when $r \ge c_*$ unless

$$|\mathfrak{t}_2 - \mathfrak{t}_1 \alpha_n|^2 \le c_* e^{-2n(r'-r)} (\exp(-c_*^{-1} e^r) + c_* \exp(-c_*^{-1} e^{(r-r')})) .$$

(A4.9)

Take $r = \frac{1}{4} r$ and $r' = \frac{3}{4} r$ in this bound to see that

$$|\mathfrak{t}_2 - \mathfrak{t}_1 \alpha_n|^2 \le c_* e^{-nr} \exp(-c_*^{-1} e^{r/4}) .$$

(A4.10)

Much the same argument proves that $|t_4 - t_3 \alpha_{n-2}|^2$ is also bounded by what is written on the right hand side of (A4.1). Using these bounds in (A4.6) gives the bound in (A4.3).

<u>Step 4</u>: Note in particular that (A4.3) with (A4.5) as a bound for $|\mathfrak{r}_1'| + |\mathfrak{r}_2'|$ gives the assertion made by the proposition for the $s \in [\frac{1}{4}r, r - \kappa_*]$ part of $S^3-K$ in the cases when $n \geq 2$ and $j = 0$. To obtain the desired bounds for the $s \leq \frac{1}{4}r$ part of $S^3-K$, use (A4.10) in the $r = \frac{1}{4}r + 1$ and $r' = \frac{3}{4}r$ version of (A4.7) and (A4.8) to see that

$$\int_{s \leq \frac{1}{4}r+1} |\mathfrak{P}_+|^2 \leq c_* \, e^{-nr} \exp(-c_*^{-1} e^{r/4}) \,.$$

(A4.11)

Use this bound with the instance of Lemma A2.2 that has V being the $s \leq \frac{1}{2}r + 1$ part of $S^3-K$ and $\mathfrak{f}$ being the function $s \to \mathfrak{f}(s) = (\frac{1}{2}r + 1 - s)$ to obtain the bound

$$\int_{s \leq \frac{1}{4}(r+1)} |\nabla\nabla\mathfrak{P}_+|^2 \leq c_* \, e^{-nr} \exp(-c_*^{-1} e^{r/4}) \,.$$

(A4.12)

This last bound with standard elliptic regularity arguments bounds the pointwise norm of $\mathfrak{P}_+$ on the $s \leq \frac{1}{4}r$ part of $S^3-K$ by $c_* \, e^{-nr} \exp(-c_*^{-1} e^{r/4})$ which is slightly stronger than the bound asserted by the proposition where $s \leq \frac{1}{4}r$. Much the same argument using $\mathfrak{P}_-$ will bound the pointwise norm of $\mathfrak{P}_-$ on the $s \leq \frac{1}{4}r$ part of $S^3-K$ by $c_* \, e^{-nr} \exp(-c_*^{-1} e^{r/4})$ also. These $\mathfrak{P}_+$ and $\mathfrak{P}_-$ bounds complete the proof of the proposition when $n \geq 3$ and $j = 0$.

*Part 2*: This part of the proof deals with the cases when $n = 2$ and $n = 0$; but only the $n = 2$ case will be discussed because the argument in the case $n = 0$ is obtained from the $n = 2$ case by switching the roles of $\mathfrak{T}_+$ and $\mathfrak{T}_-$ in what follows.

In the case $n = 2$, the top bullet of (A4.4) again depicts $\mathfrak{T}_+$ on the $s \in [\kappa_*, r-\kappa_*]$ part of $S^3-K$. Proposition A3.4 supplies a $\mathfrak{T}_+^{(n=2)}$ that obeys the top bullet of (A3.13) on the whole of $S^3-K$ and is such that $\mathfrak{P}_+ = \mathfrak{T}_+ - \mathfrak{T}_+^{(2)}$ is described by the top bullet of (A4.6). However, the depiction of $\mathfrak{T}_-$ on the $s \in [\kappa_*, r-\kappa_*]$ part of $S^3-K$ by the lower bullet of (A4.4) must be amended in this case with the addition of a term having the form

$$\alpha \, e^{3s} \begin{pmatrix} 0 & 0 & 1 \\ 0 & 0 & i \\ 1 & i & 0 \end{pmatrix}$$

(A4.13)

with $\alpha \in \mathbb{C}$. (This addition is needed because of what is said by Item b) of the second bullet of Lemma A3.3.) Proposition A3.5's top bullet now supplies a solution (denoted here by $\mathfrak{T}_-^{(2\diamond)}$) to the n = 2 version of the lower bullet in (A3.13) on the whole of $S^3-K$, and the Proposition A3.5's second bullet supplies second solution, $\mathfrak{T}_-^{(n=2)}$, to this same equation; and these have the following beneficial property: If $\mathfrak{P}_-$ is defined by the rule $\mathfrak{P}_- = \mathfrak{T}_- - t_2\mathfrak{T}_-^{(2)} - \alpha\mathfrak{T}_-^{(2\diamond)}$, then this version of $\mathfrak{P}_-$ is described by the n = 2 version of the second bullet in (A4.6).

Granted the preceding, then the arguments in Steps 3 and 4 of Part 1 can be repeated but for cosmetic changes to prove Proposition A4.1 when j = 0 and n = 2.

*Part 3*: The seven steps that follow in this part of the proof deal with the n = 1 Fourier mode when j = 0.

<u>Step 1</u>: The tensor $\mathfrak{T}_+$ on the $s \in [\kappa_*, r - \kappa_*]$ part of $S^3-K$ is depicted by the n = 1 version of the top bullet of (A4.4). This being the case, the n = 1 version of Proposition A3.4 supplies a tensor $\mathfrak{T}_+^{(n=1)}$ that solves the top bullet of (A3.13) when n = 1 and is such that $\mathfrak{P}_+ = \mathfrak{T}_+ - t_1\mathfrak{T}_+^{(1)}$ has the form depicted in the top bullet of (A4.6) on the $s \in [\kappa_*, r - \kappa_*]$ part of $S^3-K$. The arguments in Step 3 of Part 1 for $\mathfrak{P}_+$ can be repeated to see that the corresponding $(t_2 - t_1 \alpha_{n=1})$ obeys the n = 1 version of (A4.10). The arguments in Step 4 of Part 1 can also be repeated using just $\mathfrak{P}_+$ to see that $|\mathfrak{P}_+|$ obeys the bound that is asserted for $|\mathcal{X}^{(2)} - \mathcal{X}_\infty^{(2)}|$ by Proposition A4.1.

<u>Step 2</u>: The tensor $\mathfrak{T}_-$ can be written on the $s \in [\kappa_*, r - \kappa_*]$ part of $S^3-K$ as

$$\mathfrak{T}_- = \alpha\,(e^s\,t \begin{pmatrix} 1 & i & 0 \\ i & -1 & 0 \\ 0 & 0 & 0 \end{pmatrix} + e^s\,\bar{t} \begin{pmatrix} 1 & -i & 0 \\ -i & -1 & 0 \\ 0 & 0 & 0 \end{pmatrix}) + \gamma e^{3s} \begin{pmatrix} -1 & 0 & 0 \\ 0 & -1 & 0 \\ 0 & 0 & 2 \end{pmatrix} + \mathfrak{r}_2,$$

(A4.14)

with the notation as follows: The complex number $t$ is from Proposition A3.7. What is denoted in (A4.14) by $\alpha$ and $\gamma$ are complex numbers; and $\mathfrak{r}_2$ is a tensor with no $(k_1=0, k_2=0)$ Fourier mode for the T action whose norm is bounded by what is written in (A4.5). By way of an explanation for (A4.14): The appearance of the term with $\gamma$ is allowed by Item a) of the second bullet of Lemma A3.3. The form of the other term (with the two matrices appearing with the respective factors $t$ and $\bar{t}$) is dictated by the fact that $\mathfrak{T}_-$ defines a traceless Codazzi tensor on the $s \le r$ part of $S^3-K$. As explained in Part 5 of Section A3e, the appearance of the $(k_1=0, k_2=0)$ mode of a traceless Codazzi tensor on the $0 \le s \le r$ part of $S^3-K$ is constrained to have the form in (A4.14) when it extends to the $s \le 0$ part of $S^3-K$.

Step 3: Being that $\mathfrak{T}_-$ defines a Codazzi tensor on the $s \le r$ part of $S^3 - K$, it has an associated Ferus/Lafontaine class in $H^1(S^3 - K; V) \otimes_\mathbb{R} \mathbb{C}$. Denote this class by $\mathfrak{h}_-$. Let $\mathfrak{T}^{(c-)}$ denote the $\mathfrak{h} = \mathfrak{h}_-$ version of what is denoted by $\mathfrak{T}_+^{(c)}$ in the second bullet of Proposition A3.7. This Codazzi tensor can be written where $s \ge 0$ as

$$\mathfrak{T}^{(c-)} = \alpha \, ( e^s \, t \begin{pmatrix} 1 & i & 0 \\ i & -1 & 0 \\ 0 & 0 & 0 \end{pmatrix} + e^s \, \overline{t} \begin{pmatrix} 1 & -i & 0 \\ -i & -1 & 0 \\ 0 & 0 & 0 \end{pmatrix}) + \mathfrak{r}_2',$$

(A4.15)

where $\mathfrak{r}_2'$ has no $(k_1 = 0, k_2 = 0)$ Fourier mode and it is such that $|\mathfrak{r}_2'| \le c_* \exp(-e^s/c_*)|\mathfrak{h}_-|$. By way of an explanation, the appearance of the same complex number $\alpha$ in (A4.15) and (A4.14) follows from the fact that the restrictions to any constant $s$ slice of $\mathfrak{T}_-$ and $\mathfrak{T}^{(c-)}$ give the same class in $H^1(T; V)$. The bound on $|\mathfrak{r}_2'|$ comes courtesy of Proposition A3.7; but note that the bound is proportional to the $H^1(S^3 - K; V)$ norm of the class $\mathfrak{h}_-$. The next step explains why this norm is apriori bounded by $c_*$.

Step 4: To bound the size of $|\mathfrak{h}_-|$, one must first invoke a suitable instance of Lemma A2.1 and the fact that $X$ has $L^2$ norm equal at most 1 on the $s \le \kappa_*$ part of $S^3 - K$ to obtain a $c_*$ bound on the $L^2$ norm of $\nabla \nabla X$ on the $s \le 2$ part of $S^3 - K$. The latter bound and the $L^2$ norm bound for $X$ then lead via a Sobolev inequality to an $L^2$ norm bound on $\nabla X$ on this same part of $s \le 2$ part of $S^3 - K$. The resulting $L_2^2$ Sobolev bound with some standard elliptic regularity imply a $c_*$ bound for the $C^2$ norm of $X$ on the $s \le 1$ part of $S^3 - K$. The $C^2$ bound implies a $c_*$ bound on the $C^1$ norm of $\mathfrak{T}_-$ on the $s \le 1$ part of $S^3 - K$.

Let $\mathfrak{U}$ denote the locally finite cover of $S^3 - K$ from Step 4 of the proof of Proposition A3.7; and let $\mathfrak{U}_1 \subset \mathfrak{U}$ denote the subset of balls that intersect the $s \le 1$ part of $S^3 - K$. Let $B$ denote a ball from $\mathfrak{U}_1$. Since $\mathfrak{T}_-$ is a Codazzi tensor on $B$, it can be written on $B$ as $\mathfrak{T}_{-ab}|_B = \nabla_a \nabla_b f_B - \delta_{ab} f_B$. The function $f_B$ is unique up to the addition of a function that obeys (A3.46). This freedom can be used to find a version of $f_B$ with a $c_*$ bound on its $C^3$ norm. It follows as a consequence that the cocycles $\{h_{BB'} = f_B - f_{B'}\}_{B,B' \in \mathfrak{U}_1}$ have a priori $c_*$ bounds on their $C^3$ norms.

With the preceding in mind, let $\{\mathfrak{h}_i\}_{i=1,\ldots,N}$ denote the basis for $H^1(S^3 - K; V)$ from Step 6 of the proof of Proposition A3.7. For each such $\mathfrak{h}_i$, choose cocycle representatives $\{h_{iBB'}\}_{B,B' \in \mathfrak{U}_1}$. These can and should be chosen without regards to either $X$ or $R$. This guarantees that the $C^3$ norms of all of these functions are bounded by a number that is independent of $X$ and $R$; and this is to say that their $C^3$ norms are bounded by $c_*$.

Since $\mathfrak{h} = \sum_{i \in \{1,\ldots,N\}} a_i \mathfrak{h}_i$ with $\{a_i\} \in \mathbb{R}$, there exists $\{h_B \in \mathcal{V}(B)\}_{B \in \mathfrak{U}_1}$ such that

$$h_{BB'} = \sum_{i \in \{1,\ldots,N\}} a_i h_{iBB'} + h_B - h_{B'}$$

(A4.16)

for each pair B, B′ ∈ $\mathfrak{U}_1$ with B∩B′ ≠ ø. Since each function from the set $\{h_{BB'}\}_{B,B' \in \mathfrak{U}_1}$ and from each i ∈ {1, …, N} version of $\{h_{iBB'}\}_{B,B' \in \mathfrak{U}_1}$ has a $c_*$ bound on its $C^3$ norm the identity in (A4.16) can hold only if there is a corresponding $c_*$ bound on the norms of each i ∈ {1, …, N} version of $a_i$. (Keep in mind that each B ∈ $\mathfrak{U}_1$ version of $\mathcal{V}(B)$ is a fixed, 4-dimensional vector space as is each B,B′ version of $\mathcal{V}(B \cap B')$. As a consequence, (A4.16) is a linear equation between vector spaces whose dimensions are fixed a priori by the number of sets in R and $\mathcal{X}$ independent set $\mathfrak{U}_1$ and their mutual intersections.)

These $c_*$ bounds for $\{|a_i|\}_{i \in \{1,\ldots,N\}}$ leads directly to a $c_*$ bound for $\sum_{i=1,\ldots N} |a_i|$ which is the desired a priori $c_*$ bound for $|\mathfrak{h}_-|$.

<u>Step 5</u>: Use $\mathfrak{T}^{(f)}$ to denote here the traceless Codazzi tensor that is described by the first bullet of Proposition A3.7. Define $\mathfrak{P}_-$ to be $\mathfrak{T}_- - \mathfrak{T}^{(c-)} - \beta \mathfrak{T}^{(f)}$ with $\mathfrak{T}^{(c-)}$ being the Codazzi tensor in (A4.15). The tensor $\mathfrak{P}_-$ is a traceless Codazzi tensor on the s ≤ $r$ part of $S^3$−K with Ferus/Lafontaine class zero whose norm is bounded by the expression in (A4.5) on the s ∈ [$\kappa_*, r - \kappa_*$] part of $S^3$−K. The Ferus/Lafontaine class of $\mathfrak{P}_-$ being zero, this tensor can be written as $\nabla \nabla u - g_K u$ with $u$ being a function on the s < $r$ part of $S^3$−K. Since $\mathfrak{P}_-$ is traceless, the function $u$ necessarily obeys the equation $\nabla^\dagger \nabla u + 3u = 0$.

<u>Step 6</u>: The bound for $|\mathfrak{P}_-|$ on the s ∈ [$\kappa_*, r - \kappa_*$] part of $S^3$−K implies a similar bound for both $u$ and $\nabla \nabla u$. Here is why: The function $u$ on the s ∈ [0, $r$] part of $S^3$−K lacks a ($k_1$=0, $k_2$=0) Fourier mode for the T action because $\mathfrak{P}_-$ lacks this same mode. (Note that the ($k_1$=0, $k_2$=0) mode is a function of s, and since $\nabla^\dagger \nabla u + 3u = 0$, it must be a linear combination of $e^{3s}$ and $e^{-s}$.) With the vanishing of the ($k_1$=0, $k_2$=0) mode understood, take the $L^2$ norm of $\nabla \nabla u - g_K u$ on the s = r torus and use the triangle inequality to see that

$$\tfrac{1}{2} \int_{s=r} |\nabla \nabla u|^2 - 4 \int_{s=r} |u|^2 \leq \int_{s=r} |\mathfrak{P}_-|^2$$

(A4.17)

The left hand side of this equation is no smaller than $c_*^{-1} e^{2r} \int_{s=r} |u^\perp|^2$ when r ≥ $c_*$ because of the absence of the ($k_1$=0, $k_2$=0) mode. Thus, (A4.17) leads to the inequality

$$\int_{s=r} |\nabla\nabla u|^2 + \int_{s=r} |u|^2 \le c_*(\exp(-c_*^{-1} e^r) + \exp(-c_*^{-1} e^{(r-r)}))$$

(A4.18)

when $r \in [\kappa_*, r - \kappa_*]$.

Step 7: Keeping (A4.18) in mind, integrate the identity $u(\nabla^\dagger \nabla u + 3u) = 0$ over the $s \le \frac{1}{2} r$ part of $S^3 - K$ and then integrate by parts to see that

$$\int_{s \le r/2} |\nabla u|^2 + \int_{s \le r/2} |u|^2 \le c_* \exp(-c_*^{-1} e^{r/2})$$

(A4.19)

Now integrate $(\nabla^\dagger \nabla u + 3u)^2$ over the $s \le \frac{1}{2} r$ part of $S^3 - K$ and integrate by parts to bound the $L^2$ norm of $\nabla\nabla u$ on this same part of $S^3 - K$ by the expression on the left hand side of (A4.19) (with a larger version of $c_*$). This leads directly to a similar $L^2$ bound for $\mathfrak{P}_-$, and then standard elliptic regularity gives $C^0$ bound that is claimed by the proposition.

### c) Proof of Proposition A4.1 when $j \ne 0$

The proof of the proposition when $j \ne 0$ has nine steps. What follows outlines these steps for the $n \ge 0$ case. The $n < 0$ case follows from what is written below because the complex conjugate $\mathcal{X}^{(n)}$ is in the kernel of $\mathcal{L}_{g_K}^\dagger$ and it has Fourier mode number $-n$. Let $\kappa_*$ now denote the larger of the versions of $\kappa$ that appear in Proposition A3.1 and in Lemma A3.8. Assume now that $r \ge \kappa_*^2$ and that the sum of the $L^2$ norms of $\mathcal{X}^{(n)}$ on the $s \le 2\kappa_*$ and the $s \in [r - 2\kappa_*, r]$ parts of $S^3 - K$ is equal to 1.

Step 1: The $L^2$ norm of $j$ on the $s \le \kappa_*$ and the $s \in [r - \frac{7}{4}\kappa_*, r - \frac{1}{4}\kappa_*]$ parts of $S^3 - K$ is bounded by $c_*$. This can be proved using an instance of Lemma A2.2 to first bound the $L^2$ norm of $\nabla\nabla \mathcal{X}$ where $s \le \frac{7}{8}\kappa_*$ and where $s \in [r - \frac{15}{8}\kappa_*, r - \frac{1}{8}\kappa_*]$; and then using an integration by parts with a suitable cut-off function to bound the square of the $L^2$ norm of the first derivatives of $\mathcal{X}$ on a slightly smaller domain by the product of the $L^2$ norms of $\mathcal{X}$ and $\nabla\nabla \mathcal{X}$ on $s \le \frac{7}{8}\kappa_*$ and $s \in [r - \frac{15}{8}\kappa_*, r - \frac{1}{8}\kappa_*]$ parts of $S^3 - K$. (Keep in mind that the second bullet of (A3.12) writes $j$ as a linear combination of first derivatives of $\mathcal{X}$.)

Step 2: When $n \ne \pm 2$, define $\mathcal{X}^{(j)}$ and $\mathfrak{s}^{(j)}$ from $j$ using the formulas in Proposition A3.13. Let $(\mathcal{X}', \mathfrak{s}') = (\mathcal{X} - \mathcal{X}^{(j)}, \mathfrak{s} - \mathfrak{s}^{(j)})$. The pair $(\mathcal{X}', \mathfrak{s}')$ obeys (A3.12) with $j$ replaced by 0. The tensor $\mathcal{X}'$ is therefore amenable to the treatment in Parts 1-3 of the proof with $r$ replaced by $r´ = r - \frac{1}{4}\kappa_*$. In particular, the assertion of Proposition A4.1 holds for $\mathcal{X}´$. Proposition A4.1's assertions for $\mathcal{X}^{(j)}$ are proved in the subsequent Steps 3-7.

The n = ±2 case is discussed in Steps 8 and 9.

<u>Step 3</u>: The number n can be any non-negative integer in this step. Write j as $j = j_0 + j_\perp$ on the $s \in [0, r´]$ part of $S^3$–K with $j_0$ being the $(k_1 = 0, k_2 = 0)$ Fourier mode from j for the T action on this part of $S^3$–K. Lemma A3.8 is used to bound $j_\perp$ and its covariant derivative on the $s \in [\kappa_*, r´-\kappa_*]$ part of $S^3$–K by the expression in (A4.5). Meanwhile, $j_0$ is described by the integer n version of Lemma A3.10. If n > 0, denote the coefficients in this version of Lemma A3.10 by $(a_+, a_-)$. In the case n = 0, denote the coefficients from this version of Lemma A3.10 by $(a_1, a_2, a_3)$

<u>Step 4</u>: If n = 0, then the vector $(a_1, a_2) \in \mathbb{R}^2$ is a real multiple of $(\mathfrak{z}+\overline{\mathfrak{z}}, i(\mathfrak{z}-\overline{\mathfrak{z}}))$ with $\mathfrak{z}$ coming from (A3.56). This constraint on $(a_1, a_2)$ follows from cohomological considerations: The closed 1-form $j_\perp$ is d-closed on the $s \in [\kappa_*, r´-\kappa_*]$ part of $S^3$–K, but it is exact because it has no $(k_1 = 0, k_2 = 0)$ Fourier mode. Thus, the cohomology class of j is that of $j_0$. Since the j extends over the $s \leq 0$ part of $S^3$–K, the cohomology class of $j_0$ is necessarily a real multiple of the generator of the 1-dimensional image of the restriction map from $H^1(S^3$–K$; \mathbb{R})$ to any constant s slice of the $s \geq 0$ part of $S^3$–K. This image is generated by the form $\hat{\jmath}$ that is depicted in (A3.56).

If n = 0, then the coefficient $a_3$ is zero for the following reason: On the one hand, the integral of $*j$ over any constant s = r slice of $S^3$–K (for $r \in [\kappa_*, r´-\kappa_*]$) is equal to $(4\pi)^2 \det(\mathfrak{m}) a_3$ with $\mathfrak{m}$ being the matrix that appears in (3.3)-(3.6). This follows because the integral of any component of $j_\perp$ is zero and because the Hodge dual of the $(a_1, a_2)$ terms that appear in Lemma A3.10's formula restricts as zero to any constant s slice. On the other hand, the integral of $*j$ over the s = r slice is zero because Stokes' theorem writes this integral as the integral of $d*j$ over the $s \leq r$ part of $S^3$–K and $d*j = 0$.

<u>Step 5</u>: If n > 0, let $j^{(n)}$ denote the solution to (A3.11) on $S^3$–K that is described by Lemma A3.A5. If n = 0, let $j^{(0)}$ denote the solution that is described by Lemma A3.12. In any case, it follows from what is said in Steps 1-4 that j can be written as $\alpha j^{(n)} + j_- + \mathfrak{r}$ on the $s \in [\kappa_*, r´-\kappa_*]$ part of $S^3$–K with $\alpha \in \mathbb{C}-0$ and with $j_-$ and $\mathfrak{r}$ as follows: What is denoted by $\mathfrak{r}$ has no $(k_1 = 0, k_2 = 0)$ Fourier mode for the T action; and its norm and that of its covariant derivative are bounded by the expression in (A4.5). What is denoted by $j_-$ is zero when n = 0; and for n > 0 the number $\alpha$ is $a_+$ and $j_-$ has the form

$$j_- = (a_- - a_+ \hat{a}_n) \, e^{(1-n)s} \begin{pmatrix} 1 \\ -i \\ 0 \end{pmatrix}$$

(A4.20)

with $\hat{a}_n$ coming from Lemma A3.11.

<u>Step 6</u>: In the case when $n = 0$, the 1-form $j - \alpha j^{(0)}$ is harmonic with no $(k_1 = 0, k_2 = 0)$ Fourier mode for the T action on the $s \in [0, r']$. This being the case, the top bullet in Lemma A3.8 can be invoked to bound the norm of this 1-form and its covariant derivative on the $s \in [\kappa_*, r' - \kappa_*]$ part of $S^3 - K$. Since this 1-form is exact, it has the form $du$ with $u$ being a harmonic function with no $(k_1 = 0, k_2 = 0)$ Fourier mode for the T action on the $s \in [0, r]$ part of $S^3 - K$. In any event, $\nabla^\dagger \nabla u$ is equal to 0. Integrate the identity $u(\nabla^\dagger \nabla u) = 0$ over the $s \leq \frac{1}{2} r$ part of $S^3 - K$ and integrate by parts to see that this version of $u$ also obeys the bound in (A4.19). Standard elliptic regularity can then be used to obtain a similar bound for the $L^2$ norm of $\nabla \nabla u$ on the $s \leq \frac{1}{2} r$ part of $S^3 - K$ and then for the pointwise norms of $u$, $\nabla u$ and $\nabla \nabla u$ on the $s \leq \frac{1}{2} r$ part of $S^3 - K$. In particular, these bounds lead to the same sort of bounds for the norm of $\nabla(j - \alpha j^{(0)})$.

Denote by $\mathcal{X}^{(j)}$ the symmetric, traceless tensor that is given by the $n = 0$ version of Proposition A3.13 from $j$. Let $\mathcal{X}^{(j^{(0)})}$ denote the corresponding tensor that comes from the $n = 0$ version of Proposition A3.13 using $j^{(0)}$. Then $\mathcal{X}^{(j)} - \alpha \mathcal{X}^{(j^{(0)})}$ obeys the bounds that are asserted by the $n = 0$ version of the Proposition A4.1 by virtue of the bounds that were derived in the preceding paragraph step for $\nabla(j - \alpha j^{(0)})$.

<u>Step 7</u>: If $n > 0$, then an argument much like that used to derive (A4.10) from (A4.7) can be used with Lemma A3.9 playing the role played by Lemma A3.2 to see that

$$|a_- - a_+ \hat{a}_n|^2 \leq c_* \, e^{-nr} \exp(-c_*^{-1} e^{r/2}) \, .$$

(A4.21)

The version of the argument that is used here replaces $\mathfrak{P}_+$ with $j - a_+ j^{(n)}$ (and thus by $j_- + \mathfrak{r}$)) in (A4.7) and (A4.8); and it replaces $\mathfrak{t}_2$ by $a_-$ and $\mathfrak{t}_1$ by $a_+$ and $\alpha_n$ by $\hat{a}_n$ in (A4.7) and (A4.9). This bound with the aforementioned bound for $\mathfrak{r}$ and the $j - a_+ j^{(n)}$ version of (A4.8) leads via arguments much like those in Step 4 of Part 1 to an $j - a_+ j^{(n)}$ analog of (A4.11) and then one of (A4.12) that make the following assertions:

- $\displaystyle\int_{s \leq \frac{1}{4} r + 1} |j - a_+ j^{(n)}|^2 \leq c_* \, e^{-nr} \exp(-c_*^{-1} e^{r/4}) \, .$

- $$\int_{s \leq \frac{1}{4}(r+1)} |\nabla(j - a_+ j^{(n)})|^2 \leq c_* e^{-nr} \exp(-c_*^{-1} e^{r/4}) \ .$$

(A4.22)

These last two bounds with standard elliptic regularity arguments (using (A3.11)) can be used to bound the pointwise norms of $j - a_+ j^{(n)}$ and $\nabla(j - a_+ j^{(n)})$ by $c_* e^{-nr} \exp(-c_*^{-1} e^{r/4})$ on the $s \leq \frac{1}{4} r$ part of $S^3 - K$.

If $n > 0$ but $n \neq 2$, define $\mathcal{X}^{(j)}$ and $\mathfrak{s}^{(j)}$ from $j$ using the formulas in Proposition A3.13. Likewise define $\mathcal{X}^{(j^{(n)})}$ from $j^{(n)}$. The tensor $\mathcal{X}^{(j)} - a_+ \mathcal{X}^{(j^{(n)})}$ obeys the bounds that are asserted by Proposition A4.1 by virtue of what is said in the previous paragraph and in Step 5 about the norm of $\nabla(j - a_+ j^{(n)})$.

Step 8: This step and Step 9 treat the $n = 2$ case. This step explains how to find the appropriate solution to the $n = 2$ versions of (A3.11) and (A3.12) on $S^3 - K$ to subtract from $\mathcal{X}$. The next step explains how to derive what is asserted by the proposition.

It follows from what is said by Proposition A3.13 that the coefficient $a_+$ must equal zero in this case; which is to say that the $(k_1 = 0, k_2 = 0)$ Fourier mode from $j$ for the $T$ action on the $s \in [0, r')$ part of $S^3 - K$ has the form

$$a_- e^{-s} \begin{pmatrix} 1 \\ -i \\ 0 \end{pmatrix} \ .$$

(A4.23)

Here, $a_- \in \mathbb{R}$; and it follows from (A4.21) that $|a_-| \leq c_* \exp(-c_*^{-1} e^{r/2})$.

Write $j$ on the $s \in [0, r')$ part of $S^3 - K$ as $j_0 + \mathfrak{r}$ with $j_0$ as depicted in (A4.23) and with $\mathfrak{r}$ lacking the $(k_1 = 0, k_2 = 0)$ Fourier for the $T$ action. Lemma A3.8 describes $\mathfrak{r}$ and its covariant derivative. A calculation finds that the tensor given by $\nabla_a j_0^b + \nabla_b j_0^a$ is identically zero. Keeping this in mind, introduce by way of notation $(\mathcal{X}_0, \mathfrak{s}_0)$ to denote the $(k_1 = 0, k_2 = 0)$ Fourier mode from $(\mathcal{X}, \mathfrak{s})$ on the $s \in [0, r)$ part of $S^3 - K$. The $(\mathcal{X}_0, \mathfrak{s}_0, j_0)$ version of the equations in (A3.12) when written in terms of $\mathfrak{T}_{0+} = \mathcal{X}_0 + i \mathfrak{s}_0$ and $\mathfrak{T}_{0-} = \mathcal{X}_0 - i \mathfrak{s}_0$ say that

- $\varepsilon^{bcd} \nabla_c \mathfrak{T}_{0+}^{ad} + 3 i \mathfrak{T}_{0+}^{ab} = 0$ .
- $\varepsilon^{bcd} \nabla_c \mathfrak{T}_{0-}^{ad} + i \mathfrak{T}_{0-}^{ab} = -i \nabla_a j_0^b$ .

(A4.24)

Let $\mathfrak{Q}$ denote the symmetric, traceless tensor with all entries zero but for the pairs $(\mathfrak{Q}^{31}, \mathfrak{Q}^{32})$ and $(\mathfrak{Q}^{13}, \mathfrak{Q}^{23})$ which are $\mathfrak{Q}^{13} = \mathfrak{Q}^{31} = \frac{1}{2} a_- e^{-s}$ and $\mathfrak{Q}^{23} = \mathfrak{Q}^{32} = -\frac{i}{2} a_- e^{-s}$. Having defined $\mathfrak{Q}$, set $\mathfrak{T}_{0-}' = \mathfrak{T}_{0-} - \mathfrak{Q}$. This tensor $\mathfrak{T}_{0-}'$ obeys

$$\varepsilon^{bcd}\nabla_c\mathfrak{T}_{0-}{}^{ad} + i\mathfrak{T}_{0-}{}^{ab} = 0 \ .$$

(A4.25)

This last equation with the top equation in (A4.24) imply that the pair $(\mathfrak{T}_{0+}, \mathfrak{T}_{0-})$ have the following properties: The tensor $\mathfrak{T}_{0+}$ is described by the n = 2 version of the top bullet in (A4.4) with $\mathfrak{r}_1 = 0$; and the tensor $\mathfrak{T}_{0-}$ is a sum of what is depicted in the n = 2 version of the lower bullet in (A4.4) and what is depicted in (A4.13).

With the preceding understood, the arguments in Step 2 of Part 1 and those in Part 2 find a tensor $\mathfrak{T}_+^{(n=2)}$ from the n = 2 version of Proposition A3.4 so that $\mathfrak{P}_+ = \mathfrak{T}_+ - \mathfrak{t}_1\mathfrak{T}_+^{(2)}$ is described by the top bullet of (A4.6). Meanwhile, the arguments from Part 3 find tensors $\mathfrak{T}_-^{(2)}$ and $\mathfrak{T}_-^{(2\diamond)}$ from Proposition A3.5 such that $\mathfrak{P}_- = \mathfrak{T}_- - \mathfrak{t}_2\mathfrak{T}_-^{(2)} - \alpha\mathfrak{T}_-^{(2\diamond)}$ is described by the n = 2 version of the second bullet in (A4.6) except that $\mathfrak{r}_1$ in this case contains the contribution from $\mathfrak{Q}$ which is T invariant. But even so, the bound on $\mathfrak{r}_1$ holds because of the $c_*\exp(-c_*^{-1}e^{r/2})$ bound on $|a_-|$.

Step 9: Since $a_+ = 0$, the n=2 version of (A4.22) asserts a $c_*\exp(-c_*^{-1}e^{r/4})$ bound for the $L^2$ norms of both j and $\nabla j$ on the $s \leq \frac{1}{4}(r+1)$ part of $S^3 - K$. As noted directly after (A4.22), these lead to the pointwise norm bound $|j| + |\nabla j| \leq c_*\exp(-c_*^{-1}e^{r/4})$ on the $s \leq \frac{1}{4}r$ part of $S^3 - K$.

Were j = 0, then $\mathfrak{P}_+$ and $\mathfrak{P}_-$ would obey the n = 2 versions of (A3.13). With $j \neq 0$, they obey an inhomogeneous version of these equations that can be written as

- $\varepsilon^{bcd}\nabla_c\mathfrak{P}_+{}^{ad} + 3i\mathfrak{P}_+{}^{ab} = \mathfrak{r}_+$,
- $\varepsilon^{bcd}\nabla_c\mathfrak{P}_-{}^{ad} + i\mathfrak{P}_-{}^{ab} = \mathfrak{r}_-$,

(A4.26)

where $\mathfrak{r}_+$ and $\mathfrak{r}_-$ are terms that involve j and $\nabla j$. Their norms in any event are bounded by

$$c_*(e^{-s}\exp(-c_*^{-1}e^{r/4}) + \exp(-c_*^{-1}e^{(r-s)}))$$

(A4.27)

on the $s \leq r' - \kappa_*$ part of $S^3 - K$. This follows where $s \leq \frac{1}{4}r$ from what is said about j in the preceding paragraph; and it follows where $s \in [\frac{1}{4}r, r' - \kappa_*]$ from what is said about the norm of $a_-$ in the preceding paragraph and from what is said by Lemma A3.8 about the sum of the non-zero T action Fourier modes of j for the part of $S^3 - K$.

The arguments for the second bullet of Lemma A3.2 can be repeated using (A4.26) and the bounds in (A4.27) in lieu of (A3.13) to obtain an inequality that reads

$$\int_{s \leq r'} |\mathfrak{P}_\pm|^2 + c_*\exp(-c_*^{-1}e^{r/4}) \geq e^{(r'-r)}\int_{s \leq r} |\mathfrak{P}_\pm|^2$$

(A4.28)

when r´ ∈ [1, $\frac{1}{2}$ r] and r ∈ [0, r´]. This equation leads directly to an analog of (A4.28):

$$e^{(r´-r)} \int_{s \leq r} |\mathfrak{P}_\pm|^2 \leq \int_{r \leq s \leq r´} |\mathfrak{P}_\pm|^2 + c_* \exp(-c_*^{-1} e^{r/4}) \; ;$$

(A4.29)

and then to analogs of (A4.10) for the coefficients in the $\mathfrak{P}_\pm$ versions of (A4.6). Granted these analogs, then these argument used in Step 4 of Part 1 for (A4.11) and (A4.12) lead to

$$\int_{s \leq \frac{1}{4}(r+1)} |\nabla\nabla \mathfrak{P}_\pm|^2 + \int_{s \leq \frac{1}{4}r+1} |\mathfrak{P}_\pm|^2 \leq c_* \exp(-c_*^{-1} e^{r/4}) \; .$$

(A4.30)

With these bounds in hand, then standard elliptic regularity arguments give the pointwise bounds that are asserted by Proposition A4.1 for this n = 2 case.

### d) Proof of Proposition A4.2

Write $\mathcal{X} - \mathcal{X}_N$ as $\sum_{n \in \mathbb{Z}: |n| > N} \mathcal{X}^{(n)}$ with each $\mathcal{X}^{(n)}$ having Fourier mode number n. Let $\kappa_*$ now denote the larger of the versions of $\kappa$ that appear in Propositions A3.1 and A4.1 and in Lemma A3.8. Let $z_n$ denote here the sum of the $L^2$ norms of $\mathcal{X}^{(n)}$ on the $s \leq \kappa_*$ and the $s \in [r - \kappa_*, r]$ pars of $(S^3 - K) \times S^1$. It follows from Proposition A4.1 and what is said about the various n ∈ $\mathbb{Z}$ versions of $\mathcal{X}_\infty^{(n)}$ by Propositions A3.4 and A3.5 and A3.7 and A3.13, and Lemmas A3.11 and A3.12 that the $L^2$ norm on the $s \in [\frac{3}{4}r, \frac{3}{4}r+2]$ part of $S^3 - K$ of the $(k_1 = 0, k_2 = 0)$ Fourier mode from $\mathcal{X}^{(n)}$ for the T action on this part of $S^3 - K$ is at most $c_* e^{|n|c_*} e^{-|n|r/2} z_n$. It follows from these bounds and from Proposition A3.1 that the $L^2$ norm of $\mathcal{X} - \mathcal{X}_N$ on this same part of $S^3 - K$ is at most

$$c_*(e^{Nc_*} e^{-Nr/4} + \exp(-c_*^{-1} e^{r/4})) \; .$$

(A4.31)

Now invoke Lemma A2.2 with V being the $s \leq \frac{3}{4}r+2$ part of $S^3 - K$, with $\mathfrak{f}$ being the function $\frac{3}{4}r+2 - s$, and with $\mathcal{X} - \mathcal{X}_N$ used in lieu of $\mathcal{X}$. What with (A4.31), this call to Lemma A2.2 leads to the bound

$$\int_{s \leq \frac{3}{4}r+1} |\nabla\nabla(\mathcal{X} - \mathcal{X}_N)|^2 \leq c_* \int_{s \leq \frac{3}{4}r+1} |\mathcal{X} - \mathcal{X}_N|^2 + c_*(e^{Nc_*} e^{-Nr/4} + \exp(-c_*^{-1} e^{r/4})) \; .$$

(A4.32)

The integral on the left hand side of (A4.32) is no smaller than $c_*^{-1} N^4$. Therefore, supposing that $N > c_*$, then (A4.32) implies that

$$\int_{s \leq \frac{3}{4}r+1} |\nabla\nabla(\mathcal{X} - \mathcal{X}_N)|^2 + N^4 \int_{s \leq \frac{3}{4}r+1} |\mathcal{X} - \mathcal{X}_N|^2 \leq c_*(e^{Nc_*} e^{-Nr/2} + \exp(-c_*^{-1} e^{r/2})) \ .$$

(A4.33)

This last bound with some standard elliptic regularity arguments leads to the pointwise bound that is asserted by the proposition.

### A5. The proof of Propositions 5.3 and 5.4

Sections A5d and A5e supply the proofs of Propositions 5.3 and 5.4. By way of a look ahead, the arguments write an $\iota$ invariant element in the kernel of some large R version of $\mathcal{L}_{g_R}^\dagger$ on $\mathbb{T}_K$ as a pair $(\mathcal{X}_T, \mathcal{X}_K)$ in the manner of (5.4) and (5.5). Having done so, then what is said in Section A2 about $\mathcal{X}_T$ is compared with what is said in Sections A3 and A4 about $\mathcal{X}_K$. These comparison's take place where the domains of $\mathcal{X}_T$ and $\mathcal{X}_K$ overlap, the $T \times ((\mathbb{D}-D) \cup \iota(\mathbb{D}-D))$ part of $\mathcal{X}_T$'s domain which is the the $s \geq 0$ part of the domain $(S^3 - N_K) \times S^1$ for $\mathcal{X}_K$. Here is what the comparison involves: Introduce again $\mathcal{X}_{T0}$ to denote the $(k_1=0, k_2=0)$ from $\mathcal{X}_T$ for the T action on $\mathcal{X}_T$'s domain. As observed in (A2.20) and (A2.22), the components of $\mathcal{X}_{T0}$ can be written as

- $\mathcal{X}_{T0}{}^{33} = -(\mathcal{X}_{T0}{}^{11} + \mathcal{X}_{T0}{}^{22}) = s \ ,$
- $\mathcal{X}_{T0}{}^{13} - i\mathcal{X}_{T0}{}^{23} = \cdots c_1 z + c_0 + c_{-1} z^{-1} + c_{-2} z^{-2} + \cdots \ ,$
- $\mathcal{X}_{T0}{}^{11} - i\mathcal{X}_{T0}{}^{12} = (\cdots a_1 z + a_0 + a_{-2} z^{-2} + \cdots) \bar{z} + (\cdots b_1 z + b_0 + b_{-1} z^{-1} + b_{-2} z^{-2} + \cdots) \ ,$

(A5.1)

with $s$ being a real number and $\{c_k\}_{k \in \mathbb{Z}}$, $\{a_k\}_{k \in \mathbb{Z}; k \neq -1}$ and $\{b_k\}_{k \in \mathbb{Z}}$ being complex numbers. Meanwhile, Proposition A4.1 and the propositions in Section A3 lead to a depiction of $\mathcal{X}_K$ as a sum of $S^1$ action Fourier modes that has the form

$$\mathcal{X}_K = \sum_{n \in \mathbb{Z}} (\mathcal{X}_\infty^{(n)} + \mathfrak{r}_n)$$

(A5.2)

with each $n \in \mathbb{Z}$ version of $\mathcal{X}_\infty^{(n)}$ and $\mathfrak{r}_n$ having $S^1$ action Fourier mode number n, with $\mathcal{X}_\infty^{(n)}$ being in the kernel of $\mathcal{L}_{g_K}^\dagger$ on the whole of $(S^3-K) \times S^1$ and with $\mathfrak{r}_n$ being small in a suitable sense. The depiction of $\mathcal{X}_K$ in (A5.1) is first used in conjunction with (5.5) to relate the coefficients $\{c_k, a_{k-1}, b_k\}_{k \in \mathbb{Z}}$ that appear in (A5.1) to the coefficients that define the elements $\{\mathcal{X}_\infty^{(n)}\}_{n \in \mathbb{Z}}$ from (A5.2) (and vice-versa). This is done in Section A5a. Sections A5b and A5c derive some consequences of these relations. The consequence are then used in Sections A5e and A5f to prove Propositions 5.3 and 5.4.

### a) $(\mathcal{X}_K, \mathcal{X}_T)$ where $T \times (T' - (D \cup \iota(D)))$ and $(S^3 - N_K) \times S^1$ intersect

Fix some $R > 1$ to define $\mathbb{T}_K$ and its locally conformally flat metric $\mathfrak{g}_R$ as directed in Section 3b. Write a given $\iota$-invariant element in the kernel of $\mathcal{L}_{\mathfrak{g}_R}^\dagger$ as a pair $(\mathcal{X}_T, \mathcal{X}_K)$ in the manner of (5.4) and (5.5). This subsection rewrites (5.5) with the depiction of $\mathcal{X}_T$ in (A5.1) and the depiction of $\mathcal{X}_K$ in (A5.2) in mind.

To be sure of the notation, let $(t_1, t_2, t_3, t_4)$ again denote the $\mathbb{R}/2\pi\mathbb{Z}$ coordinates for $T \times T'$ with $(t_1, t_2)$ being coordinates for the T factor and $(t_3, t_4)$ being coordinates for the T' factor. Let $(e^1, e^2, e^3, e^4)$ denote the constant basis for $T^*(T \times T')$ as defined in Part 1 of Section 5c using the coordinate differentials. Use (5.7) to define the corresponding basis $\{\omega^a\}_{a=1,2,3}$ for the $T \times T'$ version of $\Lambda^+$. (This is the basis that is used to depict $\mathcal{X}_{T0}$ in (1.1) with the complex number $z$ on $\mathbb{D}$ defined so that $dz = dt_3 + i dt_4$.) The same $(t_1, t_2)$ coordinates from Part 1 with the coordinate $s$ and the $\mathbb{R}/2\pi\mathbb{Z}$ coordinate $\theta$ give coordinates for the $s \geq 0$ part of $(S^3 - K) \times S^1$. Let $\{\mathbf{e}^1, \mathbf{e}^2, \mathbf{e}^3, \mathbf{e}^4\}$ denote the orthonormal frame for $T^*((S^3 - K) \times S^1)$ with $\mathbf{e}^1 = e^{-s} e^1$, $\mathbf{e}^2 = e^{-s} e^2$, $\mathbf{e}^3 = ds$ and $\mathbf{e}^4 = d\theta$. The corresponding basis for $\Lambda^+$ is defined by (5.7) with $\{\mathbf{e}^k\}_{k=1,2,3,4}$ used in lieu of $\{e^k\}_{k=1,2,3,4}$. It is denoted in this part of the subsection by $\{w^1, w^2, w^3\}$ to distinguish it by notation from the basis for $\Lambda^+(T \times T')$ that was defined in Part 1.

The map $\psi$ whose pull-back appears in (5.5) has no effect on the pair $(t_1, t_2)$ while writing $(t_3, t_4)$ as $(t_* + e^{s-R}\cos\theta, t_* + e^{s-R}\sin\theta)$. It follows as a consequence that this map $\psi$ pulls back the basis $\{\omega^a\}_{a=1,2,3}$ as follows:

$$\psi^*(\omega^1 + i\omega^2) = e^{2(s-R)+i\theta}(w^1 + iw^2) \quad \text{and} \quad \psi^*\omega^3 = e^{2(s-R)} w^3.$$

(A5.3)

Thus, if $\mathcal{X}_T$ is written on $T \times (\mathbb{D} - D)$ as $\mathcal{X}_T^{ab} \omega^a \otimes \omega^b$ and $\psi^* \mathcal{X}_T$ is written as $\mathfrak{X}^{ab} w^a \otimes w^b$, then $\{\mathfrak{X}^{ab}\}_{a,b \in \{1,2,3\}}$ are given in terms of $\{\mathcal{X}_T^{ab}\}_{a,b \in \{1,2,3\}}$ by using the following formulae:

- $\mathfrak{X}^{33} = e^{4(s-R)} \mathcal{X}_T^{33}$ and $\mathfrak{X}^{11} + \mathfrak{X}^{22} = e^{4(s-R)}(\mathcal{X}_T^{11} + \mathcal{X}_T^{22})$.
- $\mathfrak{X}^{13} - i\mathfrak{X}^{23} = e^{4(s-R)+i\theta}(\mathcal{X}_T^{13} - i\mathcal{X}_T^{23})$.
- $\mathfrak{X}^{11} - \mathfrak{X}^{22} + 2i\mathfrak{X}^{12} = e^{4(s-R)+2i\theta}(\mathcal{X}_T^{11} - \mathcal{X}_T^{22} + 2i\mathcal{X}_T^{12})$.

(A5.4)

With (A5.4) understood, and supposing that $\mathcal{X}_K$ is written on the $s \geq 0$ part of $(S^3 - N_K) \times S^1$ as $\mathcal{X}_K^{ab} w^a \otimes w^b$, then (5.5) is asking that:

- $\mathcal{X}_K^{33} = e^{2(s-R)} \mathcal{X}_T^{33}$ and $\mathcal{X}_K^{11} + \mathcal{X}_K^{22} = e^{2(s-R)}(\mathcal{X}_T^{11} + \mathcal{X}_T^{22})$.
- $\mathcal{X}_K^{13} - i\mathcal{X}_K^{23} = e^{2(s-R)+i\theta}(\mathcal{X}_T^{13} - i\mathcal{X}_T^{23})$.

- $\mathcal{X}_K{}^{11} - \mathcal{X}_K{}^{22} + 2i\,\mathcal{X}_K{}^{12} = e^{2(s-R)+2i\theta}\,(\mathcal{X}_T{}^{11} - \mathcal{X}_T{}^{22} + 2i\,\mathcal{X}_T{}^{12})\ .$

(A5.5)

These equations can be used with the formulae in (A5.1) for the coefficients $\{\mathcal{X}_{T0}{}^{ab}\}_{a,b \in \{1,2,3\}}$ of the T-invariant part of $\mathcal{X}_T$ if each instance of z and $\bar{z}$ in (A5.1) is replaced with the respective functions $e^{(s-R)+i\theta}$ and $e^{(s-R)-i\theta}$. In particular, the substitution $z = e^{(s-R)+i\theta}$ and $\bar{z} = e^{(s-R)-i\theta}$ in (A5.1) writes (A5.5) as

- $\mathcal{X}_K{}^{33} = -(\mathcal{X}_K{}^{11} + \mathcal{X}_K{}^{22}) = r\,e^{2(s-R)} + \mathfrak{r}_0\ ,$
- $\mathcal{X}_K{}^{13} - i\mathcal{X}_K{}^{23} = \cdots + c_0\,e^{2(s-R)+i\theta} + c_{-1}\,e^{(s-R)} + c_{-2}\,e^{-i\theta} + c_{-3}\,e^{-(s-R)-2i\theta}\cdots + \mathfrak{r}_+\ ,$
- $\mathcal{X}_K{}^{11} - i\mathcal{X}_K{}^{12} = (\cdots + a_0\,e^{3(s-R)+i\theta} + a_{-2}\,e^{(s-R)-i\theta} + a_{-3}\,e^{-2i\theta} + a_{-4}\,e^{-(s-R)-3i\theta}\,)$

  $+ (\cdots + b_0\,e^{2(s-R)+2i\theta} + b_{-1}\,e^{(s-R)+i\theta} + b_{-2} + b_{-3}\,e^{-(s-R)-i\theta}\cdots\,) + \mathfrak{r}_{++}\ ,$

(A5.6)

where $\mathfrak{r}_0$, $\mathfrak{r}_+$ and $\mathfrak{r}_{++}$ denote the sum of the non-zero Fourier modes of $\mathcal{X}_K$ with respect to the action of T on the $s \geq 0$ part of $(S^3 - N_K) \times S^1$. (The first bullet of Proposition A3.1 supplies a priori pointwise bounds for $\{\mathfrak{r}_0, \mathfrak{r}_+, \mathfrak{r}_{++}\}$.)

**b) Constraints on $\{c_k, a_k, b_k\}_{k \in \mathbb{Z}}$ and the $S^1$ modes**

This section uses (A5.6) to relate the coefficients that appear in the depiction of $\mathcal{X}_{T0}$ in (A5.1) and thus $\mathcal{X}_K$ in (A5.6) with the coefficients that define the elements $\{\mathcal{X}_\infty{}^{(n)}\}_{n \in \mathbb{Z}}$ that appear in (A5.2). There are five parts to what follows. By way of notation, $\kappa_*$ is used to denote the larger of the versions of $\kappa$ from Propositions A3.1 and A4.1 and A4.2 and from Lemma A3.8.

*Part 1*: Supposing that n is a given integer, let $\mathcal{X}_K{}^{(n)}$ denote the part of $\mathcal{X}_K$ with $S^1$ action Fourier mode number n. Let $z_{K,n}$ denote the sum of its $L^2$ norms on the $s \leq \kappa_*$ part and the $s \in [R + \ln(\tfrac{1}{4} t_*) - \kappa_*, R + \ln(\tfrac{1}{4} t_*)]$ parts of $(S^3 - K) \times S^1$. Let $\mathcal{X}_\infty{}^{(n)}$ denote the mode number n element in the kernel of $\mathcal{L}_{g_K}^\dagger$ on the whole of $(S^3 - K) \times S^1$ that is supplied by Proposition A4.1 with input $\mathcal{X}^{(n)} = \mathcal{X}_K{}^{(n)}$ and with $r = R + \ln(\tfrac{1}{4} t_*)$. Having mode number n, this $\mathcal{X}_\infty{}^{(n)}$ can be written as $\mathcal{X}_\infty{}^{(n)} = e^{in\theta}\,\mathcal{X}_\infty$ where $\mathcal{X}_\infty$ comes from a data set $(\mathcal{X}_\infty, \mathfrak{s}_\infty, j_\infty)$ that solves the mode number n version of the equations in (A3.11) and (A3.12).

The element $\mathcal{X}_\infty$ for any given mode number $n \in \mathbb{Z}$ can in turn be written as $\mathcal{X}_\infty = \mathcal{X}_{\infty,0} + \mathcal{X}_{\infty,j}$ with $\mathcal{X}_{\infty,0}$ coming from a data set that solves the $j = 0$ version of (A3.12) using the given value of n, and with $\mathcal{X}_{\infty,j}$ determined from a solution to the mode number n version of (A3.13) using the formula in Proposition A3.13. The upcoming Parts 3-6 use the depiction in (A5.6) of $\mathcal{X}_K$ on the $s \geq 0$ part of $(S^3 - N_K) \times S^1$ to relate the coefficients in (A5.1) to those that define $\mathcal{X}_{\infty,0}$ and $\mathcal{X}_{\infty,j}$. This task is facilitated by the observation below:

*Let* n *denote a non-zero, $S^1$ action Fourier mode number. The respective $(k_1=0, k_2=0)$ Fourier modes from $\mathcal{X}_{\infty,0}$ and $\mathcal{X}_{\infty,j}$ for the T action on the $s \geq 0$ part of $S^3-K$ are pointwise orthogonal.*

(A5.7)

To prove this assertion, first write these tensors on the $s \geq 0$ part of $S^3-K$ as $3 \times 3$ matrices using the basis $\{w^a\}_{a=1,2,3}$. If $|n| \geq 3$, then the matrix for $\mathcal{X}_{\infty,0}$ is the left hand matrix in what follows and $\mathcal{X}_{\infty,j}$ is the right hand matrix:

$$\begin{pmatrix} * & * & 0 \\ * & * & 0 \\ 0 & 0 & 0 \end{pmatrix} \text{ and } \begin{pmatrix} 0 & 0 & * \\ 0 & 0 & * \\ * & * & 0 \end{pmatrix}.$$

(A5.8)

That this is so follows from Proposition A3.4 and Lemma A3.11. For $|n|=2$ meanwhile, Proposition A3.5 writes $\mathcal{X}_{\infty,0}$ as a sum of matrices that have both of the forms in (A5.8), but in this case, Proposition A3.13 says that there is no $\mathcal{X}_{\infty,j}$ to confuse the issue. In the case $n = \pm 1$, the matrix form of the $(k_1 = 0, k_2=0)$ Fourier mode from $\mathcal{X}_{\infty,0}$ is a sum of a matrix of the sort depicted on the left in (A5.8) and a multiple of the diagonal matrix

$$\begin{pmatrix} -1 & 0 & 0 \\ 0 & -1 & 0 \\ 0 & 0 & 2 \end{pmatrix},$$

(A5.9)

whereas $\mathcal{X}_{\infty,j}$ is as depicted on the right in (A5.8). (The form of $\mathcal{X}_{\infty,j}$ follows from what is said in Proposition A3.7 and that of $\mathcal{X}_{\infty,j}$ follows from what is said by Lemma A3.11.)

In the case when $n = 0$, the matrix form of $\mathcal{X}_{\infty,0}$ is a sum of both matrices in (A5.8), whereas $\mathcal{X}_{\infty,j}$ has the form that is depicted on the right in (A5.8). (These forms follow from what is said by Proposition A3.5 and Lemma A3.12.) Therefore, (A5.7) need not be true when $n = 0$.

More is said about the forms of the matrices in all of the cases for n in the subsequent parts of this subsection.

*Part 2*: To say more about $\mathcal{X}_{\infty,0}$, let $\mathfrak{s}_{\infty,0}$ denote the tensor that pairs with $\mathcal{X}_{\infty,0}$ to solve the $j=0$ version of (A3.12) for the given mode number n. Having defined $\mathfrak{s}_{\infty,0}$, introduce the tensors $\mathcal{T}_+ = \mathcal{X}_{\infty,0} + i\mathfrak{s}_{\infty,0}$ and $\mathcal{T}_- = \mathcal{X}_{\infty,0} - i\mathfrak{s}_{\infty,0}$ that obey the respective top and bottom equations in (A3.13) using the given mode number n.

This part of the subsection considers those mode numbers with absolute value greater than 1. (The story on $\mathcal{X}_{\infty,0}$ when n = ±1 is told in Part 4; and the n = 0 story is told in Part 5.) Start with the case where the mode number n obey $|n| \geq 3$. Proposition A3.4 says what $\mathcal{T}_+$ looks like, and the mode number n - 2 version says what $\mathcal{T}_-$ looks like. (The mode number n version of the lower equation in (A3.13) is the mode number n - 2 version of the top equation in (A3.13).) Supposing that n ≥ 3, it follows from Proposition A3.4 that the ($k_1 = 0$, $k_2 = 0$) Fourier mode from $\mathcal{X}_{\infty,0}$ for the T action on the s ≥ 0 part of $S^3 - K$ can be written schematically as

$$(a\, e^{(|n|+2)s} + b\, e^{|n|s}) \begin{pmatrix} 1 & i & 0 \\ i & -1 & 0 \\ 0 & 0 & 0 \end{pmatrix} + (a\alpha_{|n|} e^{-|n|s} + b\alpha_{|n|-2} e^{-(|n|-2)s}) \begin{pmatrix} 1 & -i & 0 \\ -i & -1 & 0 \\ 0 & 0 & 0 \end{pmatrix}$$
(A5.10)

where $a$ and $b$ are complex numbers. (The numbers $\alpha_{|n|}$ and $\alpha_{|n|-2}$ come from the respective integer n and n-2 versions of Proposition A3.4.) If n ≤ -3, then the ($k_1 = 0, k_2 = 0$) Fourier mode from $\mathcal{X}_{\infty,0}$ for the T action on the s ≥ 0 part of $S^3 - K$ is the complex conjugate of what is depicted in (A5.10), thus

$$(\bar{a}\, e^{(|n|+2)s} + \bar{b}\, e^{|n|s}) \begin{pmatrix} 1 & -i & 0 \\ -i & -1 & 0 \\ 0 & 0 & 0 \end{pmatrix} + (\bar{a}\bar{\alpha}_{|n|} e^{-|n|s} + \bar{b}\bar{\alpha}_{|n|-2} e^{-(|n|-2)s}) \begin{pmatrix} 1 & i & 0 \\ i & -1 & 0 \\ 0 & 0 & 0 \end{pmatrix}.$$
(A5.11)

This is because the mode number -n version of $\mathcal{X}_K$ for any n ∈ ℤ is the complex conjugates of the mode number n version.

Of particular note: These $|n| \geq 3$ modes in $\mathcal{X}_K$ do not contribute to the $\mathcal{X}_K$ components $\mathcal{X}_K^{33}$ and $\mathcal{X}_K^{13} - i\mathcal{X}_K^{23}$ in (A5.5). However, they do contribute to $\mathcal{X}_K^{11} - i\mathcal{X}_K^{12}$ in (A5.5); and this contribution can be written in terms of $a$ and $b$ as

- $\mathcal{X}_{\infty,0}^{11} - i\mathcal{X}_{\infty,0}^{12} = 2(a\, e^{(n+2)s} + b\, e^{ns})$      when n ≥ 3.
- $\mathcal{X}_{\infty,0}^{11} - i\mathcal{X}_{\infty,0}^{12} = 2(\bar{a}\bar{\alpha}_{|n|} e^{-|n|s} + \bar{b}\bar{\alpha}_{|n|-2} e^{-(|n|-2)s})$ when n ≤ -3.

(A5.12)

Comparing this with (A5.6) (and invoking Proposition A4.1) leads to the following observations in the n ≥ 3 case:

- $a_{n-1} = 2a\, e^{(n+2)R}$ ,    $b_{n-2} = 2b\, e^{nR}$ .

- $a_{-n-1} \approx 2\bar{b}\,\bar{\alpha}_{|n|-2}\, e^{-(|n|-2)R}$ , $b_{-n-2} \approx 2\bar{a}\,\bar{\alpha}_{|n|}\, e^{-|n|R}$ .

(A5.13)

The notation here (and subsequently) has $\approx$ signifying the following: Supposing that $e$ and $e'$ are two complex numbers, then $e \approx e'$ means that $|e - e'| \le c_* \exp(-c_*^{-1} e^{R/4})\, z_{K,n}$.

When n = 2, there is an additional term in (A5.10) that comes from $\mathfrak{T}_-$; it is a constant multiple of the tensor in the top bullet of the n=0 instance of Proposition A3.5:

$$c\, e^{3s} \begin{pmatrix} 0 & 0 & 1 \\ 0 & 0 & i \\ 1 & i & 0 \end{pmatrix} .$$

(A5.14)

This term does not affect the conclusions of (A5.13) when n = ±2 so these identities still hold. The term in (A5.14) does contribute to the depiction of $X_K^{13} - i X_K^{23}$ in (A5.6). In particular, a comparison with the second bullet in (A5.6) leads to the identity $2c = c_1 e^{-3R}$. In the n = -2 case, there is an extra contribution to (A5.11) which is the complex conjugate of (A5.14). The complex conjugate of (A5.14) makes no contribution to the depiction in (A5.6) of $X_K^{13} - i X_K^{23}$. This implies that the coefficient $c_{-3}$ in (A5.6) has norm $|c_{-3}| \approx 0$. (There is a contribution to $c_{-3}$ of size $\approx 0$ from $X_{\infty,j}$.) To summarize:

$$c_1 = 2c\, e^{3R} \quad \text{and} \quad c_{-3} \approx 0.$$

(A5.15)

*Part 3*:  This part of the subsection discusses $X_{\infty,j}$ when the mode number n obeys $|n| \ge 1$. (The n = 0 version of $X_{\infty,j}$ is discussed in the Part 5 and $X_{\infty,j} = 0$ if n = ±2.) To start, fix n to be either 1 or at least 3. It then follows from Lemma A3.11 and Proposition A3.13 that the $(k_1=0, k_2=0)$ Fourier mode from $X_{\infty,j}$ for the T action on the $s \ge 0$ part of $S^3 - K$ can be written as

$$c\, (e^{(1+|n|)s} \begin{pmatrix} 0 & 0 & 1 \\ 0 & 0 & i \\ 1 & i & 0 \end{pmatrix} + \hat{a}_n\, e^{(1-|n|)s} \begin{pmatrix} 0 & 0 & 1 \\ 0 & 0 & -i \\ 1 & -i & 0 \end{pmatrix})$$

(A5.16)

where $c$ is a complex number and $\hat{a}_n$ is from Lemma A3.11. The cases when n = -1 or when n ≤ -3 are the complex complex conjugates of the corresponding version of (A5.16).

Given the preceding depiction, it then follows that $X_{\infty j}$ contributes only to the $X_K^{13} - i X_K^{23}$ term in (A5.6). Comparing (A5.16) and its complex conjugate with what is written in the second bullet of (A5.6) leads to the following: If $n = 1$ or $n \geq 3$, then

$$c_{n-1} = 2c\, e^{(n+1)R}$$

(A5.17)

Suppose next that $n = -1$ or $n \leq -3$. Then $X_{\infty j}$ is the complex conjugate of (A5.16). In particular, looking at the complex conjugate of the (A5.16) and looking at the second bullet of (A5.6) leads to the observation

$$c_{-|n|-1} \approx 2\bar{c}\,\bar{\hat{a}}_{|n|}\, e^{(1-|n|)R} \ .$$

(A5.18)

when $n = -1$ and $n = -3$.

*Part 4*: This part of the subsection talks about the $n = \pm 1$ versions of $X_{\infty,0}$. The discussion starts with the case when $n = +1$. In this case, $\mathfrak{T}_+ = X_{\infty,0} + i\mathfrak{s}_{\infty,0}$ is described in Proposition A3.4 and $\mathfrak{T}_- = X_{\infty,0} - i\mathfrak{s}_{\infty,0}$ is described by Proposition A3.7. Their descriptions gives a depiction of the $(k_1=0, k_2=0)$ Fourier mode from $X_{\infty,0}$ for the T action on the $s \geq 0$ part of $S^3 - K$. This depiction is directly below:

$$a\,(e^{3s}\begin{pmatrix} 1 & i & 0 \\ i & -1 & 0 \\ 0 & 0 & 0 \end{pmatrix} + \alpha_1 e^{-s}\begin{pmatrix} 1 & -i & 0 \\ -i & -1 & 0 \\ 0 & 0 & 0 \end{pmatrix}) +$$

$$b(e^s t\begin{pmatrix} 1 & i & 0 \\ i & -1 & 0 \\ 0 & 0 & 0 \end{pmatrix} + e^s \bar{t}\begin{pmatrix} 1 & -i & 0 \\ -i & -1 & 0 \\ 0 & 0 & 0 \end{pmatrix}) + c e^{3s}\begin{pmatrix} -1 & 0 & 0 \\ 0 & -1 & 0 \\ 0 & 0 & 2 \end{pmatrix} \ .$$

(A5.19)

where $a$, $b$ and $c$ are complex numbers. (The complex number $\alpha_1$ comes courtesy of Proposition A3.4 and the complex number $t$ comes from Proposition A3.7). Comparing (A5.19) with (A5.6) leads to the identifications:

$$a_0 = 2a\,e^{3R} \quad and \quad b_{-1} = 2b\,t\,e^R \quad and \quad c = 0 \ .$$

(A5.20)

The $n = -1$ version of $X_{\infty,0}$ is the complex conjugate of what is written in (A5.20). Comparing the complex conjugate of (A5.20) with (A5.5) leads to additional observations:

$$a_{-2} = 2\bar{b}\,\bar{t}\,e^R \quad \text{and} \quad b_{-3} \approx \bar{a}\bar{\alpha}_1\,e^{-R}\,.$$

(A5.21)

*Part 5*: This part of the subsection talks about the n = 0 case. Let $\mathfrak{T}_+ = \mathcal{X}_{\infty,0} + i\mathfrak{s}_{\infty,0}$ and let $\mathfrak{T}_- = \mathcal{X}_{\infty,0} - i\mathfrak{s}_{\infty,0}$. The former is described by the n = 0 version of Proposition A3.5 and the latter is described by the n = -2 version of Proposition A3.5. These two versions of Proposition A3.5 lead to the following depiction of the ($k_1=0, k_2=0$) Fourier mode from $\mathcal{X}_{\infty,0}$ for the T action on the s ≥ 0 part of $S^3-K$:

$$(ae^{2s} + \bar{a}\bar{\alpha}_0)\begin{pmatrix} 1 & i & 0 \\ i & -1 & 0 \\ 0 & 0 & 0 \end{pmatrix} + (\bar{a}e^{2s} + a\alpha_0)\begin{pmatrix} 1 & -i & 0 \\ -i & -1 & 0 \\ 0 & 0 & 0 \end{pmatrix}$$

$$e^{3s}(c\begin{pmatrix} 0 & 0 & 1 \\ 0 & 0 & i \\ 1 & i & 0 \end{pmatrix} + \bar{c}\begin{pmatrix} 0 & 0 & 1 \\ 0 & 0 & -i \\ 1 & -i & 0 \end{pmatrix})\,,$$

(A5.22)

where $a$ and $c$ are complex numbers, and where $\alpha_0$ is from the second bullet of Proposition A3.5. In the case when n = 0, then the ($k_1=0, k_2=0$) Fourier mode from $\mathcal{X}_{\infty,j}$ for the T action on the s ≥ 0 part of $S^3-K$ can be written as

$$p\,e^s(\mathfrak{z}\begin{pmatrix} 1 & i & 0 \\ i & -1 & 0 \\ 0 & 0 & 0 \end{pmatrix} + \bar{\mathfrak{z}}\begin{pmatrix} 1 & -i & 0 \\ -i & -1 & 0 \\ 0 & 0 & 0 \end{pmatrix})\,,$$

(A5.23)

where $p \in \mathbb{R}$ and where $\mathfrak{z}$ is the complex number from (A3.56). Comparing (A5.22) and (A5.23) with (A5.6) leads to the identifications:

- $a = c = 0$.
- $b_{-2} \approx 2\bar{a}\bar{\alpha}_0$
- $c_{-1} = 2p\mathfrak{z}\,e^R$.

(A5.24)

(Because $a$ and $c$ are zero, it follows from Proposition A3.13 that the components of $\mathcal{X}_\infty$ when n = 0 can be written as $\frac{1}{4}p(\nabla_a j^b + \nabla_b j^a)$ with $p$ from (A5.23) and with j being the 1-form that is described in Lemma A3.12.)

**c)** $\mathcal{X}_T$ **without** $\{c_0, c_{-1}, c_{-2}\}$ **and** $\{a_{-2}, a_{-3}\}$ **and** $\{b_0, b_{-1}, b_{-2}\}$

The upcoming Lemma A5.1 makes an assertion to the effect that the relevant coefficients in (A5.1) are some subset of $\{c_0, c_{-1}, c_{-2}\}$ and $\{a_{-2}, a_{-3}\}$ and $\{b_0, b_{-1}, b_{-2}\}$. (As is explained in the next subsection, $c_{-2}$, $a_{-3}$ and $b_{-2}$ are irrelevant.) To set the stage for the lemma, fix $R > c_*$ and then an $\iota$ invariant element in the kernel of the operator $\mathcal{L}_{g_R}^\dagger$ on $\mathbb{T}_K$. Write this element as $(\mathcal{X}_T, \mathcal{X}_K)$ according to (5.4) and (5.5). Let $z_K$ denote the $L^2$ norm of $\mathcal{X}_K$ on the $s \leq \kappa_*$ part of $(S^3 - N_K) \times S^1$; and let $z_T$ denote the $L^2$ norm of $\mathcal{X}_T$ on the part of $T \times T'$ where the distance to both $p_*$ and $\iota(p_*)$ is at least $\frac{1}{4} t_* e^{-\kappa_*}$. Assume that $z_K + z_T = 1$.

Let $\mathcal{X}_{T0}$ denote the $(k_1 = 0, k_2 = 0)$ Fourier mode of $\mathcal{X}_T$ for the T action on $T \times (T' - (D \cup \iota(D)))$. Invoke Proposition A2.7 using $e^{-R} \mathcal{X}_{T0}$ for $\mathcal{X}$ and the number 2 for the integer N; and use $\mathfrak{X}_T$ in what follows to denote the resulting $\mathfrak{X}_{N=2}$ element from Proposition A2.7. By way of a reminder, $\mathfrak{X}_T$ is an $\iota$-invariant element in the kernel of $\mathcal{L}^\dagger$ on $(T' - (p_* \cup \iota(p_*)))$ that is defined by a real number, $s$, and three sets of complex numbers, $\{c_0, c_1, c_2\}$ and $\{a_0, a_1\}, \{b_0, b_1, b_2\}$. The formula for the entries of $\mathfrak{X}_T$ is as follows:

- $\mathfrak{X}_T^{33} = -(\mathfrak{X}_T^{11} + \mathfrak{X}_T^{22}) = s$ .
- $\mathfrak{X}_T^{13} - i\mathfrak{X}_T^{23} = c_0 + c_1 x + c_2 x^2$ .
- $\mathfrak{X}_T^{11} - i\mathfrak{X}_T^{12} = (a_0 + a_1 x) u + (b_0 + b_1 x + b_2 x^2)$ .

(A5.25)

As explained in Section A2d, the element $\mathfrak{X}_T$ can be viewed as an $\iota$ and T invariant element in the kernel of the operator $\mathcal{L}_{g_T}^\dagger$ on $T \times (T' - (D \cup \iota(D)))$. With $\mathfrak{X}_T$ viewed in this light, $\mathcal{X}_T - e^R \mathfrak{X}_T$ is an $\iota$ invariant element in the kernel of $\mathcal{L}_{g_T}^\dagger$ on $T \times (T' - (D \cup \iota(D)))$. Therefore, it has its own version of (A5.1); but the $\mathcal{X}_T - e^R \mathfrak{X}_T$ version of (A5.1) has neither $s$, nor $(c_0, c_{-1}, c_{-2})$ nor $(a_0, a_{-1})$ nor $(b_0, b_{-1}, b_{-2})$ because the coefficients in (A5.25) were chosen to make the $\mathcal{X}_T - e^R \mathfrak{X}_T$ version of these numbers equal to 0. However, if k is 3 or more, then the numbers $\{c_{-k}, a_{-k-1}, b_{-k}\}$ in the $\mathcal{X}_T - e^R \mathfrak{X}_T$ version of (A5.1) are the same as those in the $\mathcal{X}_T$ version. (Note that these are the coefficients that correspond to the terms in (A5.6) with positive powers of $e^{-(s-R)}$.) For the record, the $\mathcal{X}_T - e^R \mathfrak{X}_T$ version of (A5.1) is written below with $\mathcal{Z}_{T0}$ used in the formulas as shorthand for $(\mathcal{X}_{T0} - e^R \mathfrak{X}_T)$

- $\mathcal{Z}_{T0}^{33} = -(\mathcal{Z}_{T0}^{11} + \mathcal{Z}_{T0}^{22}) = 0$ ,
- $\mathcal{Z}_{T0}^{13} - i\mathcal{Z}_{T0}^{23} = \cdots c_1' z + c_{-3} z^{-3} + \cdots$ ,
- $\mathcal{Z}_{T0}^{11} - i\mathcal{Z}_{T0}^{12} = (\cdots a_1' z + a_0' + a_{-4} z^{-4} + \cdots) \bar{z} + (\cdots b_1' z + b_{-3} z^{-3} + \cdots)$ ,

(A5.26)

The promised Lemma A5.1 bounds the norm $\mathcal{X}_T - e^R \mathfrak{X}_T$.

**Lemma A5.1**: *There exists $\kappa > \kappa_*$ with the following significance: Supposing that $R \geq \kappa^2$, fix an $\iota$-invariant element in the kernel of $\mathcal{L}_{g_R}^\dagger$ on $\mathbb{T}_K$. Define $\mathcal{X}_T$ and $\mathcal{X}_K$ from this element as instructed by (5.4) and (5.5). Assume that $z_T + z_K = 1$. Having defined $\mathcal{X}_T$, then define $\mathfrak{X}_T$ as intructed above. Let $D'$ denote the disk in $T'$ with radius $\kappa e^{-R}$ with center $p_*$. Then*

$$|\mathcal{X}_T - e^R \mathfrak{X}_T| \leq \kappa((r^{-3} e^{-3R} + e^{-R}) + \exp(-c_*^{-1} e^{R/4}))$$

*on $T \times (T' - (D' \cup \iota(D')))$.*

*Proof of Lemma A5.1*: The proof has five steps.

<u>Step 1</u>: Proposition A2.7 asserts (in part) that if $r \in (2e^{-R}, \tfrac{1}{4} t_*)$, then

$$|\mathcal{X}_T - e^R \mathfrak{X}_T| \leq c_* r^{-3} e^{-R} z_K \quad on \ T \times (T' - (D_{2r} \cup \iota(D_{2r}))).$$

(A5.27)

Taking $r = r_0$ with $r_0 > c_*^{-1}$, then this bound gives the lemma's bound for the points in the complement of $T \times (D_{2r_0} - D')$ and its $\iota$ image. The formula in (A5.26) will be used to obtain the lemma's bounds for the points in these last parts of $T \times (T' - (D' \cup \iota(D')))$. This will require suitable bounds on the coefficients that appear in (A5.26).

<u>Step 2</u>: A suitable bound on the primed coefficients in (A5.26) follows from the $r = c_*^{-1}$ version of (A5.27). These bounds read

$$|c_n'| + |a_{n-1}'| + |b_n'| \leq c_*^n e^{-R} z_T$$

(A5.28)

for any given $n \geq 0$. These bounds follow from (A5.27) because $c_n'$, $a_{n-1}'$ and $b_n'$ can be written as linear combinations of components of the integral of $(\tfrac{z}{|z|})^n Z_{T0}$ on the $|z| = \tfrac{1}{16} t_*$ circle in $\mathbb{D}$ and the $|z| = \tfrac{1}{32} t_*$ circle in $\mathbb{D}$. Integrals on two circles are needed to distinguish the contributions from $a_{n-1}'$ and $b_n'$.

The bound in (A5.28) will be used in Step 5 to bound the contribution of the primed coefficients in (A5.26) to the norm of $\mathcal{X}_T - e^R \mathfrak{X}_T$.

Step 3: A suitable bound on the coefficients in (A5.26) with negative index, thus $\{c_{-k}, a_{-k-1}, b_{-k}\}_{k=3,4,...}$ requires the results of the preceding subsection. To obtain the desired bound, note first that

$$|c_n| + |a_{n-1}| + |b_n| \leq c_*^{n+1} e^R z_T \quad \text{when } n \geq 0 ,$$

(A5.29)

because anything larger than the right hand side will run afoul of the fact that $z_T$ is an upper bound for the $L^2$ norm of $\mathcal{X}_T$ on the part of $\mathbb{D}$ where $|z| \geq c_*^{-1}$. With (A5.29) in hand, fix $n \geq 2$ for the moment and let $(a, b)$ denote the coefficients from n's version of (A5.13). The top bullet of (A5.13) and (A5.29) lead to the bounds

$$|a| \leq c_*^n e^{-(n+1)R} z_T \quad \text{and} \quad |b| \leq c_*^n e^{-(n-1)R} z_T ,$$

(A5.30)

and these with the lower bullet in (A5.13) lead in turn to

- $|a_{-n-1}| \leq c_*^n e^{-(2n-3)R} z_T + \exp(-c_*^{-1} e^{R/4}) z_K$
- $|b_{-n-2}| \leq c_*^n e^{-(2n+1)R} z_T + \exp(-c_*^{-1} e^{R/4}) z_K$

(A5.31)

The lower bullet in (A5.31) also holds when $n = 1$ because of the left hand equality in (A5.20) and the right hand equality in (A5.21), and because $|a_0| \leq c_* e^R z_T$. By the same token, if $n = 1$ or if $n \geq 3$, then

$$|c_{-n-1}| \leq c_*^n e^{-(2n-1)R} z_T + \exp(-c_*^{-1} e^{R/4}) z_K .$$

(A5.32)

Indeed, this follows from (A5.29) with (A5.17) and (A5.18). This also holds for $n = 2$ because of the right hand inequality in (A5.15).

The bounds in (A5.31) and (A5.32) with what is said in Step 4 is used in Step 5 to bound the contribution of the uprimed coefficients in (A5.26) to the norm of $\mathcal{X}_T - e^R \mathfrak{X}_T$.

Step 4: Fix $r_0 \in (2e^{-R}, \frac{1}{4} t_*)$ and invoke the $r = r_0$ version of (A5.27) to prove that $|\mathcal{X}_T - e^R \mathfrak{X}_T| \leq c_* r_0^{-3} e^{-R} z_T$ on $T \times (T' - (D_{2r_0} \cup \iota(D_{2r_0})))$. To bound the norm of $\mathcal{X}_T - e^R \mathfrak{X}_T$ on $T \times (D_{2r_0} - D')$, fix $N \geq c_*$ for the moment to invoke Proposition A4.2 with $r = R + \ln(\frac{1}{4} t_*)$ and with $X$ being $\mathcal{X}_K$. As explained directly, this proposition implies that

$$e^{-4R} \sum_{n \geq 2N} c_*^{-n} e^{2nR} (|c_{-n}|^2 + |a_{-n-1}|^2 + |b_{-n}|^2) \leq c_* (e^{-NR/c_*} + \exp(-\kappa^{-1} e^{r/2})) z_K$$

(A5.33)

when $N \geq c_*$. To see this, let $\mathcal{X}_{K,N}$ denote the sum of the Fourier modes of $\mathcal{X}_K$ with mode number obeying $|n| \leq N$. The left hand side of (A5.33) is no greater than the $L^2$ norm of $\mathcal{X}_K - \mathcal{X}_{K,N}$ on the $s \in [1, 2]$ part of $(S^3 - N_K) \times S^1$. This is a consequence of (A5.6). The right hand side of (A5.33) is the bound for the $L^2$ norm of $\mathcal{X}_K - \mathcal{X}_{K,N}$ from Proposition A4.2.

The bound in (A5.33) leads directly to the following: If $|z| \in (c_* e^{-R}, c_*^{-1})$, then

$$\sum_{n \geq 2N} |z|^{-2n} (|c_{-n}|^2 + |a_{-n-1}|^2 + |b_{-n}|^2) \leq c_* (e^{-NR/c_*} + \exp(-\kappa^{-1} e^{r/2})) z_K.$$
(A5.34)

<u>Step 5</u>: Fix $N \geq c_*$ so that (A5.34) holds. Write $\mathcal{X}_T - e^R \mathfrak{X}_T$ as $A + B + C$ where $A$ is the sum of the terms in (A5.26) with primed coefficients, and $B$ is the sum of the terms with unprimed coefficients $\{c_{-n}, a_{-n-1}, b_{-n} : 3 \leq n \leq 2N\}$. Meanwhile, $C$ is the sum of the terms in (A5.26) with unprimed coefficients $\{c_{-n}, a_{-n-1}, b_{-n} : n > 2N\}$. As explained directly, if $r \in (c_* e^{-R}, c_*^{-1})$ and if $|z| \in (r, 2r)$, then the norms of $A$, $B$ and $C$ obey

- $|A| \leq \sum_{n \geq 1} c_*^n r^n e^{-R} z_T$,
- $|B| \leq \sum_{3 \leq n \leq 2N} c_*^n (r^{-n} e^{-(2n-3)R} z_T + \exp(-c_*^{-1} e^{R/4})) z_K$,
- $|C| \leq c_* (e^{-NR/c_*} + \exp(-\kappa^{-1} e^{r/2})) z_K$.

(A5.35)

The $A$ bound is a consequence of (A5.28); the $B$ bound is a consequence of (A5.31) and (A5.32); and the $C$ bound is a consequence of (A5.34). If $N \geq c_*$ and if $r \leq r_0$ and $r_0 \leq c_*^{-1}$, then the $A$ bound in (A5.35) is no greater than $c_* e^{-R} z_T$. Meanwhile, the $B$ bound is at most

$$c_* (r^{-3} e^{-3R} + e^{-R}) z_T + c_* e^{4NR} \exp(-c_*^{-1} e^{R/4}) z_K ;$$
(A5.36)

If $N \in [c_*, c_*^2]$, then the $C$ bound in (A5.35) is at most $c_* e^{-R} z_T$ also.

### d) Proof of Proposition 5.3

Proposition 5.3 is a corollary to the upcoming Proposition A5.2. The homomorphism $\mathcal{Q}_T$ in Proposition 5.3 is the composition $\mathcal{I}_T \circ \mathcal{Y}_T$ of the two homomorphism $\mathcal{I}_T$ and $\mathcal{Y}_T$ that are defined momentarily and appear in Proposition A5.2.

To set the stage and the notation for this proposition, let $\kappa_\diamond$ denote the version of $\kappa$ from Lemma A5.1. Having fixed $R \geq \kappa_\diamond^2$, let $\mathcal{H}_{K,R}$ denote the $\iota$ invariant kernel of the operator $\mathcal{L}_{g_R}^\dagger$ on the R version of $\mathbb{T}_K$. Supposing that $\mathcal{X}$ is an element in this kernel, write it as a pair $(\mathcal{X}_T, \mathcal{X}_K)$ using the rules in (5.4) and (5.5). Proposition A5.2 uses $z_T$ to denote the $L^2$ norm of $\mathcal{X}_T$ on the part of $T \times T'$ where the distance to both $p_*$ and $\iota(p_*)$ is greater than $\frac{1}{4} t_* e^{-\kappa_\diamond}$. Meanwhile, the proposition uses $z_K$ to denote the respective $L^2$ norm of $\mathcal{X}_K$ on the $s \le \kappa_\diamond$ part of its domain, $(S^3 - N_K) \times S^1$.

Proposition A5.2 refers to a homomorphism (it is called $\Upsilon_T$) from $\mathcal{H}_{K,R}$ to $\mathbb{R} \oplus \mathbb{C}^3$ that is defined as follows: Introduce $\mathcal{X}_{T0}$ again to denote the $(k_1 = 0, k_2 = 0)$ Fourier mode from $\mathcal{X}_T$ for the T action on $T \times (T' - (D \cup \iota(D)))$. Use $\mathcal{X}_{T0}$ as instructed in the previous subsection (see (A5.25)) to obtain, in particular, the complex numbers $\{c_0, c_1\}$ and $a_0$ and $\{\mathfrak{b}_0, \mathfrak{b}_1\}$. The homomorphism $\Upsilon_T$ is defined by the rule

$$X \to \Upsilon_T(X) = (\mathfrak{Re}(\bar{\mathfrak{z}} \, c_1), c_0, \mathfrak{b}_0, \mathfrak{b}_1) \, .$$

(A5.33)

(The notation here uses $\mathfrak{Re}(\cdot)$ to denote the real part of the given complex number.)

The statement of Proposition A5.2 also refers to a homomorphism (denoted by $\mathfrak{I}_T$) from $\mathbb{R} \oplus \mathbb{C}^3$ to the $\iota$ invariant kernel of the operator $\mathcal{L}_{\vartheta_T}^\dagger$ on $T \times (T' - \{p_*, \iota(p_*)\})$. To define $\mathfrak{I}_T$, agree first to depict elements in the domain of $\mathcal{L}_{\vartheta_T}^\dagger$ as symmetric, traceless matrix valued functions using the basis for $\Lambda^+$ of $T \times T'$ that from Section A5a. Now, if $(t, w_0, y_0, y_1) \in \mathbb{R} \oplus \mathbb{C}^3$, then its $\mathfrak{I}_T$ image is the symmetric, traceless matrix with entries

- $\mathfrak{I}_T^{33} = -(\mathfrak{I}_T^{11} + \mathfrak{I}_T^{22}) = 0$,
- $\mathfrak{I}_T^{13} - i \mathfrak{I}_T^{23} = e^R (w_0 + t \mathfrak{z} \, x)$,
- $\mathfrak{I}_T^{11} - i \mathfrak{I}_T^{12} = e^R (y_0 + y_1 x + \bar{y}_1 u)$.

(A5.34)

What is denoted here by $\mathfrak{z}$ is the complex number that appears in (A3.56).

**Proposition A5.2**: *There exists $\kappa > \kappa_\diamond$ such that if $R > \kappa^2$, then the map $\Upsilon_T : \mathcal{H}_{K,R} \to \mathbb{R} \oplus \mathbb{C}^3$ has the following property: Letting $X \in \mathcal{H}_{K,R}$ denote a given element write it in the manner of (5.4) and (5.5) as a pair $(\mathcal{X}_T, \mathcal{X}_K)$. Assume that $z_K + z_T = 1$. Let $D'$ denote the disk in $T'$ with center $p_*$ and radius $\kappa e^{-R}$. Then*

$$|\mathcal{X}_T - (\mathfrak{I}_T \circ \Upsilon_T)(X)| \le \kappa((r^{-2} e^{-2R} + e^{-R}) + \exp(-c_*^{-1} e^{R/4}))$$

*on the domain $T \times (T' - (D' \cup \iota(D')))$.*

Proposition 5.3 follows from Proposition A5.2 if the homomorphism $\mathcal{Q}_T$ in Proposition 5.3 is the composition $\mathfrak{I}_T \circ \Upsilon_T$.

*Proof of Proposition A5.2*: The proof has two parts.

*Part 1*: The identities in Section A5b with Proposition A4.1 constrain the data $r$, $\{c_0, c_1, c_2\}$ and $\{a_0, a_1\}$, $\{b_0, b_1, b_2\}$ from (A5.25). Here is the first set of constraints:

- $r \approx 0$.
- $|c_2| \le c_* e^{-2R}$
- $|a_1| \le c_* e^{-2R}$.
- $b_2 \approx 0$.

(A5.35)

Indeed, the observation that $r \approx 0$ follows from Proposition A4.1 because the mode number 0 terms in (A5.22) and (A5.23) have no diagonal elements. The observation that $b_2 \approx 0$ follows from the top and middle bullets of (A5.24) because $|b_2| \le c_* |b_{-2}| e^{-R}$. The derivation of the bound for the norm of $c_2$ starts with the fact that $|c_2| \le |c_{-2}| e^{-R}$. Meanwhile, the n = 1 version of (A5.18) leads to the bound $|c_{-2}| \le c_* |c|$, and (A5.17) in the case n = 1 implies that $|c| \le c_* e^{-R}$. The derivation of the bound for the norm of $a_1$ starts with the fact that $|a_1| \le c_* |a_{-3}| e^{-R}$. The second bullet of the n = 2 version of (A5.13) bounds $|a_{-3}|$ by $c_* |b|$, and the top bullet in the same version of (A5.13) bounds $|b|$ by $c_* e^{-R}$. Therefore, $|a_1| \le c_* e^{-2R}$. The bounds

$$|c_{-2}| + |a_{-3}| + |b_{-2}| \le c_* e^{-R}$$

(A5.36)

are stated here for the record.

The coefficients $c_1$ is also constrained: Let $\mathfrak{z}$ denote the unit length complex number that appears in (A3.56). Up to a small error, $c_1$ is a real multiple of $\mathfrak{z}$,

$$|c_1 - 2q\mathfrak{z}| \le c_*(e^{-2R} + \exp(-c_*^{-1} e^{R/4}))$$

(A5.37)

with $q$ being a real number. This observation follows from the third bullet of (A5.24) and the $|c_2|$ bound from the second bullet of (A5.35).

The last constraint concerns the pair $(a_0, b_1)$, they are nearly complex conjugates:

$$|a_0 - \overline{b_1}| \le c_*(e^{-2R} + \exp(-c_*^{-1} e^{R/4})) .$$

(A5.38)

Indeed, this follows from (A5.20) and (A5.21) because $a_0$ differs from $a_{-2}$ by no more than $c_* e^{-2R}$ given the bound in the second bullet of (A5.35); and $b_1 \approx b_{-1}$ given the bound in the fourth bullet of (A5.35).

*Part 2*: Define $\mathcal{X}_T$ as in (A5.25). Then $\mathcal{X}_T$ can be written as $(\mathcal{I}_T \circ \Upsilon_T)(\mathcal{X}) + \mathfrak{P}$ with $\mathfrak{P}$ denoting the traceless, symmetric matrix valued function on $T \times (T' - \{p_*, \iota(p_*)\})$ whose components are as follows:

- $\mathfrak{P}^{33} = -(\mathfrak{P}^{11} + \mathfrak{P}^{22}) = s$.
- $\mathfrak{P}^{13} - i\mathfrak{P}^{23} = i \, \mathfrak{Im}(\overline{\mathfrak{z}} c_1) \, \mathfrak{z} x + c_2 x^2$.
- $\mathfrak{P}^{11} - i\mathfrak{P}^{12} = ((a_0 - \overline{b}_1) + a_1 x) u + b_2 x^2$.

(A5.39)

Supposing that $r \in (c_* e^{-R}, \frac{1}{4} t_*)$, then what is said in Part 1 implies that

$$|\mathfrak{P}| \le c_* (r^{-2} e^{-2R} + \exp(-c_*^{-1} e^{R/4}))$$

(A5.40)

on $T \times (T' - (D_r \cup \iota(D_r)))$. Let $D'$ denote the disk in $T'$ with center $p_*$ and radius $c_* e^{-R}$. Lemma A5.1 and the bound in (A5.40) imply in turn that

$$|\mathcal{X}_T - e^R(\mathcal{I}_T \circ \Upsilon_T)(\mathcal{X})| \le c_* ((r^{-2} e^{-2R} + e^{-R}) + \exp(-c_*^{-1} e^{R/4}))$$

(A5.41)

on $T \times (T' - (D' \cup \iota(D')))$, which is what is asserted by Proposition A5.2.

**e) Proof of Proposition 5.4**

The upcoming Proposition A5.3 implies what is said by Proposition 5.4. Proposition A5.3 refers to a homomorphism (it is called $\Upsilon_K$) from $\mathcal{H}_{K,R}$ to the (complex) vector space $H^1(S^3 - K; V) \otimes_{\mathbb{R}} \mathbb{C}$ that is defined as follows: Construct the tensor $\mathcal{A}$ from $\mathcal{X}_K$ using (A3.1), and write $\mathcal{A}$ in terms of the tensor $\mathfrak{s}$ and the 1-form $\mathfrak{j}$ as done in (A3.3). Let $\mathfrak{T}$ denote the symmetric, traceless section of $\otimes^2 T^*(S^3 - N_K)$ over $(S^3 - N_K) \times S^1$ whose entries with respect to any given local orthonormal frame are $\{\mathcal{X}^{ab} - i\mathfrak{s}^{ab} - \frac{1}{2}(\nabla_a \mathfrak{j}^b + \nabla_b \mathfrak{j}^a)\}_{a,b=1,2,3}$. Write the $n = 1$ Fourier mode of $\mathfrak{T}$ as $e^{i\theta} \mathfrak{T}^{(1)}$. It follows from Propositions A3.7 and A3.13 that $\mathfrak{T}^{(1)}$ is a (complex) Codazzi tensor. As explained in Part 5 of Section A3c, any real Codazzi tensor on $S^3 - N_K$ has a corresponding Ferus/Lafontaine class in $H^1(S^3 - K; V)$. Since $\mathfrak{T}^{(1)}$ has complex coefficients, its Ferus/Lafontaine class is in $H^1(S^3 - K; V) \otimes_{\mathbb{R}} \mathbb{C}$. The element $\Upsilon_K(\mathcal{X})$ is this Ferus/Lafontaine class of $\mathfrak{T}^{(1)}$.

Proposition A5.3 also refers to a homomorphism called $\mathcal{I}_K$ that maps the vector space $H^1(S^3-K;V)\otimes_\mathbb{R}\mathbb{C}$ to the kernel of the operator $\mathcal{L}^\dagger_{g_K}$ on $(S^3-K)\times S^1$. To define $\mathcal{I}_K$, choose $\mathfrak{h}\in H^1(S^3-K;V)\otimes_\mathbb{R}\mathbb{C}$ and write $\mathfrak{h}$ as $\mathfrak{h}_1+i\mathfrak{h}_2$ with $\mathfrak{h}_1$ and $\mathfrak{h}_2$ being real valued Codazzi tensors. Assign to $\mathfrak{h}$ the complex Codazzi tensor $\mathfrak{T}_C(\mathfrak{h}_1)+i\mathfrak{T}_C(\mathfrak{h}_2)$ with $\mathfrak{T}_C(\cdot)$ from the second bullet in Proposition A3.7. Denote this complex Codazzi tensor by $\mathfrak{T}_C(\mathfrak{h})$. Now set $\mathcal{I}_K(\mathfrak{h}) = \frac{1}{2}(e^{i\theta}\mathfrak{T}_C(\mathfrak{h}) + e^{-i\theta}\overline{\mathfrak{T}_C(\mathfrak{h})})$. (It follows from (A3.11) and (A3.12) that this is in the kernel of $\mathcal{L}^\dagger_{g_K}$ on $(S^3-K)\times S^1$.)

A certain homomorphism from $(\mathbb{R}\oplus\mathbb{C}^3)\oplus(H^1(S^3-N_K)\otimes_\mathbb{R}\mathbb{C})$ to $\mathbb{C}$ appears in Proposition A5.3, it is denoted by $\lambda_K$ and it is defined as follows: Let L denote the linear extension to $H^1(S^3-N_K)\otimes_\mathbb{R}\mathbb{C}$ of the eponymous homomorphism in Proposition A3.7. This extention is a homomorphism $H^1(S^3-N_K)\otimes_\mathbb{R}\mathbb{C}$ to $\mathbb{C}$. Let $\mathfrak{t}$ denote the unit length complex number from Proposition A3.7. The linear functional $\lambda_K$ sends any given data set $((\mathfrak{t}, w_0, y_0, y_1), \mathfrak{h})$ to the complex number $y_1 - 2L(\mathfrak{h})\mathfrak{t}$.

**Proposition A5.3**: *There exists $\kappa > 1$ such that if $R > \kappa$, then the homomorphism*

$$\Upsilon_T\oplus\Upsilon_K : \mathcal{H}_{K,R} \to (\mathbb{R}\oplus\mathbb{C}^3) \oplus (H^1(S^3-N_K)\otimes_\mathbb{R}\mathbb{C})$$

*is injective. Moreover, if $X\in\mathcal{H}_{K,R}$ is written in the manner of (5.4) and (5.5) as a pair $(X_T, X_K)$ and if these are such that $z_T + z_K = 1$, then*
- $|\lambda_K(X)| \leq \kappa\exp(-\kappa^{-1}e^{R/4})$.
- $|X_K - (\mathcal{I}_K\circ\Upsilon_K)(X)| \leq \kappa\, e^{2s-R}$ *on the* $s \leq \frac{5}{8}R$ *part of* $(S^3-N_K)\times S^1$.

*Proof of Proposition A5.3*: The assertion that $\Upsilon_T\oplus\Upsilon_K$ is injective follows jointly from Proposition A5.2 and from the second bullet of this proposition which is proved momentarily. To prove the first bullet of the proposition, reintroduce the number $\mathfrak{b}$ from (A5.19). The number $\mathfrak{b}$ is $L(\Upsilon_K(X))$, this being a direct consequence of the definition of $\Upsilon_K$ and the formula in the second bullet of Proposition A3.7. The number $\mathfrak{b}$ is also related to the number $b_{-1}$ from (A5.1) by the rule $b_{-1} = 2\mathfrak{b}\mathfrak{t}\,e^R$ in (A5.20). Meanwhile, $b_{-1}e^{-R}$ differs from the number $\mathfrak{b}_1$ in (A5.25) by at most $c_*\exp(-c_*^{-1}e^{R/4})$ because of the bound in (A5.35) on the norm of $\mathfrak{b}_2$.

The proof of the second bullet of Proposition A5.2 has two parts.

*Part 1*: Return now to the milieu of Section A5b. Fix $n \geq 3$ and then reintroduce $\mathcal{X}_{\infty,0}$ from Part 2 of Section A5b. The top bullet in (A5.13) holds for $n \geq 3$, and the latter identity with (A5.29) imply the bound

$$|\mathcal{X}_{\infty,0}| \leq c_*^n \, (e^R e^{n(s-R)} + e^{-(n-1)R}(e^{-(|n|-2)s} + \exp(-c_*^{-1} e^s)))$$

(A5.43)

on the $s \leq R$ part of $S^3 - K$. This bound is also obeyed when $n = 2$ because the top bullet in (A5.13) holds when $n = 2$ and because the coefficient $c$ that appears in (A5.15) is bounded by $c_* e^{-2R}$. (The latter bound follows because $c = \frac{1}{2} c_1 e^{-3R}$ with $c_1$ coming from (A5.1).)

Supposing that $n = 1$ or that $n \geq 3$, reintroduce the corresponding version of $\mathcal{X}_{\infty,j}$ from Part 3 of Section A5b. The identity in (A5.17) and the bounds in (A5.29) imply that

$$|\mathcal{X}_{\infty,j}| \leq c_*^n \, (e^R e^{(n+1)(s-R)} + e^{-nR}(e^{-(|n|-1)s} + \exp(-c_*^{-1} e^s))) \, .$$

(A5.44)

*Part 2*: Let $(q, c_0, b_0, b_1) \in \mathbb{R} \times \mathbb{C}^3$ denote $\Upsilon_T(\mathcal{X})$. The number $q$ and $b_1$ are needed for this part of the proof. With regards to $q$, let $\mathcal{X}_{(0)}$ denote the $n = 0$ version of $\mathcal{X}_{\infty,j}$ from Part 5 of Section 11b with the number $p$ from (A5.23) equal to $\frac{1}{4} q$. To say more about $\mathcal{X}_{(0)}$, let $j$ denote for the moment the harmonic 1-form on $S^3 - K$ whose $(k_1 = 0, k_2 = 0)$ Fourier mode on the $s \geq 0$ part of $S^3 - K$ is depicted in (A3.56). The coefficients of $\mathcal{X}_{(0)}$ with respect to any given orthonormal frame are $\{\frac{1}{16} q(\nabla_a j^b + \nabla_b j^a)\}_{a,b=1,2,3}$. This $\mathcal{X}_{(0)}$ is significant because the sum of the $n = 0$ versions of $\mathcal{X}_{\infty,0}$ and $\mathcal{X}_{\infty,j}$ from Part 5 of Section A5b is very nearly $\mathcal{X}_0$:

$$|\mathcal{X}_{\infty,0} + \mathcal{X}_{\infty,j} - \mathcal{X}_0| \leq c_* \exp(-c_*^{-1} e^s)) \, ,$$

(A5.45)

which is a consequence of (A5.24) and Proposition A3.1.

Meanwhile, the number $b_1$ from $\Upsilon_T$ is determined almost entirely by $\Upsilon_K(\mathcal{X})$ via the first bullet in Proposition A5.3. It follows as a consequence of this observation (and the bound $|a_0| \leq c_* e^R z_T$ from (A5.29)) that the $n = 1$ version of $\mathcal{X}_{\infty,0}$ from Part 4 of Section A5b is nearly the same as $(\mathcal{I}_K \circ \Upsilon_K)(\mathcal{X})$ in the following sense:

$$|\mathcal{X}_{\infty,0} - (\mathcal{I}_K \circ \Upsilon_K)(\mathcal{X})| \leq c_*((e^R e^{3(s-R)} + e^{-2R} e^{-s}) + \exp(-c_*^{-1} e^{R/4})) \, .$$

(A5.46)

This bound follows from (A5.20) and (A5.19).

*Part 3*: Fix a large positive integer to be denoted by N; upper and lower bounds will be given momentarily. Let $\mathcal{X}_{K,N}$ denote the sum of the $S^1$ action Fourier modes of $\mathcal{X}_K$ with the mode number having norm at most N. Use this in the $r = R + \ln(\frac{1}{4} t_*)$ version of Proposition A4.2 to obtain the pointwise bound

$$|\mathcal{X}_K - \mathcal{X}_{K,N}| \leq c_*(e^{Nc_*} e^{-NR/4} + \exp(-c_*^{-1} e^{r/4}))$$
(A5.47)

on the $s \leq \frac{5}{8} R$ part of $(S^3 - N_K) \times S^1$. Write $\mathcal{X}_{K,N}$ as the Fourier mode sum $\sum_{|n| \leq N} \mathcal{X}^{(n)}$ with any given version of $\mathcal{X}^{(n)}$ being the Fourier mode number n part of $\mathcal{X}$ for the $S^1$ action. Invoke Proposition A4.1 using each $|n| \leq N$ version of $\mathcal{X}^{(n)}$ and $r = R + \ln(\frac{1}{4} t_*)$ to obtain the corresponding $\mathcal{X}_\infty^{(n)}$. Let $\mathcal{X}_{N\infty} = \sum_{|n| \leq N} \mathcal{X}_\infty^{(n)}$. This is close to $\mathcal{X}_N$ in the sense that

$$|\mathcal{X}_{K,N} - \mathcal{X}_{N\infty}| \leq c_* \sum_{|n| \leq N} (e^{-|n|s} e^{-|n|r/2} \exp(-c_*^{-1} e^{r/4}) + \exp(-c_*^{-1} e^{(r-s)}))$$
(A5.48)

on the $s \leq R + \ln(\frac{1}{4} t_*) - c_*$ part of $(S^3 - N_K) \times S^1$. This sum implies the bound

$$|\mathcal{X}_{K,N} - \mathcal{X}_{N\infty}| \leq c_*(\exp(-c_*^{-1} e^{r/4}) + N\exp(-c_*^{-1} e^{(r-s)}))$$
(A5.49)

on this same part of $(S^3 - N_K) \times S^1$. In turn, $\mathcal{X}_{N\infty}$ can be compared with $\mathcal{X}_0$ and $(\mathcal{I}_K \circ \mathcal{Y}_K)(\mathcal{X})$ using the bounds in (A5.43)-(A5.46). The latter lead directly to the following bound:

$$|\mathcal{X}_{N\infty} - \mathcal{X}_{(0)} - (\mathcal{I}_K \circ \mathcal{Y}_K)(\mathcal{X})| \leq c_*((e^R e^{3(s-R)} + e^{-2R} e^{-s}) + N\exp(-c_*^{-1} e^{R/4})).$$
(A5.50)

Taken together, the bounds in (A5.47) and (A5.49) and (A5.50) and those from Parts 1 and 2 imply what is asserted by the second bullet of Proposition A5.3 if $R \geq c_*$ and $N = e^{R/32}$. Note in this regard that the n = 2 case of (A5.43) and the n = 1 case of (A5.44) are the largest contributions to the right hand side of the inequality in the second bullet of Proposition A5.3; and these account for the specific powers of $e^s$ (two) and $e^{-R}$ (one) that appear in this inequality.